\input rztex.qel

\beginOptionen
\Brotschrift{\cmrom 11 pt}
\Rubriken{\lang}
\Verweise
\Verzeichnisse {\standard}
\endOptionen

\uKopfzeile{}{\Seitennummer}{}
\uKopflinie 0 pt
\Rand = 2.75 cm
\Kolumnenbreite = 16.25 cm

\redefine\tagform#1{\Abschnittswurzel\Zahl\aktAbschnitt -#1}
\redefine\tagstyle#1{\rm(#1)}

\def\Seq#1{ \{\mkern -4 mu [ #1 ]\mkern -4  mu \} }
\def\Seqn#1{ \{\mkern -4 mu [ #1_n ]\mkern -4 mu \} }

\def\Ent#1{ [\mkern - 2.5 mu  [ #1 ] \mkern - 2.5 mu ] }

\def\bigSeq#1{ \bigl\{\mkern -6 mu \bigl[
#1 \bigr]\mkern -6  mu \bigr\} }

\def\F{\text{I} \mkern - 3.5 mu \text{F} }
\def\D{\text{I} \mkern - 3.5 mu \text{D} }
\def\H{\text{I} \mkern - 3.5 mu \text{H} }

\def\Bar{\, \vert \,}
\def\bigBar{\, \bigl\vert \bigr. \,}

\def\biggBar{\, \biggl\vert \biggr. \,}

\Titelblatt

\bigskip

\bigskip

\bigskip

{\LARGE \bf

\centreline {NONLINEAR SEQUENCE TRANSFORMATIONS}

\medskip

\centreline {FOR THE}

\medskip

\centreline {ACCELERATION OF CONVERGENCE}

\medskip

\centreline {AND THE}

\medskip

\centreline {SUMMATION OF DIVERGENT SERIES}
}

\bigskip

\bigskip

\bigskip

\bigskip

{\large
\centreline {Ernst Joachim Weniger}
\centreline {Institut f\"ur Physikalische und Theoretische Chemie}
\centreline {Universit\"at Regensburg}
\centreline {D-8400 Regensburg}
\centreline {Federal Republic of Germany}
}

\bigskip

\bigskip

\bigskip

\bigskip

\bigskip

\bigskip

\vfill

\Ueberschrift \large Abstract

\bigskip

\noindent Slowly convergent series and sequences as well as divergent
series occur quite frequently in the mathematical treatment of
scientific problems. In this report, a large number of mainly nonlinear
sequence transformations for the acceleration of convergence and the
summation of divergent series are discussed. Some of the sequence
transformations of this report  as for instance Wynn's $\epsilon$
algorithm or Levin's sequence transformation are well established in
the literature on convergence acceleration, but the majority of them is
new. Efficient algorithms for the evaluation of these transformations
are derived. The theoretical properties of the sequence transformations
in convergence acceleration and summation processes are analyzed.
Finally, the performance of the sequence transformations of this report
are tested by applying them to certain slowly convergent and divergent
series, which are hopefully realistic models for a large part of the
slowly convergent or divergent series that can occur in scientific
problems and in applied
mathematics.

\bigskip
\bigskip
\bigskip
\bigskip
\bigskip

\centreline {\bf Published as Computer Physics Reports Vol.\ 10, 189 -
371 (1989)}

\neueSeite
\Titelblatt

\ungeradeSeite

\neueSeite

\Pageno = 1

\redefine \Seitennummer {\roemisch \pageno}

\Ueberschrift \large Contents

\medskip

\Inhaltsverzeichnis

\keinTitelblatt\neueSeite

\redefine \Seitennummer {\Zahl \pageno}
\Pageno = 1

\keinTitelblatt\neueSeite

\beginAbschnittsebene
\aktAbschnitt = 0

\Abschnitt Introduction

\vskip - 2 \jot

\beginAbschnittsebene

\medskip

\Abschnitt Infinite series and their evaluation

\smallskip

\aktTag = 0

Infinite series are ubiquitous in the mathematical analysis of
scientific problems. They naturally appear in the evaluation of
integrals, in the solutions of differential and integral equations, or
as Fourier series. They are also used for both the definition and the
evaluation of many of the special functions of mathematical physics.
The conventional approach for the evaluation of an infinite series
consists in computing a finite sequence of partial sums
$$
s_n = \sum_{k=0}^{n} \, a_k
\tag
$$
by adding up one term after the other. Then, the magnitude of the
truncation error is estimated. If the sequence of partial sums $s_0,
\ldots , s_n$ has not converged yet to the desired accuracy, additional
terms must be added until convergence has finally been achieved. With
this approach it is at least in principle possible to determine the
value of an infinite series as accurately as one likes provided that
one is able to compute a sufficiently large number of terms accurately
enough to overcome eventual numerical instabilities.

However, in many scientific problems one will only be able to compute a
relatively small number of terms. In addition, particularly the series
terms with higher summation indices are often affected by serious
inaccuracies which may lead to a catastrophic accumulation of round-off
errors.

Consequently, if an infinite series is to be evaluated by adding one
term after the other, an infinite series will be of practical use only
if it converges after a sufficiently small number of terms.
Unfortunately, many counterexamples are known in which alternative
methods for the evaluation of infinite series must be used since in
these cases the conventional approach of evaluating an infinite series
does not suffice.

For instance, when Haywood and Morgan [1] performed a discrete
basis-set calculation of the Bethe logarithm of the $1s$ state of the
hydrogen atom, they found that even 120 basis functions gave no more
than 2 -- 3 decimal digits and they estimated that approximately
$10^{10}$ basis functions would be needed to obtain an accuracy of more
than 10 decimal digits. Haywood and Morgan also showed that with the
help of a suitable convergence acceleration method an accuracy of more
than 13 decimal digits can be extracted from their data.

A good mathematical model for the convergence problems which Haywood
and Morgan [1] encountered in their calculation of the Bethe logarithm
is the following series expansion for the Riemann zeta function:
$$
\zeta (z) \; = \; \sum _{n=0}^{\infty} \; (n+1)^{-z} \, .
\tag
$$

It is well known that this infinite series converges if $Re(z) > 1$
holds. However, if $Re(z)$ is only slightly larger than one, the rate
of convergence becomes extremely slow. For instance, Bender and Orszag
remark in their book (see p. 379 of ref. [2]) that about $10^{20}$
terms of the above series expansion would be needed to compute $\zeta
(1.1)$ accurate to one percent. Bender and Orszag also show that only
10 terms of the series in connection with a specially designed
acceleration method are needed to compute $\zeta (1.1)$ to 26 decimal
digits (see table 8.7 on p. 380 of ref. [2]).

Even more striking examples for the inadequacy of the conventional
approach towards the evaluation of infinite series are some
Rayleigh-Schr\"odinger perturbation expansions of elementary quantum
mechanical systems. For instance, if the following normalization for
the Hamiltonian of the quartic anharmonic oscillator is used,
$$
{\hat H} \; = \; {\hat p^2} + {\hat x^2} +
\beta {\hat x^4} \, ,
\tag
$$
then it follows from the results obtained by Bender and Wu (see eq.
(1.8) of ref. [3]) that the coefficients $c_n$ of the power series in
the coupling constant $\beta$ for the ground state energy eigenvalue
$E_0 (\beta )$ of the quartic anharmonic oscillator,
$$
E_0 (\beta ) \; = \;
\sum_{n=0}^{\infty} \; c_n \; {\beta}^n \, ,
\tag
$$
possess the following asymptotic behaviour:
$$
c_n \; \sim \; (-1)^{n+1} \, (3/2)^n \, \Gamma(n+1/2) \, ,
\qquad n \to \infty \, .
\tag
$$

The radius of convergence of the above Rayleigh-Schr\"odinger
perturbation series is obviously zero, i.e., it diverges for all
nonzero values of $\beta$ and summation techniques have to be applied
to give this series any meaning beyond a mere formal expansion.

A good mathematical model for the kind of divergence, which occurs in
the perturbation series of the quartic anharmonic oscillator, is the
so-called Euler integral
$$
E(z) \; = \; \int\nolimits^{\infty}_{0} \frac {\e^{-t} \, \d t}
{1 + z t } \; ,
\tag
$$
and its associated asymptotic series, the so-called Euler series
$$
E(z) \; \sim \; \sum _{n = 0}^{\infty} \; (-1)^{n} n! z^n \; = \;
{}_2 F_0 (1,1;-z) \, , \qquad z \to 0 \, .
\tag
$$

The radius of convergence of the Euler series is obviously zero.
Consequently, this series diverges quite wildly for all $z \neq 0$ and
appropriate summation techniques have to be applied if numerical values
for the Euler integral are to be computed with the help of this
asymptotic series. In fact, the Euler series (1.1-7) will be used quite
frequently in this report to test the ability of a sequence
transformation to sum wildly divergent series.

\medskip

\Abschnitt A short history of sequence transformations

\smallskip

\aktTag = 0

In this section, a short sketch of the historical development of
sequence transformations will be given. A more complete treatment of
the history would be beyond the scope of this report. Consequently, the
emphasis will be on those developments which laid the foundations for
the sequence transformations which are discussed in this report.

The idea of applying suitable transformations for the acceleration of
the convergence of a series or for the summation of a divergent series
is almost as old as analysis itself. According to Knopp (see p. 249 of
ref. [4]) the first series transformation was published by Stirling [5]
already in 1730, and in 1755 Euler [6] published the series
transformation which now bears his name.

These early ideas were later extended and refined as well as
supplemented by convergence proofs, and they finally led to the theory
of regular matrix transformations. Let $\Seq {s_n}$ be a sequence of
partial sums of a series according to eq. (1.1-1). Then, a new sequence
$\Seq {s^{\prime}_ n}$ with hopefully better convergence is obtained by
forming weighted means of the elements of the original sequence,
$$
s^{\prime}_ n \; = \; \sum _{k = 0}^{n} \; {\mu }_{n k} s_k \, .
\tag
$$

The main appeal of these matrix transformations lies in the fact that
for the weights ${\mu}_{n k}$ with $k,n \in \N_0$, which define such a
transformation, some necessary and sufficient conditions could be
formulated which ensure the regularity of the transformation. This
implies that such a regular matrix transformation can safely be applied
to any convergent sequence $\Seq {s_n}$ since the transformed sequence
$\Seq {s^{\prime}_n}$ will also converge to the same limit. A fairly
complete survey of the properties of such matrix transformations can be
found in books by Knopp [4], Hardy [7], Petersen [8], Peyerimhoff [9],
Zeller and Beekmann [10], and Powell and Shah [11].

This general applicability of regular matrix transformations to all
convergent sequences is undoubtedly quite advantageous from a
theoretical point of view. However, for the practical usefulness of a
transformation in actual computations this general applicability is
more likely a hindrance rather than an advantage. This may sound
paradoxical. But one cannot realistically expect that a given method
will be particularly efficient in a special case if it is
simultaneously required that this method should also be able to work in
all cases.

Consequently, in recent years emphasis has shifted towards the more
powerful but also more specialized nonlinear transformations.
Theoretically, nonlinear transformations are much more complicated than
matrix transformations and their properties are by no means completely
understood. In addition, nonlinear transformations are usually
nonregular, i.e., it is not guaranteed that the transformed sequence
$\Seq {s^{\prime}_n}$ will converge at all, let alone to the same limit
as the original sequence $\Seq {s_n}$. Hence, unless additional
information about the sequence to be transformed is available, the use
of a nonlinear sequence transformation may be risky. However, this
undeniable disadvantage is compensated by the empirical fact that if a
nonlinear transformation works, it frequently does so in a spectacular
fashion.

The probably oldest nonlinear sequence transformation is the famous
$\Delta^2$ process,
$$
s^{\prime}_n \; = \; s_n \, - \, \frac
{[ \Delta s_n ]^2} {\Delta^2 s_n} \; = \; s_n \, - \, \frac
{[s_{n+1} - s_n]^2} {s_{n+2} - 2 s_{n+1} + s_n}
\; , \qquad n \in \N_0 \, .
\tag
$$

This sequence transformation is named after Aitken [12] who published
this transformation in 1926 but there are indications that it is in
fact much older. For instance, Todd (see p. 5 of ref. [13]) claims that
this transformation was in principle already known to Kummer [14].

It is generally accepted that the current interest in nonlinear
transformations is due to two articles by Shanks [15] and Wynn [16],
respectively. Shanks rediscovered in 1955 a sequence transformation
which was originally derived in 1941 by Schmidt [17]. Wynn showed only
one year later how this sequence transformation, which was originally
defined as the ratio of two determinants, can be computed quite
efficiently by a nonlinear recursive scheme which is now commonly
called the $\epsilon $ algorithm. These two articles by Shanks [15] and
Wynn [16] had an enormous impact since they stimulated a large amount
of research not only in various branches of mathematics but also in
theoretical physics and in other sciences. This is amply demonstrated
by the long lists of references in books by Baker [18], Brezinski [19
-- 21], Baker and Graves-Morris [22], and Wimp [23], and also by a
recent review article by Brezinski [24].

This active research on nonlinear transformations contributed
significantly not only to the understanding of Pad\'e approximants or the
$\epsilon $ algorithm, but also led to the discovery of several other
sequence transformations. For instance, in 1956 Wynn [25] introduced
the so-called $\rho $ algorithm which is essentially an intelligent and
efficient way of computing and extrapolating even-order convergents of
an interpolating continued fraction.

In 1971 Brezinski [26] introduced his so-called $\theta $ algorithm
which may be interpreted to be some kind of improved and accelerated
$\epsilon $ algorithm. Brezinski's derivation of this powerful
algorithm was purely heuristic. It was emphasised by Brezinski [27]
that this heuristic approach is not restricted to Wynn's $\epsilon $
algorithm and can also be used in the case of other sequence
transformations. Some examples of new sequence transformations, which
were derived in that way, will be given later.

Another class of sequence transformations was introduced by Levin [28]
in 1973. According to Smith and Ford [29,30] who compared the
performances of several linear as well as nonlinear sequence
transformations, some variants of the Levin transformation are probably
the most powerful and most versatile convergence accelerators
currently known and they are also able to sum even wildly divergent
series. The sequence transformations introduced by Levin are also the
basis of a large part of this article since they are the starting point
for the derivation of several new sequence transformations which offer
in some cases computational advantages, in particular if wildly
divergent series are to be summed. In addition, a theoretical analysis
of the properties of the new transformations can also often be done
more easily than in the case of the Levin transformation.

A general extrapolation algorithm, which encompasses the majority of
the currently known extrapolation methods and also many of the new
sequence transformations of this report as special cases, was developed
independently by Brezinski [31] and H{\aa}vie [32].

Finally, Germain-Bonne [33] developed in 1973 a formal theory of
convergence acceleration which is of considerable importance not only
for this report. By means of Germain-Bonne's theory it can in some
cases be decided whether a given transformation is regular, i.e.,
whether the convergence of a sequence $\Seqn s$ to some limit s implies
the convergence of the transformed sequence $\Seq {s^{\prime}_n}$ to
the same limit $s$. Also, in some cases it can be decided by a priori
considerations whether the transformed sequence $\Seq {s^{\prime}_n}$
will converge faster than the original sequence $\Seqn s$.

\medskip

\Abschnitt Organization of this report

\smallskip

\aktTag = 0

No attempt is made to treat all aspects of the acceleration of
convergence and the summation of divergent series. The emphasis of this
report is on convergence acceleration and summation by means of
nonlinear sequence transformations. Linear sequence transformations are
only treated if they are special cases of nonlinear sequence
transformations. The nonlinear sequence transformations, which occur in
this report, are designed to handle convergent or divergent sequences
of partial sums of infinite series as they occur in scientific
applications or in the theory of special functions. However, the
specific problems, which arise in connection with the acceleration of
the convergence of Fourier series or of orthogonal expansions, are not
treated. Also, the acceleration or summation of multidimensional
sequences -- or vector sequences, as they are called in the literature
on convergence acceleration -- is not considered in this report.

Several nonlinear sequence transformations as for instance Aitken's
$\Delta^2$ process [12], Wynn's $\epsilon$ algorithm [16], Wynn's
$\rho$ algorithm [25], Brezinski's $\theta$ algorithm [26], and Levin's
sequence transformation [28] are now relatively well known and many
applications of these transformations have been reported in the more
recent literature. The properties of these nonlinear sequence
transformations are reviewed shortly in this report and efficient
algorithms for their computation are discussed. However, the emphasis
of this report is on the derivation of new nonlinear sequence
transformations, on the construction of efficient algorithms for their
computation, and on the analysis of their properties in convergence
acceleration and summation processes.

In this report, the sequence transformations are always computed with
the help of linear or nonlinear 2-dimensional recurrence formulas.
Also, it is always tried to find computational schemes for these
recursions which are optimal with respect to storage requirements. Such
an optimization is actually not necessary if the sequence
transformations are programmed in FORTRAN because then storage space
would not be a problem even if less efficient computational algorithms
would be used. If, however, sequence transformations are programmed in
a formal manipulation language such as REDUCE, MACSYMA or MAPLE, it is
probably a good idea to use such an optimized computational scheme
since storage restrictions may then be much more severe.

Some listings of FORTRAN 77 programs are included in the text. In order
to save space, all comments and also all IF statements, which check the
validity of the input data, were removed from the programs.
Consequently, these FORTRAN 77 programs are not ``good'' programs which
comply with the recommendations of books on programming style. The sole
purpose of these program listings is to facilitate the understanding of
the sometimes relatively intricate computational algorithms which are
described in this report.

In order to make this report more selfcontained, in section 2 the
mathematical terminology, which is specific for this report, as well as
the most important mathematical concepts and techniques, which are
needed for the derivation and understanding of sequence
transformations, are introduced.

In section 3 general properties of nonlinear sequence transformations
are discussed. In addition, it has been attempted to give a motivation
for some of the most important concepts and assumptions, which are the
basis for the construction of a large class of nonlinear sequence
transformations.

Section 4 deals with Wynn's $\epsilon$ algorithm [16], which is an
efficient algorithm for the computation of the Shanks transformation
[15] or -- if the elements of the sequence to be accelerated or summed
are the partial sums of a power series -- of Pad\'e approximants. Section
5 deals with Aitken's $\Delta^2$ process [12] and its iteration, which
are both close relatives of Wynn's $\epsilon$ algorithm. Section 6
deals with Wynn's $\rho$ algorithm [25], which is structurally almost
identical with Wynn's $\epsilon$ algorithm. Also, a new sequence
transformation is constructed by iterating the explicit expression for
$\rho_2^{(n)}$ along the lines of Aitken's iterated $\Delta^2$ process.

Section 7 deals with with Levin's sequence transformation [28], and
several other sequence transformations which are either special cases
or generalizations of Levin's sequence transformation. Levin's sequence
transformation is also the starting point for two new classes of
sequence transformations which are treated in sections 8 and 9. The
difference between Levin's sequence transformation and the new sequence
transformations is that Levin's sequence transformation is based upon
the assumption that the remainders of the partial sums can be
approximated by truncated Poincar\'e-type asymptotic expansions whereas
the new sequence transformations assume that the remainders can be
approximated by truncations of factorial series and related expansions
which are also based upon Pochhammer symbols.

Section 10 deals with Brezinski's $\theta$ algorithm [26], which was
derived by modifying the recursive scheme for Wynn's $\epsilon$
algorithm, and a closely related sequence transformation which is
obtained by iterating the expression for $\theta_2^{(n)}$. In section
11, the recursive schemes of several other linear and nonlinear
sequence transformations are modified along the lines of Brezinski's
$\theta$ algorithm and several new nonlinear sequence transformations
are derived.

The practical usefulness of the original version of Germain-Bonne's
formal theory of convergence acceleration [33] is quite limited since
it can only analyze the properties of a sequence transformation if its
recursive scheme satisfies some very restrictive conditions. Many
sequence transformations of this report do not satisfy these
conditions. Consequently, in section 12 Germain-Bonne's formal theory
of convergence acceleration is modified in such a way that the
properties of the sequence transformations of this report can also be
analyzed.

Unfortunately, Germain-Bonne's theory cannot be applied in all cases of
interest. In particular, it cannot be used for the analysis of the
summation of wildly divergent Stieltjes series as they for instance
occur in the Rayleigh-Schr\"odinger perturbation expansion for the energy
eigenvalues of the quartic anharmonic oscillator. In section 13 the
transformation of sequences of partial sums of convergent and divergent
Stieltjes series is analyzed. The estimates, which are obtained in this
way, indicate that some variants of the new sequence transformations,
which are discussed in sections 8 and 9, should sum a divergent
Stieltjes series as good or even somewhat better than the analogous
variants of Levin's sequence transformation, and that they should all
be far more efficient than Pad\'e approximants. These theoretical
estimates are supported quite convincingly by some numerical examples.

One of the most complicated computational problems, which can occur in
this context, is the acceleration of the convergence of infinite
series with terms $a_n$ that all have the same sign and that decay
like a fixed power $n^{- \alpha}$ with $\alpha > 1$ as $n \to \infty$.
A good example of such an extremely slowly convergent infinite series
with positive terms is the series (1.1-2) for the Riemann $\zeta$
function. In the case of such slowly convergent series with positive
terms, Germain-Bonne's formal theory of convergence acceleration [33]
also does not help. Section 14 deals with the acceleration of the
convergence of series of that type. Some exactness results and also
some error estimates are derived. However, it is relatively difficult
to obtain theoretical results. Consequently, the emphasis in section 14
is on numerical testing.

Finally, section 15 contains a condensed review of the properties of
the sequence transformations which are treated in this report.

\endAbschnittsebene

\neueSeite

\Abschnitt Terminology

\vskip - 2 \jot

\beginAbschnittsebene

\medskip

\Abschnitt Special mathematical symbols and special functions

\smallskip

In this report, essentially standard mathematical terminology will be
used. In particular, $\N$ stands for the set of positive integers $n =
1, 2, 3, \ldots$, and $\N_0$ stands for the set of nonnegative integers
$n = 0, 1, 2, \ldots $ . Also, $\R$ and $\C$ denote the sets of real
and complex numbers, respectively. The following notations are,
however, nonstandard:

\smallskip

\beginBeschreibung \zu \Laenge{$\Seqn s$ :}

\item{$\Seqn s$ \hfill:} Sequence of elements $s_n$ with $n \in \N_0$.
It is always tacitly assumed that the sequence elements with negative
indices, $s_{-1}, s_{-2}, s_{-3}, \ldots$, are zero.

\item{$\Ent x$ \hfill :} Integral part of $x \in \R$, i.e., the largest
integer $m$ satisfying the inequality $m \le x$.

\item{$\F^n$ \hfill :} Set of all vectors $(x_1, \ldots , x_n) \in
\R^n$ with all components being different from zero, i.e., $x_j \neq 0$
for all $j = 1, 2, \ldots, n$.

\item{$\D^n$ \hfill :} Set of all vectors $(x_1, \ldots , x_n) \in
\R^n$ with all components being distinct, i.e., $i \neq j$ implies
$x_i \neq x_j$ for all $i,j = 1,2,\ldots, n$.

\item{$\H^n$ \hfill :} Intersection of $\F^n$ and $\D^n$, i.e., the set
of vectors $(x_1, \ldots , x_n) \in \R^n$ with all components being
nonzero and distinct.

\endBeschreibung

\smallskip

Sometimes sums or products will occur in which the lower limit is
greater than the upper limit. In this report, we shall always use the
convention that such an {\it empty sum} will be interpreted as zero,
i.e.,
$$
\sum_{k=m}^n \, a_k \; = \; 0 \, ,
\qquad \text{if} \quad m > n \, ,
\tag
$$
and that such an {\it empty product} will be interpreted as one, i.e.,
$$
\prod_{k=m}^n \, a_k \; = \; 1 \, ,
\qquad \text{if} \quad m > n \, .
\tag
$$

For the commonly occurring special functions of mathematical physics
the notation and the conventions of Magnus, Oberhettinger, and Soni
[34] will be used in this report unless explicitly stated.

\medskip

\Abschnitt Order symbols

\smallskip

\aktTag = 0

Let $f(z)$ and $g(z)$ be two functions defined on some domain $D$ in
the complex plane and let $z_0$ be a limit point of $D$, possibly the
point at infinity. Then,
$$
f(z) \; = \; O(g(z)) \, , \qquad z \to z_0
\tag
$$
means that there is a positive constant $A$ and a neighbourhood $U$ of
$z_0$ such that
$$
\vert f(z) \vert \; \le \; A \vert g(z) \vert
\tag
$$
for all $z \in U \cap D$. If $g(z)$ does not vanish on $U \cap D$ this
simply means that $f(z)/g(z)$ is bounded on $U \cap D$. Also,
$$
f(z) \; = \; o(g(z)) \, , \qquad z \to z_0
\tag
$$
means that for any positive number $\epsilon \in \R$ there exists a
neighbourhood $U$ of $z_0$ such that
$$
\vert f(z) \vert \; \le \; \epsilon \vert g(z) \vert
\tag
$$
for all $z \in U \cap D$. If $g(z)$ does not vanish on $U \cap D$ this
simply means that $f(z)/g(z)$ approaches zero as $z \to z_0$.

\medskip

\Abschnitt Asymptotic sequences and asymptotic expansions

\smallskip

\aktTag = 0

A finite or infinite sequence of functions $\Seq {\Phi_n (z)}$ with $n
\in \N_0$, which are defined on some domain D of complex numbers on
which all $\Phi_n (z)$ are nonzero except possibly at $z_0$, is called
an {\it asymptotic sequence} as $z \to z_0$ if, for all $n \in \N_0$,
$$
\Phi_{n+1} (z) \; = \; o(\Phi_n (z)) \, , \qquad z \to z_0 \, .
\tag
$$

Examples for asymptotic sequences with $n \in \N_0$ are $\Seq
{(z-z_0)^n}$ as $z \to z_0$ or $\Seq {(\log z )^{-n}}$ as $z \to
\infty$.

The formal series
$$
f(z) \sim \sum_{n=0}^{\infty} \; c_n \Phi_n (z) \, ,
\tag
$$
which need not be convergent, is called an {\it asymptotic expansion}
of $f(z)$ with respect to the asymptotic sequence $\Seq {\Phi_n (z)}$
in the {\it sense of Poincar\'e} if, for every $m \in \N_0$,
$$
f(z) - \sum_{n=0}^m \;c_n \Phi_n (z) \; = \; o( \Phi_m (z)) \, ,
\qquad z \to z_0 \, .
\tag
$$

If such a Poincar\'e-type asymptotic expansion exists, it is unique, and
its coefficients $c_n$ can be computed recursively,
$$
c_m \; = \;
\lim_{z \to z_0} \; \Bigl\{ \bigl[ f(z) - \sum_{n=0}^{m-1} \; c_n
\Phi_n (z) \bigr] \; / \; \Phi_m (z) \Bigr\}
\, , \qquad m \in \N_0 \, .
\tag
$$

The first term of the asymptotic expansion (2.3-2) is usually called
the {\it dominant} or {\it leading term} and one frequently writes
$$
f(z) \; \sim \; \Phi_0 (z) \, ,
\tag
$$
indicating that $f(z)/\Phi_0 (z)$ tends to $c_0$ as $z \to z_0$.

A particularly simple asymptotic sequence as $z \to \infty$ is the set
$\Seq {\Psi (z) / z^n}$, $n \in \N_0$, with $\Psi (z)$ being a suitable
function. If a given function $f(z)$ possesses an asymptotic expansion
with respect to this sequence,
$$
f(z) \; \sim \; \Psi (z) \sum_{n=0}^{\infty} \; c_n / z^n \, ,
\qquad z \to \infty \, ,
\tag
$$
then the ratio $f (z)/\Psi (z)$ can be expressed as an {\it asymptotic
power series} in $1/z$,
$$
f(z)/ \Psi (z) \; \sim \; \sum_{n=0}^{\infty} \; c_n / z^n \, ,
\qquad z \to \infty \, .
\tag
$$

\medskip

\Abschnitt Finite differences

\smallskip

\aktTag = 0

Let $f$ be a function defined on the set of integers $\N_0$. Then, the
{\it forward difference} $\Delta f(n)$ is defined by the relationship
$$
\Delta f(n) \; = \; f(n+1) - f(n), \qquad n \in \N_0 \, .
\tag
$$

Higher powers of the difference operator $\Delta$ can be defined
recursively, i.e.,
$$
\beginAligntags
" \Delta^k f(n) \; " = \; " \Delta [ \Delta^{k-1} f(n) ] \, ,
 \qquad k \in \N, \\
\tag
" \Delta^0 f(n) \; " = \; " f(n) \, . \\
\tag
\endAligntags
$$

The {\it shift operator} $E$ is defined by the relationship
$$
E f(n) \; = \; f(n+1) \, .
\tag
$$

Higher powers of $E$ can again be defined recursively. Obviously, we
have
$$
\beginAligntags
" E^k f(n) \; " = \; " f(n+k), \qquad k \in \N \, , \\
\tag
" E^0 f(n) \; " = \; " f(n) \, . \\
\tag
\endAligntags
$$

It follows at once from their definitions that the operators $\Delta$
and $E$ are connected by the relationship
$$
\Delta \; = \; E - 1 \, .
\tag
$$

This relationship can be combined with the binomial theorem to give
$$
\Delta^k f(n) \; = \; (-1)^k \, \sum_{j=0}^k (-1)^j \binom k j
f(n+j) \, , \qquad k \in \N_0 \, .
\tag
$$

In the following text it will always be tacitly assumed that in the
case of several indices the difference operator $\Delta$ and the shift
operator $E$ will only act upon $n$ and not on other indices.

\medskip

\Abschnitt Special sequences

\smallskip

\aktTag = 0

Let us assume that the sequence $\Seqn s$ either converges to some
limit $s$, or, if it diverges, can be summed by an appropriate
summation method to give $s$. In the case of divergence $s$ is
frequently called {\it antilimit}. Then, the partition of a sequence
element $s_n$ into the {\it limit} or {\it antilimit} $s$ and the {\it
remainder} $r_n$ according to
$$
s_n \; = \; s \; + \; r_n
\tag
$$
makes sense for all $n \in \N_0$. If $s_n$ is the partial sum of a
series,
$$
s_n \; = \; \sum_{k=0}^n \, a_k \, ,
\tag
$$
the remainder $r_n$ obviously satisfies
$$
r_n \; = \; - \sum_{k=n+1}^{\infty} \; a_k \, .
\tag
$$

For the ratio of two consecutive terms of an infinite series we write
$$
\rho_n \; = \; a_{n+1} / a_n, \qquad n \in \N_0 \, .
\tag
$$

The magnitude of the remainder $r_n$ is a natural measure for the
convergence of a sequence or series. Often, it is also of considerable
interest to analyze the asymptotics of the sequence of remainders
$\Seqn r$ as $n \to \infty $. Let $\Seq {\phi_k (n)}$, $k \in \N_0$, be
a suitable asymptotic sequence as $n \to \infty$ with $\phi_0 (n) = 1$.
Then, $\omega_n$ denotes the dominant part of $r_n$ with respect to the
asymptotic sequence $\phi_k (n)$, i.e.,
$$
r_n / \omega_n \; \sim \; \sum_{k=0}^{\infty} \; c_k \phi_k (n)
\, , \qquad n \to \infty \, .
\tag
$$

Sequences of {\it remainder estimates} $\Seq {\omega_n}$ will be of
considerable importance in this report. The reason is that it is often
possible to obtain at least some structural information about the
behaviour of the dominant term of a remainder $r_n$ as $n \to \infty $.
It will become clear later that those convergence acceleration or
summation methods, which explicitly utilize the information contained
in the remainder estimates $\Seq {\omega_n}$, are frequently
particularly efficient.

Many sequence transformations do not only require the input of the
sequences $\Seqn s$ and $\Seq {\omega_n}$ but also the input of an
additional sequence of auxiliary quantities as for instance
interpolation points. In this report, such an {\it auxiliary sequence}
will usually be denoted by $\Seqn x$.

\medskip

\Abschnitt Types of convergence

\smallskip

\aktTag = 0

It is neither possible nor desirable to set up a complete
classification scheme which is able to cover all types of convergence.
However, in the majority of all practical applications only a few types
of convergence occur. Consequently, special names were given to them in
the literature.

Let us assume that the sequence $\Seqn s$, which converges to some limit
$s$, satisfies
$$
\lim_{n \to \infty} \; \frac {s_{n+1} - s} {s_n - s} \; = \; \lim_{n
\to \infty} \; \frac {r_{n+1}}{r_n} \; =\; \rho \, .
\tag
$$

If $ 0 < |\rho | < 1$ holds, we say that the sequence $\Seqn s$
converges {\it linearly}, if $\rho = 1$ holds, we say that $\Seqn s$
converges {\it logarithmically}, and if $\rho = 0$ holds, we say that
$\Seqn s$ converges {\it hyperlinearly}. Of course, $|\rho | > 1$
implies that the sequence $\Seqn s$ diverges.

The standard example for linear convergence is the sequence of partial
sums of the geometric series,
$$
s_n (z) \; = \; \sum_{k=0}^n \; z^k \; = \;
\frac {1-z^{n+1}} {1-z} \, ,
\qquad 0 < |z| < 1 \, , \quad n \in \N_0 \, .
\tag
$$

The sequence of partial sums of the series (1.1-2) for the Riemann zeta
function is a good example for logarithmic convergence. Also, it can be
shown quite easily that the partial sums of the power series for the
exponential function form a sequence which converges hyperlinearly.

The above definitions for hyperlinear, linear and logarithmic
convergence do not seem to be particularly well suited for the
classification of infinite series because normally only the terms $a_k$
of a series but not the remainders $r_n$ are known. However, Wimp
showed on p. 6 of his book [23] that if $0 < \vert \rho \vert < 1$
holds, the two statements
$$
\lim_{n \to \infty} \; (r_{n+1} / r_n) \; = \; \rho
\tag
$$
and
$$
\lim_{n \to \infty} \; (a_{n+1} / a_n) \; = \; \rho
\tag
$$
are equivalent. In addition, Clark, Gray, and Adams [35] showed that if
the terms $a_k$ of a convergent series are all real and have the same
sign, then
$$
\lim_{n \to \infty} \; (r_{n+1} / r_n) \; = \;
\lim_{n \to \infty} \; (a_{n+1} / a_n) \; = \; \rho \, .
\tag
$$

In eq. (2.6-5), the case $\rho = 1$, which corresponds to logarithmic
convergence, is not excluded. Hence, it is at least possible to
classify linearly and logarithmically convergent series -- which are of
particular interest in connection with convergence acceleration methods
-- according to the behaviour of their terms $a_n$ as $n \to \infty$.

Sequences and series, which converge hyperlinearly, often converge so
rapidly that not much can be gained by convergence acceleration
methods. Consequently, hyperlinear convergence is more or less
neglected in the literature on convergence acceleration. This is not
entirely justified because in some situations the use of convergence
acceleration methods can indeed be quite helpful. A simple and
nevertheless striking example, which shows that convergence
acceleration methods may be quite useful even in the case of
hyperlinear convergence, is the power series for the exponential
function. It cannot be used for the computation of $\e^{-x}$ if $x$ is
positive and large because then large terms with alternating signs
would lead to cancellation and to severe numerical instabilities.
However, if suitable convergence acceleration methods are used,
remarkably accurate results can be obtained after a relatively small
number of terms.

Let us assume that two sequences $\Seqn s$ and $\Seq {s^{\prime}_n}$
both converge to the same limit $s$. We shall say that the sequence
$\Seq {s^{\prime}_n}$ {\it converges more rapidly} than $\Seqn s$ if
$$
\lim_{n \to \infty} \; \frac { s^{\prime}_n - s }{ s_n - s}
\; = \; 0 \, .
\tag
$$

In convergence acceleration processes, this definition is somewhat
inconvenient since it requires the knowledge of the limit $s$ which is
usually not known. Consequently, it would be desirable to replace eq.
(2.6-6) by the alternative condition
$$
\lim_{n \to \infty} \; \frac { s^{\prime}_{n+1} - s^{\prime}_n }
{ s_{n+1} - s_n} \; = \;
\lim_{n \to \infty} \; \frac { \Delta s^{\prime}_n }
{ \Delta s_n} \; = \; 0 \, .
\tag
$$

However, it seems that it is not possible to prove the equivalence of
eqs. (2.6-6) and (2.6-7) without making explicit assumptions about how
fast the sequences $\Seqn s$ and $\Seq {s^{\prime}_n}$ approach their
common limit $s$.

If the sequence $\Seqn s$ converges linearly, the transformed sequence
$\Seq {s^{\prime}_n}$ can only converge more rapidly than $\Seqn s$ if
it converges at least linearly or even faster. In this case, the
equivalence of the two conditions (2.6-6) and (2.6-7) follows at once
from the relationship
$$
\frac {s^{\prime}_{n+1} - s^{\prime}_n} {s_{n+1} - s_n} \; = \;
\frac { s^{\prime}_n - s } { s_n - s} \;
\frac
{[(s^{\prime}_{n+1} - s) \, / \, (s^{\prime}_n - s)] \, - \, 1}
{[(s_{n+1} - s) \, / \, (s_n - s)] \, - \, 1} \, .
\tag
$$

However, if $\Seqn s$ converges logarithmically, the denominator of the
second term on the right-hand side of eq. (2.6-8) approaches zero as $n
\to \infty$. In this case, some additional assumptions about the rate
of convergence of $\Seq {s^{\prime}_n}$ to $s$ have to be made in order
to be able to show that the two conditions (2.6-6) and (2.6-7) are
indeed equivalent.

\medskip

\Abschnitt Sequence transformations

\smallskip

\aktTag = 0

In this report a sequence transformation ${\cal T}$ will always be a
rule which transforms a given sequence $\Seqn s$ into a new sequence
$\Seq {s^{\prime}_n}$,
$$
{\cal T} : {\Seqn s} \; \mapsto \; {\Seq {s^{\prime}_n}} \, ,
\qquad n \in \N_0 \, .
\tag
$$

Since a computational algorithm can only involve a finite number of
operations, only finite subsets of a sequence $\Seqn s$ can be
associated to a new sequence element $s^{\prime}_m$ by a sequence
transformation ${\cal T}$. In addition, it will always be assumed that
the finite subset, which is to be transformed, will entirely consist of
consecutive elements. This means that only subsets of the type $\{ s_n,
s_{n+1}, \ldots , s_{n+\ell} \}$ with $n,\ell \in \N_0$ will be
considered. Since the subsets, which are to be transformed, contain
$\ell + 1$ elements, $\ell$ will frequently be called the {\it order}
of the transformation ${\cal T}$. Hence, if all elements of the
sequence $\Seqn s$ are real and if no sequences of interpolation points
$\Seqn x$ or remainder estimates $\Seq {\omega_n}$ are needed, a
sequence transformation ${\cal T}$ of order $\ell$ is a map of the
following type:
$$
{\cal T} : {\R}^{\ell + 1} \; \to \; \R \, .
\tag
$$

In this report, a sequence transformation ${\cal T}$ can always be
represented by an infinite set of doubly indexed quantities
$T_k^{(n)}$ with $k,n \in \N_0$. The superscript $n$ always indicates
the minimal index occurring in the finite subset of sequence elements
which are used for the computation of the transform $T_k^{(n)}$, and
the subscript $k$ is a measure for the complexity of such a $T_k^{(n)}$.

The quantities $T_k^{(n)}$ are gauged in such a way that $T_0^{(n)}$
will always correspond to an untransformed sequence element, i.e.,
$$
T_0^{(n)} \; = \; s_n \, , \qquad n \in \N_0 \, .
\tag
$$

Increasing values of $k$ imply that the order $\ell$ of the transform
$T_k^{(n)}$ also increases. This means that for every $k,n \in \N_0$
the sequence transformation ${\cal T}$ will produce a new transform
for which we shall write
$$
T_k^{(n)} \; = \; {\cal T}(s_n, s_{n+1}, \ldots , s_{n+\ell}) \, .
\tag
$$

Here, the order $\ell$ is of course a function of $k$. The exact
relationship, which connects the subscript $k$ and the order $\ell$, is
specific for every sequence transformation ${\cal T}$. In this report
we shall encounter a variety of different relationships such as $\ell =
k$, $\ell = k + 1$, $\ell = 2 k$ or even $\ell = 3 k$.

The transforms $T_k^{(n)}$ with $k,n \in \N_0$ can be displayed in a
2-dimensional array which is called the {\it table} of the sequence
transformation {\cal T}. In this report, the transforms $T_k^{(n)}$
will always be ordered in a rectangular scheme in such a way that the
superscript $n$ indicates the row and the subscript $k$ the column of
the array. Hence, in this report the table of a transformation ${\cal
T}$ will always be displayed in the following way:

$$
\matrix{
T_0^{(0)} " T_1^{(0)} " T_2^{(0)} " \ldots " T_n^{(0)} " \ldots \\
T_0^{(1)} " T_1^{(1)} " T_2^{(1)} " \ldots " T_n^{(1)} " \ldots \\
T_0^{(2)} " T_1^{(2)} " T_2^{(2)} " \ldots " T_n^{(2)} " \ldots \\
T_0^{(3)} " T_1^{(3)} " T_2^{(3)} " \ldots " T_n^{(3)} " \ldots \\
\vdots    " \vdots    " \vdots    " \ddots " \vdots    " \ddots \\
T_0^{(n)} " T_1^{(n)} " T_2^{(n)} " \ldots " T_n^{(n)} " \ldots \\
\vdots  " \vdots " \vdots " \ddots " \vdots  " \ddots } \tag
$$

In the process of convergence acceleration or summation only those
elements of the table of a sequence transformation ${\cal T}$ should be
computed which will be needed to obtain convergence up to a certain
accuracy. Of course, this implies that one has to decide in advance
which transforms $T_k^{(n)}$ should be used for that purpose. In this
context it will be advantageous to introduce the following terminology.

A sequence $\Seq {(n_j,k_j)}$ of ordered pairs of integers $n_j,k_j \in
\N_0$ is called a {\it path} if $n_0 = k_0 = 0$ and if for all integers
$j\in \N_0$ we have $n_{j+1} \ge n_j$ and $k_{j+1} \ge k_j$ and if for
each $j\in \N_0$ either one or both of the two relations $n_{j+1} = n_j
+ 1$ and $k_{j+1} = k_j + 1$ are true. Obviously, $n_j + k_j \to
\infty$ as $j \to \infty$.

Paths where $k_j$ is ultimately constant are called {\it vertical
paths}, and paths where $n_j$ is ultimately constant are called {\it
horizontal paths}.

We shall say that a sequence transformation ${\cal T}$ is {\it regular}
on a given path ${\cal P} = \Seq {(n_j,k_j)}$ if for every convergent
sequence $\Seqn s$ we have
$$
\lim_{j \to \infty} \; T_{k_j}^{(n_j)} \; = \; s \, .
\tag
$$

Next, we want to define what we mean by saying that a transformation
${\cal T}$ is called {\it accelerative on a path} ${\cal P}$ for a
sequence $\Seqn s$. In the literature on convergence acceleration the
following definition is the most common one:
$$
\lim_{j \to \infty} \; \frac
{T_{k_j}^{(n_j)} - s} { s_{n_j} - s } \; = \; 0 \, .
\tag
$$

However, since for a given subscript $k$ a sequence transformation
$\cal T$ always acts on $\ell + 1$ consecutive sequence elements $s_n,
s_{n+1}, \ldots , s_{n+\ell}$ with $\ell$ being a function of $k$, it
would actually be better to say that ${\cal T}$ is accelerative on a
path ${\cal P}$ for a sequence $\Seqn s$ if
$$
\lim_{j \to \infty} \; \frac
{T_{k_j}^{(n_j)} - s} { s_{n_j + \ell_j} - s } \; = \; 0 \, .
\tag
$$

Hence, if the second definition is used a transformation ${\cal T}$
will be called accelerative on a path ${\cal P} = \Seq {(n_j, k_j)}$ if
the transforms $T_{k_j}^{(n_j)}$ converge faster than the last elements
$s_{n_j + \ell_j}$ of the strings $s_{n_j}, s_{n_j + 1}, \ldots ,
s_{n_j + \ell_j}$ which are used for the computation of the transforms
$T_{k_j}^{(n_j)}$ .

In this report both definitions (2.7-7) and (2.7-8) will be used.
However, it will always be stated explicitly which of the two different
definitions is actually meant.

Let us again assume that a sequence $\Seqn s$ converges to some limit
$s$. A sequence transformation ${\cal T}$ will be called {\it exact}
for the sequence $\Seqn s$ if some integer $\ell_0 \in \N_0$ exists
such that the application of ${\cal T}$ to every {\it finite} string
$s_n, \ldots, s_{n+\ell}$ of sequence elements with $\ell \ge \ell_0$
yields the exact limit $s$ of this sequence.

\endAbschnittsebene

\neueSeite

\Abschnitt On the derivation of sequence transformations

\vskip - 2 \jot

\beginAbschnittsebene

\medskip

\Abschnitt General properties of nonlinear sequence transformations

\smallskip

\aktTag = 0

It was remarked earlier that in many cases of physical interest the
conventional approach of evaluating an infinite series by adding one
term after the other does not suffice. Examples are logarithmically
convergent series which may converge so slowly that an evaluation by
adding one term after the other would overstep even the potential of
modern supercomputers, or divergent series as they for instance occur
in Rayleigh-Schr\"odinger perturbation theory.

In such cases it is necessary to replace the conventional process of
evaluating a series by a {\it generalized summation process} which is
able to associate a numerical value even to a prohibitively slowly
convergent or divergent series.

The generalized summation processes of this report are transformations
which are defined on finite subsets of the sequence $\Seqn s$ of
partial sums. Let ${\cal T}_{\ell}$ be such a generalized summation
process which acts upon $\ell + 1$ partial sums $s_n, \ldots,
s_{n+\ell}$. In view of the fact that in the case of convergent series
with real terms we have
$$
\sum_{n=0}^{\infty} \; ( \alpha a_n + \beta b_n ) \; = \;
\alpha \sum_{n=0}^{\infty} \; a_n \, + \,
\beta \sum_{n=0}^{\infty} \; b_n \, ,
\qquad \alpha , \beta \in \R \, ,
\tag
$$
it seems natural to require that such a generalized summation process
${\cal T}_{\ell}$ should also be {\it linear}. Therefore, let us assume
that $\Seqn s$ and $\Seqn t$ are two sequences of partial sums of real
terms which converge to $s$ and $t$, respectively. Thus, a generalized
summation process ${\cal T}_{\ell}$ should satisfy
$$
\beginAligntags
" {\cal T}_{\ell} (\alpha s_n + \beta t_n, \ldots ,
\alpha s_{n+\ell} + \beta t_{n+\ell}) \; = \;
\alpha {\cal T}_{\ell} (s_n, \ldots , s_{n+\ell}) \; + \;
\beta {\cal T}_{\ell} (t_n, \ldots , t_{n+\ell}) \, , \\
" \alpha,\beta \in \R \, , \qquad \ell , n \in \N_0 \, .
\\ \tag
\endAligntags
$$

Also, such a generalized summation process should preserve the limit of
a convergent sequence, i.e., it should be {\it regular}. Hence, a {\it
linear and regular generalized summation process} ${\cal T}_{\ell}$
should satisfy for all sequences $\Seqn s$ and $\Seqn t$ which converge
to $s$ and $t$, respectively,
$$
\lim_{n \to \infty} {\cal T}_{\ell}
(\alpha s_n + \beta t_n, \ldots ,
\alpha s_{n+\ell} + \beta t_{n+\ell}) \; = \;
\alpha s \; + \; \beta t \, .
\tag
$$

Unfortunately, the generalized summation processes considered in this
report will in general be neither linear nor regular and we have to
content ourselves with a weaker requirement. Let $\Seqn s$ be a
sequence and let $\alpha$ and $\tau$ be two constants. We may only
assume that a generalized summation process ${\cal T}_{\ell}$ is {\it
invariant under translation}, i.e., that for all admissible $\ell, n
\in \N_0$
$$
{\cal T}_{\ell}
(\alpha s_n + \tau, \ldots , \alpha s_{n+\ell} + \tau) \; = \;
\alpha {\cal T}_{\ell} (s_n, \ldots , s_{n+\ell})
\; + \; \tau \, .
\tag
$$

It must be emphasized that because of the nonregularity of the
generalized summation processes ${\cal T}_{\ell}$ of this report we
also cannot assume that either
$$
\lim_{n \to \infty} {\cal T}_{\ell} (s_n, \ldots , s_{n+\ell})
\; = \; s
\tag
$$
or -- if this limit is defined -- that
$$
\lim_{\ell \to \infty} {\cal T}_{\ell} (s_n, \ldots , s_{n+\ell})
\; = \; s
\tag
$$
holds for arbitrary convergent sequences $\Seqn s$.

Nonlinearity and nonregularity are undeniably unpleasant complications
which one would like to avoid instinctively. However, they are
essential and indispensable since the power of the sequence
transformations, which are discussed in this report, stems from their
nonlinearity and nonregularity.

\medskip

\Abschnitt An example: Convergence acceleration of alternating series

\smallskip

\aktTag = 0

If we want to construct a generalized summation process, which is able
to accelerate the convergence of an infinite series, we are confronted
with the practical problem that the information contained in a finite
string of partial sums $s_0, s_1, \ldots, s_m$ has to be extracted and
utilized in a way which is more efficient than the conventional
process of adding up one term after the other. If we again assume that
for all $n \in \N_0$ a sequence element $s_n$ can be partitioned into
the limit $s$ and the remainder $r_n$ according to
$$
s_n \; = \; s + r_n \, ,
\tag
$$
then this essentially means that we have to find a way of eliminating
the remainder $r_n$ and determining the limit $s$ at least
approximately by exploiting the information stored in the finite string
$s_0, s_1, \ldots, s_m$ of partial sums.

Essentially the same problem of eliminating the remainder $r_n$ and
determining the antilimit $s$ at least approximately arises if we try
to sum a divergent series. The only difference is that in the case of a
divergent series the remainder $r_n$ does not vanish as $n \to \infty$,
and that the antilimit $s$ cannot be obtained by simply adding up the
terms of the series. Instead, the antilimit $s$ of a sequence can only
be determined with the help of a suitable summation method.

Since we cannot assume that it will be possible to obtain the numerical
values of the remainders $r_n$ directly, the best we can hope for is
something which may be called {\it structural information}. In order to
clarify this concept we will consider a simple example. Let us assume
that the sequence elements $s_n$ are partial sums of a series with real
and strictly alternating terms,
$$
s_n \; = \; \sum_{k=0}^n \; (-1)^k b_k \, ,
\tag
$$
which means that all $b_n$ with $n \in \N_0$ have the same sign. The
remainder $r_n$ of $s_n$ is then given by
$$
r_n \; = \; - \sum_{k=n+1}^{\infty} \; (-1)^k b_k \, .
\tag
$$

Let us now also assume that all $b_n$ with $n \in \N_0$ are positive
and strictly decreasing with increasing $n$ and that they vanish as $n
\to \infty$. This implies that the series converges to some limit $s$.
In addition, it can be shown that the sequence of remainders $\Seqn r$
is also strictly alternating and that the magnitude of a remainder $r_n$
is bounded by the first term which was not included in the partial sum
$s_n$ (see p. 259 of ref. [4]),
$$
\vert r_n \vert < b_{n+1} \, , \qquad n \in \N_0.
\tag
$$

Now, we have to find a way of utilizing this structural information
about the behaviour of the sequence of remainders $\Seqn r$. Simply
adding the next term $(-1)^{n+1} b_{n+1}$ to $s_n$ would only produce
$s_{n+1}$ and we would not gain anything substantial. Consequently, we
need an additional assumption which will help us to construct a
sequence transformation for alternating series.

It is a relatively natural idea to assume that the ratio $r_n /
[(-1)^{n+1} b_{n+1}]$ can be expanded in a Poincar\'e-type asymptotic
power series in the variable $1/(n+1)$, i.e.,
$$
r_n \; \sim \; (-1)^{n+1} b_{n+1} \,
\sum_{j=0}^{\infty} \; c_j \, (n+1)^{-j} \, ,
\qquad n \to \infty \, .
\tag
$$

The assumption, that such a Poincar\'e-type asymptotic expansion exists,
will enable us to derive a sequence transformation which is capable of
accelerating the convergence of alternating series.

However, the complete elimination of such a remainder $r_n$ on the
basis of the asymptotic expansion (3.2-5) will not be possible since
any computational algorithm can only determine a finite number of the
unknown linear coefficients $c_j$ in eq. (3.2-5). Consequently, we can
only construct a sequence transformation which is able to eliminate
model remainders of the following type:
$$
{\tilde r_n} \; = \; (-1)^{n+1} b_{n+1} \>
\sum_{j=0}^{k-1} \, c_j \, (n+1)^{-j}
\, , \qquad n \in \N_0 \, .
\tag
$$

Model remainders of this type are obtained by truncating the infinite
series in eq. (3.2-5) after the first $k$ terms. Hence, at least for
sufficiently large values of $k$ and $n$ the model remainders ${\tilde
r_n}$ should approximate the actual remainders $r_n$ very well. This
implies that the partial sums $s_n$ can also be approximated very well
by the elements of the following model sequence, which contain only
finitely many terms:
$$
{\tilde s}_n \; = \; s \> + \> (-1)^{n+1} b_{n+1} \>
\sum_{j=0}^{k-1} \, c_j \, (n+1)^{-j}
\, , \qquad n \in \N_0 \, .
\tag
$$

In eq. (3.2-7) there occur $k+1$ unknowns, the limit $s$ and the $k$
coefficients $c_0, \ldots , c_{k-1}$. Since all unknowns occur {\it
linearly} their determination poses in principle no problems. All that
is needed are the numerical values of $k+1$ sequence elements, e.g.,
the string ${\tilde s}_n, {\tilde s}_{n+1}, \ldots , {\tilde s}_{n+k}$,
and it is possible to determine the limit $s$ of the model sequence
(3.2-7).

Hence, we only have to use Cramer's rule in order to see that the limit
$s$ is given by the following ratio of determinants:

\smallskip

$$
s \; = \; \frac
{
\vmatrix
{
{\tilde s}_n " \ldots " {\tilde s}_{n+k} \\ [1\jot]
(-1)^{n+1} b_{n+1} " \ldots " (-1)^{n+k+1} b_{n+k+1} \\
\vdots " \ddots " \vdots \\ [1\jot]
{\displaystyle \frac {(-1)^{n+1} b_{n+1} } {(n+1)^{k-1}} } " \ldots "
{\displaystyle \frac {(-1)^{n+k+1} b_{n+k+1}} {(n+k+1)^{k-1}}}
}
}
{
\vmatrix
{
1 " \ldots " 1       \\ [1\jot]
(-1)^{n+1} b_{n+1} " \ldots " (-1)^{n+k+1} b_{n+k+1} \\
\vdots " \ddots " \vdots \\ [1\jot]
{\displaystyle \frac {(-1)^{n+1} b_{n+1} } {(n+1)^{k-1}} } " \ldots "
{\displaystyle \frac {(-1)^{n+k+1} b_{n+k+1}} {(n+k+1)^{k-1}} }
}
} \, .
\tag
$$

\smallskip

Now we could try to replace the elements ${\tilde s}_n$ of the model
sequence (3.2-7) in the first determinant in eq. (3.2-8) by the
partial sums $s_n$. This would certainly not produce the exact limit
$s$ of the alternating series since the partial sums $s_n$ satisfy eq.
(3.2-7) only approximately. However, if the sequence elements ${\tilde
s}_n$ are able to approximate the partial sums $s_n$ with sufficient
accuracy then we can hope that the ratio of determinants, in which now
the partial sums $s_n, s_{n+1}, \ldots , s_{n+k}$ occur, will be a
better approximation to the limit $s$ than the last partial sum
$s_{n+k}$ which occurs in the ratio of determinants.

We shall see later that this is indeed the case. Actually, with the
help of our simple arguments we found the determinantal representation
of the sequence transformation $d_k^{(n)} (\beta, s_n)$ with $\beta =
1$. This transformation, which will be defined later in eq. (7.3-9),
is a special case of a very powerful sequence transformation which was
introduced by Levin [28]. Levin's general sequence transformation and
its numerous variants are discussed quite extensively in section 7 of
this report. In section 13 it is also shown that $d_k^{(n)} (\beta,
s_n)$, eq. (7.3-9), is able to sum even wildly divergent series and to
accelerate the convergence of linearly convergent series.

Our approach did not lead to a representation of this sequence
transformation $d_k^{(n)} (\beta, s_n)$, eq. (7.3-9), which is
completely satisfactory from a computational point of view.
Determinantal representations of sequence transformations are
computationally quite unattractive since the reliable and economical
evaluation of determinants is a more or less unsolved problem of
numerical mathematics. Consequently, it is important to find other
methods for the computation of a sequence transformation. In section 7,
it will be shown how the determinantal representation (3.2-8) can be
replaced by other representations which are better suited for numerical
work. However, the concepts and principles, which will be used for the
derivation of a large part of the sequence transformations of this
report, should now be clear. They can be summarized as follows:

\medskip

\noindent (1): Consider a model sequence with elements ${\tilde s}_n =
s + {\tilde r}_n$ and assume that their remainders ${\tilde r}_n$ can
be partitioned into a remainder estimate $\omega_n$ multiplied by some
other quantity $z_n$. This implies that the elements ${\tilde s}_n$ of
the model sequence satisfy:
$$
{\tilde s}_n \; = \; s + \omega_n \, z_n , \qquad n \in \N_0 \,.
\tag
$$
(2): Assume that an operator ${\hat T}$ , which is defined on finite
subsets of sequences and which is linear, annihilates the quantities
$z_n$ defined in eq. (3.2-9), i.e., ${\hat T} (z_n) = 0$. If we rewrite
eq. (3.2-9) in the following way
$$
( {\tilde s}_n \, - \, s ) / \omega_n \; = \, z_n \, ,
\tag
$$
we see that a sequence transformation ${\cal T} ({\tilde s}_n,
\omega_n)$, which is exact for the model sequence eq. (3.2-9), i.e.,
which satisfies ${\cal T} ({\tilde s}_n, \omega_n) = s$, is given by
the following ratio:
$$
{\cal T} ({\tilde s}_n, \omega_n) \; = \; \frac
{ {\hat T} ( {\tilde s}_n / \omega_n ) }
{ {\hat T} ( 1 / \omega_n ) } \, .
\tag
$$
(3): Replace the elements ${\tilde s}_n$ of the model sequence (3.2-9)
in the expression defining the sequence transformation ${\cal T}
({\tilde s}_n, \omega_n)$ -- in this report either the ratio of two
determinants, an explicit expression, or a recursive scheme --  by the
elements of the sequence $\Seqn s$ which is to be transformed.

The crucial step in this approach is the choice of an appropriate
sequence of model remainders $\Seq {{\tilde r}_n}$ since the ${\tilde
r}_n$ should have a mathematical structure which permits the
construction of a manageable annihilation operator ${\hat T}$. In
addition, the ${\tilde r}_n$  should also be capable of producing
good approximations for remainders $r_n$ which occur in actual
problems, because only then we may hope that the sequence of transforms
will converge more rapidly than the original sequence $\Seqn s$. These
aims are usually accomplished by partitioning ${\tilde r}_n$ into a
remainder estimate $\omega_n$ multiplied by a finite sum,
$$
{\tilde r}_n \; = \; \omega_n \> \sum_{j=0}^m \, c_j \, \phi_j (n) \, ,
\qquad m,n \in \N_0 \, .
\tag
$$

Here, the $\Seq {\phi_j (n)}$ with $j,n \in \N_0$ are a suitable set of
functions -- usually an asymptotic sequence as for instance
$(n+1)^{-j}$ with $j,n \in\N_0$ -- for which a sufficiently simple
annihilation operator can be found and the $c_j$ are so far completely
unspecified coefficients which are responsible for the flexibility of
this ansatz.

Once the asymptotic sequence $\Seq {\phi_j (n)}$ is chosen, the crucial
problem is the determination of a suitable sequence of remainder
estimates $\Seq {\omega_n}$.

Although it may not be obvious at first sight, the determinantal
expression (3.2-8) is exactly of the form of eq. (3.2-11). To see this,
one only has to divide both determinants in eq. (3.2-8) by the product
$(-1)^{n+1} b_{n+1}$ $\cdots$ $(-1)^{n+k+1} b_{n+k+1}$.

For the sake of simplicity no distinction between the elements of a
model sequence $\Seqn {\tilde s}$ and the elements of a sequence $\Seqn
s$, which is to be transformed, will be made from here on. This means
that if a sequence transformation is constructed on the basis of a
model sequence, then in the explicit expression or the recursive
scheme, which defines this transformation, the sequence elements $s_n$
and not the elements ${\tilde s}_n$ of the model sequence will occur.

\medskip

\Abschnitt The general extrapolation algorithm by Brezinski and
H{\aa}vie

\smallskip

\aktTag = 0

It is a typical feature of a large part of the modern nonlinear
sequence transformations that they are by construction exact for
special model sequences. The most general ansatz described in the
literature was introduced independently by Brezinski [31] and H{\aa}vie
[32]. They assume a model sequence of the following type:
$$
s_n \; = \; s + \sum_{j=0}^{k-1} \> c_j \, f_j (n) , \qquad k,n \in
\N_0 \, .
\tag
$$

Concerning the set $\Seq {f_j (n)}$ with $j,n \in \N_0$ it is assumed
that the $f_j (n)$ are known functions of $n$, but otherwise, they are
essentially completely arbitrary. Hence, the ansatz (3.3-1)
incorporates convergent as well as divergent sequences, depending upon
the behaviour of the functions $f_j (n)$ as $n \to \infty$.

In eq. (3.3-1), there occur $k+1$ unknowns, the limit or antilimit $s$
and the $k$ coefficients $c_0, \ldots , c_{k-1}$. Since all the
unknowns in eq. (3.3-1) occur linearly, the numerical values of $k+1$
sequence elements $s_n, s_{n+1}, \ldots , s_{n+k}$ have to be known in
order to be able to determine the limit or antilimit $s$ with the help
of Cramer's rule. Consequently, the general extrapolation algorithm
$E_k (s_n)$ by Brezinski and H{\aa}vie, which is by design exact for
sequences of the type of eq. (3.3-1), can be formulated as the ratio of
two determinants,

\smallskip

$$
E_k ( s_n ) \; = \; \frac
{
\vmatrix{
s_n "  \ldots " s_{n+k} \\ [1\jot]
f_0 (n) " \ldots " f_0 (n+k) \\
\vdots " \ddots " \vdots \\
f_{k-1} (n) " \ldots " f_{k-1} (n+k) }
}
{
\vmatrix{
1 " \ldots " 1       \\ [1\jot]
f_0 (n) " \ldots " f_0 (n+k) \\
\vdots " \ddots " \vdots \\
f_{k-1} (n) " \ldots " f_{k-1} (n+k) }
} \, .
\tag
$$

\smallskip

Brezinski and H{\aa}vie were also able to derive a recursive scheme for
the computation of the transforms $E_k (s_n)$, which is, however,
relatively complicated. A description of a FORTRAN IV program, which
computes the transforms $E_k (s_n)$ via this recursive scheme, can be
found in ref. [36].

Brezinski [31] showed that the general extrapolation algorithm $E_k
(s_n)$ contains the majority of the currently known sequence
transformations as special cases, among them Levin's sequence
transformation [28]. It will be seen later that many of the new
sequence transformations, which will be discussed in this report, are
actually special cases of the general extrapolation algorithm $E_k
(s_n)$.

Consequently, it may seem that it is sufficient to consider only the
general extrapolation algorithm $E_k (s_n)$ and not its numerous
special cases. However, the complicated structure of its recursive
scheme [31,32,36] makes the general extrapolation algorithm $E_k (s_n)$
computationally much less efficient than its special cases. Of
considerable importance is also the following aspect: In practical
applications it is certainly helpful to know that for arbitrary
functions $f_j (n)$ the sequence transformation $E_k (s_n)$ can be
computed recursively. But it is clearly of greater practical relevance
to find out which set $\Seq {f_j (n)}$ produces the best results for a
given sequence $\Seqn s$. Questions of that kind can only be answered
by studying special transformations and by exploiting their specific
properties.

\medskip

\Abschnitt Iterated sequence transformations

\smallskip

\aktTag = 0

In section 3.3, it was mentioned that a large part of the modern
nonlinear sequence transformations are constructed in such a way that
they are exact for certain model sequences. However, it is also
possible to find new sequence transformations which are not constructed
on the basis of model sequences.

Let us assume that a sequence $\Seqn s$ is to be transformed by a
sequence transformation $T_k^{(n)}$ with $k, n \in \N_0$ and that for
some $\kappa \in \N$, which is usually a relatively small number, a
transform $T_{\kappa}^{(n)}$ can be expressed explicitly in terms of
the sequence elements $s_n, s_{n+1}, \ldots, s_{n + \lambda}$, i.e.,
$$
T_{\kappa}^{(n)} \; = \;
F (s_n, s_{n+1}, \ldots, s_{n+\lambda}) \, .
\tag
$$

Then, a new sequence transformation ${\mit \Theta}_k^{(n)}$ can be
obtained by iterating the expression for $T_{\kappa}^{(n)}$. This means
that we define
$$
{\mit \Theta}_0^{(n)} \; = \; s_n \, , \qquad n \in \N_0 \, ,
\tag
$$
and that eq. (3.4-1) is rewritten in the following way:
$$
{\mit \Theta}_1^{(n)} \; = \; F_0 \left(
{\mit \Theta}_0^{(n)}, {\mit \Theta}_0^{(n+1)}, \ldots,
{\mit \Theta}_0^{(n + \lambda)} \right) \, ,
\qquad n \in \N_0 \, .
\tag
$$

This relationship can now be used to construct a recursive scheme by
means of which the transforms ${\mit \Theta}_k^{(n)}$ with $k \ge 2$
can be computed. One simple possibility of obtaining a recursive scheme
would be to assume that eq. (3.4-3) corresponds to the special case $k
= 0$ of the following, more general recursive scheme:
$$
{\mit \Theta}_{k+1}^{(n)} \; = \; F_k \left(
{\mit \Theta}_k^{(n)}, {\mit \Theta}_k^{(n+1)}, \ldots,
{\mit \Theta}_k^{(n + \lambda)} \right) \, ,
\qquad k, n \in \N_0 \, .
\tag
$$

Later, we shall encounter several very powerful sequence
transformations which are derived by iterating explicit expressions for
other sequence transformations. A well known example is Aitken's
iterated $\Delta^2$ process which is obtained by iterating the explicit
expression for Aitken's $\Delta^2$ process, eq. (1.2-2). Interestingly,
it often happens that the properties of the new sequence transformation
differ significantly from the properties of the transformation from
which it was derived.

\endAbschnittsebene

\neueSeite

\Abschnitt The epsilon algorithm and related topics

\vskip - 2 \jot

\beginAbschnittsebene

\medskip

\Abschnitt The Shanks transformation

\smallskip

\aktTag = 0

In his article on nonlinear sequence transformations Shanks [15]
considered the following model sequence:
$$
s_n \; = \; s + \sum_{j=0}^{k-1} \> c_j \, \Delta s_{n+j} \, ,\qquad n
\in \N_0 \, .
\tag
$$
If the sequence elements $s_n$ are partial sums of an infinite series,
$$
s_n \; = \; \sum_{\nu = 0}^n \, a_{\nu} \, ,
\tag
$$
the above model sequence can also be rewritten in the following way:
$$
s_n \; = \; s + \sum_{j=0}^{k-1} \> c_j \, a_{n+j+1} \, , \qquad n \in
\N_0 \, .
\tag
$$

Essentially this means that the limit $s$ of the infinite series is
approximated by the partial sum $s_n$ plus a weighted sum of the next
$k$ terms $a_{n+1}, a_{n+2}, \ldots , a_{n+k}$. As in the previous
examples the model sequence $s_n$ contains $k+1$ unknowns -- the limit
or antilimit $s$ and the $k$ linear coefficients $c_0, \ldots ,
c_{k-1}$ -- which all occur linearly. Consequently, according to
Cramer's rule the sequence transformation $e_k (s_n)$, which is by
construction exact for the model sequence (4.1-1), can be defined by
the following ratio of determinants:

\smallskip

$$
e_k ( s_n ) \; = \; \frac
{
\vmatrix{
s_n "  \ldots " s_{n+k} \\
\Delta s_n " \ldots " \Delta s_{n+k} \\
\vdots " \ddots " \vdots \\
\Delta s_{n+k-1} " \ldots " \Delta s_{n + 2 k - 1} }
}
{
\vmatrix{
1 " \ldots " 1       \\
\Delta s_n " \ldots " \Delta s_{n+k} \\
\vdots " \ddots " \vdots \\
\Delta s_{n+k-1} " \ldots " \Delta s_{n + 2 k - 1} }
} \; .
\tag
$$

\smallskip

For the computation of the transform $e_k ( s_n )$ the sequence
elements $s_n, \ldots , s_{n + 2 k}$ are needed. This implies that $e_k
( s_n )$ is a transformation of order $2 k$.

The sequence transformation $e_k ( s_n )$ is a special case of the
general extrapolation algorithm $E_k (s_n)$ introduced by Brezinski
[31] and H{\aa}vie [32]. To see this we only have to replace $f_j (n)$
in eq. (3.3-2) by $\Delta s_{n+j}$. The transformation (4.1-4) was
originally introduced in 1941 by Schmidt [17] who used it for the
iterative solution of linear systems. In 1955 it was rediscovered by
Shanks [15] who also analyzed some of the mathematical properties of
the sequence transformation $e_k (s_n)$ and derived several interesting
results. For instance, he was able to show that this transformation is
also exact for model sequences with remainders that are sums of
exponentials:
$$
s_n \; = \; s + \sum_{j=0}^{k-1} \, c_j {\lambda_j}^{\! n} \, , \qquad
n \in \N_0 \, .
\tag
$$

Concerning the $\lambda_j$ it is assumed that they are ordered in such
a way that their magnitudes are strictly decreasing, i.e.,
$$
| \lambda_0 | > | \lambda_1 | > \cdots > | \lambda_{k-1} | \, .
\tag
$$

If the condition $| \lambda_0 | < 1$ is satisfied the model sequence
(4.1-5) converges. In analogy with so-called {\it physical transients},
which disappear after a sufficiently long time, Shanks called the terms
on the right-hand side of eq. (4.1-5) {\it mathematical transients}
since all ${\lambda_j}^{\! n}$ with $| \lambda_j | < 1$ vanish as $n
\to \infty$. However, this concept of a mathematical transient has to
be used here in a broader sense since the Shanks transformation $e_k
(s_n)$ can also be used for the summation of divergent sequences and
series. The model sequence (4.1-5) diverges if at least one of the
$\lambda_j$ satisfies $| \lambda_k | \ge 1$, because then such a term
${\lambda_k}^{\! n}$ will not vanish as $n \to \infty$.

Shanks [15] also showed that the transformation $e_k (s_n)$ and Pad\'e
approximants are closely related. Let us assume that $f (z)$ is
analytic in a neighbourhood of $z = 0$,
$$
f (z) \; = \;
\sum_{\nu=0}^{\infty} \> \alpha_{\nu} \, z^{\nu} \, .
\tag
$$

Following the notation of Baker and Graves-Morris [22] we say that the
Pad\'e approximant of $f (z)$ is the ratio of two polynomials $p_{\ell}
(z)$ and $q_m (z)$ of degrees $\ell$ and $m$, respectively, and write
$$
[ \>\ell \> / \> m \>]_f \>(z) \; = \;
p_{\ell} (z) \> / \> q_m (z) \, .
\tag
$$

The Pad\'e approximants $[ \, \ell \, / \, m \, ]$ with $\ell , m \in
\N_0$ are displayed in a 2-dimensional rectangular scheme called the
{\it Pad\'e table} in such a way that the first index $\ell$ indicates
the column and the second index $m$ the row of the array. The
coefficients of the two polynomials $p_{\ell} (z)$ and $q_m (z)$ are
defined by the relationship
$$
f (z) \, - \, p_{\ell} (z) \, / \, q_m (z) \; = \; O(z^{\ell + m
+1}) \, .
\tag
$$

This implies that the coefficients of the Taylor expansion of the Pad\'e
approximant $p_{\ell} (z) \> / \> q_m (z)$ have to agree with the
series coefficients $\alpha_k$ up to the coefficient $\alpha_{\ell +
m}$. Let $f_n (z)$ stand for a partial sum of the power series (4.1-7),
$$
f_n (z) \; = \; \sum_{\nu=0}^n \> \alpha_{\nu} \, z^{\nu} \> .
\tag
$$

Shanks [15] could show that the application of his transformation to
the sequence $\Seq {f_n (z)}$ produces the following elements of the
Pad\'e table:
$$
e_k ( f_n (z) ) \; = \; [ \, n \, + \, k \, / \, k \, ]_f \, (z) \, ,
\qquad k,n \in \N_0 \, .
\tag
$$

\medskip

\Abschnitt Wynn's epsilon algorithm

\smallskip

\aktTag = 0

As it stands, the Shanks transformation $e_k (s_n)$, eq. (4.1-4), is
not particularly useful because of its definition as the ratio of two
determinants. Fortunately, only one year after the publication of
Shanks' article [15] Wynn [16] found a nonlinear recursive scheme which
is now commonly called the $\epsilon$ algorithm:
$$
\beginAligntags
" \epsilon_{-1}^{(n)} " \; = \; " 0 \, , \hfill \epsilon_0^{(n)}
\; = \; s_n \, , \erhoehe\aktTag \\ \tag*{\tagnr a}
" \epsilon_{k+1}^{(n)} " \; = \; " \epsilon_{k-1}^{(n+1)} \, + \,
1 / [\epsilon_{k}^{(n+1)} - \epsilon_{k}^{(n)} ]
\, , \qquad k,n \in \N_0 \, . \\ \tag*{\tagform\aktTagnr b}
\endAligntags
$$

Wynn [16] was able to show that the elements of the $\epsilon$ table
with even subscripts give the Shanks transformation,
$$
\epsilon_{2 k}^{(n)} \; = \; e_k (s_n) \, , \qquad k,n \in \N_0 \, ,
\tag
$$
whereas the elements of the $\epsilon$ table with odd subscripts are
only auxiliary quantities satisfying
$$
\epsilon_{2 k + 1}^{(n)} \, = \; 1 / e_k ( \Delta s_n ) \, , \qquad k,n
\in \N_0 \, .
\tag
$$

The publication of the $\epsilon$ algorithm, which allows a simple and
efficient evaluation of the Shanks transformation, stimulated an
enormous amount of research. According to Wimp (see p. 120 of ref.
[23]) over 50 articles on the $\epsilon$ algorithm were published by
Wynn alone, and at least 30 articles by Brezinski. As a fairly complete
source of references Wimp recommends Brezinski's first book [19]. Since
the main concern of this report are other sequence transformations and
not Wynn's $\epsilon$ algorithm, only those properties of the
$\epsilon$ algorithm will be discussed which are relevant for an
understanding of its power as well as its limitations as a convergence
acceleration and summation method.

In a later article, Wynn [37] analyzed the convergence properties of
the $\epsilon$ algorithm by applying it to several model sequences. For
instance, in the case of sequences, which have strictly alternating
remainders $r_n$ of the following type,
$$
s_n \; \sim \; s \, + \, (-1)^n \sum_{j=0}^{\infty} \, c_j \, / \, (n +
\beta )^{j+1} \, , \qquad \beta \in \R_{+} \, ,
\quad n \to \infty \, ,
\tag
$$
Wynn [37] obtained assuming $c_0 \ne 0$ for fixed $k$ an estimate which
shows that the $\epsilon$ algorithm accelerates convergence:
$$
\epsilon_{2 k}^{(n)} \, \sim \, s \, + \, \frac
{(-1)^n (k!)^2} {2^{2 k} (n + \beta)^{2 k + 1}} \> c_0 \, , \qquad n
\to \infty \; .
\tag
$$

Wynn [37] also considered sequences which generalize the model sequence
(4.1-5),
$$
s_n \; \sim \; s + \sum_{j=0}^{\infty} \, c_j {\lambda_j}^{\! n} \, ,
\qquad 1 > \lambda_0 > \lambda_1 > \cdots > 0 \, ,
\quad n \to \infty \, .
\tag
$$

For fixed $k$ Wynn [37] obtained the following estimate which shows
that the $\epsilon$ algorithm accelerates convergence:
$$
\epsilon_{2 k}^{(n)} \, \sim \, s \, + \, c_k \, \frac
{ \{ (\lambda_{k} - \lambda_{0}) (\lambda_{k} - \lambda_{1}) \ldots
(\lambda_{k} - \lambda_{k-1}) \}^2}
{\{ (1 - \lambda_{0}) (1 - \lambda_{1}) \ldots (1 - \lambda_{k-1})\}^2
} \, \lambda_{k}^n \, , \qquad n \to \infty \; .
\tag
$$

Wynn [37] also applied the $\epsilon$ algorithm to logarithmically
convergent sequences of the following type:
$$
s_n \; \sim \; s \, + \, \sum_{j=0}^{\infty} \, c_j \, / \, (n + \beta
)^{j+1} \, , \qquad \beta \in \R_{+} \, ,
\quad n \to \infty \, .
\tag
$$

Assuming $c_0 \ne 0$, Wynn [37] obtained for fixed $k$ an estimate
which shows that the $\epsilon$ algorithm -- or equivalently the Shanks
transformation -- is not able to accelerate the convergence of
sequences of that type:
$$
\epsilon_{2 k}^{(n)} \, \sim \, s \, + \, \frac
{c_0} {(k+1) (n + \beta )} \, , \qquad n \to\infty \; .
\tag
$$

This inability of accelerating logarithmic convergence is one of the
major defects of the otherwise very powerful $\epsilon$ algorithm.

If the elements $s_n$ of the sequence to be transformed are partial
sums of a power series as in eq. (4.1-10), the $\epsilon$ algorithm
produces Pad\'e approximants according to eqs. (4.1-11) and (4.2-2),
$$
\epsilon_{2 k}^{(n)} \; = \; [ \, n \, + \, k \, / \, k \, ]_f \, (z)
\, , \qquad k,n \in \N_0 \, ,
\tag
$$
and the convergence theory of Pad\'e approximants can be applied in this
case. This convergence theory can be found in the standard references
on Pad\'e approximants, as for instance the books by Baker [18] or by
Baker and Graves-Morris [22].

\medskip

\Abschnitt Programming the epsilon algorithm

\smallskip

\aktTag = 0

In this section it will be discussed how the $\epsilon$ algorithm can
be programmed efficiently. But first, the objectives of such a program,
which performs the transformation of a given sequence $\Seqn s$ with
the help of Wynn's $\epsilon$ algorithm, should be stated.

Since it is in general not possible to predict by a priori
considerations how many sequence elements will be needed until
convergence has finally been achieved, such a program should be
input-directed. This means it should read in the sequence elements
$s_0, s_1, \ldots , s_m, \ldots $ successively starting with $s_0$.
After the input of each new sequence element $s_m$ as many new elements
$\epsilon_k^{(n)}$ should be computed as is permitted by the recurrence
formula (4.2-1b). That new element $\epsilon_k^{(n)}$, which has the
largest even subscript $k$, should be used as the new approximation to
the limit of the sequence.

Let us arrange the elements $\epsilon_k^{(n)}$ of the $\epsilon$ table
in a rectangular scheme in such a way that the superscript $n$
indicates the row and the subscript $k$ the column of the 2-dimensional
array:

$$
\matrix{
\epsilon_0^{(0)} " \epsilon_1^{(0)} " \epsilon_2^{(0)} " \ldots "
\epsilon_n^{(0)} " \ldots \\
\epsilon_0^{(1)} " \epsilon_1^{(1)} " \epsilon_2^{(1)} " \ldots "
\epsilon_n^{(1)} " \ldots \\
\epsilon_0^{(2)} " \epsilon_1^{(2)} " \epsilon_2^{(2)} " \ldots "
\epsilon_n^{(2)} " \ldots \\
\epsilon_0^{(3)} " \epsilon_1^{(3)} " \epsilon_2^{(3)} " \ldots "
\epsilon_n^{(3)} " \ldots \\
\vdots " \vdots " \vdots " \ddots " \vdots " \ddots \\
\epsilon_0^{(n)} " \epsilon_1^{(n)} " \epsilon_2^{(n)} " \ldots "
\epsilon_n^{(n)} " \ldots \\
\vdots " \vdots " \vdots " \ddots " \vdots " \ddots }
\tag
$$

The entries in the first column of the array are the starting values
$\epsilon_0^{(n)} = s_n$ of the recursion according to eq. (4.2-1a).
The remaining elements of the $\epsilon$ table can be computed with
the help of the recurrence formula (4.2-1b). This nonlinear 4-term
recursion connects four elements of the $\epsilon$ table which are
located at the vertices of a rhombus:
$$
\beginMatrix
\beginFormat &\Formel\links \endFormat
\+" " \epsilon_{k}^{(n)} \qquad " \epsilon_{k+1}^{(n)} "
\\ [1\jot]
\+" \epsilon_{k-1}^{(n+1)} \qquad " \epsilon_{k}^{(n+1)} \qquad " " \\
\endMatrix
\tag
$$

\medskip

If the sequence elements $s_0, s_1, \ldots , s_m$ are used as initial
conditions according to eq. (4.2-1a), the recurrence formula (4.2-1b)
is able to compute all elements $\epsilon_j^{(\mu-j)}$ with $0 \le \mu
\le m$ and $0 \le j \le \mu$. Obviously, these elements form an
equilateral triangle located in the upper left corner of the $\epsilon$
table. If the next sequence element $s_{m+1}$ is also used as a
starting value for the recursion (4.2-1b), the triangle will be
enlarged by the neighbouring counterdiagonal $\epsilon_j^{(m-j+1)}$
with $0 \le j \le m+1$. In this context it will be advantageous to
rewrite Wynn's $\epsilon$ algorithm, eq. (4.2-1), in the following way:
$$
\beginAligntags
" \epsilon_0^{(n)} " \, = \, " s_n \, , \qquad n \ge 0 \, ,
\erhoehe\aktTag \\ \tag*{\tagnr a}
" \epsilon_1^{(n-1)} " \, = \,
" 1 / [ s_{n} - s_{n-1} ] \, ,
\qquad n \ge 1 \, , \\ \tag*{\tagform\aktTagnr b}
" \epsilon_{j}^{(n-j)} " \, = \, " \epsilon_{j-2}^{(n-j+1)} \, + \, 1 /
[\epsilon_{j-1}^{(n-j+1)} - \epsilon_{j-1}^{(n-j)}] \, ,
\quad n \ge 2 \, , \quad 2 \le j \le n \, .
\\ \tag*{\tagform\aktTagnr c}
\endAligntags
$$

Concerning the approximations to the limit it follows from eqs. (4.2-2)
and (4.2-3) that one has to distinguish between even and odd values of
the index $m$ of the last sequence element $s_m$ which was used in the
recursion. If $m$ is even, $m = 2 \mu$, our approximation to the limit
of the sequence is the transformation
$$
\left\{s_0, s_1, \ldots , s_{2 \mu} \right\} \; \to \;
\epsilon_{2 \mu}^{(0)} \; = \; e_{\mu} (s_0) \, ,
\tag
$$
and if $m$ is odd, $m = 2 \mu + 1$, we use the transformation
$$
\left\{s_1, s_2, \ldots , s_{2 \mu + 1} \right\} \; \to \;
\epsilon_{2 \mu}^{(1)} \; = \; e_{\mu} (s_1) \, .
\tag
$$

With the help of the notation $\Ent x$ for the integral part of $x$,
i.e., the largest integer $\nu$ satisfying $\nu \le x$, these two
relationships can be combined into a single equation yielding
$$
\left\{ s_{m - 2 \Ent {m/2}}, s_{m - 2 \Ent {m/2} + 1}, \ldots , s_m
\right\} \; \to \;
\epsilon_{2 \Ent {m/2}}^{(m - 2 \Ent {m/2})} \; = \;
e_{\Ent {m/2}} \bigl(s_{m - 2 \Ent {m/2}} \bigr) \, .
\tag
$$

If the sequence elements $s_n$ are partial sums of a power series as in
eq. (4.1-10), our approximations to the limit correspond according to
eq. (4.2-10) to the following staircase sequence in the Pad\'e table:
$$
[0/0]\; , \; [1/0]\; , \; [1/1]\; , \ldots ,\; [\nu / \nu]\; , \; [\nu
+ 1/ \nu] \; , \; [\nu +1/ \nu +1]\; , \ldots \; .
\tag
$$

Because of the rhombus structure (4.3-2) of the 4-term recursion in
Wynn's $\epsilon$ algorithm it appears that a program would either need
a single 2-dimensional or at least two 1-dimensional arrays. However,
Wynn [38] could show that a single 1-dimensional array is sufficient.
Wynn's algorithm, which is called {\it moving lozenge technique}, is
based upon the observation that for the computation of a new element
$\epsilon_j^{(m-j+1)}$ only the two neighbouring elements
$\epsilon_{j-1}^{(m-j+1)}$ and $\epsilon_{j-2}^{(m-j+2)}$ have to be
known but not the whole upper counterdiagonal
$\epsilon_{\mu}^{(m-\mu)}$ with $0 \le \mu \le m$. Hence, in Wynn's
moving lozenge technique these quantities are stored in auxiliary
variables while the recursion (4.3-3) moves along the current
counterdiagonal $\epsilon_j^{(m-j+1)}$ with $0 \le j \le m+1$ and
overwrites the previous entries $\epsilon_{\mu}^{(m - \mu)}$ with $0
\le \mu \le m$. A good description of Wynn's moving lozenge technique
[38] can also be found in Brezinski's second book (see pp. 326 - 327 of
ref. [20]).

Wynn [38] performed the recursion in a 1-dimensional array $E$ in which
he stored the elements of the current counterdiagonal of the $\epsilon$
table in such a way that the index of the array element coincides with
the subscript of the corresponding element of the $\epsilon$ table,
$$
\epsilon_j^{(m-j)} \; \to \; E(j) \, ,
\qquad m \ge 0 \, , \quad 0 \le j \le m \, .
\tag
$$

If the above convention is used three auxiliary variables will be
needed. But Wynn's moving lozenge technique can be improved if the
elements of the current counterdiagonal of the $\epsilon$ table are
stored in a 1-dimensional array $E$ in such a way that the superscript
of the corresponding element of the $\epsilon$ table coincides with the
index of the array element,
$$
\epsilon_j^{(m-j)} \; \to \; E(m-j) \, ,
\qquad m \ge 0 \, , \quad 0 \le j \le m \, .
\tag
$$

If this convention is used only two auxiliary variables will be needed
and the structure of the resulting computer program will also be
simpler and more elegant. The recursive scheme (4.3-3) can then be
reformulated in terms of the elements of the array $E$ in the following
way:
$$
\beginAligntags
" E(n) " \, \gets \, " s_n \, , \qquad n \ge 0 \, ,
\erhoehe\aktTag \\ \tag*{\tagnr a}
" E(n-1)" \, \gets \, " 1 / [ E(n) - E'(n-1) ] \, ,
\qquad n \ge 1 \, , \\ \tag*{\tagform\aktTagnr b}
" E(n-j) " \, \gets \, " E'(n-j+1) \, + \, 1 / [ E(n-j+1) - E'(n-j) ]
\, , \hfill \\
" n \ge 2 \, , " 2 \le j \le n \, .
\\ \tag*{\tagform\aktTagnr c}
\endAligntags
$$

The primed array elements $E'(n-j)$ and $E'(n-j+1)$ have to be stored
in auxiliary variables since they will be overwritten during the
computation of the current counterdiagonal $\epsilon_j^{(n-j)}$ with $0
\le j \le n$. The primes also indicate that the indices of the array
elements $E'(n-j)$ and $E'(n-j+1)$ refer to the occupation of $E$ after
the previous run, i.e., after the computation of the counterdiagonal
$\epsilon_j^{(n-j-1)}$ with $0 \le j \le n-1$.

If a counterdiagonal $\epsilon_{\mu}^{(m-\mu)}$ with $0 \le \mu \le m$
is computed with the help of the recursive scheme (4.3-3) and if the
elements $\epsilon_{2 k}^{(n)}$ converge -- which means that the whole
process is successful -- the computation of the elements $\epsilon_{2 k
+ 1}^{(n)}$ will necessarily involve divisions by the small quantities
$\epsilon_{2 k}^{(n+1)} - \epsilon_{2 k}^{(n)}$. This may easily lead
to an intolerable magnification of the inevitable rounding errors.
Hence, it looks as if Wynn's $\epsilon$ algorithm should be extremely
susceptible to rounding errors.

Fortunately, this is normally not the case although the elements
$\epsilon_{2 k +1}^{(n)}$ may become quite large in magnitude and may
have a very low relative accuracy due to the numerical problems
described above. But in the next step of the recursion the elements
$\epsilon_{2 k +1}^{(n)}$ serve as divisors which will dampen the
rounding errors again. Consequently, it is not clear what the overall
effect will be. However, numerical experience indicates that in most
cases of practical interest Wynn's $\epsilon$ algorithm is remarkably
stable. This experimental evidence is supported by a theoretical
analysis of the numerical stability of the $\epsilon$ algorithm which
was performed by Wynn [37] in the case of several model sequences .

In some cases -- for instance if the elements of the sequence to be
transformed are the partial sums of the Taylor series of a rational
function -- it may happen that the difference $\epsilon_{2 k}^{(n+1)} -
\epsilon_{2 k}^{(n)}$ vanishes. If pathologies of that kind occur, the
so-called {\it singular rules} of the $\epsilon$ algorithm can be used
which were also derived by Wynn [39]. A good discussion of these
singular rules and of related problems can also be found in section
4.1.2 of Brezinski's second book [20]. There, one can also find
listings of FORTRAN IV programs for Wynn's $\epsilon$ algorithm which
partly use the singular rules mentioned above (see pp. 338 - 352 of
ref. [20]).

According to the limited experience of the author pathologies of that
kind occur only rarely in scientific applications. Consequently, Wynn's
singular rules are not used in the following FORTRAN 77 program EPSAL
which computes the Shanks transformation of a given sequence by means
of Wynn's $\epsilon$ algorithm. However, a good program should take
some precautions against an approximate equality of the elements
$\epsilon_{k}^{(n+1)}$ and $\epsilon_{k}^{(n)}$ since in this case the
reciprocal of the difference of these two elements could exceed the
largest floating point number representable on the computer. This would
lead to overflow and to an error termination of the program.

This safeguard against overflow can be accomplished by defining two
variables HUGE and TINY. Their values should be close to but not
identical with the largest and smallest floating point numbers
representable on the computer. If the difference $\epsilon_{k}^{(n+1)}
- \epsilon_{k}^{(n)}$ is smaller in magnitude than TINY, then
$\epsilon_{k+1}^{(n)}$ will be set equal to HUGE. If this approximate
equality of the elements $\epsilon_{k}^{(n+1)}$ and
$\epsilon_{k}^{(n)}$ was accidental the program can continue with the
computation of the other elements of the $\epsilon$ table producing
numbers which are normally not noticeably affected. It is also possible
that this approximate equality of the elements $\epsilon_{k}^{(n+1)}$
and $\epsilon_{k}^{(n)}$ was not accidental but due to convergence.
However, in this case the program should have been stopped before.

In order to monitor the exact or approximate vanishing of the
denominators, it may also be a good idea to define in a FORTRAN program
for Wynn's $\epsilon$ algorithm an error variable, for instance IFAIL,
whose value is changed if one of the differences $\epsilon_{k}^{(n+1)}
- \epsilon_{k}^{(n)}$ is smaller in magnitude than TINY.

The following FORTRAN 77 program EPSAL uses the modification (4.3-10)
of Wynn's moving lozenge technique. It is safeguarded against
approximate equality of the elements $\epsilon_{k}^{(n+1)}$ and
$\epsilon_{k}^{(n)}$ by using two variables HUGE and TINY as described
above. The elements $s_n$ with $n = 0, 1, 2, \ldots $ of the sequence
to be transformed have to be computed in a DO loop in the calling
program. Whenever a new sequence element $s_n$ is computed in the outer
DO loop this subroutine EPSAL has to be called again and a new
counterdiagonal of the $\epsilon$ table will be calculated. The new
sequence element $s_n$ is read in via the variable SOFN and the
approximation to the limit is returned via the variable ESTLIM.

Finally, it is important to note that this subroutine EPSAL only
calculates the approximations to the limit according to eqs. (4.3-3)
and (4.3-4) and does not analyze the convergence of the whole process.
This has to be done in the calling program.

\bigskip

\listenvon{epsal.for}

\endAbschnittsebene

\neueSeite

\Abschnitt The iteration of Aitken's $\Delta^2$ process

\vskip - 2 \jot

\beginAbschnittsebene

\medskip

\Abschnitt Aitken's $\Delta^2$ transformation and its iteration

\smallskip

\aktTag = 0

Let us consider the following model sequence which is obtained by
setting $k=1$ in eq. (4.1-5):
$$
s_n \; = \; s \, + \, c \lambda^n \; , \qquad c \ne 0 \, , \;
|\lambda| < 1 \, , \; n \in \N_0 \, .
\tag
$$

Each sequence element $s_n$ contains the three unknowns $c$,
$\lambda$, and the limit $s$. Consequently, a sequence transformation
will at least require three elements of the above model sequence for
the determination of the limit $s$. In order to derive such a
transformation we form the first and second differences of $s_n$:
$$
\beginAligntags
" \Delta s_n \; " = \; " c \lambda^n ( \lambda - 1 ) \, , \\
\tag
" \Delta^2 s_n \; " = \; " c \lambda^n ( \lambda - 1 )^2 \, . \\
\tag
\endAligntags
$$

A short computation shows that the following sequence transformation is
exact for the model sequence (5.1-1):
$$
{\cal A}_1^{(n)} \; = \; s_n \; - \; \frac
{ [ \Delta s_n ]^2 } { \Delta^2 s_n } \, ,
\qquad n \in\N_0 \; .
\tag
$$

This sequence transformation is Aitken's well-known $\Delta^2$ process
[12]. The structure of this transformation explains quite clearly why
it bears this name.

It follows at once from the derivation of the sequence transformation
${\cal A}_1^{(n)}$ via the model sequence (5.1-1) that it is a special
case of the Shanks transformation, eq. (4.1-4), or Wynn's $\epsilon$
algorithm, eq. (4.2-1),
$$
{\cal A}_1^{(n)} \; = \; e_1 ( s_n ) \; = \; \epsilon_2^{(n)}\; .
\tag
$$

Many other representations for Aitken's $\Delta^2$ process can be
derived by suitable manipulations of eq. (5.1-4). Examples are:
$$
\beginAligntags
{\cal A}_1^{(n)}\; " = \; s_{n+1} \> - \> \frac
{ [\Delta s_n ] [\Delta s_{n+1}] } { \Delta^2 s_n} \\
\tag
" = \; s_{n+2} \> - \> \frac
{ [ \Delta s_{n+1} ]^2 } { \Delta^2 s_n} \\
\tag
" = \; \frac { s_{n+2} s_n - [ s_{n+1} ]^2 } { \Delta^2 s_n} \\
\tag
" = \; \frac
{ [\Delta s_{n+1} ] s_{n+1} - [ \Delta s_n ] s_{n+2} }
{ \Delta^2 s_n} \\
\tag
" = \; \frac
{ [ \Delta s_{n+1} ] s_{n} - [ \Delta s_n ] s_{n+1} }
{ \Delta^2 s_n} \\
\tag
" = \; s_{n+1} \> + \> \frac {1} { \Delta [ 1/ \Delta s_n ] } \; , \\
\tag
" = \; \frac
{\Delta \, [s_{n+1} / \Delta s_n ]} {\Delta \, [1/ \Delta s_n ]} \; . \\
\tag
\endAligntags
$$

Aitken's $\Delta^2$ process was studied in articles by Shanks [15],
Clark, Gray, and Adams [35], Lubkin [40], Tucker [41,42], Cordellier
[43], and Bell and Phillips [44]. A multidimensional generalization of
Aitken's transformation to vector sequences was discussed by MacLeod
[45]. Modifications of Aitken's $\Delta^2$ process were proposed by
Drummond [46] and by Bj{\o}rstad, Dahlquist, and Grosse [47].

The properties of Aitken's $\Delta^2$ process are discussed in books by
Brezinski (see pp. 37 - 40 of ref. [19] and pp. 43 - 45 of ref. [20])
and Wimp (see pp. 149 - 152 of ref. [23]). Those properties which are
particularly important for our purposes can be summarized as follows:

\beginBeschreibung \zu \Laenge{(ii):}

\item {(i) \hfill :} The $\Delta^2$ process accelerates linear
convergence.

\item {(ii):} The $\Delta^2$ process is regular but not accelerative
for logarithmically convergent sequences of the type of eq. (4.2-8).

\endBeschreibung

This shows that Aitken's $\Delta^2$ process has similar properties as
Wynn's $\epsilon$ algorithm. In view of eq. (5.1-5) this is not
surprising. However, one cannot expect that Aitken's $\Delta^2$ process
will be as powerful as Wynn's $\epsilon$ algorithm. The reason is that
the transform ${\cal A}_1^{(n)}$ is produced by only three sequence
elements $s_n, s_{n+1}$, and $s_{n+2}$ which implies that ${\cal
A}_1^{(n)}$ is a transformation of order $\ell = 2$. This will
certainly limit the power as well as practical usefulness of this
transformation.

If the accelerative power of Aitken's transformation turns out to be
insufficient and if it is necessary to use a more powerful sequence
transformation one could of course use Wynn's $\epsilon$ algorithm
which because of eq. (5.1-5) can be considered to be a more complex and
also more powerful generalization of Aitken's $\Delta^2$ process.
Another alternative, which also produces sequence transformations with
higher transformation orders, would be to iterate the $\Delta^2$
process. This means that Aitken's $\Delta^2$ process will be applied to
the transformed sequence $\bigSeq {{\cal A}_1^{(n)}}$ yielding a new
sequence $\bigSeq {{\cal A}_2^{(n)}}$. This process can in principle be
repeated indefinitely.

In order to obtain some heuristic motivation for this iteration, let us
apply Aitken's $\Delta^2$ process to the following model sequence which
generalizes the sequence (5.1-1):
$$
s_n \; = \; s \, + \, a x^n \, + \, b y^n \; , \qquad \, 0 < | y | < |
x | < 1 \; , \quad a,b \ne 0 \, .
\tag
$$

A short calculation shows that Aitken's $\Delta^2$ process eliminates
the dominating term $a x^n$ from the model sequence (5.1-13):
$$
{\cal A}_1^{(n)} \; = \; s \, + \ \frac
{b [(x-y)/(x-1)]^2 y^n} {1 + (b/a) [(y-1)/(x-1)]^2 (y/x)^n} \; .
\tag
$$

Since we have by assumption $0 < | y | < | x | < 1$, the transformed
sequence (5.1-14) converges faster than the original sequence (5.1-13).
Also, since $(y/x)^n$ vanishes as $n \to \infty$, at least for large
values of $n$ the elements of the resulting sequence (5.1-14) have
essentially the same structure as the elements of the sequence (5.1-1).

For the iteration of Aitken's $\Delta^2$ process each of the numerous
representations for ${\cal A}_1^{(n)}$ given above can be used since
they are all mathematically equivalent. However, the various
representations for ${\cal A}_1^{(n)}$ differ considerably in their
numerical stability. In the book by Press, Flannery, Teukolsky, and
Vetterling (see p. 133 of ref. [48]) it is remarked that Aitken's
$\Delta^2$ process should be computed with the help of eq. (5.1-4)
since the other equivalent representations are numerically less
reliable. Numerical studies performed by the author confirmed this
statement. Consequently, in this report an iteration of Aitken's
$\Delta^2$ process will always be based upon eq. (5.1-4). If we
identify the sequence elements $s_n$ with the initial values ${\cal
A}_0^{(n)}$ of the recursion we obtain the following nonlinear
recursive scheme:
$$
\beginAligntags
" {\cal A}_0^{(n)} \, " = " \, s_n \; ,
\hfill \erhoehe\aktTag \\ \tag*{\tagnr a}
" {\cal A}_{k+1}^{(n)} \, " = " \, {\cal A}_{k}^{(n)} \, - \, \frac
{\bigl[\Delta {\cal A}_{k}^{(n)}\bigr]^2}
{\Delta^2 {\cal A}_{k}^{(n)}}
\; , \qquad k,n \in \N_0 \, . \\
\tag*{\tagform\aktTagnr b}
\endAligntags
$$

As usual, the difference operator $\Delta$ acts upon the superscript
$n$ and not upon the subscript $k$. It follows from this recurrence
formula that the computation of ${\cal A}_k^{(n)}$ requires the
sequence elements $s_n, s_{n+1}, \ldots , s_{n + 2 k}$. Consequently,
${\cal A}_k^{(n)}$ is a transformation of order $2 k$. In this respect
${\cal A}_k^{(n)}$ is equivalent to $\epsilon_{2 k}^{(n)}$ which needs
the same set $s_n, s_{n+1}, \ldots , s_{n + 2 k}$ of sequence elements
for its computation. However, we shall see later that the numerical
properties of Wynn's $\epsilon$ algorithm and Aitken's iterated
$\Delta^2$ process often differ considerably although they are both
generalizations of the same sequence transformation ${\cal A}_1^{(n)}$,
eq. (5.1-4).

The numerical properties of Aitken's iterated $\Delta^2$ process were
studied by Smith and Ford [30]. Concerning the theoretical properties
of Aitken's iterated $\Delta^2$ process, very little seems to be known.
Apparently, there is only one article by Hillion [49] in which the
theoretical properties of Aitken's iterated $\Delta^2$ process were
studied. Hillion was able to find a model sequence for which the
iterated $\Delta^2$ process is exact. He also derived a determinantal
representation for the transforms ${\cal A}_k^{(n)}$. However,
Hillion's expressions contain in both cases explicitly the lower
order transforms ${\cal A}_0^{(n)}, \ldots , {\cal A}_{k-1}^{(n)},
\ldots , {\cal A}_0^{(n+k)}, \ldots , {\cal A}_{k-1}^{(n+k)}$.
Consequently, it seems that not much insight about the properties of
Aitken's iterated $\Delta^2$ process can be gained by these results.

\medskip

\Abschnitt Programming Aitken's iterated $\Delta^2$ process

\smallskip

\aktTag = 0

A program for Aitken's iterated $\Delta^2$ process should have the same
features as the subroutine EPSAL which transforms a given sequence
$\Seqn s$ with the help of Wynn's $\epsilon$ algorithm. This means it
should read in the sequence elements $s_0, s_1, \ldots , s_m , \ldots$
successively starting with $s_0$. After the input of each new sequence
element $s_m$ as many new elements ${\cal A}_k^{(n)}$ of the Aitken
table should be computed as it is permitted by the recursive scheme
(5.1-15). That element ${\cal A}_k^{(n)}$, which has the largest
subscript $k$, should be used as the new approximation to the limit of
the sequence $\Seqn s$.

Let us arrange the elements ${\cal A}_k^{(n)}$ of the Aitken table in
rectangular scheme in such a way that the superscript $n$ indicates the
row and the subscript $k$ the column of the 2-dimensional array:

$$
\matrix{
{\cal A}_0^{(0)} " {\cal A}_1^{(0)} " {\cal A}_2^{(0)} " \ldots " {\cal
A}_n^{(0)} " \ldots \\
{\cal A}_0^{(1)} " {\cal A}_1^{(1)} " {\cal A}_2^{(1)} " \ldots " {\cal
A}_n^{(1)} " \ldots \\
{\cal A}_0^{(2)} " {\cal A}_1^{(2)} " {\cal A}_2^{(2)} " \ldots " {\cal
A}_n^{(2)} " \ldots \\
{\cal A}_0^{(3)} " {\cal A}_1^{(3)} " {\cal A}_2^{(3)} " \ldots " {\cal
A}_n^{(3)} " \ldots \\
\vdots " \vdots  " \vdots  " \ddots " \vdots " \ddots \\
{\cal A}_0^{(n)} " {\cal A}_1^{(n)} " {\cal A}_2^{(n)} " \ldots " {\cal
A}_n^{(n)} " \ldots \\
\vdots  " \vdots " \vdots " \ddots " \vdots  " \ddots }
\tag
$$

The entries in the first column of the array are the starting values
${\cal A}_0^{(n)} = s_n$ of the recursion according to eq. (5.1-15a).
The remaining elements of the Aitken table can be computed with the
help of the recurrence formula (5.1-15b). The 4 elements, which are
connected by this nonlinear recursion, form a pattern in the Aitken
table which looks like the move of a knight on the chessboard:

$$
\beginMatrix
\beginFormat &\Formel\links \endFormat
\+" {\cal A}_k^{(n)} \qquad " {\cal A}_{k+1}^{(n)} " \\ [1\jot]
\+" {\cal A}_k^{(n+1)} \qquad " " \\ [1\jot]
\+" {\cal A}_k^{(n+2)} \qquad " " \\
\endMatrix
\tag
$$

\medskip

This pattern implies that the recursion (5.1-15b) has to proceed along
a relatively complicated path in the Aitken table if the elements $s_0,
s_1, \ldots , s_m , \ldots $ are read in successively and if one tries
to increase the subscript $k$ as much as possible. In this context, it
is advantageous to rewrite the recursive scheme (5.1-15) in the
following way:
$$
\beginAligntags
" {\cal A}_0^{(n)} \, " = " \, s_n \, , \qquad n \ge 0 \, ,
\erhoehe\aktTag \\ \tag*{\tagnr a}
" {\cal A}_{j}^{(n - 2 j)} \, " = " \, {\cal A}_{j-1}^{(n - 2 j)} \, -
\, \frac
{ \bigl[ \Delta {\cal A}_{j - 1}^{(n - 2 j)} \bigr]^2 }
{\Delta^2 {\cal A}_{j - 1}^{(n - 2 j)}} \; ,
\qquad n \ge 2 \, , \quad 1 \le j \le \Ent {n/2} \, . \\
\tag*{\tagform\aktTagnr b}
\endAligntags
$$

Here, $\Ent {n/2}$ denotes the integral part of $n/2$, i.e., the
largest integer $\nu$ satisfying $\nu \le n/2$. If the sequence
elements $s_0, s_1, \ldots , s_m$ are used as starting values, the
recursion (5.2-3b) is able to compute all elements ${\cal
A}_{j}^{(\mu-2 j)}$ with $0 \le \mu \le m$ and $0 \le j \le \Ent
{\mu/2}$. If the next sequence element $s_{m+1}$ is also used as a
starting value for the recursion, this set of elements of the Aitken
table will be enlarged by the string ${\cal A}_{j}^{(m-2 j+1)}$ with $0
\le j \le \Ent {(m+1)/2}$.

As in Wynn's $\epsilon$ algorithm the approximation to the limit
depends upon the index $m$ of the last sequence element $s_m$ which
was used in the recursion. If $m$ is even, $m = 2 \mu$, our
approximations to the limit of the sequence are the transformations
$$
\{ s_0, s_1, \ldots , s_{2 \mu}\} \; \to \;
{\cal A}_{\mu}^{(0)} \, ,
\tag
$$
and if $m$ is odd, $m = 2 \mu + 1$, the approximation to the limit will
be
$$
\{ s_1, s_2, \ldots , s_{2 \mu + 1}\} \; \to \;
{\cal A}_{\mu}^{(1)} \, .
\tag
$$

As in the case of Wynn's $\epsilon$ algorithm, these two relationships
can be combined into a single equation,
$$
\left\{ s_{m - 2 \Ent {m/2}}, s_{m - 2 \Ent {m/2} + 1},
\ldots , s_m \right\} \, \to \;
{\cal A}_{\Ent {m/2}}^{(m - 2 \Ent {m/2})} \, .
\tag
$$

Because of the relatively complicated geometrical structure (5.2-2) of
the recursion (5.2-3b) it appears that a program, which computes
Aitken's iterated $\Delta^2$ process, would need a 2-dimensional array.
However, a single 1-dimensional array $A$ is sufficient if the elements
of the Aitken table are stored according to the following rule:
$$
{\cal A}_{\Ent {\nu / 2}}^{(n - \nu)} \; \to \; A(n - \nu) \, ,
\qquad n \ge 0 \, , \quad 0 \le \nu \le n \, .
\tag
$$

With this convention the recurrence formula (5.2-3) can be reformulated
in terms of the elements of the 1-dimensional array $A$:
$$
\beginAligntags
" A(n) \, " \gets " \, s_n \, , \qquad n \ge 0 \, ,
\erhoehe\aktTag \\ \tag*{\tagnr a}
" A(n - 2 j) \, " \gets " \, A(n - 2 j) \, - \, \frac
{ [ \Delta A(n - 2 j) ]^2 }
{\Delta^2 A(n - 2 j)} \, ,
\qquad n \ge 2 \, , \quad 1 \le j \le \Ent {n/2} \, . \\
\tag*{\tagform\aktTagnr b}
\endAligntags
$$

Aitken's iterated $\Delta^2$ transformation makes sense only if the
second differences $\Delta^2 {\cal A}_{k-1}^{(n)}$ do not vanish for
sufficiently large values of $k$. This will certainly be guaranteed if
for fixed $k$ all elements of the sequence $\bigSeq {{\cal
A}_{k-1}^{(n)}}$ are different from zero and strictly alternating in
sign.

Unfortunately, the above statement is not particularly helpful since
only very little is known about the theoretical properties of Aitken's
iterated $\Delta^2$ process. In addition, it does not help at all if
the initial sequence $\Seqn s$ is not alternating. A related problem,
which may easily arise in this context, is that some second differences
$\Delta^2 {\cal A}_{k-1}^{(n)}$ may become so small that division would
lead to overflow. Consequently, a good program should be protected
against the exact or approximate vanishing of the second differences
$\Delta^2 {\cal A}_{k-1}^{(n)}$.

As in the case of Wynn's $\epsilon$ algorithm this can be accomplished
by introducing two variables HUGE and TINY which have values that are
close to but not identical with the largest and smallest floating point
number representable on the computer. If $\Delta^2 {\cal
A}_{k-1}^{(n)}$ is smaller in magnitude than TINY, ${\cal A}_k^{(n)}$
will be set equal to HUGE and the recursion is continued.

The following FORTRAN 77 subroutine AITKEN performs the recursive
computation of the Aitken table in a single 1-dimensional array $A$
according to eq. (5.2-7). It is safeguarded against an exact or
approximate vanishing of the second differences $\Delta^2 {\cal
A}_{k-1}^{(n)}$ by using two variables HUGE and TINY. The elements
$s_n$ with $n = 0, 1, 2, \ldots $ of the sequence to be transformed
have to be computed in a DO loop in the calling program. Whenever a new
sequence element $s_n$ is computed in the outer DO loop this subroutine
AITKEN has to be called again and a new string ${\cal A}_{j}^{(n-2j)}$
with $0 \le j \le \Ent {n/2}$ will be calculated. The new sequence
element $s_n$ is read in via the variable SOFN and the approximation to
the limit is returned via the variable ESTLIM.

It is important to note that this subroutine AITKEN only calculates the
approximations to the limit according to eqs. (5.2-4) and (5.2-5). The
convergence of the whole process has to be analyzed in the calling
program.

\bigskip

\listenvon{aitken.for}

\medskip

\endAbschnittsebene

\neueSeite

\Abschnitt Wynn's rho algorithm and related topics

\vskip - 2 \jot

\beginAbschnittsebene

\medskip

\Abschnitt Polynomial and rational extrapolation

\smallskip

\aktTag = 0

Assume that the values of a function $f(x)$ are only known at some
discrete points $x_0 < x_1 < \cdots < x_m$. It is one of the classical
problems of numerical analysis to estimate the value of $f$ at some
point $\xi \notin \{ x_0, x_1, \ldots , x_m \}$. If $x_0 < \xi < x_m$,
this problem is called {\it interpolation}, and if either $\xi < x_0$
or $x_m < \xi$, this problem is called {\it extrapolation}. These
problems and their solution are discussed in any book on numerical
analysis. More specialized treatments of these topics can be found in a
monograph on interpolation by Davis [50] or in a review article on
extrapolation processes by Joyce [51].

Extrapolation techniques can be used for the construction of
convergence acceleration methods. In this approach, the existence of a
function ${\cal S}$ of a continuous variable is postulated which
coincides on a discrete set of arguments $\Seqn x$ with the elements of
the sequence $\Seqn s$ to be transformed,
$$
{\cal S} (x_n) \; = \; s_n \, , \qquad n \in \N_0 \, .
\tag
$$

This ansatz reduces the problem of accelerating the convergence of a
sequence $\Seqn s$ to an extrapolation problem. If a finite string
$s_m, s_{m+1}, \ldots , s_{m+k}$ of $k+1$ sequence elements is known
one can construct an approximation ${\cal S}_k (x)$ to ${\cal S} (x)$
which satisfies the $k+1$ interpolation conditions
$$
{\cal S}_k (x_{m+j}) \; = \; s_{m+j} \, ,
\qquad 0 \le j \le k \, .
\tag
$$

In the next step one has to determine the value of the approximant
${\cal S}_k (x)$ for $x \to x_{\infty}$. If this can be done and if the
assumptions which are implicitly contained in this ansatz -- the
existence of a function ${\cal S} (x)$ which can be approximated at
least locally by a suitable set of interpolating functions -- are
justified, one can expect that the extrapolated value ${\cal S}_k
(x_{\infty})$ will provide a better approximation to the limit $s$ of
the sequence $\Seqn s$ than the last sequence element $s_{m+k}$ which
was used for the construction of ${\cal S}_k (x)$.

In interpolation and extrapolation problems the function under
consideration has to be modelled either in between or beyond a finite
set $x_0, x_1, \ldots, x_n$ of interpolation points by a suitable set
of interpolating functions. These interpolating functions should be
flexible and general enough to produce good approximations for large
classes of functions which can occur in practice. In addition, they
should also be simple enough to be manageable. The most common
interpolating functions are either polynomials or rational functions.
These two sets will also lead to different convergence acceleration
methods.

If interpolation by polynomials is used as the basis of a convergence
acceleration method it is implicitly assumed that the $k$-th order
approximant ${\cal S}_k (x)$ is a polynomial of degree $k$ in $x$,
$$
{\cal S}_k (x) \; = \, c_0 \, + \, c_1 x \, + \, \cdots \, + \, c_k x^k
\; .
\tag
$$

For polynomials, the most natural extrapolation point is $x = 0$.
Consequently, the interpolation points $x_n$ have to satisfy the
conditions
$$
\beginMultiline
x_0 > x_1 > x_2 > \cdots > x_m > x_{m+1} > \cdots > 0 \, ,
\erhoehe\aktTag \\ \tag*{\tagnr a}
\lim_{n \to \infty} \; x_n \; = \; 0 \, .
\\ \tag*{\tagform\aktTagnr b}
\endMultiline
$$

The choice $x=0$ as the extrapolation point implies that the
approximation to the limit is to be identified with the constant term
$c_0$ of the polynomial (6.1-3).

Several different methods for the computation of interpolating
polynomials ${\cal S}_k (x)$ are described in the mathematical
literature. Since only the constant term of a polynomial ${\cal S}_k$
has to be computed and since in most applications it is desirable to
compute simultaneously a whole string of approximants ${\cal S}_0 (0),
{\cal S}_1 (0), {\cal S}_2 (0), \ldots$, the most economical choice is
Neville's scheme [52] for the recursive computation of interpolating
polynomials. In the case $x=0$ Neville's algorithm reduces to the
following 2-dimensional linear recursive scheme (see p. 6 of ref. [20]):
$$
\beginAligntags
" {\cal N}_0^{(n)} (s_n , x_n ) \, " \; = \; " \, s_n
\, , \qquad n \in \N_0 \, ,
\erhoehe\aktTag \\ \tag*{\tagnr a}
" {\cal N}_{k+1}^{(n)} (s_n , x_n ) \, " \; = \; " \,
\frac
{x_n {\cal N}_{k}^{(n+1)} (s_{n+1} , x_{n+1} ) \, - \,
x_{n+k+1} {\cal N}_{k}^{(n)} (s_n , x_n ) }
{x_n \, - \, x_{n+k+1}} \; ,
\qquad k,n \in \N_0 \, . \\ \tag*{\tagform\aktTagnr b}
\endAligntags
$$

In the literature on convergence acceleration this variant of
Neville's recursive scheme is usually called Richardson extrapolation
[53]. Obviously, the linear transformation ${\cal N}_k^{(n)} (s_n , x_n
)$ is exact if the sequence elements $s_n$ are polynomials of degree
$k$ in the interpolation points $x_n$, i.e., for model sequences of the
following type:
$$
s_n \; = \; s \, + \, \sum_{j=0}^{k-1} \; c_j x_n^{j+1} \, ,
\qquad k,n \in \N_0 \, .
\tag
$$

The most obvious interpolation points for the Richardson extrapolation
scheme, eq. (6.1-5), are $x_n = 1 / (n + \beta)$ with $\beta > 0$ or
also $x_n = 1 / (n + \beta)^2$. These two choices are known to work
quite well in a variety of cases. However, if either one of these sets
of interpolation points $\Seqn x$ is used, the Richardson extrapolation
scheme (6.1-5) is not regular, i.e., the convergence of a sequence
$\Seqn s$ to some limit $s$ does not imply the convergence of the
transformed sequence to the same limit. In Brezinski's second book (see
pp. 37 - 38 of ref. [20]) it is shown that the regularity of the
Richardson extrapolation scheme is guaranteed only if some $a > 1$
exists such that the interpolation points $\Seqn x$ satisfy for all $n
\in \N_0$:
$$
x_n / x_{n+1} \; \ge \; a \; .
\tag
$$

This condition (6.1-7) is obviously fulfilled if the interpolation
points $\Seqn x$ satisfy $x_n = b^n$ with $0 < b < 1$ for all $n \in
\N_0$. A good discussion of the properties of the Richardson
extrapolation scheme as well as a list of various different sets of
interpolation points $\Seqn x$ can be found in Brezinski's second book
(see pp. 36 - 42 of ref. [20]).

It is well known that some functions can be approximated by polynomials
only quite poorly but by rational functions they can be approximated
very well. Consequently, it is likely that at least for some sequences
$\Seqn s$ rational extrapolation will give better results than
polynomial extrapolation. Let us therefore assume that the approximant
${\cal S}_k (x)$ can be written as the ratio of two polynomials of
degrees $\ell$ and $m$, respectively,
$$
{\cal S}_k (x) \; = \; \frac
{a_0 + a_1 x + a_2 x^2 + \cdots + a_{\ell} x^{\ell}}
{b_0 + b_1 x + b_2 x^2 + \cdots + b_m x^m} \; ,
\qquad k, \ell, m \in \N_0 \, .
\tag
$$

This rational function contains $\ell + m + 2$ coefficients $a_0,
\ldots , a_{\ell}$ and $b_0, \ldots , b_m$. However, only $\ell + m +1$
coefficients are independent since they are determined only up to a
common nonvanishing factor. Usually, one requires either $b_0 = 1$ or
$b_m = 1$. Consequently, the $k+1$ interpolation conditions (6.1-2)
will determine the coefficients $a_0, \ldots , a_{\ell}$ and $b_0,
\ldots , b_m$ provided that $k = \ell + m$ holds.

The extrapolation point $x=0$ is also the most obvious choice in the
case of rational extrapolation. Extrapolation to $x=0$ implies that
the interpolation points $\Seqn x$ have to satisfy eq. (6.1-4) and
that the approximation to the limit is to be identified with the ratio
$a_0/b_0$ of the constant terms of the polynomials in eq. (6.1-8).

However, if in eq. (6.1-8) $\ell = m$ holds, extrapolation to infinity
is also possible. In that case the interpolation points $\Seqn x$ would
have to satisfy
$$
\beginMultiline
0 < x_0 < x_1 < x_2 < \cdots < x_m < x_{m+1} < \cdots \, ,
\erhoehe\aktTag \\ \tag*{\tagnr a}
\lim_{n \to \infty} \; x_n \; = \; \infty \, .
\\ \tag*{\tagform\aktTagnr b}
\endMultiline
$$

In the case of extrapolation to infinity only the coefficients
$a_{\ell}$ and $b_{\ell}$ of the polynomials in eq. (6.1-8), which are
proportional to highest power $x^{\ell}$, contribute. Consequently, the
approximation to the limit has to be identified with the ratio
$a_{\ell}/b_{\ell}$.

As in the case of polynomial interpolation several different algorithms
for the computation of rational interpolants are described in the
literature. A discussion of the relative merits of these algorithms as
well as a survey of the relevant literature can be found in chapter III
of a book by Cuyt and Wuytack [54] which was recently published .

\medskip

\Abschnitt Wynn's rho algorithm

\smallskip

\aktTag = 0

Wynn's $\rho$ algorithm [25] is designed to compute even-order
convergents of Thiele's interpolating continued fraction [55] and to
extrapolate them to infinity. The even-order convergents are rational
functions of the following type:
$$
{\cal S}_{2 k} (x) \; = \; \frac
{a_k x^k + a_{k-1} x^{k-1} + \cdots + a_1 x + a_0 }
{b_k x^k + b_{k-1} x^{k-1} + \cdots + b_1 x + b_0} \, ,
\qquad k \in \N_0 \, .
\tag
$$

This means that the ratio $a_k/b_k$ is to be identified with the
approximation to the limit. According to Cuyt and Wuytack (see p. 214
of ref. [54]) Wynn's $\rho$ algorithm performs the computation of the
interpolating rational function (6.2-1) and its extrapolation to
infinity with a smaller number of arithmetic operations than similar
recursive algorithms.

Wynn's $\rho$ algorithm [25] is the following nonlinear recursive
scheme which is formally almost identical with Wynn's $\epsilon$
algorithm, eq. (4.2-1):
$$
\beginAligntags
" \rho_{-1}^{(n)} " \; = \; " 0 \, , \hfill \rho_0^{(n)} \; = \; s_n \,
, \erhoehe\aktTag \\ \tag*{\tagnr a}
" \rho_{k+1}^{(n)} " \; = \; " \rho_{k-1}^{(n+1)} \, + \, \frac
{ x_{n+k+1} - x_n } {\rho_{k}^{(n+1)} - \rho_{k}^{(n)} } \; , \qquad
k,n \in \N_0 \, . \\ \tag*{\tagform\aktTagnr b}
\endAligntags
$$

The only difference between Wynn's $\epsilon$ and Wynn's $\rho$
algorithm is that the $\rho$ algorithm also involves a sequence of
interpolation points $\Seqn x$ which have to satisfy eq. (6.1-9). As in
the case of Wynn's $\epsilon$ algorithm only the elements $\rho_{2
k}^{(n)}$ with even orders serve as approximations to the limit. The
elements $\rho_{2 k + 1}^{(n)}$ with odd orders are only auxiliary
quantities which diverge if the whole process converges.

Despite their formal similarity, the $\epsilon$ and $\rho$ algorithm
differ significantly in their ability of accelerating convergence. For
instance, the $\epsilon$ algorithm is exact for the model sequence
(4.1-5), and is known to be a very efficient accelerator for linearly
convergent sequences. In many cases the $\epsilon$ algorithm is also
able to sum divergent series. However, the otherwise very powerful
$\epsilon$ algorithm fails to accelerate logarithmic convergence.

The properties of Wynn's $\epsilon$ algorithm and Wynn's $\rho$
algorithm are in some sense complementary. Wynn's $\rho$ algorithm
fails to accelerate linear convergence and is not able to sum divergent
series. However, it is very powerful for some logarithmically
convergent sequences. This can easily be understood on the basis of the
following model sequence for which the transform $\rho_{2 k}^{(n)}$ is
exact:
$$
s_n \; = \, \frac
{ s x_n^k + a_1 x_n^{k-1} + \cdots + a_k }
{ x_n^k + b_1 x_n^{k-1} + \cdots + b_k } \, ,
\qquad k, n \in \N_0 \, .
\tag
$$

Since for fixed $k \in \N_0$ the zeros of the numerator and denominator
polynomials in eq. (6.2-3) are contained in a compact set and since the
interpolation points $\Seqn x$ diverge as $n \to \infty$, at least for
sufficiently large values of $n$ a rational function of that kind will
change only relatively slowly and monotonously with increasing $n$.
Certainly, such an expression will not oscillate or even diverge. This
should explain why the $\rho$ algorithm normally works well in the case
of logarithmic convergence but fails in the case of oscillating or
divergent sequences.

The properties of Wynn's $\rho$ algorithm are discussed in books by
Brezinski (see pp. 102 - 106 of ref. [19] and pp. 96 - 102 of ref.
[20]) and Wimp (see pp. 168 - 169 of ref. [23]). In these books the
connection of the $\rho$ algorithm with interpolating continued
fractions is emphasized and it is also shown that the transforms
$\rho_{2 k}^{(n)}$ can be represented as the ratio of two determinants.
But otherwise, relatively little seems to be known about this sequence
transformation.

The most obvious interpolation points $\Seqn x$ are $x_n = n + \beta$
with $\beta > 0$. With this choice, Wynn's $\rho$ algorithm assumes its
standard form:
$$
\beginAligntags
" \rho_{-1}^{(n)} " \; = \; " 0 \, , \hfill \rho_0^{(n)} \; = \; s_n \,
, \erhoehe\aktTag \\ \tag*{\tagnr a}
" \rho_{k+1}^{(n)} " \; = \; " \rho_{k-1}^{(n+1)} \, + \, \frac
{k+1} {\rho_{k}^{(n+1)} - \rho_{k}^{(n)} }, \qquad k,n \in \N_0 \, . \\
\tag*{\tagform\aktTagnr b}
\endAligntags
$$

Other possible sequences of interpolation points as for instance $x_n =
(n +\beta )^2$ with $\beta > 0$ are discussed in Brezinski's second
book [20].

As in Wynn's $\epsilon$ algorithm the approximation to the limit
depends upon the index $m$ of the last sequence element $s_m$ which
was used in the recursion. If $m$ is even, $m = 2 \mu$, our
approximation to the limit of the sequence is the transformation
$$
\{s_0, x_0; s_1, x_1; \ldots ; s_{2 \mu}, x_{2 \mu} \} \,
\to \; \rho_{2 \mu}^{(0)} \; ,
\tag
$$
and if $m$ is odd, $m = 2 \mu + 1$, we use the transformation
$$
\{s_1, x_1; s_2, x_2; \ldots ; s_{2 \mu + 1}, x_{2 \mu + 1} \} \,
\to \; \rho_{2 \mu}^{(1)} \; .
\tag
$$

With the help of the notation $\Ent x$ for the integral part of $x$,
i.e., the largest integer $\nu$ satisfying $\nu \le x$, these two
relationships can be combined into a single equation yielding
$$
\left\{ s_{m - 2 \Ent {m/2}}, x_{m - 2 \Ent {m/2}}; s_{m - 2 \Ent {m/2}
+ 1}, x_{m - 2 \Ent {m/2} + 1}; \ldots ; s_m, x_m \right\}
\; \to \;
\rho_{2 \Ent {m/2}}^{(m - 2 \Ent {m/2})} \, .\tag
$$

The elements of the $\rho$ table can be arranged in the same
rectangular scheme as the elements of the $\epsilon$ table in (4.3-1).
Since the recurrence relationships for Wynn's $\epsilon$ algorithm and
Wynn's $\rho$ algorithm are structurally identical, those elements of
the $\rho$ table which are connected by the 4-term recurrence formula
(6.2-2b), are also located in the $\rho$ table at the vertices of a
rhombus:

$$
\beginMatrix
\beginFormat &\Formel\links \endFormat
\+" " \rho_{k}^{(n)} \qquad " \rho_{k+1}^{(n)} " \\ [1\jot]
\+" \rho_{k-1}^{(n+1)} \qquad " \rho_{k}^{(n+1)} \qquad " " \\
\endMatrix
\tag
$$

\medskip

Consequently, Wynn's $\rho$ algorithm can be programmed in exactly the
same way as Wynn's $\epsilon$ algorithm. For that purpose we rewrite
the recursive scheme (6.2-2) in the following way:
$$
\beginAligntags
" \rho_0^{(n)} " \, = \, " s_n \, , \qquad n \ge 0 \, ,
\erhoehe\aktTag \\ \tag*{\tagnr a}
" \rho_1^{(n-1)} " \, = \, " \frac
{ x_n - x_{n-1} } {s_{n} - s_{n-1} } \, ,
\qquad n \ge 1 \, , \\ \tag*{\tagform\aktTagnr b}
" \rho_{j}^{(n-j)} " \, = \, " \rho_{j-2}^{(n-j+1)}
\, + \, \frac { x_n - x_{n-j} }
{ \rho_{j-1}^{(n-j+1)} - \rho_{j-1}^{(n-j)} } \, ,
\quad n \ge 2 \, , \quad 2 \le j \le n \, .
\\ \tag*{\tagform\aktTagnr c}
\endAligntags
$$

As in the case of the $\epsilon$ algorithm the modification (4.3-9) of
Wynn's moving lozenge technique (4.3-8) can be used. This means that
only a single 1-dimensional array ${\it R}$ will be needed if the
elements of the current counterdiagonal $\rho_j^{(m-j)}$ with $m \in
\N_0$ and $0 \le j \le m$ are stored in $R$ in such a way that the
superscript of the element of the $\rho$ table coincides with the index
of the corresponding array element,
$$
\rho_j^{(m-j)} \; \to \; {\it R} (m-j) \, .
\tag
$$

The only difference with Wynn's $\epsilon$ algorithm is that also a
second 1-dimensional array $\xi$ will be needed in which the
interpolation points $x_n$ are stored according to the rule
$$
x_n \; \to \; \xi (n) \, .
\tag
$$

With these two conventions the recursive scheme (6.2-9) can be
reformulated in terms of the elements of the 1-dimensional arrays ${\it
R}$ and $\xi$,
$$
\beginAligntags
" {\it R} (n) " \, \gets \, " s_n \, , \qquad n \ge 0 \, ,
\erhoehe\aktTag \\ \tag*{\tagnr a}
" {\it R} (n-1)" \, \gets \,
" \frac { \xi (n) - \xi (n-1) }
{ {\it R} (n) - {\it R}^{\prime} (n-1) } \, ,
\qquad n \ge 1 \, , \\ \tag*{\tagform\aktTagnr b}
" {\it R} (n-j) " \, \gets \, " {\it R}^{\prime} (n-j+1)
\, + \, \frac { \xi (n) - \xi (n-j) }
{ {\it R} (n-j+1) - {\it R}^{\prime} (n-j) } \, , \hfill \\
" n \ge 2 \, , \qquad 2 \le j \le n \, . \span\omit \span\omit
\\ \tag*{\tagform\aktTagnr c}
\endAligntags
$$

As in the case of Wynn's $\epsilon$ algorithm, the primed array
elements ${\it R}^{\prime} (n-j)$ and ${\it R}^{\prime} (n-j+1)$ have
to be stored in auxiliary variables. The primes also indicate that the
array elements ${\it R}^{\prime} (n-j)$ and ${\it R}^{\prime} (n-j+1)$
refer to the occupation of ${\it R}$ after the previous run, i.e.,
after the computation of the counterdiagonal $\rho_j^{n-j-1}$ with $0
\le j \le n - 1$. The listings of some FORTRAN IV programs, which
compute Wynn's $\rho$ algorithm, can be found in Brezinski's second
book (see pp. 361 - 365 of ref. [20]).

\medskip

\Abschnitt The iteration of Wynn's rho algorithm

\smallskip

\aktTag = 0

According to eq. (5.1-5), Aitken's $\Delta^2$ process is identical with
the transform $\epsilon_2^{(n)}$,
$$
{\cal A}_1^{(n)} \; = \; s_n \; - \; \frac
{ [ \Delta s_n ]^2 } { \Delta^2 s_n } \; = \; \epsilon_2^{(n)}
\, , \qquad n \in \N_0 \; .
\tag
$$

If Aitken's $\Delta^2$ process is iterated, a new sequence
transformation ${\cal A}_k^{(n)}$ results which has a similar ability
of accelerating convergence as Wynn's $\epsilon$ algorithm. However, in
section 13 we shall encounter some examples in which Aitken's iterated
$\Delta^2$ process clearly outperforms Wynn's $\epsilon$ algorithm.
This alone should justify an analysis of Aitken's iterated $\Delta^2$
transformation ${\cal A}_k^{(n)}$.

Since Wynn's $\epsilon$ algorithm (4.2-1) and Wynn's $\rho$ algorithm
(6.2-2) are formally almost identical, one can construct a new sequence
transformation by proceeding as in the case of Aitken's iterated
$\Delta^2$ process. This means that first the transform $\rho_2^{(n)}$
is expressed in terms of some sequence elements $s_n$ and interpolation
points $x_n$ The resulting expression for $\rho_2^{(n)}$ will then be
iterated.

From eqs. (6.2-2a) and (6.2-2b) we obtain the following expression for
the $\rho$ analogue of Aitken's $\Delta^2$ process:
$$
\rho_2^{(n)} \; = \; s_{n+1} \, + \, \frac
{(x_{n+2} - x_n) [ \Delta s_{n+1} ] [ \Delta s_n ] }
{ [ \Delta x_{n+1} ] [ \Delta s_{n} ] -
[ \Delta x_n ] [ \Delta s_{n+1} ] }
\, , \qquad n \in \N_0 \, .
\tag
$$

An iteration of this expression, which may be considered to be a kind
of weighted $\Delta^2$ process, can be done in a variety of ways. The
problem is that there is no unique way of choosing the indices of the
interpolation points $x_n$. However, if we take into account that in
Wynn's $\rho$ algorithm, eq. (6.2-2), the differences of the indices of
the interpolation points $x_n$ increase with increasing $k$, we see
that the following nonlinear recursive scheme should be the most
natural iteration of the transform (6.3-2):
$$
\beginAligntags
" {\cal W}_{0}^{(n)} \; " = \; " s_n \, ,
\hfill n \in \N_0 \, , \erhoehe\aktTag \\ \tag*{\tagnr a}
" {\cal W}_{k+1}^{(n)} \; " = \; " {\cal W}_{k}^{(n+1)} \,
+ \frac
{ (x_{n + 2 k + 2} - x_n ) \bigl[ \Delta {\cal W}_k^{(n+1)}
\bigr] \bigl[ \Delta {\cal W}_k^{(n)} \bigr] }
{(x_{n + 2 k + 2} - x_{n+1}) \bigl[ \Delta {\cal W}_k^{(n)} \bigr] -
(x_{n + 2 k + 1} - x_n) \bigl[ \Delta {\cal W}_k^{(n+1)} \bigr] } \, ,
\\
" k,n \in \N_0 \, . \span\omit \span\omit
\\ \tag*{\tagform\aktTagnr b}
\endAligntags
$$

As usual it is assumed that the difference operator $\Delta$ acts upon
$n$ and not upon $k$. The most obvious interpolation points are as in
Wynn's $\rho$ algorithm $x_n = n + \beta$ with $\beta > 0$. With this
choice, the iterated $\rho_2$ transformation assumes its standard form:
$$
\beginAligntags
" {\cal W}_{0}^{(n)} \; " = \; " s_n \, ,
\hfill n \in \N_0 \, , \erhoehe\aktTag \\ \tag*{\tagnr a}
" {\cal W}_{k+1}^{(n)} \; " = \; " {\cal W}_{k}^{(n+1)} \, -
\frac
{ (2 k + 2) \bigl[ \Delta {\cal W}_k^{(n+1)} \bigr]
\bigl[ \Delta {\cal W}_k^{(n)} \bigr] }
{ (2 k + 1) \Delta^2 {\cal W}_k^{(n)} }
\, , \hfill \qquad k,n \in \N_0 \, . \\
\tag*{\tagform\aktTagnr b}
\endAligntags
$$

The table of this transformation can be arranged in the same
rectangular scheme as the Aitken table (5.2-1). Also, the recurrence
formulas for Aitken's iterated $\Delta^2$ process, eq. (5.1-14), and
the recurrence formulas for the transforms ${\cal W}_k^{(n)}$ are
structurally identical. This implies that the four elements, which are
connected by the nonlinear recurrence formulas (4.3-2) or (4.3-3), also
form a pattern in the ${\cal W}$ table which looks like the move of a
knight on the chessboard:

$$
\beginMatrix
\beginFormat &\Formel\links \endFormat
\+" {\cal W}_k^{(n)} \qquad " {\cal W}_{k+1}^{(n)} " \\ [1\jot]
\+" {\cal W}_k^{(n+1)} \qquad " " \\ [1\jot]
\+" {\cal W}_k^{(n+2)} \qquad " " \\
\endMatrix
\tag
$$

\medskip

Consequently, this iterated $\rho_2$ process can be programmed in
exactly the same way as Aitken's iterated $\Delta^2$ process and only a
few minor alterations have to be done in the subroutine AITKEN. In this
context, it is advantageous to reformulate the recurrence scheme
(6.3-4) in the following way:
$$
\beginAligntags
" {\cal W}_0^{(n)} \, " = \, " s_n \, , \qquad n \ge 0 \, ,
\erhoehe\aktTag \\ \tag*{\tagnr a}
" {\cal W}_{j}^{(n - 2 j)} \, " = \,
" {\cal W}_{j-1}^{(n - 2 j + 1)} + \frac
{
( x_n - x_{n - 2 j} )
\bigl[ \Delta {\cal W}_{j - 1}^{(n - 2 j + 1)} \bigr]
\bigl[ \Delta {\cal W}_{j - 1}^{(n - 2 j)} \bigr]
}
{
( x_n - x_{n - 2 j + 1} )
\bigl[ \Delta {\cal W}_{j - 1}^{(n - 2 j)} \bigr] -
( x_{n-1} - x_{n - 2 j} )
\bigl[ \Delta {\cal W}_{j - 1}^{(n - 2 j + 1)} \bigr]
} \, , \\
" n \ge 2 \; , \quad 1 \le j \le \Ent {n/2} \, .
\span\omit \span\omit
\\ \tag*{\tagform\aktTagnr b}
\endAligntags
$$

As usual, $\Ent {n/2}$ stands for the integral part of $n/2$, i.e., the
largest integer $\nu$ satisfying $\nu \le n/2$. If the sequence
elements $s_0, s_1, \ldots , s_m$ are used as starting values, the
recursive scheme (6.3-6) is able to compute all elements ${\cal
W}_{j}^{(\mu-2 j)}$ with $0 \le \mu \le m$ and $0 \le j \le \Ent
{\mu/2}$.

As in the case of Aitken's iterated $\Delta^2$ process the
approximation to the limit depends upon the index $m$ of the last
sequence element $s_m$ which was used in the recursion. If $m$ is even,
$m = 2 \mu$, our approximations to the limit of the sequence are the
transformations
$$
\{ s_0, x_0; s_1, x_1; \ldots , s_{2 \mu}, x_{2 \mu}\}
\; \to \; {\cal W}_{\mu}^{(0)} \, ,
\tag
$$
and if $m$ is odd, $m = 2 \mu + 1$, the approximation to the limit will
be
$$
\{ s_1, x_1; s_2, x_2; \ldots , s_{2 \mu + 1}, x_{2 \mu + 1}\}
\; \to \; {\cal W}_{\mu}^{(1)} \, .
\tag
$$

These two relationships can be combined into a single equation,
$$
\left\{ s_{m - 2 \Ent {m/2}}, x_{m - 2 \Ent {m/2}}; s_{m - 2 \Ent {m/2}
+ 1}, x_{m - 2 \Ent {m/2} + 1};\ldots ; s_m, x_m \right\}
\; \to \; {\cal W}_{\Ent {m/2}}^{(m - 2 \Ent {m/2})} \, .
\tag
$$

Only two 1-dimensional arrays $w$ and $\xi$ are needed if the
interpolation points $x_n$ are stored in $\xi$ according to eq.
(6.2-11) and if the elements ${\cal W}_k^{(n)}$ are stored in $w$
according to the following rule:
$$
{\cal W}_{\Ent {\nu / 2}}^{(n - \nu)} \; \to \; w (n - \nu)
\, , \qquad n \ge 0 \, , \quad 0 \le \nu \le n \, .
\tag
$$

With this convention the recursive scheme (6.3-6) can be reformulated
in terms of the elements of the 1-dimensional arrays $w$ and $\xi$:
$$
\beginAligntags
" w (n) \, " \gets \, " s_n \, , \qquad n \ge 0 \; , \hfill
\erhoehe\aktTag \\ \tag*{\tagnr a}
" w (\ell) \, " \gets \, " w (\ell + 1) \, + \, \frac
{ [ \xi (n) - \xi (\ell + 1) ]
[ \Delta w (\ell + 1) ] [ \Delta w (\ell ) ] }
{
[ \xi (n) - \xi (\ell + 1) ] [ \Delta w (\ell) ] -
[ \xi (n-1) - \xi (\ell ) ] [ \Delta w (\ell + 1) ]
} \; , \openup 3pt \\
" \ell = n - 2 j \; , \qquad n \ge 2 \; ,
\qquad 1 \le j \le \Ent {n/2} \, . \span\omit \span\omit \\
\tag*{\tagform\aktTagnr b}
\endAligntags
$$

\endAbschnittsebene

\neueSeite

\Abschnitt The Levin transformation

\vskip - 2 \jot

\beginAbschnittsebene

\medskip

\Abschnitt The derivation of Levin's sequence transformation

\smallskip

\aktTag = 0

Levin's sequence transformation [28] is designed to be exact for model
sequences of the following type:
$$
s_n \; = \; s \, + \, \omega_n \, \sum_{j=0}^{k-1} \, c_j / (n +
\beta)^j \, , \qquad k, n \in \N_0 \, .
\tag
$$

Here, the remainder estimates $\omega_n$ are essentially arbitrary
functions of $n$. It is only assumed that they are different from zero
for all finite values of $n$. In addition, it would not make much sense
to consider in convergence acceleration and summation processes
remainder estimates which are constant. Consequently, we shall also
assume that for all finite values of $n$ the elements of the sequence
$\Seqn {\omega}$ are all distinct. Depending upon the behaviour of the
remainder estimates $\Seq {\omega_n}$ as $n \to \infty$, the sequence
$\Seqn s$ may either converge or diverge. In eq. (7.1-1) it has to be
required that $\beta + n$ must not be zero. This implies that the
parameter $\beta$ must not be zero or a negative integer. However, the
elements of the model sequence (7.1-1) will serve as finite
approximations to Poincar\'e-type asymptotic expansions of the following
kind:
$$
s_n \; \sim \; s \, + \, \omega_n \, \sum_{j=0}^{\infty} \, c_j / (n +
\beta)^j \, , \qquad n \to \infty \, .
\tag
$$

In expansions of that kind negative values of $\beta$ will lead to
different signs of the terms if either $n + \beta < 0$ or $n + \beta >
0$ holds. Since model sequences of the type of eq. (7.1-1) will be used
as approximations for asymptotic expansions of the above type and since
these approximations should be uniformly valid over a wide range of
admissible values of $n$ -- preferably for all $n \ge 0$ -- it is
necessary to require that the sign pattern of the terms of the sum in
eq. (7.1-1) must not depend upon $n$. This rules out $\beta < 0$ and we
have to require $\beta > 0$. But otherwise, the parameter $\beta$ is in
principle completely arbitrary. In the literature on Levin's sequence
transformation, only the case $\beta = 1$, which is the most obvious
choice, has been considered so far.

In eq. (7.1-1) there occur $k+1$ unknown quantities, the limit or
antilimit $s$ and the $k$ linear coefficients $c_0, \ldots , c_{k-1}$.
Hence, $k+1$ sequence elements $s_n, \ldots, s_{n+k}$ are needed for
the determination of $s$, and according to Cramer's rule the general
Levin transformation ${\cal L}_k^{(n)} (\beta , s_n , \omega_n)$ can
be defined by the following ratio of determinants:

\smallskip

$$
{\cal L}_k^{(n)} (\beta , s_n , \omega_n) \; = \; \frac
{
\vmatrix
{
s_n "  \ldots " s_{n+k} \\ [1\jot]
\omega_{n} " \ldots " \omega_{n+k} \\
\vdots " \ddots " \vdots \\ [1\jot]
\omega_{n} / (\beta + n)^{k-1} " \ldots "
\omega_{n+k} / (\beta + n + k)^{k-1}
}
}
{
\vmatrix
{
1 " \ldots " 1       \\ [1\jot]
\omega_{n} " \ldots " \omega_{n+k} \\
\vdots " \ddots " \vdots \\ [1\jot]
\omega_{n} / (\beta + n)^{k-1} " \ldots "
\omega_{n+k} / (\beta + n + k)^{k-1}
}
} .
\tag
$$

\smallskip

If the sequence elements $s_n, s_{n+1}, \ldots, s_{n+k}$ satisfy eq.
(7.1-1) then Levin's general sequence transformation is exact by
construction, i.e.,
$$
{\cal L}_k^{(n)} (\beta , s_n , \omega_n) \; = \; s \, .
\tag
$$

The representation of the general Levin transformation as the ratio of
two determinants is not well suited for practical applications because
the fast and reliable evaluation of large order determinants is a not
yet satisfactorily solved problem of numerical analysis. Thus,
alternative expressions for the general Levin transformation are highly
desirable. Fortunately, they can be derived quite easily.

Levin's original derivation [28] of nondeterminantal expressions for
his sequence transformation was based upon the observation that the
determinants in eq. (7.1-3) may be expressed in terms of Vandermonde
determinants. However, for our purposes it is advantageous to follow
Sidi's approach [56] which exploits properties of the difference
operator $\Delta$, since this approach can easily be extended to other
sequence transformations which will be treated later in this report.
For that purpose, eq. (7.1-1) is rewritten in the following way:
$$
(n + \beta)^{k-1} [s_n - s] / \omega_n \; = \; \sum_{j=0}^{k-1} \, c_j
\, (n + \beta)^{k-j-1} \, .
\tag
$$

The highest power of n, which occurs on the right-hand side of eq.
(7.1-5), is $n^{k-1}$. We now utilize the well-known fact that any
polynomial of degree $k-1$ in $n$ will be annihilated by the difference
operator $\Delta^k$. Since the difference operator $\Delta^k$ is
linear, we may conclude from eqs. (7.1-4) and (7.1-5) that the general
Levin transformation ${\cal L}_k^{(n)} (\beta , s_n , \omega_n)$ is
given by the following ratio:
$$
{\cal L}_k^{(n)} (\beta , s_n , \omega_n) \; = \; \frac
{ \Delta^k \, \{ (n + \beta)^{k-1} \> s_n / \omega_n\} }
{ \Delta^k \, \{ (n + \beta)^{k-1} / \omega_n \} } \, .
\tag
$$

With the help of eq. (2.4-8) the action of the difference operators in
eq. (7.1-6) can be expressed in closed form. This leads to a
representation of the general Levin transformation as the ratio of two
finite sums:
$$
{\cal L}_{k}^{(n)} (\beta , s_n, \omega_n) \; = \;
\frac
{\displaystyle
\sum_{j=0}^{k} \; ( - 1)^{j} \; \binom {k} {j} \;
\frac {(\beta + n +j )^{k-1}} {(\beta + n + k )^{k-1}} \;
\frac {s_{n+j}} {\omega_{n+j}} }
{\displaystyle
\sum_{j=0}^{k} \; ( - 1)^{j} \; \binom {k} {j} \;
\frac {(\beta + n +j )^{k-1}} {(\beta + n + k )^{k-1}} \;
\frac {1} {\omega_{n+j}} }
\; , \qquad k, n \in \N_0 \, .
\tag
$$

The common factor $(\beta + n + k)^{k-1}$ in eq. (7.1-7) was introduced
in order to decrease the magnitude of the terms of the numerator and
denominator sums, because otherwise overflow may happen too easily for
larger values of $k$.

A mild extension of the general Levin transformation, eq. (7.1-7), will
also be considered:
$$
{\cal L}_{k,\ell}^{(n)} (\beta , s_n, \omega_n) \; = \;
\frac
{\displaystyle
\sum_{j=0}^{k} \; ( - 1)^{j} \; \binom {k} {j} \;
\frac {(\beta + n +j )^{k-\ell-1}} {(\beta + n + k )^{k-1}} \;
\frac {s_{n+j}} {\omega_{n+j}} }
{\displaystyle
\sum_{j=0}^{k} \; ( - 1)^{j} \; \binom {k} {j} \;
\frac {(\beta + n +j )^{k-\ell-1}} {(\beta + n + k )^{k-1}} \;
\frac {1} {\omega_{n+j}} } \; ,
\qquad k,\ell , n \in \N_0 \, .
\tag
$$

For $\ell = 0$ this transformation reduces to the general Levin
transformation, eq. (7.1-7). An alternative representation for this
generalization of Levin's sequence transformation can be derived with
the help of eq. (2.4-8),
$$
{\cal L}_{k,\ell}^{(n)} (\beta , s_n , \omega_n) \; = \; \frac
{ \Delta^k \, \{ (n + \beta)^{k-\ell-1} \> s_n / \omega_n \} }
{ \Delta^k \, \{ (n + \beta)^{k-\ell-1} \> / \omega_n \} } \, .
\tag
$$

This relationship implies that this generalization of Levin's sequence
transformation is exact for sequences of the following type:
$$
s_n \; = \; s \, + \, (\beta + n)^{\ell} \> \omega_n \,
\sum_{j=0}^{k-1} \, c_j / ( \beta + n)^j
\; , \qquad k,\ell , n \in \N_0 \, .
\tag
$$

\medskip

\Abschnitt Recursive computation of the Levin transformation

\smallskip

\aktTag = 0

Another advantage of Sidi's approach [56] for the derivation of
nondeterminantal expressions for Levin's sequence transformation is
that starting from eq. (7.1-6) Fessler, Ford, and Smith [57] could
derive a recursive scheme which allows a convenient computation of both
the numerator and the denominator of the general Levin transformation,
eq. (7.1-7). In eq. (7.1-6), both numerator and denominator are of the
general form
$$
P_k^{(n)} (\beta) \; = \; \Delta^k \> X_k^{(n)} (\beta) \, ,
\qquad k, n \in \N_0 \, .
\tag
$$

As usual, it is assumed here that the difference operator $\Delta$ as
well as the shift operator $E$, which is defined in eq. (2.4-4), act
only upon $n$ and not upon $k$. Comparison with eq. (7.1-6) shows that
the quantities $X_k^{(n)} (\beta)$ satisfy the following 2-term
recursion in $k$:
$$
X_k^{(n)} (\beta) \; = \; (\beta + n) X_{k-1}^{(n)} (\beta) \, ,
\qquad k \ge 1 \, , \quad n \ge 0 \, .
\tag
$$

The following commutator relationship can be proved by complete
induction with respect to $k$ [57],
$$
\Delta^k (\beta + n) \, - \, (\beta + n) \Delta^k \; = \;
k E \Delta^{k-1} \, .
\tag
$$

Combination of eqs. (2.4-4), (7.2-1), (7.2-2), and (7.2-3) yields:
$$
\beginAligntags
P_k^{(n)} (\beta) \; " = " \;
\{ k E + (\beta + n) \Delta \} \Delta^{k-1} X_{k-1}^{(n)} (\beta)
\hfill \\ \tag
" = " \; \{ k E + (\beta + n) \Delta \} P_{k-1}^{(n)} (\beta)
\hfill \\ \tag
" = " \; (\beta + n +k) P_{k-1}^{(n+1)} (\beta) \, - \,
(\beta + n) P_{k-1}^{(n)} (\beta) \, .
\hfill \\ \tag
\endAligntags
$$

With the help of the 3-term recurrence formula (7.2-6) the numerator as
well as the denominator of the general Levin transformation ${\cal
L}_k^{(n)} (\beta , s_n , \omega_n)$ can be computed for $k \ge 1$.
However, for the sake of numerical stability and in order to make
overflow less likely it is preferable to scale the quantities
$P_k^{(n)} (\beta)$ by defining
$$
L_k^{(n)} (\beta) \; = \;
P_k^{(n)} (\beta) \, / \, (\beta + n + k)^{k-1} \, .
\tag
$$

Inserting this into eq. (7.2-6) yields the following 3-term recurrence
formula for the scaled quantities $L_k^{(n)} (\beta)$:
$$
L_{k+1}^{(n)} (\beta) \; = \; L_k^{(n+1)} (\beta) \, - \, \frac
{ (\beta + n) (\beta + n + k)^{k-1} } { (\beta + n + k + 1)^k } \>
L_k^{(n)} (\beta) \, , \qquad k, n \in \N_0 \, .
\tag
$$

If we use the starting values
$$
L_0^{(n)} (\beta) \; = \; s_n / \omega_n \, ,
\qquad n \in \N_0 \, ,
\tag
$$
the recurrence formula (7.2-8) produces the numerator of the general
Levin transformation, eq. (7.1-7) and if we use the starting values
$$
L_0^{(n)} (\beta) \; = \; 1 / \omega_n \, ,
\qquad n \in \N_0 \, ,
\tag
$$
we obtain the denominator of the general Levin transformation.

With the help of the 3-term recurrence formula (7.2-8) it is also
possible to compute both the numerator and the denominator of the
generalized Levin transformation, eq. (7.1-8). If the starting values
$$
L_0^{(n)} (\beta) \; = \; s_n / [ (\beta + n)^{\ell} \omega_n ]
\, , \qquad \ell, n \in \N_0 \, ,
\tag
$$
are used, eq. (7.2-8) produces the numerator of the transformation
(7.1-8), and the starting values
$$
L_0^{(n)} (\beta) \; = \; 1 / [ (\beta + n)^{\ell} \omega_n ]
\, , \qquad \ell, n \in \N_0 \, ,
\tag
$$
produce the denominator of the transformation (7.1-8).

The 3-term recurrence formula (7.2-6) was according to the knowledge of
the author first published by Longman [58]. However, Longman's
derivation of the recurrence formula (7.2-6) is based upon Sister
Celine's technique [59] and not on properties of the difference
operator $\Delta$ as the derivation by Fessler, Ford, and Smith [57].

\medskip

\Abschnitt Remainder estimates for the Levin transformation

\smallskip

\aktTag = 0

Until now, we have completely ignored the sequence $\Seq {\omega_n}$ of
remainder estimates and its r{\^ o}le in the process of convergence
acceleration or summation. In order to deal with this question we take
into account that the elements of the model sequence (7.1-1), for which
the general Levin transformation is exact, can be obtained from the
elements of the sequence (7.1-2) by truncating the asymptotic power
series in $1/(\beta + n)$ after the first $k$ terms.

This indicates that the Levin transformation (7.1-7) should work very
well for a given sequence $\Seqn s$ if the sequence $\Seq {\omega_n}$
of remainder estimates is chosen in such a way that $\omega_n$ is
proportional to the dominant term of an asymptotic expansion of the
remainder $r_n$,
$$
r_n \; = \; s_n - s \; = \; \omega_n [ c + O(n^{-1}) ] \, ,
\qquad n \to \infty \, .
\tag
$$

Now, one is confronted with the practical problem of finding such a
sequence $\Seq {\omega_n}$ of remainder estimates for a given sequence
$\Seqn s$. Here, it must be emphasized that a sequence $\Seq
{\omega_n}$ of remainder estimates is not determined uniquely by the
asymptotic condition (7.3-1). Consequently, it should at least in
principle always be possible to find a variety of different sequences
of remainder estimates which all satisfy eq. (7.3-1).

In some exceptional cases it is possible to derive explicit analytical
expressions for the remainder estimates $\omega_n$ which satisfy eq.
(7.3-1) -- for instance if the $s_n$ are partial sums of a series and
if the series terms $a_k$ have a sufficiently simple analytical
structure. If such an explicit expression for $\omega_n$ is used in eq.
(7.1-7), the general Levin transformation is a {\it linear} sequence
transformation.

However, in most practical applications no information about the
analytical structure of the sequence of remainders $\Seqn r$ will be
available and only the numerical values of a relatively small number of
sequence elements $s_m, s_{m+1}, \ldots , s_{m + \ell}$ will be known.
Consequently, it is necessary to find a way of obtaining the sequence
of remainder estimates $\Seq {\omega_n}$ directly from the numerical
values of the elements of the sequence $\Seqn s$. If such a sequence of
remainder estimates is used in eq. (7.1-7), the Levin transformation is
a {\it nonlinear} sequence transformation because each remainder
estimate $\omega_m$ depends explicitly upon at least one element of
$\Seqn s$.

On the basis of purely heuristic arguments Levin [28] suggested for
sequences of partial sums
$$
s_n \; = \; \sum_{\nu = 0}^n a_{\nu} \, , \qquad n \in \N_0 \, ,
\tag
$$
some simple remainder estimates which according to experience
nevertheless work remarkably well. In the case of logarithmic
convergence, i.e., if the elements of the sequence of partial sums
$s_n$ satisfy
$$
\lim_{n \to \infty} \, \frac
{s_{n+1} - s} {s_n - s} \; = \; 1 \, ,
\tag
$$
Levin [28] suggested the remainder estimate
$$
\omega_n \; = \; (\beta + n) a_n \, , \qquad n \in \N_0 \, .
\tag
$$

The use of this remainder estimate in eq. (7.1-7) yields Levin's $u$
transformation:
$$
u_{k}^{(n)} (\beta , s_n) \; = \;
\frac
{\displaystyle
\sum_{j=0}^{k} \; ( - 1)^{j} \; \binom {k} {j} \;
\frac {(\beta + n +j )^{k-2}} {(\beta + n + k )^{k-1}} \;
\frac {s_{n+j}} {a_{n+j}} }
{\displaystyle
\sum_{j=0}^{k} \; ( - 1)^{j} \; \binom {k} {j} \;
\frac {(\beta + n +j )^{k-2}} {(\beta + n + k )^{k-1}} \;
\frac {1} {a_{n+j}} }
\; .
\tag
$$

In the case of alternating series Levin [28] suggested the remainder
estimate
$$
\omega_n \; = \; a_n \, , \qquad n \in \N_0 \, .
\tag
$$

This gives Levin's $t$ transformation:
$$
t_{k}^{(n)} (\beta , s_n) \; = \;
\frac
{\displaystyle
\sum_{j=0}^{k} \; ( - 1)^{j} \; \binom {k} {j} \;
\frac {(\beta + n +j )^{k-1}} {(\beta + n + k )^{k-1}} \;
\frac {s_{n+j}} {a_{n+j}} }
{\displaystyle
\sum_{j=0}^{k} \; ( - 1)^{j} \; \binom {k} {j} \;
\frac {(\beta + n +j )^{k-1}} {(\beta + n + k )^{k-1}} \;
\frac {1} {a_{n+j}} }
\; .
\tag
$$

However, Smith and Ford [29] rightly remarked that the best simple
remainder estimate for a convergent series with strictly alternating
terms $a_{\nu}$ would be
$$
\omega_n \; = \; a_{n+1} \, , \qquad n \in \N_0 \, .
\tag
$$

Using this in eq. (7.1-7) gives Smith and Ford's [29] modification of
Levin's $t$ transformation:
$$
d_{k}^{(n)} (\beta , s_n) \; = \;
\frac
{\displaystyle
\sum_{j=0}^{k} \; ( - 1)^{j} \; \binom {k} {j} \;
\frac {(\beta + n +j )^{k-1}} {(\beta + n + k )^{k-1}} \;
\frac {s_{n+j}} {a_{n+j+1}} }
{\displaystyle
\sum_{j=0}^{k} \; ( - 1)^{j} \; \binom {k} {j} \;
\frac {(\beta + n +j )^{k-1}} {(\beta + n + k )^{k-1}} \;
\frac {1} {a_{n+j+1}} }
\; .
\tag
$$

As a third simple remainder estimate Levin [28] suggested
$$
\omega_n \; = \; \frac
{a_n a_{n+1}} {a_n - a_{n+1}} \, , \qquad n \in \N_0 \, .
\tag
$$

Comparison with eq. (5.1-6) shows that this remainder estimate is based
upon Aitken's $\Delta^2$ process. It gives Levin's $v$ transformation:
$$
v_{k}^{(n)} (\beta , s_n) \; = \;
\frac
{\displaystyle
\sum_{j=0}^{k} \; ( - 1)^{j} \; \binom {k} {j} \;
\frac {(\beta + n +j )^{k-1}} {(\beta + n + k )^{k-1}} \;
\frac {a_{n+j} - a_{n+j+1}} {a_{n+j} a_{n+j+1}} s_{n+j} }
{\displaystyle
\sum_{j=0}^{k} \; ( - 1)^{j} \; \binom {k} {j} \;
\frac {(\beta + n +j )^{k-1}} {(\beta + n + k )^{k-1}} \;
\frac {a_{n+j} - a_{n+j+1}} {a_{n+j} a_{n+j+1}} }
\; .
\tag
$$

The remainder estimates (7.3-4), (7.3-6), (7.3-8), and (7.3-10) can
also be used if the sequence $\Seqn s$, which is to be transformed, is
not a sequence of partial sums. It is only necessary to replace in eqs.
(7.3-5), (7.3-7), (7.3-9), and (7.3-11) $a_0$ by $s_0$ and $a_n$ with
$n \ge 1$ by $\Delta s_{n-1}$.

Levin's remainder estimates (7.3-4), (7.3-6), and (7.3-10) as well as
Smith and Ford's modification (7.3-8) were derived using simple
heuristic arguments. However, experience indicates that these remainder
estimates nevertheless give rise to very powerful sequence
transformations [29,30,57,60-64].

In some cases a more rigorous derivation of the remainder estimates
mentioned above can be given. For instance, in Wimp's book (see p. 19
of ref. [23]) it is shown that if the terms $a_n$ of a series satisfy
$$
a_n \; \sim \; \lambda^n n^{\Theta} \left\{ \alpha_0 +
\frac {\alpha_1}{n} + \frac {\alpha_2}{n^2} + \cdots \; \right\}
\, , \qquad n \to \infty \, ,
\tag
$$
with $\alpha_0 \ne 0$, then there exist constants $\beta_j$ and
$\gamma_j$ such that the remainders $r_n$ of the partial sums of this
series satisfy
$$
r_n \; \sim \; \frac {\lambda^{n+1} n^{\Theta}} {\lambda - 1}
\left\{ \alpha_0 + \frac {\beta_1}{n} + \frac {\beta_2}{n^2}
+ \cdots \; \right\} \, , \qquad n \to \infty \, ,
\tag
$$
if $| \lambda | < 1$, and
$$
r_n \; \sim \; - \frac {n^{\Theta + 1}} {\Theta + 1}
\left\{ \alpha_0 + \frac {\gamma_1}{n} + \frac {\gamma_2}{n^2} + \cdots
\; \right\} \, , \qquad n \to \infty \, ,
\tag
$$
if $\lambda = 1$ and $Re(\Theta) < 1$.

A comparison of eqs. (7.3-12) and (7.3-13) gives us essentially the
remainder estimate (7.3-6) which leads to the $t$ transformation, eq.
(7.3-7), and a comparison of eqs. (7.3-12) and (7.3-14) gives us
essentially the remainder estimate (7.3-4) which leads to the $u$
transformation, eq. (7.3-5).

If we replace in the model sequence (6.1-6), for which the Richardson
extrapolation scheme, eq. (6.1-5), is exact, $x_n$ by $1/(n+\beta)$, we
obtain the model sequence
$$
s_n \; = \; s \, + \, \sum_{j=0}^{k-1} \, c_j / (n+\beta)^{j+1}
\, , \qquad k, n \in \N_0 \, .
\tag
$$

This model sequence can be obtained from the model sequence (7.1-2),
for which Levin's sequence transformation, eq. (7.1-7) is exact, by
choosing $\omega_n = 1/(n+\beta)$. Hence, if we define
$$
\Lambda_k^{(n)} (\beta, s_n) \; = \; {\cal L}_k^{(n)} (\beta, s_n ,
1/(n + \beta)) \, , \qquad k, n \in \N_0 \, ,
\tag
$$
we see that that the sequence transformation $\Lambda_k^{(n)} (\beta,
s_n)$ is obviously exact for the model sequence (7.3-15). If we now use
eq. (7.1-6), we see that the transformation $\Lambda_k^{(n)} (\beta,
s_n)$ satisfies
$$
\Lambda_k^{(n)} (\beta, s_n) \; = \; \frac
{\Delta^k \{ (n+\beta)^k \> s_n\} }
{\Delta^k (n+\beta)^k }
\, , \qquad k, n \in \N_0 \, .
\tag
$$

The denominator in eq. (7.3-17) can be expressed in closed form. If we
use (see p. 4 of ref. [65])
$$
\sum_{j=0}^k \, (-1)^j \, \binom {k} {j} \, (\beta + n + j)^k
\; = \; (-1)^k \, k ! \, , \qquad k, n \in \N_0 \, ,
\tag
$$
together with eq. (2.4-8), we find:
$$
\Delta^k (n+\beta)^k \; = \; k! \, .
\tag
$$

Combination of eqs. (7.3-17) and (7.3-19) with eq. (2.4-8) gives us an
explicit expression for this sequence transformation:
$$
\Lambda_k^{(n)} (\beta, s_n) \; = \; (-1)^k \> \sum_{j=0}^k \>
(-1)^j \> \frac {(\beta+n+j)^k} {j! \> (k-j)!} \> s_{n+j} \, ,
\qquad k, n \in \N_0 \, .
\tag
$$

If we combine eq. (7.2-6) with eqs. (7.3-17) and (7.3-19), we can
derive the following recursive scheme for the sequence transformation
$\Lambda_k^{(n)} (\beta, s_n)$:
$$
\beginAligntags
" \Lambda_0^{(n)} (\beta, s_n)\; " = \; " s_n
\, , \qquad n \in \N_0 \, ,
\erhoehe\aktTag \\ \tag*{\tagnr a}
" \Lambda_{k+1}^{(n)} (\beta, s_n) \; " = \;
" \Lambda_k^{(n+1)} (\beta, s_{n+1}) \, + \,
\frac {\beta+n} {k+1} \> \Delta \Lambda_k^{(n)} (\beta, s_n)
\, , \qquad k, n \in \N_0 \, . \\ \tag*{\tagform\aktTagnr b}
\endAligntags
$$

This recursive scheme for the sequence transformation $\Lambda_k^{(n)}
(\beta, s_n)$ corresponds to the special case $x_n = 1 / (\beta + n)$
in the recursive scheme (6.1-5) which computes the Richardson
extrapolation scheme with arbitrary interpolation points $\Seqn x$.

A more complete discussion of the properties of the linear but
nonregular sequence transformation $\Lambda_k^{(n)} (\beta, s_n)$ can
be found in articles by Salzer [66, 67], Salzer and Kimbro [68], and
Wimp [69] as well as in Wimp's book (see pp. 35 - 38 of ref. [23]).

The sequence transformations $u_k^{(n)} (\beta , s_n)$, eq. (7.3-5),
and $t_k^{(n)} (\beta , s_n)$, eq. (7.3-7), require the sequence
elements $s_{n-1}, s_n, s_{n+1}, \ldots, s_{n+k}$ for their
computation, whereas $d_k^{(n)} (\beta , s_n)$, eq. (7.3-9), requires
the sequence elements $s_n, s_{n+1}, \ldots , s_{n+k+1}$. Hence, they
are all transformations of order $k+1$. The sequence transformation
$v_k^{(n)} (\beta , s_n)$, eq. (7.3-11), requires the sequence elements
$s_{n-1}, s_n, s_{n+1}, \ldots , s_{n+k+1}$ which implies that it is a
transformation of order $k+2$. The linear sequence transformation
$\Lambda_k^{(n)} (\beta, s_n)$, eq. (7.3-20), requires the sequence
elements $s_n, s_{n+1}, \ldots , s_{n+k}$, i.e., it is a transformation
of order $k$.

The situation is somewhat different if the transforms with superscript
$n = 0$ are computed because then $u_k^{(0)} (\beta , s_0)$ and
$t_k^{(0)} (\beta , s_0)$ are transformations of order $k$, whereas
$d_k^{(0)} (\beta , s_0)$ and $v_k^{(0)} (\beta , s_0)$ are
transformations of order $k+1$.

\medskip

\Abschnitt Sidi's generalization of Levin's sequence transformation

\smallskip

\aktTag = 0

As discussed in section 6.1, the Richardson extrapolation scheme, eq.
(6.1-5), is exact for model sequences of the following type:
$$
s_n \; = \; s \, + \, \sum_{j=0}^{k-1} \; c_j x_n^{j+1} \, ,
\qquad k, n \in \N_0 \, .
\tag
$$

The interpolation points $\Seqn x$ have to satisfy eq. (6.1-4). Very
natural interpolation points for the Richardson extrapolation scheme
are $x_n = 1/(n+\beta)$ with $\beta > 0$. If this set of extrapolation
points is used in eq. (7.4-1) we obtain the model sequence (7.3-15) for
which the sequence transformation $\Lambda_k^{(n)} (\beta, s_n)$, eq.
(7.3-17), is exact. This implies that the Richardson extrapolation
scheme (6.1-5) with the interpolation points $x_n = 1/(n+\beta)$ is a
special case of Levin's sequence transformation with $\omega_n = 1/(n +
\beta)$.

However, the Richardson extrapolation scheme is also in some sense more
general than Levin's sequence transformation since the interpolation
points $\Seqn x$ only have to satisfy eq. (6.1-4) but otherwise they
are completely arbitrary. In Levin's sequence transformation it is
tacitly assumed that the interpolation points $\Seqn x$ always satisfy
$x_n = 1/(n+\beta)$ with $\beta > 0$.

Now, one can try to construct a sequence transformation which combines
the advantageous features of the Levin transformation and the
Richardson extrapolation scheme. This was accomplished by Sidi [70] who
constructed a sequence transformation, which he called {\it generalized
Richardson extrapolation process}, on the basis of the following model
sequence:
$$
s_n \; = \; s \, + \, \omega_n \, \sum_{j=0}^{k-1} \, c_j x_n^j \, ,
\qquad k, n \in \N_0 \, .
\tag
$$

This model sequence combines the features of the model sequence (7.4-1)
for the Richardson extrapolation scheme, eq. (6.1-5), since it contains
arbitrary interpolation points $\Seqn x$, and of the model sequence
(7.1-1) for the Levin transformation, eq. (7.1-7), since it contains
arbitrary remainder estimates $\Seq {\omega_n}$.

For the construction of such a sequence transformation, which is exact
for the above model sequence, eq. (7.4-2) is rewritten in the following
way:
$$
[s_n - s] / \omega_n \; = \; \sum_{j=0}^{k-1} \, c_j x_n^j \, .
\tag
$$

Obviously, the right-hand side of eq. (7.4-3) is a polynomial of degree
$k-1$ in the variable $x_n$. Hence, the limit or antilimit $s$ of this
sequence can be determined if a linear operator can be found which
annihilates the polynomial on the right-hand side of eq. (7.4-3).

This annihilation of a polynomial can be accomplished with the help of
{\it divided differences} which for instance occur in Newton's
interpolation formula. A discussion of divided differences and their
properties can be found in any book on numerical analysis or also in
books on finite differences such as N\"orlund [71] or Milne-Thomson [72].

Let $\Seqn x$ with $n \in \N_0$ be a set of distinct interpolation
points. Then the divided differences of orders $0, 1, \ldots, k, k+1,
\ldots $ of a given function $f$ are defined recursively by the
relations
$$
\beginAligntags
" f [x_n] \> " = " \> f (x_n)\; ,
\hfill \erhoehe\aktTag \\ \tag*{\tagnr a}
" f [x_n, \ldots, x_{n+k+1}] \> " = " \> \frac
{f [x_{n+1}, \ldots, x_{n+k+1}] - f [x_n, \ldots, x_{n+k}] }
{ x_{n+k+1} - x_n } \, , \qquad k, n \in \N_0 \> . \\
\tag*{\tagform\aktTagnr b}
\endAligntags
$$

The divided differences $f [x_n, \ldots , x_{n+k}]$ can also be
expressed in closed form:
$$
f [x_n, \ldots , x_{n+k}] \; = \; \sum_{j=0}^k \, f (x_{n+j}) \,
\prod^k \Sb{i=0 \\ i\ne j} \, \frac {1} {x_{n+j} - x_{n+i} } \, ,
\qquad k, n \in \N_0 \, .
\tag
$$

It follows either from this expression or from the recursive scheme
(7.4-4) that the divided differences $f [x_n, \ldots ,x_{n+k}]$ are
linear functions of the initial values $f (x_n)$, $\ldots$, $f
(x_{n+k})$. In addition, it can be shown that if $p_m (x)$ is a
polynomial of degree $m$ in $x$,
$$
p_m (x) \; = \; c_0 + c_1 x + c_2 x^2 + \cdots + c_m x^m \, ,
\tag
$$
then all its divided differences with $k>m$ satisfy
$$
p_m [x_n, \ldots, x_{n+k}] \, = \, 0 \, , \qquad k > m \, .
\tag
$$

The divided differences $f [x_n, \ldots, x_{n+k}]$ with $k, n \in \N_0$
can be identified with the elements ${\it \Delta}_k^{(n)} (f)$ of a
2-dimensional rectangular array. With this convention, the recursive
scheme (7.4-4) for the computation of the divided differences of a
given function $f$ can be rewritten in the following way:
$$
\beginAligntags
" {\mit \Delta}_0^{(n)} (f) \, " = " \, f (x_n)\, ,
\qquad n \in \N_0 \, ,
\hfill \erhoehe\aktTag \\ \tag*{\tagnr a}
" {\mit \Delta}_{k+1}^{(n)} (f) \, " = " \, \frac
{ {\mit \Delta}_{k}^{(n+1)} (f) - {\mit \Delta}_{k}^{(n)} (f)}
{ x_{n+k+1} - x_n } \; , \qquad k, n \in \N_0 \, . \\
\tag*{\tagform\aktTagnr b}
\endAligntags
$$

If we assume that two functions ${\cal S} (x)$ and ${\it \Omega} (x)$
of a continuous variable $x$ exist, which coincide at the interpolation
points $x_n$ with $s_n$ and $\omega_n$, respectively,
$$
{\cal S} (x_n) \; = \; s_n \, , \qquad
{\it \Omega} (x_n) \; = \; \omega_n \, ,
\tag
$$
and which satisfy
$$
[ {\cal S} (x) - s ] / {\it \Omega} (x) \; = \;
\sum_{j=0}^{k-1} \, c_j x^j \, ,
\tag
$$
we see that we only have to compute the divided difference of order $k$
for the left-hand side of eq. (7.4-10) at the interpolation points
$x_n, \ldots , x_{n+k}$ in order to obtain the limit or antilimit $s$.
Hence, Sidi's generalized Richardson extrapolation process ${\cal
R}_k^{(n)} (s_n, \omega_n, x_n)$ can be defined in terms of divided
differences in the following way:
$$
{\cal R}_k^{(n)} (s_n,\omega_n,x_n) \; = \; \frac
{\{ {\cal S} (x) / {\it \Omega} (x) \} [x_n, \ldots , x_{n+k}] }
{\{ 1 / {\it \Omega} (x) \} [x_n, \ldots , x_{n+k}] }
\, , \qquad k, n \in \N_0 \, .
\tag
$$

It follows from eq. (7.4-8) that both numerator and denominator of this
transformation can be computed with the help of the same 3-term
recurrence formula:
$$
R_{k+1}^{(n)} \, = \, \frac
{ R_{k}^{(n+1)} - R_{k}^{(n)} }
{ x_{n+k+1} - x_n } \; , \qquad k, n \in \N_0 \, .
\tag
$$

If we use in eq. (7.4-12) the starting values
$$
R_0^{(n)} \; = \; s_n / \omega_n \, , \qquad n \in \N_0 \, ,
\tag
$$
we obtain the numerator of Sidi's generalized Richardson extrapolation
process, eq. (7.4-11), and if we use the starting values
$$
R_0^{(n)} \; = \; 1 / \omega_n \, , \qquad n \in \N_0 \, ,
\tag
$$
we obtain the denominator in eq. (7.4-11).

Obviously, the Richardson extrapolation scheme (6.1-5) is a special
case of Sidi's generalized Richardson extrapolation process ${\cal
R}_k^{(n)} (s_n,\omega_n,x_n)$. To see this one only has to specialize
${\it \Omega} (x) \; = \; x$ in eq. (7.4-11) which also implies
$\omega_n \; = \; x_n$.

With this specialization, the Richardson extrapolation scheme (6.1-5)
follows from eqs. (7.4-11) and (7.4-12). It is only necessary to
compute the divided differences for $1/x$ at the interpolation points
which can be done quite easily. The result is (see p. 8 of ref. [72]):
$$
\frac {1} {x} \> [x_n, \ldots , x_{n+k}] \, = \;
\frac {(-1)^k} {x_n \cdots x_{n+k} } \, .
\tag
$$

This implies that for $\omega_n = x_n$ the denominator of the
generalized Richardson extrapolation process ${\cal R}_k^{(n)}
(s_n,\omega_n,x_n)$ can be expressed in closed form. If we now set
$$
R_0^{(n)} \, = \, s_n / x_n \; ,
\tag
$$
and
$$
{\cal N}_k^{(n)} \, = \, (-1)^k x_n \cdots x_{n+k} R_k^{(n)}
\; , \qquad k, n \in \N_0 \, ,
\tag
$$
we immediately obtain from eqs. (7.4-11) and (7.4-12) the Richardson
extrapolation scheme, eq. (6.1-5).

Further generalizations of Sidi's generalized Richardson extrapolation
are possible. For instance, let us consider the following model
sequence:
$$
s_n \; = \; s \, + \, \sum_{\mu = 0}^{m} \, \omega_n^{(\mu)} \;
\sum_{j=0}^{k_{\mu}-1} \, c_j^{(\mu)} \> x_n^j \, ,
\qquad m, n, k_{\mu} \in \N_0 \, .
\tag
$$

This model sequence contains $m+1$ different sets of remainder
estimates $\Seq {\omega_n^{(0)}}$, $\ldots$, $\Seq {\omega_n^{(m)}}$.
If we set $m=0$ in eq. (7.4-18), we obtain the model sequence (7.4-2)
as a special case.

Sidi [73] constructed a sequence transformation which is exact for the
above model sequence. Originally, this sequence transformation was
defined as the ratio of determinants, which is computationally quite
unattractive. But recently, Ford and Sidi [74] could derive an
algorithm which permits a recursive computation of Sidi's sequence
transformation.

From Sidi's very general sequence transformation [73], which is exact
for the model sequence (7.4-18), other sequence transformations can be
obtained by specializing the interpolation points $\Seqn x$. For
instance, in earlier articles Levin and Sidi [75,76] had chosen the
interpolation points according to the rule $x_n = 1/(n+1)$ and had
obtained a generalization of Levin's sequence transformation with
several sets of remainder estimates.

The typical feature of these sequence transformations is that they
involve several sets of remainder estimates $\Seq {\omega_n^{(0)}},
\ldots , \Seq {\omega_n^{(m)}} $. Consequently, it is to be expected
that these sequence transformations should be particularly powerful if
sequences are to be accelerated which are superpositions of sequences
with different convergence types. The numerical examples presented in
the literature [74-76] confirm this opinion and it seems that the
sequence transformations, which are based upon variants of the model
sequence (7.4-18), are able to accelerate convergence even in cases in
which other transformations fail.

The power of these transformations stems from the occurrence of several
different sets of remainder estimates. This is at the same time also
the major disadvantage of these transformations. This may sound
paradoxical. However, one should take into consideration that the
popularity and the success of Levin's sequence transformation, eq.
(7.1-7), is largely due to the fact that the simple remainder estimates
(7.3-4), (7.3-6), (7.3-8), and (7.3-10) work remarkably well even in
situations in which only the numerical values of a relatively small
number of sequence elements $s_m, \ldots , s_{m + \ell}$ are known. If
we would try to use in such a situation a sequence transformation,
which is based upon a variant of the model sequence (7.4-18), we would
first have to find a way of determining numerically more than a single
set of remainder estimates. Unfortunately, no simple and manageable
theory is in sight which would yield more than a single set of
remainder estimates if only the numerical values of a few sequence
elements are known.

Consequently, if little or no information about the analytical
behaviour of the remainders $\Seqn r$ is available and if only a
relatively small number of sequence elements $s_m, \ldots , s_{m +
\ell}$ are known, it is normally not possible to use sequence
transformations, which are based upon a variant of the model sequence
(7.4-18), although they are potentially very powerful. The reason is
that such a sequence transformation has too many degrees of freedom
since it requires the input of $m+1$ different sets of remainder
estimates $ \Seq {\omega_n^{(0)}}, \ldots , \Seq {\omega_n^{(m)}}$ as
well as a set of interpolation points $\Seqn x$.

In such a situation, one is more or less forced to use a simpler and
probably also less efficient sequence transformation which, however,
does not require such a detailed knowledge about the sequence to be
transformed.

\medskip

\Abschnitt Programming the Levin transformation

\smallskip

\aktTag = 0

In this section it will be discussed how the general Levin
transformation (7.1-7) and its variants (7.1-8), (7.3-5), (7.3-7),
(7.3-9), and (7.3-11) can be programmed efficiently. It is a typical
feature of the general Levin transformation ${\cal L}_k^{(n)} (\beta,
s_n, \omega_n)$ and its variants that they can all be expressed as the
ratio of two finite sums and that both the numerator and the
denominator sum can be computed by the same 3-term recursion (7.2-8).

Consequently, a program for the general Levin transformation or any of
its variants has to compute simultaneously the numerator and
denominator sums of the transformation. In the case of the general
Levin transformation a program requires not only the input of the
sequence elements $s_n$, but also the remainder estimates $\omega_n$.
But otherwise, such a program should have essentially the same features
as the subroutines EPSAL and AITKEN, which were described in sections
4.3 and 5.2, respectively.

This means that such a program for the general Levin transformation
should read in the sequence elements $s_0$, $s_1$, $ \ldots$, $s_m$,
$\ldots$ and remainder estimates $\omega_0$, $\omega_1$, $\ldots$,
$\omega_m$, $\ldots$ successively, starting with $s_0$ and $\omega_0$.
After the input of each new pair $s_m$ and $\omega_m$ as many new
numerator and denominator sums of the Levin transformation (7.1-7) or
its variants should be computed as possible.

The elements $L_k^{(n)} (\beta)$, which either represent numerators or
denominators of the general Levin transformation and which are computed
with the help of the 3-term recurrence formula (7.2-8), can be arranged
in a rectangular scheme in such a way that the superscript $n$
indicates the row and the subscript $k$ the column of the 2-dimensional
array:

$$
\matrix{
L_0^{(0)} (\beta) " L_1^{(0)} (\beta) " L_2^{(0)} (\beta) " \ldots "
L_n^{(0)} (\beta) " \ldots \\
L_0^{(1)} (\beta) " L_1^{(1)} (\beta) " L_2^{(1)} (\beta) " \ldots "
L_n^{(1)} (\beta) " \ldots \\
L_0^{(2)} (\beta) " L_1^{(2)} (\beta) " L_2^{(2)} (\beta) " \ldots "
L_n^{(2)} (\beta) " \ldots \\
L_0^{(3)} (\beta) " L_1^{(3)} (\beta) " L_2^{(3)} (\beta) " \ldots "
L_n^{(3)} (\beta) " \ldots \\
\vdots   " \vdots    " \vdots    " \ddots " \vdots    " \ddots \\
L_0^{(n)} (\beta) " L_1^{(n)} (\beta) " L_2^{(n)} (\beta) " \ldots "
L_n^{(n)} (\beta) " \ldots \\
\vdots  " \vdots " \vdots " \ddots " \vdots  " \ddots }
\tag
$$

\medskip

The entries in the first column of the array are the starting values of
the recursion. If the starting values $L_0^{(n)} (\beta)$ are chosen
according to eq. (7.2-9), the 3-term recurrence formula (7.2-8) will
produce a table $L_k^{(n)} (\beta)$ of numerators of the general Levin
transformation, and if the starting values are chosen according to eq.
(7.2-10), a table of denominators will be computed. The 3 elements,
which are connected by the linear 3-term recursion (7.2-8), form a
triangle in the Levin table:

$$
\beginMatrix
\beginFormat &\Formel\links \endFormat
\+" L_k^{(n)} (\beta) \qquad " L_{k+1}^{(n)} (\beta) " \\ [1\jot]
\+" L_k^{(n+1)} (\beta) \qquad " " \\
\endMatrix
\tag
$$

\medskip

This pattern implies that the recursion (7.2-8) has to proceed along
counterdiagonals in the Levin table. Consequently, it is advantageous
to rewrite the recurrence formula (7.2-8) in the following way:
$$
\beginAligntags
" L_{j}^{(n-j)} (\beta) \, = \, L_{j-1}^{(n-j+1)} (\beta)\, - \, \frac
{ (\beta + n - j) (\beta + n - 1)^{j-2} }
{ (\beta + n)^{j-1} }
\> L_{j - 1}^{(n-j)} (\beta) \; , \\
" n \ge 1 \; , \qquad 1 \le j \le n \, . \\
\tag
\endAligntags
$$

It also follows from the triangular structure of this recursion that
the sequence elements $s_0, s_1, \ldots , s_m$ and the remainder
estimates $\omega_0, \omega_1, \ldots , \omega_m$ have to be known for
the computation of all elements $L_{\nu}^{(n - \nu)}$ with $0 \le n \le
m$ and $0 \le \nu \le n$ of the numerator and denominator tables. Since
the transforms with the highest values of the subscript normally give
the best results, our approximation to the limit $s$ of the sequence
$\Seqn s$ to be transformed will be:
$$
\{ s_0, \omega_0; s_1, \omega_1; \ldots ; s_m, \omega_m \} \; \to \;
{\cal L}_m^{(0)} (\beta, s_0, \omega_0) \, ,
\qquad m \in \N_0 \, .
\tag
$$

Essentially this means that we shall use the following sequence of
Levin transforms with minimal superscripts and maximal subscripts as
approximations to the limit $s$:
$$
{\cal L}_0^{(0)} (\beta, s_0, \omega_0), \;
{\cal L}_1^{(0)} (\beta, s_0, \omega_0), \; \ldots , \;
{\cal L}_m^{(0)} (\beta, s_0, \omega_0), \; \ldots \; .
\tag
$$

Because of the triangular structure (7.5-2) of the recurrence formula
(7.5-3) and since the computation proceeds along a counterdiagonal of
the Levin table, a single 1-dimensional array will be sufficient for
the computation of the $L_j^{(n-j)} (\beta)$ which are either numerator
or denominator sums of the general Levin transformation. For that
purpose the $L_j^{(n-j)} (\beta)$ are stored in a 1-dimensional array
$\ell$ in such a way that the superscript coincides with the index of
the corresponding array element:
$$
L_{\nu}^{(n-\nu)} (\beta) \; \to \; \ell (n-\nu) \, ,
\qquad n \ge 0 \, , \quad 0 \le \nu \le n \, .
\tag
$$

With this convention the recurrence formula (7.5-3) can be reformulated
in terms of the elements of the array $\ell$:
$$
\beginAligntags
" \ell (n-j) \, \gets \, \ell (n-j+1) \, - \, \frac
{ (\beta + n - j) (\beta + n - 1)^{j-2} }
{ (\beta + n)^{j-1} } \, \ell (n-j) \; , \\
" n \ge 1 \; , \qquad 1 \le j \le n \, . \\
\tag
\endAligntags
$$

This computational scheme is simpler than the corresponding scheme for
Wynn's $\epsilon$ algorithm, eq. (4.3-10), since no auxiliary variables
are needed here.

Essentially the same approach works also in the case of the Richardson
extrapolation scheme, eq. (6.1-5), or in the case of Sidi's generalized
Richardson extrapolation process, eq. (7.4-11). For instance, the
recurrence formula (7.4-12) for Sidi's generalized Richardson
extrapolation process can be rewritten in the following way:
$$
R_j^{(n-j)} \; = \; \frac
{R_{j-1}^{(n-j+1)} - R_{j-1}^{(n-j)} } {x_n - x_{n-j}} \, ,
\qquad n \ge 1 \, , \quad 1 \le j \le n \, .
\tag
$$

If the elements $R_j^{(n-j)}$ and the interpolation points $x_n$ are
stored in 1-dimensional arrays $r$ and $\xi$ according to the rules
$$
\beginAligntags
" R_j^{(n-j)} \; " \to \; " r (n-j) \, ,
\qquad " n \ge 0 \, , \quad 0 \le j \le n \, ,
\\ \tag
" x_n \; " \to \; " \xi (n) \, , \qquad " n \ge 0 \, ,
\\ \tag
\endAligntags
$$
the recurrence formula (7.5-8) can be reformulated in terms of the
elements of the arrays $r$ and $\xi$:
$$
r (n-j) \, \gets \, \frac {r (n-j+1) \, - \, r (n-j)}
{\xi (n) \, - \, \xi (n-j)} \, ,
\qquad n \ge 1 \; , \quad 1 \le j \le n \, .
\tag
$$

Similarly, the Richardson extrapolation scheme, eq. (6.1-5), can be
rewritten in the following way:
$$
\beginAligntags
" {\cal N}_0^{(n)} (s_n , x_n ) \; " = \; " s_n \, ,
\qquad n \ge 0 \, , \erhoehe\aktTag \\ \tag*{\tagnr a}
" {\cal N}_j^{(n-j)} (s_{n-j} , x_{n-j}) \; " = \; " \frac
{x_{n-j} {\cal N}_{j-1}^{(n-j+1)} (s_{n-j+1} , x_{n-j+1})
- x_n {\cal N}_{j-1}^{(n-j)} (s_{n-j} , x_{n-j}) }
{x_{n-j} - x_n } \, , \\
" n \ge 1 \, , \qquad 1 \le j \le n \, . \span\omit \span\omit
\\ \tag*{\tagform\aktTagnr b}
\endAligntags
$$

If the elements ${\cal N}_j^{(n-j)} (s_{n-j} , x_{n-j})$ and the
interpolation points $x_n$ are stored in 1-dimensional arrays ${\mit
N}$ and $\xi$ according to the rules
$$
\beginAligntags
" {\cal N}_j^{(n-j)} (s_{n-j} , x_{n-j}) \; " \to \;
" {\mit N} (n-j) \, ,
\qquad " n \ge 0 \, , \quad 0 \le j \le n \, ,
\erhoehe\aktTag \\ \tag*{\tagnr a}
" x_n \; " \to \; " \xi (n) \, , \qquad " n \ge 0 \, ,
\\ \tag*{\tagform\aktTagnr b}
\endAligntags
$$
the recurrence formula (7.5-12) can be reformulated in terms of the
elements of the arrays ${\mit N}$ and $\xi$:
$$
\beginAligntags
" {\mit N} (n) \, " \gets \; " s_n \, , \qquad n \ge 0 \, ,
\erhoehe\aktTag \\ \tag*{\tagnr a}
" {\mit N} (n-j) \, " \gets \, " \frac
{ \xi (n-j) {\mit N} (n-j+1) \, - \, \xi (n) {\mit N} (n-j)}
{\xi (n-j) \, - \, \xi (n)} \, ,
\qquad n \ge 1 \; , \quad 1 \le j \le n \, .
\\ \tag*{\tagform\aktTagnr b}
\endAligntags
$$

Obviously, the computational scheme (7.5-11) for Sidi's generalized
Richardson extrapolation process and the computational scheme (7.5-14)
for the Richardson extrapolation process are structurally identical
with the computational scheme (7.5-7) for Levin's general sequence
transformation. This implies that programs for the Richardson
extrapolation process and Sidi's generalization of the Richardson
extrapolation process would have the same features as a program for
Levin's general sequence transformation.

A program which computes the general Levin transformation (7.1-7) or
any of its variants has to take precautions against an exact or
approximate vanishing of the denominator sum. Again, this can be
accomplished by introducing two variables HUGE and TINY which have
values that are close to but not identical with the largest and
smallest floating point number representable on the computer. If the
denominator sum of the transform ${\cal L}_m^{(0)} (\beta , s_0 ,
\omega_0)$ is smaller in magnitude than TINY, then ${\cal L}_m^{(0)}
(\beta , s_0 , \omega_0)$ will be set equal to HUGE. This check is only
necessary if the approximation to the limit according to eq. (7.5-5) is
computed.

The following FORTRAN 77 subroutine GLEVIN computes the numerator and
denominator sum of the general Levin transformation ${\cal L}_m^{(0)}
(\beta , s_0 , \omega_0)$ with the help of the recurrence formula
(7.5-3) in two 1-dimensional arrays ARUP and ARLO. It is safeguarded
against an exact or approximate vanishing of the denominator sum by
using the variables HUGE and TINY described above. The sequence
elements $s_n$ and the remainder estimates $\omega_n$ with $n = 0, 1,
2, \ldots $ have to be computed in a DO loop in the calling program.
Whenever a new pair $s_n$ and $\omega_n$ is computed in the outer DO
loop this subroutine GLEVIN has to be called again and a new string of
transforms ${\cal L}_j^{(n-j)} (\beta , s_{n-j} , \omega_{n-j})$ with
$0 \le j \le n$ is computed. The new sequence element $s_n$ is read in
via the variable SOFN, and the new remainder estimate $\omega_n$ is
read in via the variable ROFN. The approximation to the limit, which is
given by the ratio ARUP(0) / ARLO(0), is returned via the variable
ESTLIM.

Again, it should be noted that GLEVIN only calculates the approximation
to the limit according to eq. (7.5-4). The convergence of the whole
process has to be analyzed in the calling program.

Finally, it should be noted that the description of a FORTRAN IV
program, which computes Levin's $u$ transformation, eq. (7.3-5), with
$\beta = 1$, can be found in ref. [57].

\bigskip

\listenvon{glevin.for}

\medskip

\endAbschnittsebene

\neueSeite

\Abschnitt Sequence transformations based upon factorial series

\vskip - 2 \jot

\beginAbschnittsebene

\medskip

\Abschnitt Factorial series

\smallskip

\aktTag = 0

In extensive numerical studies performed by Smith and Ford [29,30] and
also in other articles [57, 60 - 64] it was demonstrated that the
general Levin transformation (7.1-7) and its variants (7.3-5), (7.3-7),
(7.3-9), and (7.3-11) are remarkably powerful sequence transformations.
Consequently, if one tries to derive alternative sequence
transformations it should definitely be worthwhile to try to retain as
many of the advantageous features of the Levin transformation as
possible.

It is the conviction of the author that the power of the Levin
transformation is due to the fact that a sequence $\Seq {\omega_n}$ of
remainder estimates is explicitly included in the transformation. This
is not necessarily an advantage because if the remainder estimates
$\omega_n$, which are used, are poor approximants of the actual
remainders $r_n$, the Levin transformation will lose much of its
efficiency. However, if the remainder estimates are good approximants,
it is likely that the Levin transformation will produce excellent
results.

In section 7.4, it was shown that Sidi's generalized Richardson
extrapolation process (7.4-11) -- which is also a generalization of
Levin's sequence transformation -- is by construction exact if the
remainders $r_n$ of the sequence to be transformed can be written as an
remainder estimate $\omega_n$ multiplied by a polynomial of degree
$k-1$ in $x_n$,
$$
r_n \; = \; \omega_n \, \sum_{j=0}^{k-1} \, c_j x_n^j \, ,
\qquad k, n \in \N_0 \, .
\tag
$$

The interpolation points $x_n$ have to satisfy eq. (6.1-4) which means
that they have to approach zero as $n \to \infty$. If we choose in eq.
(8.1-1) $x_n = 1/(n+\beta)$ we obtain the remainder of the model
sequence (7.1-1) which is the basis for the construction of Levin's
sequence transformation. Model sequences with remainders of the above
type can be viewed to be finite approximations of sequence elements
$s_n$ which can be written as Poincar\'e-type asymptotic expansions with
respect to the asymptotic sequence $\Seq {\omega_n x_n^j}$ with $n,j
\in \N_0$,
$$
s_n \; \sim \; s \, + \, \omega_n \,
\sum_{j=0}^{\infty} \, c_j x_n^j \, , \qquad n \to \infty \, .
\tag
$$

Essentially this means that the sequence of remainder estimates $\Seq
{\omega_n}$ should be chosen in such a way that the ratio
$(s_n-s)/\omega_n$ can be written as an asymptotic power series in the
interpolation points $\Seqn x$,
$$
(s_n - s) / \omega_n \; \sim \;
\sum_{j=0}^{\infty} \, c_j x_n^j \, , \qquad n \to \infty \, .
\tag
$$

If one tries to construct alternative sequence transformations, which
also incorporate explicit remainder estimates via the auxiliary
sequence $\Seq {\omega_n}$, the simplest approach would be to replace
the asymptotic power series on the right-hand side of eq. (8.1-3) by
some other kind of expansion. This means that in eq. (8.1-3) instead of
the powers $\Seq {x_n^j}$ some other asymptotic sequence $\Seq {\phi_j
(n)}$ with $n,j \in \N_0$ would have to be used. Consequently, it would
be necessary to construct a transformation which is exact for the
following class of model sequences:
$$
s_n \; = \; s \, + \, \omega_n \, \sum_{j=0}^{k-1}
\, c_j \phi_j (n) \, , \qquad k, n \in \N_0 \, .
\tag
$$

Such a transformation would also be a special case of the general
extrapolation algorithm $E_k (s_n)$, eq. (3.3-2), which was introduced
by Brezinski [31] and H{\aa}vie [32]. This follows immediately if $f_j
(n)$ in eq. (3.3-1) is replaced by $\omega_n \phi_j (n)$.

In principle, every set $\Seq {\phi_j (n)}$ of functions of $n$ could
be used in eq. (8.1-4) which satisfies
$$
\beginAligntags
" \phi_0 (n) \; " = \; " 1 \, , " " n \in \N_0 \, ,
\erhoehe\aktTag \\ \tag*{\tagnr a}
" \phi_{j+1} (n) \; " = \; " o ( \phi_j (n) ) \, ,
\qquad " " j \in \N_0 \, , \quad n \to \infty \, .
\\ \tag*{\tagform\aktTagnr b}
\endAligntags
$$

However, such a minimal requirement on the set $\Seq {\phi_j (n)}$
would not suffice to make a new transformation practically useful, let
alone to give it any advantage over already existing transformations.

In order to be practically useful, a new sequence transformation should
produce excellent numerical results in convergence acceleration and
summation processes. Preferably, it should be as good as the Levin
transformation or maybe even better. However, this would not be enough.
Since the evaluation of large order determinants, as they for instance
occur in eq. (3.3-2), is computationally very unattractive, a
comparatively simple recursive scheme, which allows a fast and reliable
computation of the transformation, would also be of considerable
importance. The derivation of an explicit expression of the type of eq.
(7.1-7) for the new transformation would also be desirable since this
would give us a better chance of understanding the mechanism as well as
the shortcomings of the new transformation.

It is not a simple task to find an alternative asymptotic sequence
other than powers $\Seq {x_n^j}$ with $n,j \in \N_0$ which leads to a
sequence transformation satisfying the requirements mentioned above.
However, it will become clear later that a new class of sequence
transformations with most of the advantageous features of the Levin
transformation and some new ones can be derived quite easily if it is
assumed that the ratio $(s_n - s)/\omega_n$ is expressed as a factorial
series and not as an asymptotic power series as in eq. (8.1-3).

Let $\Omega (z)$ be a function which vanishes as $| z | \to \infty$.
Then, a factorial series for $\Omega (z)$ is an expansion of the
following type,
$$
\Omega (z) \; = \; \frac {c_0} {z} \, + \, \frac {c_1} {z (z+1)}
\, + \, \frac {c_2} {z (z+1) (z+2)} \, + \, \cdots \;
= \; \sum_{\nu=0}^{\infty} \frac {c_{\nu}} {(z)_{\nu+1}} \, .
\tag
$$

Here, $(z)_{\nu+1}$ is a Pochhammer symbol which is commonly defined as
the ratio of two gamma functions (see p. 3 of ref. [34]),
$$
(z)_{\nu+1} \; = \; \Gamma (z+\nu+1) / \Gamma (z) \; = \;
z (z+1) \ldots (z+\nu) \, , \qquad \nu \in \N_0 \, .
\tag
$$

Factorial series have a long tradition in mathematics. For instance, a
large part of Stirling's book [5], which was published in 1730, deals
with factorial series. In the nineteenth century the theory of
factorial series was developed and refined by a variety of authors. A
fairly complete survey of the older literature on this subject can be
found in books by Nielsen [77] and N\"orlund [71]. In these two books
good treatments of the fundamental properties of factorial series can
be found.

Factorial series have a remarkable property which will also be utilized
quite profitably in this report: It is extremely simple to apply higher
powers of the difference operator $\Delta$ to a factorial series.
Consequently, factorial series play a similar r{\^ o}le in the theory
of difference equations as power series in the theory of differential
equations. This explains why factorial series were often treated in the
classical literature on finite differences, e.g., in books by N\"orlund
[71,78] and Milne-Thomson [72].

Quite interesting in the context of convergence acceleration and
summation is also Borel's book on divergent series [79] in which the
connection between factorial series and summability is emphasized.

However, it seems that in recent years mathematicians have lost
interest in factorial series. This can be concluded from the fact that
only quite rarely references dealing with factorial series can be found
in the more modern mathematical literature. Notable exceptions are a
book by Wasow [80], which contains a chapter on factorial series, and
an article by Iseki and Iseki [81] on remainder estimates of truncated
factorial series. In the opinion of the author this declining interest
in factorial series is quite deplorable because the numerical potential
of factorial series has not yet been fully exploited.

The fact, that the argument $z$ of a factorial series occurs in
Pochhammer symbols and not in the form of inverse powers as in
asymptotic power series, has some far-reaching consequences for the
convergence properties of factorial series.

A power series converges in the interior of a circle which may coincide
with the whole complex plane $\C$ or which may shrink to a single point
as in the case of divergent asymptotic series. However, if a factorial
series converges then according to Landau [82] it converges in a
half-plane. This means that if a factorial series converges for some
$z_0 \in \C$ it also converges with the possible exception of the
points $z = 0, -1, -2,\ldots$ for all $z \in \C$ with $Re(z) >
Re(z_0)$.

The different convergence properties of power series and factorial
series are demonstrated quite drastically by the following two infinite
series which both have the same numerical coefficients $c_m = (-1)^m
m!$:
$$
\beginAligntags
" \frac {1} {x} \, - \, \frac {1} {x^2} \, + \, \frac
{2} {x^3} \, - \, \frac {6} {x^4} \, + \cdots \; " = \;
" \sum_{m=0}^{\infty} \frac {(-1)^m m!} {x^{m+1}} \, ,
\\ \tag
" \frac {1} {x} \, - \, \frac {1} {(x)_2} \, + \, \frac
{2} {(x)_3} \, - \, \frac {6} {(x)_4} \, + \cdots \;
" = \; " \sum_{m=0}^{\infty} \frac {(-1)^m m!} {(x)_{m+1}} \, .
\\ \tag
\endAligntags
$$

The power series diverges for all finite values of $x \in \R$, whereas
the factorial series converges for all $x > 0$.

Because of the different convergence properties of factorial and power
series it may happen that a given function $\Omega (z)$, which
possesses a representation as a divergent asymptotic power series,
$$
\Omega (z) \; \sim \; \frac {c'_0}{z} \, + \, \frac {c'_1}{z^2}
\, + \, \frac {c'_3}{z^3} \, + \, \cdots \; ,
\qquad z \to \infty \, ,
\tag
$$
possesses also a representation as a convergent factorial series
according to eq. (8.1-6).

The algebraic processes, by means of which the two series expansions
(8.1-6) and (8.1-10) can be transformed into each other, were already
described by Stirling [5] in 1730. A more modern description of
Stirling's method can be found in Nielsen's book (see pp. 272 - 282 of
ref. [77]). A detailed investigation of the problems associated with
the transformation of an asymptotic series into a convergent factorial
series can be found in a long article by Watson [83].

\medskip

\Abschnitt A factorial series analogue of Levin's transformation

\smallskip

\aktTag = 0

The following model sequence will be the basis for the new class of
sequence transformations which will be discussed in this section:
$$
s_n \; = \; s \, + \, \omega_n \, \sum_{j=0}^{k-1} \, c_j / (n +
\beta)_j \, , \qquad k, n \in \N_0 \, .
\tag
$$

This sequence is formally almost identical with the model sequence
(7.1-1) which is the basis of the Levin transformation. The only
difference is that the powers $(n+\beta)^j$ in eq. (7.1-1) are replaced
by Pochhammer symbols $(n+\beta)_j$. Concerning the sequence $\Seq
{\omega_n}$ of remainder estimates it is again assumed that the
$\omega_n$ are known functions of $n$ which are different from zero and
distinct for all finite values of $n$. But otherwise, the $\omega_n$
are in principle completely arbitrary.

The parameter $\beta$ in eq. (8.2-1) is subject to the restriction that
the Pochhammer symbols $(n+\beta)_j$ must not be zero for all $n,j \in
\N_0$. This is certainly guaranteed if $\beta$ is not a negative
integer or zero. However, the elements of the model sequence (8.2-1)
will serve as finite approximations to factorial series of the
following kind:
$$
s_n \; \sim \; s \, + \, \omega_n \, \sum_{j=0}^{\infty} \, c_j / (n +
\beta)_j \, , \qquad n \to \infty \, .
\tag
$$

In expansions of that kind negative values of $\beta$ will lead to
different signs of the terms of this factorial series if either $n +
\beta < 0$ or $n + \beta > 0$ holds. Since the sign pattern of the
terms of such a factorial series should not change as $n$ increases, we
see that as in the case of the Levin transformation the additional
restriction $\beta > 0$ is necessary. But otherwise, $\beta$ is in
principle completely arbitrary.

In eq. (8.2-1) there occur $k+1$ unknown quantities, the limit or
antilimit $s$ and the $k$ linear coefficients $c_0, \ldots , c_{k-1}$.
Hence, if $k+1$ sequence elements $s_n, \ldots , s_{n+k}$ are known,
the sequence transformation ${\cal S}_k^{(n)} (\beta, s_n, \omega_n)$
can be defined according to Cramer's rule by the following ratio of
determinants:

\smallskip

$$
{\cal S}_k^{(n)} (\beta , s_n , \omega_n) \; = \; \frac
{
\vmatrix
{
s_n " \ldots " s_{n+k} \\ [1\jot]
\omega_{n} " \ldots " \omega_{n+k} \\
\vdots " \ddots " \vdots \\ [1\jot]
\omega_{n} / (\beta + n)_{k-1} " \ldots "
\omega_{n+k} / (\beta + n + k)_{k-1}
}
}
{
\vmatrix
{
1 " \ldots " 1       \\ [1\jot]
\omega_{n} " \ldots " \omega_{n+k} \\
\vdots " \ddots " \vdots \\ [1\jot]
\omega_{n} / (\beta + n)_{k-1} " \ldots "
\omega_{n+k} / (\beta + n + k)_{k-1}
}
} .
\tag
$$

\smallskip

If the sequence elements $s_n, \ldots , s_{n+k}$ satisfy eq. (8.2-1),
then obviously the sequence transformation ${\cal S}_k^{(n)} (\beta ,
s_n , \omega_n)$ is exact by construction, i.e.,
$$
{\cal S}_k^{(n)} (\beta , s_n , \omega_n) \; = \; s \, .
\tag
$$

As in the case of the Levin transformation it would be desirable to
find some alternative representation for the transformation ${\cal
S}_k^{(n)} (\beta , s_n , \omega_n)$. Fortunately, this can be
accomplished as easily as in the case of the Levin transformation. For
that purpose eq. (8.2-1) is rewritten in the following way:
$$
(\beta + n)_{k-1} [s_n - s] / \omega_n \; = \; \sum_{j=0}^{k-1} \, c_j
\> (\beta + n + j)_{k-j-1} \, .
\tag
$$

The highest power of $n$, which occurs on the right-hand side of eq.
(8.2-5), is $n^{k-1}$. Hence, if we apply the difference operator
$\Delta^k$ to eq. (8.2-5), the sum on the right-hand side, which is a
polynomial of degree $k-1$ in $n$, will be annihilated and we may
conclude from eqs. (8.2-4) and (8.2-5) that the sequence transformation
${\cal S}_k^{(n)} (\beta , s_n , \omega_n)$ is given by the following
ratio:
$$
{\cal S}_k^{(n)} (\beta , s_n , \omega_n) \; = \; \frac
{\Delta^k \, \{(\beta + n)_{k-1} \> s_n / \omega_n\} }
{\Delta^k \, \{(\beta + n)_{k-1} / \omega_n\} } \, .
\tag
$$

With the help of eq. (2.4-8) we obtain a representation of this
transformation as the ratio of two finite sums:
$$
{\cal S}_{k}^{(n)} (\beta , s_n, \omega_n) \; = \;
\frac
{\displaystyle
\sum_{j=0}^{k} \; ( - 1)^{j} \; \binom {k} {j} \;
\frac {(\beta + n +j )_{k-1}} {(\beta + n + k )_{k-1}} \;
\frac {s_{n+j}} {\omega_{n+j}} }
{\displaystyle
\sum_{j=0}^{k} \; ( - 1)^{j} \; \binom {k} {j} \;
\frac {(\beta + n +j )_{k-1}} {(\beta + n + k )_{k-1}} \;
\frac {1} {\omega_{n+j}} }
\; , \qquad k, n \in \N_0 \, .
\tag
$$

The common factor $(\beta+n+k)_{k-1}$ in eq. (8.2-7) was introduced in
order to decrease the magnitude of the terms of the numerator and
denominator sums, because otherwise overflow may happen too easily for
larger values of $k$.

The transformation (8.2-7) had already been treated by Sidi (see eq.
(1.9) of ref. [84]) who used this as well as some other transformations
for the derivation of explicit expressions for Pad\'e approximants of
some special hypergeometric series. However, it seems that Sidi did not
consider the transformation (8.2-7) to be a sequence transformation in
its own right. This is certainly an undeserved neglect. It will become
clear later that the transformation (8.2-7) is very powerful. We shall
see in section 13 that for divergent Stieltjes series, as they for
instance occur in the perturbation expansion of the quartic anharmonic
oscillator [3,85 - 88], it is certainly one of the most efficient
summation methods which is currently known.

As in the case of the Levin transformation, we also consider the
following mild extension of the sequence transformation ${\cal
S}_{k}^{(n)} (\beta , s_n, \omega_n)$:
$$
{\cal S}_{k,\ell}^{(n)} (\beta , s_n, \omega_n) \; = \;
\frac
{\displaystyle
\sum_{j=0}^{k} \; ( - 1)^{j} \; \binom {k} {j} \;
\frac {(\beta + n + \ell + j)_{k-\ell-1}} {(\beta + n + k )_{k-1}} \;
\frac {s_{n+j}} {\omega_{n+j}} }
{\displaystyle
\sum_{j=0}^{k} \; ( - 1)^{j} \; \binom {k} {j} \;
\frac {(\beta + n + \ell + j)_{k-\ell-1}} {(\beta + n + k )_{k-1}} \;
\frac {1} {\omega_{n+j}} } \; ,
\qquad k, \ell , n \in \N_0 \, .
\tag
$$

For $\ell = 0$ this transformation reduces to the transformation ${\cal
S}_k^{(n)} (\beta , s_n, \omega_n)$, eq. (8.2-7). An alternative
representation for the generalized transformation (8.2-8) can be
derived with the help of eq. (2.4-8),
$$
{\cal S}_{k,\ell}^{(n)} (\beta , s_n , \omega_n) \; = \; \frac
{\Delta^k \, \{(\beta + n + \ell)_{k-\ell-1} \> s_n / \omega_n\}}
{\Delta^k \, \{(\beta + n + \ell)_{k-\ell-1} \> / \omega_n\}}
\, .
\tag
$$

From this relationship we may deduce immediately that the
transformation (8.2-8) is exact for model sequences of the following
type:
$$
s_n \; = \; s \, + \, (\beta + n)_{\ell} \> \omega_n \,
\sum_{j=0}^{k-1} \, c_j / ( \beta + n)_j
\; , \qquad k, \ell , n \in \N_0 \, .
\tag
$$

\medskip

\Abschnitt Recurrence formulas

\smallskip

\aktTag = 0

Next, it will be shown how the numerator and denominator of the
transformation (8.2-7) can be computed recursively. It will turn out
that virtually the same technique can be used as in the case of the
Levin transformation. In eq. (8.2-7) both numerator and denominator are
of the general form
$$
Q_k^{(n)} (\beta) \; = \; \Delta^k \> Y_k^{(n)} (\beta) \, ,
\qquad k, n \in \N_0 \, .
\tag
$$

As usual, it is assumed that the difference operator $\Delta$ as well
as the shift operator $E$, which is defined in eq. (2.4-4), act only
upon $n$ and not upon $k$. The quantities $Y_k^{(n)} (\beta)$ satisfy
the following 2-term recursion in $k$:
$$
Y_k^{(n)} (\beta) \; = \; (\beta + n + k -2) \>
Y_{k-1}^{(n)} (\beta) \, ,
\qquad k \ge 1 \, , \quad n \ge 0 \, .
\tag
$$

Combination of eqs. (2.4-4), (7.2-3), (8.3-1) and (8.3-2) yields: $$
\beginAligntags
Q_k^{(n)} (\beta) \; " = " \;
\{ k E + (\beta + n + k - 2 ) \Delta \} \Delta^{k-1}
Y_{k-1}^{(n)} (\beta) \hfill \\ \tag
" = " \; \{ k E + (\beta + n + k - 2) \Delta \}
Q_{k-1}^{(n)} (\beta) \hfill \\ \tag
" = " \; (\beta + n +2 k - 2) Q_{k-1}^{(n+1)} (\beta) \, - \,
(\beta + n + k - 2) Q_{k-1}^{(n)} (\beta) \, .
\hfill \\ \tag
\endAligntags
$$

With the help of the 3-term recurrence formula (8.3-5) both the
numerator as well as the denominator of ${\cal S}_k^{(n)} (\beta , s_n,
\omega_n)$ can be computed for $k \ge 1$. However, as in the case of
the Levin transformation it is preferable to compute instead the scaled
quantities
$$
S_k^{(n)} (\beta) \; =
\; Q_k^{(n)} (\beta) \> / \> (\beta + n + k)_{k-1} \, .
\tag
$$

If we insert this into eq. (8.3-5), we obtain the following recurrence
formula for the scaled quantities $S_k^{(n)} (\beta)$:
$$
S_{k+1}^{(n)} (\beta) \; = \; S_k^{(n+1)} (\beta)
\, - \, \frac
{ (\beta + n + k ) (\beta + n + k - 1) }
{ (\beta + n + 2 k ) (\beta + n + 2 k - 1) } \>
S_k^{(n)} (\beta) \, , \qquad k, n \ge 0 \, .
\tag
$$

If we use the starting values
$$
S_0^{(n)} (\beta) \; = \; s_n / \omega_n \, ,
\qquad n \in \N_0 \, ,
\tag
$$
the 3-term recursion (8.3-7) produces the numerator of the
transformation (8.2-7), and if we use the starting values
$$
S_0^{(n)} (\beta) \; = \; 1 / \omega_n \, ,
\qquad n \in \N_0 \, ,
\tag
$$
we obtain the denominator of the transformation (8.2-7).

With the help of the 3-term recursion (8.3-7) it is also possible to
compute both the numerator and the denominator of the generalized
transformation (8.2-8). If the starting values
$$
S_0^{(n)} (\beta) \; = \; s_n / [(\beta + n)_{\ell} \>\omega_n]
\, , \qquad \ell, n \in \N_0 \, ,
\tag
$$
are used, eq. (8.3-7) produces the numerator of the generalized
transformation (8.2-8), and the starting values
$$
S_0^{(n)} (\beta) \; = \; 1 / [(\beta + n)_{\ell} \> \omega_n]
\, , \qquad \ell, n \in \N_0 \, ,
\tag
$$
give the denominator of the transformation (8.2-8).

Since the transforms with the highest values of the subscript normally
give the best results, our approximation to the limit $s$ of the
sequence $\Seqn s$ to be transformed will be the same as in the case of
the Levin transformation,
$$
\{ s_0, \omega_0; s_1, \omega_1; \ldots ; s_m, \omega_m \} \; \to \;
{\cal S}_m^{(0)} (\beta, s_0, \omega_0) \, ,
\qquad m \in \N_0 \, .
\tag
$$

Essentially this means that we shall use the following sequence of
transforms with minimal superscripts and maximal subscripts as
approximations to the limit $s$:
$$
{\cal S}_0^{(0)} (\beta, s_0, \omega_0), \;
{\cal S}_1^{(0)} (\beta, s_0, \omega_0), \; \ldots , \;
{\cal S}_m^{(0)} (\beta, s_0, \omega_0), \; \ldots \; .
\tag
$$

The recursive computation of the sequence transformation ${\cal
S}_k^{(n)} (\beta , s_n, \omega_n)$ can be done in virtually the same
way as in the case of the Levin transformation. For that purpose it is
recommendable to reformulate the 3-term recursion (8.3-7) in the
following way:
$$
\beginAligntags
" S_j^{(n-j)} (\beta) \; = \; S_{j-1}^{(n-j+1)} (\beta)
\, - \, \frac
{ (\beta + n - 1 ) (\beta + n - 2) }
{ (\beta + n + j - 2 ) (\beta + n + j - 3) } \>
S_{j-1}^{(n-j)} (\beta) \, , \\
" n \ge 1 \, , \qquad 1 \le j \le n \, .
\\ \tag
\endAligntags
$$

As in the case of the Levin transformation, only a single 1-dimensional
array will be needed for the computation of the $S_j^{(n-j)} (\beta)$
which are either numerator or denominator sums of the transformation
(8.2-7). For that purpose the $S_j^{(n-j)} (\beta)$ are stored in a
1-dimensional array ${\mit s}$ in such a way that the superscript
coincides with the index of the corresponding array element:
$$
S_{\nu}^{(n-\nu)} (\beta) \; \to \; {\mit s} (n-\nu) \, ,
\qquad n \ge 0 \, , \quad 0 \le \nu \le n \, .
\tag
$$

With this convention the recursive scheme (8.3-14) can be reformulated
in terms of the elements of the array ${\mit s}$:
$$
\beginAligntags
" {\mit s} (n-j) \; \gets \; {\mit s} (n-j+1) \, - \, \frac
{ (\beta + n - 1 ) (\beta + n - 2) }
{ (\beta + n + j - 2 ) (\beta + n + j - 3) } \>
{\mit s} (n-j) \, , \\
" n \ge 1 \, , \qquad 1 \le j \le n \, .
\\ \tag
\endAligntags
$$

\medskip

\Abschnitt Explicit remainder estimates

\smallskip

\aktTag = 0

It still has to be discussed how the auxiliary sequence $\Seq
{\omega_n}$ in eq. (8.2-7) should be chosen. The simplest approach
would be to proceed as in the case of the Levin transformation. There,
it was argued that the auxiliary sequence $\Seq {\omega_n}$ should be
chosen in such a way that $\omega_n$ is proportional to the dominant
term of the asymptotic expansion of the remainder $r_n$,
$$
r_n \; = \; s_n - s \; = \; \omega_n [ c + O(n^{-1}) ] \, ,
\qquad n \to \infty \, .
\tag
$$

Since the dominant term will not be affected if an asymptotic expansion
is transformed into a factorial series or vice versa, it should be
possible to use the same simple remainder estimates for sequences of
partial sums as in the case of the Levin transformation.

Hence, the remainder estimate (7.3-4) will be used in eq. (8.2-7). This
gives an analogue of Levin's $u$ transformation, eq. (7.3-5):
$$
y_{k}^{(n)} (\beta , s_n) \; = \;
\frac
{\displaystyle
\sum_{j=0}^{k} \; ( - 1)^{j} \; \binom {k} {j} \;
\frac {(\beta + n +j + 1 )_{k-2}} {(\beta + n + k )_{k-1}} \;
\frac {s_{n+j}} {a_{n+j}} }
{\displaystyle
\sum_{j=0}^{k} \; ( - 1)^{j} \; \binom {k} {j} \;
\frac {(\beta + n +j + 1 )_{k-2}} {(\beta + n + k )_{k-1}} \;
\frac {1} {a_{n+j}} }
\; .
\tag
$$

In the same way, the remainder estimate (7.3-6) can be used in eq.
(8.2-7). This gives an analogue of Levin's $t$ transformation, eq.
(7.3-7):
$$
{\tau}_{k}^{(n)} (\beta , s_n) \; = \;
\frac
{\displaystyle
\sum_{j=0}^{k} \; ( - 1)^{j} \; \binom {k} {j} \;
\frac {(\beta + n +j )_{k-1}} {(\beta + n + k )_{k-1}} \;
\frac {s_{n+j}} {a_{n+j}} }
{\displaystyle
\sum_{j=0}^{k} \; ( - 1)^{j} \; \binom {k} {j} \;
\frac {(\beta + n +j )_{k-1}} {(\beta + n + k )_{k-1}} \;
\frac {1} {a_{n+j}} }
\; .
\tag
$$

The use of the remainder estimate (7.3-8) in eq. (8.2-7) gives an
analogue of Levin's $d$ transformation, eq. (7.3-9):
$$
{\delta}_{k}^{(n)} (\beta , s_n) \; = \;
\frac
{\displaystyle
\sum_{j=0}^{k} \; ( - 1)^{j} \; \binom {k} {j} \;
\frac {(\beta + n +j )_{k-1}} {(\beta + n + k )_{k-1}} \;
\frac {s_{n+j}} {a_{n+j+1}} }
{\displaystyle
\sum_{j=0}^{k} \; ( - 1)^{j} \; \binom {k} {j} \;
\frac {(\beta + n +j )_{k-1}} {(\beta + n + k )_{k-1}} \;
\frac {1} {a_{n+j+1}} }
\; .
\tag
$$

Finally, the remainder estimate (7.3-10) gives an analogue of Levin's
$v$ transformation, eq. (7.3-11):
$$
{\phi}_{k}^{(n)} (\beta , s_n) \; = \;
\frac
{\displaystyle
\sum_{j=0}^{k} \; ( - 1)^{j} \; \binom {k} {j} \;
\frac {(\beta + n +j )_{k-1}} {(\beta + n + k )_{k-1}} \;
\frac {a_{n+j} - a_{n+j+1}} {a_{n+j} a_{n+j+1}} s_{n+j} }
{\displaystyle
\sum_{j=0}^{k} \; ( - 1)^{j} \; \binom {k} {j} \;
\frac {(\beta + n +j )_{k-1}} {(\beta + n + k )_{k-1}} \;
\frac {a_{n+j} - a_{n+j+1}} {a_{n+j} a_{n+j+1}} }
\; .
\tag
$$

If one of the remainder estimates (7.3-4), (7.3-6), (7.3-8), and
(7.3-10) is used in eq. (8.2-7), ${\cal S}_{k}^{(n)} (\beta , s_n,
\omega_n)$ is a {\it nonlinear} sequence transformation. If, however,
remainder estimates $\Seq {\omega_n}$ are used that do not depend
explicitly upon the elements of the sequence $\Seqn s$, ${\cal
S}_{k}^{(n)} (\beta , s_n, \omega_n)$ is a {\it linear} sequence
transformation.

Next, a factorial series analogue of the linear sequence transformation
$\Lambda_k^{(n)} (\beta, s_n)$, eq. (7.3-20), will be constructed. A
factorial series analogue of the model sequence (7.3-15), for which the
sequence transformation $\Lambda_k^{(n)} (\beta, s_n)$ is exact, would
be
$$
s_n \; = \; s \, + \, \sum_{j=0}^{k-1} \, c_j / (n+\alpha)_{j+1}
\, , \qquad k, n \in \N_0 \, .
\tag
$$

This model sequence can be obtained from the model sequence (8.2-1),
for which the sequence transformation (8.2-7) is exact, by choosing
$\beta = \alpha + 1$ and $\omega_n = 1/(n+\alpha)$. Hence, if we define
$$
{\cal F}_k^{(n)} (\alpha, s_n) \; = \; {\cal S}_k^{(n)} (\alpha +1, s_n
, 1/(n + \alpha)) \, , \qquad k, n \in \N_0 \, ,
\tag
$$
we see that that the sequence transformation ${\cal F}_k^{(n)} (\alpha,
s_n)$ is obviously exact for the model sequence (8.4-6). If we now use
eq. (8.2-6), we see that the transformation ${\cal F}_k^{(n)} (\alpha,
s_n)$ satisfies
$$
{\cal F}_k^{(n)} (\alpha, s_n) \; = \; \frac
{\Delta^k \{ (n+\alpha)_k \> s_n\} }
{\Delta^k (n+\alpha)_k }
\, , \qquad k, n \in \N_0 \, .
\tag
$$

The denominator in eq. (8.4-8) can be expressed in closed form. We only
have to use
$$
\Delta^k (n+\alpha)_k \; = \; k! \, .
\tag
$$

This is a special case of the following general relationship which can
be proved by complete induction in $k$,
$$
\Delta^k \> \frac {\Gamma (a+n)} {\Gamma (b+n)} \; = \;
(-1)^k \, (b-a)_k \> \frac {\Gamma (a+n)} {\Gamma (b+n+k)} \, .
\tag
$$

Combining eqs. (8.4-8) and ( 8.4-9) with eq. (2.4-8) gives us
$$
{\cal F}_k^{(n)} (\alpha, s_n) \; = \; (-1)^k \, \sum_{j=0}^k \>
(-1)^j \> \frac {(\alpha+n+j)_k} {j! \> (k-j)!} \> s_{n+j} \, ,
\qquad k, n \in \N_0 \, .
\tag
$$

If we combine eq. (8.3-5) with eqs. (8.4-8) and (8.4-9), we obtain the
following recursive scheme for the sequence transformation ${\cal
F}_k^{(n)} (\alpha, s_n)$:
$$
\beginAligntags
" {\cal F}_0^{(n)} (\alpha, s_n)\; " = \; " s_n
\, , \qquad n \in \N_0 \, ,
\erhoehe\aktTag \\ \tag*{\tagnr a}
" {\cal F}_{k+1}^{(n)} (\alpha, s_n) \; " = \; "
{\cal F}_k^{(n+1)} (\alpha, s_{n+1}) \, + \,
\frac {\alpha+n+k}{k+1} \> \Delta {\cal F}_k^{(n)} (\alpha, s_n)
\, , \qquad k, n \in \N_0 \, . \\ \tag*{\tagform\aktTagnr b}
\endAligntags
$$

The transformation ${\cal F}_k^{(n)} (\alpha, s_n)$ can be computed in
essentially the same way as the Richardson extrapolation process. For
that purpose it is recommendable to rewrite the above recursive scheme
in the following way:
$$
\beginAligntags
" {\cal F}_0^{(n)} (\alpha, s_n)\; " = \; " s_n
\, , \qquad n \ge 0 \, ,
\erhoehe\aktTag \\ \tag*{\tagnr a}
" {\cal F}_j^{(n-j)} (\alpha, s_{n-j}) \; " = \;
" {\cal F}_{j-1}^{(n-j+1)} (\alpha, s_{n-j+1}) \, + \,
\frac {\alpha+n-1} {j} \>
\Delta {\cal F}_{j-1}^{(n-j)} (\alpha, s_{n-j}) \, , \\
" n \ge 1 \, , \quad 1 \le j \le n \, . \span\omit \span\omit
\\ \tag*{\tagform\aktTagnr b}
\endAligntags
$$

If the ${\cal F}_j^{(n)} (\alpha, s_n)$ are stored in a 1-dimensional
array ${\mit f}$ according to the rule
$$
{\cal F}_j^{(n-j)} (\alpha, s_{n-j}) \;
\to \; {\mit f} (n-j) \, ,
\qquad n \ge 0 \, , \quad 0 \le j \le n \, ,
\tag
$$
we see that the recursive scheme (8.4-13) can be reformulated in terms
of the elements of the array ${\mit f}$:
$$
\beginAligntags
" {\mit f} (n) \; " = \; " s_n \, , \qquad n \ge 0 \, ,
\erhoehe\aktTag \\ \tag*{\tagnr a}
" {\mit f} (n-j) \; " = \; " {\mit f} (n-j+1) \, + \,
\frac {\alpha+n-1}{j} \> \Delta {\mit f} (n-j) \, ,
\qquad n \ge 1 \, , \quad 1 \le j \le n \, .
\\ \tag*{\tagform\aktTagnr b}
\endAligntags
$$

A discussion of the linear but nonregular sequence transformation
${\cal F}_k^{(n)} (\alpha, s_n)$ can be found in Wimp's book (see pp.
38 - 40 of ref. [23]). However, the recursive scheme (8.4-12) for the
computation of this transformation seems to be new.

The sequence transformations $y_k^{(n)} (\beta , s_n)$, eq. (8.4-2),
and $\tau_k^{(n)} (\beta , s_n)$, eq. (8.4-3), require the sequence
elements $s_{n-1}, s_n, s_{n+1}, \ldots, s_{n+k}$ for their
computation, whereas $\delta_k^{(n)} (\beta , s_n)$, eq. (8.4-4),
requires the sequence elements $s_n, s_{n+1}, \ldots , s_{n+k+1}$.
Hence, they are all transformations of order $k+1$. The sequence
transformation ${\phi}_k^{(n)} (\beta , s_n)$, eq. (8.4-5), requires
the sequence elements $s_{n-1}, s_n, s_{n+1}, \ldots , s_{n+k+1}$ which
implies that it is a transformation of order $k+2$. The linear
sequence transformation ${\cal F}_k^{(n)} (\alpha, s_n)$, eq. (8.4-11),
requires the sequence elements $s_n, s_{n+1}, \ldots , s_{n+k}$, i.e.,
it is a transformation of order $k$.

The situation is somewhat different if the transforms with superscript
$n = 0$ are computed because then $y_k^{(0)} (\beta , s_0)$ and
${\tau}_k^{(0)} (\beta , s_0)$ are transformations of order $k$,
whereas ${\delta}_k^{(0)} (\beta , s_0)$ and ${\phi}_k^{(0)} (\beta ,
s_0)$ are transformations of order $k+1$.

\endAbschnittsebene

\neueSeite

\Abschnitt Other generalizations of Levin's sequence transformation

\vskip - 2 \jot

\beginAbschnittsebene

\medskip

\Abschnitt Asymptotic approximations based upon Pochhammer symbols

\smallskip

\aktTag = 0

In the last section it was demonstrated how a new class of sequence
transformations can be derived in exactly the same way as the Levin
transformation which is generally accepted to be a very powerful
convergence acceleration and summation method [29,30,57, 60-64].

The only difference between the Levin transformation and the new
transformation discussed in the last section is that the Levin
transformation assumes that the ratio $r_n / \omega_n$ can be expressed
as an asymptotic power series whereas the new transformation assumes
that $r_n / \omega_n$ can be expressed as a factorial series.
Consequently, the analytical expressions for the various Levin
transformations and those for the analogous variants of the new
transformation can easily be transformed into each other. For instance,
one only has to replace the powers $(\beta+n+j)^{k-1}$ in the
expression for the general Levin transformation ${\cal L}_k^{(n)}
(\beta, s_n , \omega_n)$, eq. (7.1-7), by Pochhammer symbols
$(\beta+n+j)_{k-1}$ in order to obtain the analogous new transformation
${\cal S}_k^{(n)} (\beta, s_n , \omega_n)$, eq. (8.2-7).

However, these new transformations, which were discussed in the last
section, do not yet exhaust all possibilities of constructing other
simple generalizations of the Levin transformation, which nevertheless
retain most of the advantages of the Levin transformation. For
instance, in recent articles on large order perturbation theory
asymptotic approximations of the following general type were considered
[89-91]:
$$
\beginAligntags
" f (z) \; " \sim \;
" \frac {c_0} {z} \, + \, \frac {c_1} {z (z-1)}
\, + \, \frac {c_2} {z (z-1) (z-2)} \, + \, \cdots \, + \,
\frac {c_n} {z (z-1) (z-2) \ldots (z-n)} \\
" " = \; " \sum_{\nu=0}^n \> (-1)^{\nu+1} \;
\frac {c_{\nu}} {(-z)_{\nu+1} }
\, , \qquad \vert z \vert \to \infty \, .
\\ \tag
\endAligntags
$$

Superficially, such an expression looks very much like a truncated
factorial series since the argument $z$ occurs also in Pochhammer
symbols. However, the fact that the Pochhammer symbols in eq. (9.1-1)
are of the type
$$
\beginAligntags
" (-z)_{\nu+1} \;" = \; " (-z)(-z+1) \ldots (-z+\nu) \\
" "= \; " (-1)^{\nu+1} z (z-1) \ldots (z-\nu) \, ,
\qquad \nu \in \N_0 \, ,
\\ \tag
\endAligntags
$$
has some far-reaching consequences. For instance, if $z$ is a positive
real number, eq. (9.1-1) makes sense only if $z > n$ holds. For $n >
z$, either the later terms in the sum will show irregular sign
patterns, or, if $z$ happens to be a positive integer, some Pochhammer
symbols will then be zero. Consequently, for a fixed value of $z$ the
summation limit $n$ in eq. (9.1-1) cannot be extended to infinity, and
such an expression cannot be considered to be the truncation of an
asymptotic series after a finite number of terms. Instead, an expression
such as eq. (9.1-1) has to be interpreted to be some kind of asymptotic
approximation involving only a finite number of terms.

The author is not aware of any reference in the mathematical
literature, in which expressions like the one in eq. (9.1-1) are
treated and their properties are analyzed. Consequently, the material
in this section is somewhat experimental and its mathematical basis is
not as solid as in the other sections.

However, it must be emphasized that these objections do not exclude the
possibility that finite sums of the type of eq. (9.1-1) may yield
excellent approximations if suitable restrictions on $z$ and $n$ are
made. We shall see later in section 13 that sequence transformations
which are based upon asymptotic approximations of the type of eq.
(9.1-1) are indeed able to produce excellent results in convergence
acceleration and summation processes.

\medskip

\Abschnitt New sequence transformations based upon Pochhammer symbols

\smallskip

\aktTag = 0

The following model sequence will be the basis for the new class of
sequence transformations which will be derived in this section:
$$
s_n \; = \; s \, + \, \omega_n \, \sum_{j=0}^{k-1} \, c_j / (- \gamma -
n)_j \, , \qquad k,n \in \N_0 \, .
\tag
$$

This sequence is formally almost identical with the model sequence
(8.2-1). The only difference is that the Pochhammer symbols
$(n+\beta)_j$ in eq. (8.2-1) are replaced by Pochhammer symbols $(-
\gamma - n)_j$. Concerning the sequence $\Seq {\omega_n}$ of remainder
estimates it is again assumed that the $\omega_n$ are known functions
of $n$ which have to be different from zero and distinct for all finite
values of $n$. But otherwise, the $\omega_n$ are in principle
completely arbitrary.

The parameter $\gamma$ in eq. (9.2-1) is not only subject to the
restriction that the Pochhammer symbols $(- \gamma - n)_j$ must not be
zero for all admissible values of $n$ and $j$. Also, the regular sign
pattern of the Pochhammer symbols in eq. (9.2-1) must not be destroyed.
These two restrictions suggest that $\gamma$ should be a positive
number satisfying $\gamma \ge k - 1$.

In eq. (9.2-1) there occur $k+1$ unknown quantities, the limit or
antilimit $s$ and the $k$ linear coefficients $c_0, \ldots , c_{k-1}$.
Hence, if $k+1$ sequence elements $s_n, \ldots , s_{n+k}$ are known,
the sequence transformation ${\cal M}_k^{(n)} (\gamma, s_n, \omega_n)$
can be defined according to Cramer's rule by the following ratio of
determinants:

\smallskip

$$
{\cal M}_k^{(n)} (\gamma , s_n , \omega_n) \; = \; \frac
{
\vmatrix
{
s_n " \ldots " s_{n+k} \\ [1\jot]
\omega_{n} " \ldots " \omega_{n+k} \\
\vdots " \ddots " \vdots \\ [1\jot]
\omega_{n} / (- \gamma - n)_{k-1} " \ldots "
\omega_{n+k} / (- \gamma - n - k)_{k-1}
}
}
{
\vmatrix
{
1 " \ldots " 1       \\ [1\jot]
\omega_{n} " \ldots " \omega_{n+k} \\
\vdots " \ddots " \vdots \\ [1\jot]
\omega_{n} / (- \gamma - n)_{k-1} " \ldots "
\omega_{n+k} / (- \gamma - n - k)_{k-1}
}
} .
\tag
$$

\smallskip

If the sequence elements $s_n, \ldots ,s_{n+k}$ satisfy eq. (9.2-1),
then obviously the sequence transformation ${\cal M}_k^{(n)} (\gamma ,
s_n , \omega_n)$ is exact by construction, i.e.,
$$
{\cal M}_k^{(n)} (\gamma , s_n , \omega_n) \; = \; s \, .
\tag
$$

Again it would be desirable to have some alternative representation for
the transformation ${\cal M}_k^{(n)} (\gamma , s_n , \omega_n)$.
Fortunately, this can be accomplished quite easily. For that purpose
eq. (9.2-1) is rewritten in the following way:
$$
(- \gamma - n)_{k-1} [s_n - s] / \omega_n \; = \;
\sum_{j=0}^{k-1} \, c_j \> (- \gamma - n + j)_{k-j-1} \, .
\tag
$$

The highest power of $n$, which occurs on the right-hand side of eq.
(9.2-4), is $n^{k-1}$. Hence, if we apply the difference operator
$\Delta^k$ to eq. (9.2-4), the sum on the right-hand side, which is a
polynomial of degree $k-1$ in $n$, will be annihilated and we may
conclude from eqs. (9.2-3) and (9.2-4) that the sequence transformation
${\cal M}_k^{(n)} (\beta , s_n , \omega_n)$ is given by the following
ratio:
$$
{\cal M}_k^{(n)} (\gamma , s_n , \omega_n) \; = \; \frac
{\Delta^k \, \{(- \gamma - n)_{k-1} \> s_n / \omega_n\} }
{\Delta^k \, \{(- \gamma - n)_{k-1} / \omega_n\} } \, .
\tag
$$

If we use eq. (2.4-8), we see that this transformation can be
represented as the ratio of two finite sums:
$$
{\cal M}_{k}^{(n)} (\gamma ,s_n, \omega_n) \; = \;
\frac
{\displaystyle
\sum_{j=0}^{k} \; ( - 1)^{j} \; \binom {k} {j} \;
\frac {(- \gamma - n - j )_{k-1}} {(- \gamma - n - k )_{k-1}} \;
\frac {s_{n+j}} {\omega_{n+j}} }
{\displaystyle
\sum_{j=0}^{k} \; ( - 1)^{j} \; \binom {k} {j} \;
\frac {(- \gamma - n - j )_{k-1}} {(- \gamma - n - k )_{k-1}} \;
\frac {1} {\omega_{n+j}} } \; , \qquad k,n \in \N_0 \, .
\tag
$$

The common factor $(- \gamma-n-k)_{k-1}$ in eq. (9.2-6) was introduced
in order to decrease the magnitude of the terms of the numerator and
denominator sums, because otherwise overflow may happen too easily for
larger values of $k$.

As in the previous sections, we also consider a mild extension of the
sequence transformation ${\cal M}_{k}^{(n)} (\gamma ,s_n, \omega_n)$:
$$
{\cal M}_{k,\ell}^{(n)} (\gamma ,s_n, \omega_n) \; = \;
\frac
{\displaystyle
\sum_{j=0}^{k} \; ( - 1)^{j} \; \binom {k} {j} \;
\frac {(- \gamma - n - j + \ell)_{k-\ell-1}}
{(- \gamma - n - k )_{k-1}} \;
\frac {s_{n+j}} {\omega_{n+j}} }
{\displaystyle
\sum_{j=0}^{k} \; ( - 1)^{j} \; \binom {k} {j} \;
\frac {(- \gamma - n - j + \ell)_{k-\ell-1}}
{(- \gamma - n - k )_{k-1}} \;
\frac {1} {\omega_{n+j}} } \; ,
\qquad k,\ell ,n \in \N_0 \, .
\tag
$$

For $\ell = 0$ this transformation reduces to the transformation ${\cal
M}_k^{(n)} (\gamma ,s_n, \omega_n)$, eq. (9.2-6). An alternative
representation for the generalized transformation (9.2-7) can be
derived with the help of eq. (2.4-8),
$$
{\cal M}_{k,\ell}^{(n)} (\gamma , s_n , \omega_n) \; = \; \frac
{\Delta^k \, \{(- \gamma - n + \ell)_{k-\ell-1}
\> s_n / \omega_n\}}
{\Delta^k \, \{(- \gamma - n + \ell)_{k-\ell-1}
\> / \omega_n\}} \, .
\tag
$$

From this relationship we may deduce immediately that the
transformation (9.2-7) is exact for model sequences of the following
type:
$$
s_n \; = \; s \, + \, (- \gamma - n)_{\ell} \> \omega_n \,
\sum_{j=0}^{k-1} \, c_j / ( - \gamma - n)_j
\; , \qquad k,\ell ,n \in \N_0 \, .
\tag
$$

\medskip

\Abschnitt Recurrence formulas

\smallskip

\aktTag = 0

Next, it will be shown how the numerator and denominator of the
transformation (9.2-6) can be computed recursively. It will turn out
that virtually the same technique can be used as in sections 7.2 and
8.3. In eq. (9.2-6) both numerator and denominator are of the general
form
$$
{\mit R}_k^{(n)} (\gamma) \; = \; \Delta^k \> Z_k^{(n)} (\gamma) \, ,
\qquad k,n \in \N_0 \, .
\tag
$$

As usual, it is assumed that the difference operator $\Delta$ as well
as the shift operator $E$, which is defined in eq. (2.4-4), act only
upon $n$ and not upon $k$. The quantities $Z_k^{(n)} (\gamma)$ satisfy
the following 2-term recursion in $k$:
$$
Z_k^{(n)} (\gamma) \; = \; (- \gamma - n + k - 2) \>
Z_{k-1}^{(n)} (\gamma) \, ,
\qquad k \ge 1 \, , \quad n \ge 0 \, .
\tag
$$

Combination of eqs. (2.4-4), (7.2-3), (9.3-1) and (9.3-2) yields: $$
\beginAligntags
{\mit R}_k^{(n)} (\gamma) \; " = " \;
\{ (- \gamma - n + k - 2 ) \Delta - k E \} \Delta^{k-1}
Z_{k-1}^{(n)} (\gamma) \hfill \\ \tag
" = " \; \{ (- \gamma - n + k - 2) \Delta - k E \}
{\mit R}_{k-1}^{(n)} (\gamma) \hfill \\ \tag
" = " \; (\gamma + n - k + 2) {\mit R}_{k-1}^{(n)} (\gamma) \, - \,
(\gamma + n + 2) {\mit R}_{k-1}^{(n+1)} (\gamma) \, .
\hfill \\ \tag
\endAligntags
$$

With the help of the 3-term recurrence formula (9.3-5) both the
numerator as well as the denominator of the transformation ${\cal
M}_k^{(n)} (\gamma ,s_n, \omega_n)$ can be computed for $k \ge 1$.
However, it is again preferable to compute instead the scaled quantities
$$
M_k^{(n)} (\gamma) \; =
\; {\mit R}_k^{(n)} (\gamma) \> / \> (- \gamma - n - k)_{k-1} \, .
\tag
$$

If we insert this into eq. (9.3-5), we obtain the following recurrence
formula for the scaled quantities $M_k^{(n)} (\gamma)$:
$$
M_{k+1}^{(n)} (\gamma) \; = \; M_k^{(n+1)} (\gamma)
\, - \, \frac
{ \gamma + n - k + 1 } { \gamma + n + k + 1 } \>
M_k^{(n)} (\gamma) \, , \qquad k, n \ge 0 \, .
\tag
$$

If we use the starting values
$$
M_0^{(n)} (\gamma) \; = \; s_n / \omega_n \, ,
\qquad n \in \N_0 \, ,
\tag
$$
the 3-term recursion (9.3-7) produces the numerator of the
transformation (9.2-6), and if we use the starting values
$$
M_0^{(n)} (\gamma) \; = \; 1 / \omega_n \, ,
\qquad n \in \N_0 \, ,
\tag
$$
we obtain the denominator of the transformation (9.2-7).

With the help of the 3-term recursion (9.3-7) it is also possible to
compute both the numerator and the denominator of the generalized
transformation (9.2-7). If the starting values
$$
M_0^{(n)} (\gamma) \; = \;
s_n / [(- \gamma - n)_{\ell} \> \omega_n] \, ,
\qquad \ell, n \in \N_0 \, ,
\tag
$$
are used, the 3-term recursion (9.3-7) produces the numerator of the
generalized transformation (9.2-7), and the starting values
$$
M_0^{(n)} (\gamma) \; = \; 1 / [(- \gamma - n)_{\ell} \> \omega_n]
\, , \qquad \ell, n \in \N_0 \, ,
\tag
$$
produce the denominator of the transformation (9.2-7).

Since the transforms with the highest values of the subscript normally
give the best results, our approximation to the limit $s$ of the
sequence $\Seqn s$ to be transformed will be the same as in the
previous two sections, i.e.,
$$
\{ s_0, \omega_0; s_1, \omega_1; \ldots ; s_m, \omega_m \} \; \to \;
{\cal M}_m^{(0)} (\gamma, s_0, \omega_0) \, ,
\qquad m \in \N_0 \, .
\tag
$$

Essentially this means that we shall use the following sequence of
transforms with minimal superscripts and maximal subscripts as
approximations to the limit $s$:
$$
{\cal M}_0^{(0)} (\gamma, s_0, \omega_0), \;
{\cal M}_1^{(0)} (\gamma, s_0, \omega_0), \; \ldots , \;
{\cal M}_m^{(0)} (\gamma, s_0, \omega_0), \; \ldots \; .
\tag
$$

The recursive computation of the sequence transformations ${\cal
M}_k^{(n)} (\gamma ,s_n, \omega_n)$ can be done in virtually the same
way as in the previous two sections. For that purpose it is
recommendable to reformulate the 3-term recursion (9.3-7) in the
following way:
$$
M_j^{(n-j)} (\gamma) \; = \; M_{j-1}^{(n-j+1)} (\gamma)
\, - \, \frac
{ \gamma + n - 2 j + 2 } { \gamma + n } \>
M_{j-1}^{(n-j)} (\gamma) \, ,
\qquad n \ge 1 \, , \quad 1 \le j \le n \, .
\tag
$$

Again, only a single 1-dimensional array will be needed for the
computation of the $M_j^{(n-j)} (\gamma)$ which are either numerator or
denominator sums of the transformation (9.2-6). For that purpose the
$M_j^{(n-j)} (\gamma)$ are stored in a 1-dimensional array ${\mit m}$
in such a way that the superscript coincides with the index of the
corresponding array element:
$$
M_{\nu}^{(n-\nu)} (\gamma) \; \to \; {\mit m} (n-\nu) \, ,
\qquad n \ge 0 \, , \quad 0 \le \nu \le n \, .
\tag
$$

With this convention the recursive scheme (9.3-14) can be reformulated
in terms of the elements of the array ${\mit m}$:
$$
{\mit m} (n-j) \; \gets \; {\mit m} (n-j+1) \, - \, \frac
{ \gamma + n - 2 j + 2 } { \gamma + n } \>
{\mit m} (n-j) \, , \qquad n \ge 1 \, , \quad 1 \le j \le n \, .
\tag
$$

\medskip

\Abschnitt Explicit remainder estimates

\smallskip

\aktTag = 0

It still must be discussed how the auxiliary sequence $\Seq {\omega_n}$
in eq. (9.2-6) should be chosen. The simplest approach would again
consist of using the same simple remainder estimates for sequences of
partial sums as in the case of the Levin transformation.

Hence, if we use in eq. (9.2-6) the remainder estimate
$$
\omega_n \; = \; (- \gamma - n) a_n \, ,
\tag
$$
we obtain an analogue of Levin's $u$ transformation, eq. (7.3-5):
$$
Y_{k}^{(n)} (\gamma ,s_n) \; = \;
\frac
{\displaystyle
\sum_{j=0}^{k} \; ( - 1)^{j} \; \binom {k} {j} \;
\frac
{(- \gamma - n - j + 1)_{k-2}} {(- \gamma - n - k )_{k-1}} \;
\frac {s_{n+j}} {a_{n+j}} }
{\displaystyle
\sum_{j=0}^{k} \; ( - 1)^{j} \; \binom {k} {j} \;
\frac
{(- \gamma - n - j + 1)_{k-2}} {(- \gamma - n - k )_{k-1}} \;
\frac {1} {a_{n+j}} }
\; .
\tag
$$

In the same way, the remainder estimate (7.3-6) can be used in eq.
(9.2-6). This gives an analogue of Levin's $t$ transformation, eq.
(7.3-7):
$$
T_{k}^{(n)} (\gamma ,s_n) \; = \;
\frac
{\displaystyle
\sum_{j=0}^{k} \; ( - 1)^{j} \; \binom {k} {j} \;
\frac
{(- \gamma - n - j )_{k-1}} {(- \gamma - n - k )_{k-1}} \;
\frac {s_{n+j}} {a_{n+j}} }
{\displaystyle
\sum_{j=0}^{k} \; ( - 1)^{j} \; \binom {k} {j} \;
\frac
{(- \gamma - n - j )_{k-1}} {(- \gamma - n - k )_{k-1}} \;
\frac {1} {a_{n+j}} }
\; .
\tag
$$

The use of the remainder estimate (7.3-8) in eq. (9.2-6) gives an
analogue of Levin's $d$ transformation, eq. (7.3-9):
$$
{\Delta}_{k}^{(n)} (\gamma ,s_n) \; = \;
\frac
{\displaystyle
\sum_{j=0}^{k} \; ( - 1)^{j} \; \binom {k} {j} \;
\frac
{(- \gamma - n - j )_{k-1}} {(- \gamma - n - k )_{k-1}} \;
\frac {s_{n+j}} {a_{n+j+1}} }
{\displaystyle
\sum_{j=0}^{k} \; ( - 1)^{j} \; \binom {k} {j} \;
\frac
{(- \gamma - n - j )_{k-1}} {(- \gamma - n - k )_{k-1}} \;
\frac {1} {a_{n+j+1}} }
\; .
\tag
$$

Finally, the remainder estimate (7.3-10) gives an analogue of Levin's
$v$ transformation, eq. (7.3-11):
$$
{\Phi}_{k}^{(n)} (\gamma ,s_n) \; = \;
\frac
{\displaystyle
\sum_{j=0}^{k} \; ( - 1)^{j} \; \binom {k} {j} \;
\frac
{(- \gamma - n - j )_{k-1}} {(- \gamma - n - k )_{k-1}} \;
\frac {a_{n+j} - a_{n+j+1}} {a_{n+j} a_{n+j+1}} s_{n+j} }
{\displaystyle
\sum_{j=0}^{k} \; ( - 1)^{j} \; \binom {k} {j} \;
\frac
{(- \gamma - n - j )_{k-1}} {(- \gamma - n - k )_{k-1}} \;
\frac {a_{n+j} - a_{n+j+1}} {a_{n+j} a_{n+j+1}} }
\; .
\tag
$$

If one of the remainder estimates (7.3-6), (7.3-8), (7.3-10), and
(9.4-1) is used in eq. (9.2-6), ${\cal M}_{k}^{(n)} (\gamma ,s_n,
\omega_n)$ is a {\it nonlinear} sequence transformation. If, however,
remainder estimates $\Seq {\omega_n}$ are used that do not depend
explicitly upon the elements of the sequence $\Seqn s$, ${\cal
M}_{k}^{(n)} (\gamma ,s_n, \omega_n)$ is a {\it linear} sequence
transformation.

Next, an analogue of the Richardson extrapolation scheme (6.1-5) will
be introduced which is based upon asymptotic approximations involving
Pochhammer symbols of the type of eq. (9.1-1). For that purpose we
consider the following model sequence:
$$
s_n \; = \; s \, + \, \sum_{j=0}^{k-1} \, c_j / (-\zeta-n)_{j+1}
\, , \qquad k,n \in \N_0 \, .
\tag
$$

This model sequence can be obtained from the model sequence (9.2-1),
for which the sequence transformation (9.2-6) is exact, by choosing
$\zeta = \gamma + 1$ and $\omega_n = - 1/(\zeta + n)$. Hence, if we
define
$$
{\cal P}_k^{(n)} (\zeta, s_n) \; = \; {\cal M}_k^{(n)} (\zeta -1, s_n ,
- 1/(n + \zeta)) \, , \qquad k,n \in \N_0 \, ,
\tag
$$
we see that that the sequence transformation ${\cal P}_k^{(n)} (\zeta,
s_n)$ is obviously exact for the model sequence (9.4-6). If we now use
eq. (9.2-5), we see that the transformation ${\cal P}_k^{(n)} (\zeta,
s_n)$ satisfies
$$
{\cal P}_k^{(n)} (\zeta, s_n) \; = \; \frac
{\Delta^k \{ (- \zeta-n)_k \> s_n\} }
{\Delta^k (- \zeta-n)_k }
\, , \qquad k,n \in \N_0 \, .
\tag
$$

The denominator in eq. (9.4-8) can be expressed in closed form. We only
have to use
$$
\Delta^k (- \zeta-n)_k \; = \; (-1)^k \, k! \, .
\tag
$$

This is a special case of the following general relationship which can
be proved by complete induction in $k$,
$$
\Delta^k \> \frac {\Gamma (a-n)} {\Gamma (b-n)} \; = \;
(b-a)_k \> \frac {\Gamma (a-n-k)} {\Gamma (b-n)} \, .
\tag
$$

Combination of eqs. (9.4-8) and ( 9.4-9) with eq. (2.4-8) yields
$$
{\cal P}_k^{(n)} (\zeta, s_n) \; = \; \sum_{j=0}^k \> (-1)^j
\> \frac {(- \zeta-n-j)_k} {j! \> (k-j)!} \> s_{n+j} \, ,
\qquad k,n \in \N_0 \, .
\tag
$$

If we combine eq. (9.3-5) with eqs. (9.4-8) and (9.4-9), we obtain the
following recursive scheme for the sequence transformation ${\cal
P}_k^{(n)} (\zeta, s_n)$:
$$
\beginAligntags
" {\cal P}_0^{(n)} (\zeta, s_n)\; " = \; " s_n \, ,
\qquad n \in \N_0 \, , \erhoehe\aktTag \\ \tag*{\tagnr a}
" {\cal P}_{k+1}^{(n)} (\zeta, s_n) \; " = \; "
{\cal P}_k^{(n+1)} (\zeta, s_{n+1}) \, + \,
\frac {\zeta+n-k} {k+1}
\> \Delta {\cal P}_k^{(n)} (\zeta, s_n) \, ,
\qquad k,n \in \N_0 \, . \\ \tag*{\tagform\aktTagnr b}
\endAligntags
$$

The transformation ${\cal P}_k^{(n)} (\zeta, s_n)$ can be computed in
essentially the same way as the Richardson extrapolation process. For
that purpose it is recommendable to rewrite the above recursive scheme
in the following way:
$$
\beginAligntags
" {\cal P}_0^{(n)} (\zeta, s_n)\; " = \; " s_n \, ,
\qquad n \ge 0 \, , \erhoehe\aktTag \\ \tag*{\tagnr a}
" {\cal P}_j^{(n-j)} (\zeta, s_{n-j}) \; " = \; "
{\cal P}_{j-1}^{(n-j+1)} (\zeta, s_{n-j+1}) \, + \,
\frac {\zeta + n - 2 j + 1} {j} \>
\Delta {\cal P}_{j-1}^{(n-j)} (\zeta, s_{n-j}) \, , \\
" n \ge 1 \, , \qquad 1 \le j \le n \, . \span\omit \span\omit
\\ \tag*{\tagform\aktTagnr b}
\endAligntags
$$

If the ${\cal P}_j^{(n)} (\zeta, s_n)$ are stored in a 1-dimensional
array ${\mit p}$ according to the rule
$$
{\cal P}_j^{(n-j)} (\zeta, s_{n-j}) \; \to \; {\mit p} (n-j) \, ,
\qquad n \ge 0 \, , \quad 0 \le j \le n \, ,
\tag
$$
we see that the recursive scheme (9.4-14) can be reformulated in terms
of the elements of the array ${\mit p}$:
$$
\beginAligntags
" {\mit p} (n) \; " \gets \; " s_n \, , \qquad n \ge 0 \, ,
\erhoehe\aktTag \\ \tag*{\tagnr a}
" {\mit p} (n-j) \; " \gets \; " {\mit p} (n-j+1) \, + \,
\frac {\zeta + n - 2 j +1 } {j} \> \Delta {\mit p} (n-j) \, , \qquad n
\ge 1 \, , \quad 1 \le j \le n \, . \qquad
\\ \tag*{\tagform\aktTagnr b}
\endAligntags
$$

The sequence transformations $Y_k^{(n)} (\gamma , s_n)$, eq. (9.4-2),
and $T_k^{(n)} (\gamma , s_n)$, eq. (9.4-3), require the sequence
elements $s_{n-1}, s_n, s_{n+1}, \ldots, s_{n+k}$ for their
computation, whereas ${\Delta}_k^{(n)} (\gamma , s_n)$, eq. (9.4-4),
requires the sequence elements $s_n, s_{n+1}, \ldots , s_{n+k+1}$.
Hence, they are all transformations of order $k+1$. The sequence
transformation ${\Phi}_k^{(n)} (\gamma , s_n)$, eq. (9.4-5), requires
the sequence elements $s_{n-1}, s_n, s_{n+1}, \ldots , s_{n+k+1}$ which
implies that it is a transformation of order $k+2$. The linear
sequence transformation ${\cal P}_k^{(n)} (\zeta, s_n)$, eq. (9.4-11),
requires the sequence elements $s_n, s_{n+1}, \ldots , s_{n+k}$, i.e.,
it is a transformation of order $k$.

The situation is somewhat different if the transforms with superscript
$n = 0$ are computed because then $Y_k^{(0)} (\gamma , s_0)$ and
$T_k^{(0)} (\gamma, s_0)$ are transformations of order $k$, whereas
${\Delta}_k^{(0)} (\gamma, s_0)$ and ${\Phi}_k^{(0)} (\gamma, s_0)$ are
transformations of order $k+1$.

\medskip

\Abschnitt Drummond's sequence transformation

\smallskip

\aktTag = 0

Let us assume that $P_{k-1} (n)$ is a polynomial of degree $k-1$ in
$n$. We want to derive a sequence transformation ${\cal D}_k^{(n)}
(s_n, \omega_n)$, which is by construction exact for the following
model sequence:
$$
s_n \; = \; s \, + \, \omega_n \, P_{k-1} (n) \, ,
\qquad k,n \in \N_0 \, .
\tag
$$

Concerning the sequence $\Seq {\omega_n}$ of remainder estimates it is
again assumed that the $\omega_n$ are different from zero for all
finite values of $n$, but otherwise they are in principle completely
arbitrary. For the derivation of a sequence transformation, which is
exact for the above model sequence, we rewrite eq. (9.5-1) in the
following way:
$$
[s_n - s]/\omega_n \; = \; P_{k-1} (n) \, ,
\qquad k,n \in \N_0 \, .
\tag
$$

Since $P_{k-1} (n)$ is by assumption a polynomial of degree $k-1$ in
$n$, the right-hand side of eq. (9.5-2) will be annihilated by the
difference operator $\Delta^k$. Hence, we can define the sequence
transformation ${\cal D}_k^{(n)} (s_n, \omega_n)$ by the following
ratio:
$$
{\cal D}_k^{(n)} (s_n , \omega_n) \; = \; \frac
{\Delta^k \, \{ s_n / \omega_n\} }
{\Delta^k \, \{ 1 / \omega_n\} } \, .
\tag
$$

With the help of eq. (2.4-8) we obtain a representation of this
transformation as the ratio of two finite sums:
$$
{\cal D}_{k}^{(n)} (s_n, \omega_n) \; = \;
\frac
{\displaystyle
\sum_{j=0}^{k} \; ( - 1)^{j} \; \binom {k} {j} \;
\frac {s_{n+j}} {\omega_{n+j}} }
{\displaystyle
\sum_{j=0}^{k} \; ( - 1)^{j} \; \binom {k} {j} \;
\frac {1} {\omega_{n+j}} }
\; , \qquad k,n \in \N_0 \, .
\tag
$$

The special case $\omega_n = a_n$ of this sequence transformation was
originally derived by Drummond [92]. Later, it was rederived by Sidi
[84].

Both numerator and denominator of the sequence transformation (9.5-4)
can be computed with the help of the following 3-term recurrence
formula, which is an immediate consequence of eq. (9.5-3):
$$
{\mit D}_{k+1}^{(n)} \; = \;
{\mit D}_k^{(n+1)} \, - \, {\mit D}_k^{(n)} \, ,
\qquad k,n \ge 0 \, .
\tag
$$

If we use the starting values
$$
D_0^{(n)} \; = \; s_n / \omega_n \, ,
\qquad n \in \N_0 \, ,
\tag
$$
the recurrence formula (9.5-5) produces the numerator of Drummond's
transformation (9.5-4), and if we use the starting values
$$
D_0^{(n)} \; = \; 1 / \omega_n \, ,
\qquad n \in \N_0 \, ,
\tag
$$
we obtain the denominator of Drummond's transformation (9.5-4).

Essentially the same computational algorithm, which was used for the
other sequence transformations of sections 7, 8, and 9 can also be used
in the case of Drummond's sequence transformation ${\cal D}_k^{(n)}
(s_n, \omega_n)$. Consequently, our approximation to the limit $s$ of
the sequence $\Seqn s$ to be transformed will be
$$
\{ s_0, \omega_0; s_1, \omega_1; \ldots ; s_m, \omega_m \} \; \to \;
{\cal D}_m^{(0)} (s_0, \omega_0) \, ,
\qquad m \in \N_0 \, .
\tag
$$

This means that we shall again use a sequence of transforms with
minimal superscripts and maximal subscripts as approximations to the
limit $s$:
$$
{\cal D}_0^{(0)} (s_0, \omega_0), \;
{\cal D}_1^{(0)} (s_0, \omega_0), \; \ldots , \;
{\cal D}_m^{(0)} (s_0, \omega_0), \; \ldots \; .
\tag
$$

Since the recursive computation of the sequence transformation ${\cal
D}_k^{(n)} (s_n, \omega_n)$ can be done in virtually the same way as in
the case of the Levin transformation, it is recommendable to
reformulate the recursive scheme (9.5-5) in the following way:
$$
D_j^{(n-j)} \; = \; D_{j-1}^{(n-j+1)} \, - \,
D_{j-1}^{(n-j)} \, , \qquad n \ge 1 \, , \quad 1 \le j \le n \, .
\tag
$$

Again, a single 1-dimensional array will be sufficient for the
computation of the quantities $D_j^{(n-j)}$, which are either numerator
or denominator sums of the sequence transformation (9.5-4), if the
$D_j^{(n-j)}$ are stored in a 1-dimensional array ${\mit d}$ according
to the following rule:
$$
D_{\nu}^{(n-\nu)} \; \to \; {\mit d} (n-\nu) \, ,
\qquad n \ge 0 \, , \quad 0 \le \nu \le n \, .
\tag
$$

With this convention the recurrence formula (9.5-10) can be
reformulated in terms of the elements of the array ${\mit d}$:
$$
{\mit d} (n-j) \; \gets \; {\mit d} (n-j+1) \, - \,
{\mit d} (n-j) \, , \qquad n \ge 1 \, , \quad 1 \le j \le n \, .
\tag
$$

If we compare Drummond's sequence transformation ${\cal D}_{k}^{(n)}
(s_n, \omega_n)$, eq. (9.5-4), with the ana\-logous sequence
transformations ${\cal L}_k^{(n)} (\beta, s_n, \omega_n)$, eq. (7.1-7),
${\cal S}_k^{(n)} (\beta, s_n, \omega_n)$, eq. (8.2-7), and ${\cal
M}_k^{(n)} (\gamma, s_n, \omega_n)$, eq. (9.2-6), we see that the
numerator and denominator sums of these transformations contain
additional $n$-dependent coefficients such as $(\beta+n+j)^{k-1}$,
$(\beta+n+j)_{k-1}$, or $(-\gamma - n - j)_{k-1}$, which are all of
order $O (n^{k-1})$.

Similarly, in the difference operator representation (9.5-3) for
Drummond's sequence transformation $\Delta^k$ acts only upon
$s_n/\omega_n$ and $1/\omega_n$ whereas in the analogous difference
operator representations (7.1-6), (8.2-5), and (9.2-5) for the other
transformations mentioned above weighted differences of $s_n/\omega_n$
and $1/\omega_n$ are formed.

This implies that in the case of Drummond's sequence transformation,
eq. (9.5-4), the sequence elements $s_n, \ldots, s_{n+k}$ and the
remainder estimates $\omega_n, \ldots, \omega_{n+k}$, which are needed
for the computation of ${\cal D}_{k}^{(n)} (s_n, \omega_n)$, all
contribute equally in the numerator and denominator sums.

This is not true in the case of the other sequence transformations
mentioned above since they contain additional weights. Essentially,
this means that the information contained in the sequence elements and
remainder estimates with larger indices will be emphasized more
strongly in the computation of these sequence transformations.

Since the later elements of a convergent sequence $\Seqn s$ are usually
closer to the limit $s$ than the earlier elements, it seems plausible
to expect that Drummond's sequence transformation, eq. (9.5-4), which
does not give special weight to the sequence elements and remainder
estimates with higher indices, should normally be somewhat less
powerful than the transformations ${\cal L}_k^{(n)} (\beta, s_n,
\omega_n)$, ${\cal S}_k^{(n)} (\beta, s_n, \omega_n)$, and ${\cal
M}_k^{(n)} (\gamma, s_n, \omega_n)$. We shall see later that this
assumption is indeed normally true.

\endAbschnittsebene

\neueSeite

\Abschnitt Brezinski's theta algorithm and related topics

\vskip - 2 \jot

\beginAbschnittsebene

\medskip

\Abschnitt The derivation of Brezinski's theta algorithm

\smallskip

\aktTag = 0

It is well known that Wynn's $\epsilon$ algorithm accelerates linear
convergence quite efficiently and that it is also able to sum even
wildly divergent series. However, the $\epsilon$ algorithm is not able
to accelerate logarithmic convergence. In the same way, Wynn's $\rho$
algorithm is certainly one of the better accelerators for logarithmic
convergence but fails to accelerate linear convergence and to sum
divergent series. It would certainly be desirable to modify either the
$\epsilon$ or the $\rho$ algorithm in such a way that the advantageous
features of the $\epsilon$ and the $\rho$ algorithm could be combined.
For that purpose, let us consider a recursive scheme of the following
type:
$$
\beginAligntags
" T_{-1}^{(n)} \; " = \; " 0 \, , \hfill
" T_0^{(n)} \; = \; s_n \, ,
\erhoehe\aktTag \\ \tag*{\tagnr a}
" T_{k+1}^{(n)} \; " = \; " T_{k-1}^{(n+1)} \, + \,
{\mit w}_k \, D_k^{(n)} \, , \qquad " k,n \in \N_0 \, .
\\ \tag*{\tagform\aktTagnr b}
\endAligntags
$$

$D_k^{(n)}$ is a quantity which depends upon one or several other
elements $T_{\kappa}^{(\nu)}$ of the table of this transformation. It
is assumed that the functional form of $D_k^{(n)}$ is known. The
quantity ${\mit w}_k$ is for the moment unspecified. Later, we will try
to derive an expression for ${\mit w}_k$ which will guarantee that the
above recursive scheme will lead to an acceleration of convergence.

The recursive scheme (10.1-1) contains the $\epsilon$ and the $\rho$
algorithm as special cases. If we choose ${\mit w}_k = 1$ together with
$$
D_k^{(n)} \; = \; \frac {1} {T_k^{(n+1)} - T_k^{(n)}} \, , \qquad k,n
\in \N_0 \, ,
\tag
$$
the recursive scheme (10.1-1) corresponds to Wynn's $\epsilon$
algorithm, eq. (4.2-1), and if we choose ${\mit w}_k = 1$ together with
$$
D_k^{(n)} \; = \; \frac {x_{n+k+1} - x_n} {T_k^{(n+1)} -
T_k^{(n)}} \, , \qquad k,n \in \N_0 \, ,
\tag
$$
the recursive scheme (10.1-1) corresponds to Wynn's $\rho$ algorithm,
eq. (6.2-2).

We now want to analyze how the quantity ${\mit w}_k$ has to be chosen
in order to guarantee that the sequence transformation $T_k^{(n)}$ will
lead to an acceleration of convergence.

In Wynn's $\epsilon$ or $\rho$ algorithm, only the transforms with even
subscripts are used as approximations to the limit. The transforms with
odd subscripts are only auxiliary quantities which diverge if the whole
process converges. Since either the $\epsilon$ or the $\rho$ algorithm
will be our starting point for the construction of a new sequence
transformation, we assume that $T_k^{(n)}$ behaves in the same way.
This means that only the transforms with even subscripts will be used
as approximations to the limit whereas the transforms with odd
subscripts are only auxiliary quantities which diverge if the
transforms with even subscripts converge.

Brezinski [26] argued that the exact numerical values of the transforms
with odd subscripts do not really matter as long as they diverge if the
whole process converges. Consequently, the most convenient choice for
${\mit w}_{2 k}$ in eq. (10.1-1b) would be to proceed as in Wynn's
$\epsilon$ or $\rho$ algorithm, i.e.,
$$
{\mit w}_{2 k} = 1 \, , \qquad k \in \N_0 \, .
\tag
$$

The parameters ${\mit w}_{2 k + 1}$ can be determined by requiring that
for fixed $k \in \N_0$ the sequence $\bigSeq {T_{2 k + 2}^{(n)}}$
should converge more rapidly than the sequence $\bigSeq {T_{2
k}^{(n+1)}}$ in the following sense:
$$
\lim_{n \to \infty} \; \frac
{\Delta T_{2 k + 2}^{(n)}}
{\Delta T_{2 k}^{(n+1)}} \; = \; 0 \, , \qquad k \in \N_0 \, .
\tag
$$

If we form in eq. (10.1-1b) the first difference with respect to $n$,
we see that condition (10.1-5) is automatically fulfilled if we choose
$$
{\mit w}_{2 k + 1} \; = \; - \lim_{n \to \infty} \, \frac
{\Delta T_{2 k}^{(n+1)}} {\Delta D_{2 k + 1}^{(n)}} \, ,
\qquad k \in \N_0 \, .
\tag
$$

Unfortunately, in situations of practical interest it will normally not
be be possible to compute this limit $n \to \infty$. As a manageable
alternative, Brezinski [26] suggested to use instead:
$$
{\mit w}_{2 k + 1}^{(n)} \; = \; - \, \frac
{\Delta T_{2 k}^{(n+1)}} {\Delta D_{2 k + 1}^{(n)}} \, ,
\qquad k,n \in \N_0 \, .\tag
$$

This choice together with eq. (10.1-4) leads to the following recursive
scheme for the sequence transformation $T_k^{(n)}$:
$$
\beginAligntags
" T_{-1}^{(n)} \; " = \; " 0 \, , \hfill
" T_0^{(n)} \; = \; s_n \, ,
\erhoehe\aktTag \\ \tag*{\tagnr a}
" T_{2 k + 1}^{(n)} \; " = \; " T_{2 k-1}^{(n+1)}
\, + \, D_{2 k}^{(n)} \, , \hfill "
\\ \tag*{\tagform\aktTagnr b}
" T_{2 k+2}^{(n)} \; " = \; " T_{2 k}^{(n+1)} \, - \,
\frac {\Delta T_{2 k}^{(n+1)}} {\Delta D_{2 k + 1}^{(n)}}
\, D_{2 k+1}^{(n)} \, , \qquad " k,n \in \N_0 \, .
\\ \tag*{\tagform\aktTagnr c}
\endAligntags
$$

If we choose in this recursive scheme $D_k^{(n)}$ according to eq.
(10.1-2) -- which corresponds to Wynn's $\epsilon$ algorithm -- we
obtain Brezinski's $\theta$ algorithm [26]:
$$
\beginAligntags
" \theta_{-1}^{(n)} \; " = \; " 0 \, , \hfill
" \theta_0^{(n)} \; = \; s_n \, ,
\erhoehe\aktTag \\ \tag*{\tagnr a}
" \theta_{2 k + 1}^{(n)} \; " = \; " \theta_{2 k-1}^{(n+1)}
\, + \, 1 / [\Delta \theta_{2 k}^{(n)}] \, , \hfill "
\\ \tag*{\tagform\aktTagnr b}
" \theta_{2 k+2}^{(n)} \; " = \; " \theta_{2 k}^{(n+1)} \, +
\, \frac
{[\Delta \theta_{2 k}^{(n+1)}] \,
[\Delta \theta_{2 k + 1}^{(n+1)}]}
{\Delta^2 \theta_{2 k+1}^{(n)}}
\, , \qquad " k,n \in \N_0 \, .
\\ \tag*{\tagform\aktTagnr c}
\endAligntags
$$

As usual it is assumed that the difference operator $\Delta$ acts only
upon the superscript $n$ and not upon the subscript $k$.

Brezinski's derivation of his $\theta$ algorithm, which is based upon
the the somewhat arbitrary choice (10.1-6), was purely experimental.
However, it was certainly a very successful experiment. In numerical
studies performed by Smith and Ford [29,30] it was demonstrated that
Brezinski's $\theta$ algorithm is a very powerful as well as a very
versatile sequence transformation since it is able to accelerate both
linear and logarithmic convergence quite efficiently and to sum even
wildly divergent series.

Unlike most of the other sequence transformations in this report,
Brezinski's $\theta$ algorithm was not derived via a model sequence. In
addition, the recursive scheme (10.1-9) is significantly more
complicated than the recursive schemes of most other nonlinear sequence
transformations. This should explain why only relatively little is
known about the theoretical properties of Brezinski's $\theta$
algorithm. In his second book, Brezinski showed that the $\theta$
algorithm is invariant under translation according to eq. (3.1-4) (see
p. 106 of ref. [20]). Also, Smith and Ford could prove that
$\theta_2^{(n)}$ accelerates linear convergence (see pp. 225 - 226 of
ref. [29]). Short discussions of the properties of the $\theta$
algorithm can be found in books by Brezinski [19,20] and Wimp [23].

Inspired by the success of his $\theta$ algorithm, Brezinski [27]
suggested to use the approach, which led to the derivation of the
$\theta$ algorithm, also in the case of other sequence transformations.
Since Wynn's $\epsilon$ algorithm -- which is the starting point for
Brezinski's $\theta$ algorithm -- and Wynn's $\rho$ algorithm are
structurally almost identical, it is a relatively obvious idea to try
to use Brezinski's concept for the construction of a new sequence
transformation which would be based upon Wynn's $\rho$ algorithm. This
can be accomplished quite easily. One only has to insert eq. (10.1-3)
instead of eq. (10.1-2) into eqs. (10.1-8b) and (10.1-8c) to obtain the
following recursive scheme:
$$
\beginAligntags
" \Theta_{-1}^{(n)} " " = \; " 0 \, , \hfill
\Theta_0^{(n)} \; = \; s_n \, , \qquad
\erhoehe\aktTag \\ \tag*{\tagnr a}
" \Theta_{2 k + 1}^{(n)} " " = \; " \Theta_{2 k-1}^{(n+1)}
\, + \, \frac {x_{n + 2 k + 1} - x_n}
{\Delta \Theta_{2 k}^{(n)}} \, , \hfill " \;
\\ \tag*{\tagform\aktTagnr b}
" \Theta_{2 k+2}^{(n)} " " = \; " \Theta_{2 k}^{(n+1)} \, -
\, \frac
{[x_{n + 2 k + 2} - x_n] \, [\Delta \Theta_{2 k}^{(n+1)}] \,
[\Delta \Theta_{2 k + 1}^{(n+1)}]}
{[x_{n + 2 k + 2} - x_{n+1}] \, [\Delta \Theta_{2 k+1}^{(n)}] -
[x_{n + 2 k + 1} - x_n] \, [\Delta \Theta_{2 k+1}^{(n+1)}]}
\, , \qquad \\
" k,n \; \in \; \N_0 \, . \span\omit \span\omit \span\omit
\\ \tag*{\tagform\aktTagnr c}
\endAligntags
$$

Numerical tests showed that this sequence transformation
$\Theta_k^{(n)}$ is more versatile than Wynn's $\rho$ algorithm, from
which it was derived, since it is able to accelerate linear convergence
and to sum some divergent series. However, it is unfortunately much
less efficient than Wynn's $\rho$ algorithm in the case of logarithmic
convergence, and it is also not particularly powerful in the case of
linear convergence or divergence. This shows that in the case of Wynn's
$\rho$ algorithm Brezinski's experimental choice (10.1-6) does not lead
to the same spectacular success as in the case of Wynn's $\epsilon$
algorithm since the resulting sequence transformation $\Theta_k^{(n)}$
-- although clearly more versatile than Wynn's $\rho$ algorithm -- is
not able to compete with other, more specialized sequence
transformations.

\medskip

\Abschnitt Programming Brezinski's theta algorithm

\smallskip

\aktTag = 0

A program for Brezinski's $\theta$ algorithm should have the same
features as the other programs described in this report. This means it
should read in the sequence elements $s_0, s_1, \ldots ,s_m , \ldots$
successively starting with $s_0$. After the input of each new sequence
element $s_m$ as many new elements $\theta_k^{(n)}$ should be computed
as is permitted by the recursive scheme (10.1-9). That new element
$\theta_k^{(n)}$, which has the largest even subscript $k$, should be
used as the new approximation to the limit of the sequence $\Seqn s$.

Let us arrange the elements of the $\theta$ table in a rectangular
scheme in such a way that the superscript $n$ indicates the row and the
subscript $k$ the column of the 2-dimensional array:

$$
\matrix{
\theta_0^{(0)} " \theta_1^{(0)} " \theta_2^{(0)} " \ldots "
\theta_n^{(0)} " \ldots \\
\theta_0^{(1)} " \theta_1^{(1)} " \theta_2^{(1)} " \ldots "
\theta_n^{(1)} " \ldots \\
\theta_0^{(2)} " \theta_1^{(2)} " \theta_2^{(2)} " \ldots "
\theta_n^{(2)} " \ldots \\
\theta_0^{(3)} " \theta_1^{(3)} " \theta_2^{(3)} " \ldots "
\theta_n^{(3)} " \ldots \\
\vdots    " \vdots    " \vdots    " \ddots " \vdots    " \ddots \\
\theta_0^{(n)} " \theta_1^{(n)} " \theta_2^{(n)} " \ldots "
\theta_n^{(n)} " \ldots \\
\vdots  " \vdots " \vdots " \ddots " \vdots  " \ddots }
\tag
$$

The entries in the first column of the array are the starting values
$\theta_0^{(n)} = s_n$ of the recursion according to eq. (10.1-9a). The
remaining elements of the $\theta$ table can be computed with the help
of the recurrence formulas (10.1-9b) and (10.1-9c). The 4 elements,
which are connected by the nonlinear recursion (10.1-9b), form the same
pattern in the $\theta$ table as the 4 elements of the $\epsilon$ table
which are connected by eq. (4.2-1b):

$$
\beginMatrix
\beginFormat &\Formel\links \endFormat
\+" " \theta_{2 k}^{(n)} \qquad " \theta_{2 k+1}^{(n)} "
\\ [1\jot]
\+" \theta_{2 k-1}^{(n+1)} \qquad " \theta_{2 k}^{(n+1)} \qquad " " \\
\endMatrix
\tag
$$

\medskip

The 6 elements, which are connected by the nonlinear recursion
(10.1-9c), form the following pattern in the $\theta$ table:
$$
\beginMatrix
\beginFormat &\Formel\links \endFormat
\+" " \theta_{2 k+1}^{(n)} \qquad " \theta_{2 k+2}^{(n)} "
\\ [1\jot]
\+" \theta_{2 k}^{(n+1)} \qquad
" \theta_{2 k+1}^{(n+1)} \qquad " " \\ [1\jot]
\+" \theta_{2 k}^{(n+2)} \qquad
" \theta_{2 k+1}^{(n+2)} \qquad " " \\
\endMatrix
\tag
$$

\medskip

These two patterns show that the recursions (10.1-9b) and (10.1-9c)
have to proceed along a relatively complicated path in the $\theta$
table. Let us assume that the sequence elements $s_0, s_1, \ldots ,
s_{m-1}$ had been read in and as many elements of the $\theta$ table
had been computed as it is permitted by the recursive scheme (10.1-9).
After the input of the next sequence element $s_m$ the string
$\theta_j^{(m - \Ent {3 j / 2})}$ with $0 \le j \le \Ent {(2 m + 1)/3}$
can be computed. Again, $\Ent {x}$ stands for the integral part of $x$,
i.e., the largest integer $\nu$ satisfying $\nu \le x$. In this
context, it is recommendable to rewrite the recursive scheme (10.1-9)
in the following way:
$$
\beginAligntags
" \theta_0^{(n)} \; " = \; " s_n \, , \qquad n \ge 0 \, , "
\erhoehe\aktTag \\ \tag*{\tagnr a}
" \theta_1^{(n-1)} \; " = \; "
1/[\theta_0^{(n)} - \theta_0^{(n-1)}]
\, , \qquad n \ge 1 \, , "
\\ \tag*{\tagform\aktTagnr b}
" \theta_{2 j}^{(n-3 j)} \; " = \; " \theta_{2 j -2}^{(n-3 j+1)}
\, + \, \frac
{ [\Delta \theta_{2 j - 2}^{(n-3 j+1)}] \,
[\Delta \theta_{2 j - 1}^{(n-3 j+1)}] }
{ \Delta^2 \theta_{2 j - 1}^{(n-3 j)} }
\, , \qquad n \ge 3 \, , \; 1 \le j \le \, \Ent {n/3} \, , " \quad
\\ \tag*{\tagform\aktTagnr c}
" \theta_{2 j + 1}^{(n-3 j-1)} \; " = \;
" \theta_{2 j-1}^{(n-3 j)}
\, + \, 1 / [\Delta \theta_{2 j}^{(n-3 j-1)}] \, , \qquad
n \ge 4 \, , \; j \le \Ent {(n-1)/3} \, . "
\\ \tag*{\tagform\aktTagnr d}
\endAligntags
$$

It follows either from this recursive scheme or equivalently from the
two geometric patterns (10.2-2) and (10.2-3) that for the computation
of the transform $\theta_{2 k}^{(n)}$ the sequence elements $s_n,
s_{n+1}, \ldots , s_{n + 3 k}$ have to be known. Thus, $\theta_{2
k}^{(n)}$ is a transformation of order $\ell = 3 k$.

The approximation to the limit will depend upon the index $m$ of the
last sequence element $s_m$ which was used in the recursive scheme
(10.2-4). If $m$ is a multiple of 3, $m = 3 \mu$, our approximation to
the limit will be the transformation
$$
\{ s_0, s_1, \ldots ,s_{3 \mu}\} \; \to \; \theta_{2 \mu}^{(0)} \, ,
\tag
$$
if we have $m = 3 \mu + 1$, our approximation to the limit will be
$$
\{ s_1, s_2, \ldots ,s_{3 \mu + 1}\} \; \to \; \theta_{2 \mu}^{(1)} \, ,
\tag
$$
and if we have $m = 3 \mu + 2$, our approximation to the limit will be
$$
\{ s_2, s_3, \ldots ,s_{3 \mu + 2}\} \; \to \; \theta_{2 \mu}^{(2)} \, .
\tag
$$

These three relationships can be combined into a single equation
yielding
$$
\left\{ s_{m - 3 \Ent {m/3}}, s_{m - 3 \Ent {m/3} + 1}, \ldots , s_m
\right\} \; \to \;
\theta_{2 \Ent {m/3}}^{(m - 3 \Ent {m/3})} \, .
\tag
$$

Because of the complicated structure of the nonlinear recursive scheme
(10.2-4) a single 1-dimensional array will not suffice for the
computation of the new string $\theta_j^{(n - \Ent {3 j /2})}$ with $0
\le j \le \Ent {(2 n +1)/ 3}$ after the input of the last sequence
element $s_n$. Two 1-dimensional arrays $A$ and $B$ together with three
auxiliary variables will be needed.

We shall use the convention that if the index $n$ of the last sequence
element $s_n$, which was read in, is even, $n = 2 m$, the new string
will be stored in the array $A$ according to the rule
$$
\theta_j^{(2 m - \Ent {3 j / 2})} \; \to \; A (j) \, ,
\qquad 0 \le j \le \Ent {(4 m + 1) / 3} \, ,
\tag
$$
and if the index $n$ of the last sequence element $s_n$ is odd, $n = 2
m + 1$, the new string will be stored in the array $B$ according to the
rule
$$
\theta_j^{(2 m - \Ent {3 j / 2} + 1)} \; \to \; B (j) \, ,
\qquad 0 \le j \le \Ent {(4 m + 3) / 3} \, .
\tag
$$

Let us now assume that the index $n$ of the last sequence element
$s_n$, which was read in, is even, i.e., $n = 2 m$. Before the
computation of the new string (10.2-9), which is to be stored in $A$,
the array $B$ contains the elements $\theta_j^{(2 m -\Ent {3 j/2}-1)}$
with $0 \le j \le \Ent {(4 m - 1)/3}$, whereas in $A$ the elements
$\theta_j^{(2 m -\Ent {3 j/2}-2)}$ with $0 \le j \le \Ent {(4 m -
3)/3}$ are stored. The recursive scheme (10.2-4) can then be expressed
in terms of the elements of the arrays $A$ and $B$ in the following
way:
$$
\beginAligntags
" A(0) \; " \gets \; " s_{2 m} \, , \qquad m \in \N_0 \, ,
\erhoehe\aktTag \\ \tag*{\tagnr a}
" A(1) \; " \gets \; " 1 / [ A(0) - B(0)] \, ,
\\ \tag*{\tagform\aktTagnr b}
" A(2 j) \; " \gets \; " A'(2 j - 2) \, + \, \frac
{[B(2 j - 2) - A'(2 j - 2)] \, [A(2 j -1) - B(2 j - 1)]}
{A(2 j - 1) - 2 B(2 j -1) + A'(2 j - 1)} \, , \\
" j \le \Ent {2 m / 3} \, , \span\omit \span\omit
\\ \tag*{\tagform\aktTagnr c}
" A(2 j + 1) \; " \gets \; " A'(2 j - 1) \, + \,
1/[B(2 j) - A(2 j)] \, ,
\qquad j \le \Ent {(2 m - 1) / 3} \, .
\\ \tag*{\tagform\aktTagnr d}
\endAligntags
$$

The primed array elements $A'(2 j - 2)$ and $A'(2 j - 1)$ refer to the
occupation of $A$ after the computation of the string $\theta_j^{(2 m
-\Ent {3 j/2}-2)}$ with $0 \le j \le \Ent {(4 m - 3)/3}$. Since these
elements are overwritten during the computation of the new string, they
have to be stored in auxiliary variables.

Let us now assume that the index $n$ of the last sequence element
$s_n$, which was read in, is odd, i.e., $n = 2 m+1$. This implies that
the r{\^ o}le of the two arrays $A$ and $B$ has to be interchanged.
Before the computation of the new string (10.2-10), which is to be
stored in $B$, the array $A$ contains the elements $\theta_j^{(2 m
-\Ent {3 j/2})}$ with $0 \le j \le \Ent {(4 m + 1)/3}$, and in $B$ the
elements $\theta_j^{(2 m -\Ent {3 j/2}-1)}$ with $0 \le j \le \Ent {(4
m - 1)/3}$ are stored. The recursive scheme (10.2-4) can then be
expressed in terms of the elements of the arrays $A$ and $B$ in the
following way:
$$
\beginAligntags
" B(0) \; " \gets \; " s_{2 m + 1} \, , \qquad m \in \N_0 \, ,
\erhoehe\aktTag \\ \tag*{\tagnr a}
" B(1) \; " \gets \; " 1 / [ B(0) - A(0)] \, ,
\\ \tag*{\tagform\aktTagnr b}
" B(2 j) " \; \gets \; " B'(2 j - 2) \, + \, \frac
{[A(2 j - 2) - B'(2 j - 2)] \, [B(2 j -1) - A(2 j - 1)]}
{B(2 j - 1) - 2 A(2 j -1) + B'(2 j - 1)} \, , \\
" j \le \Ent {(2 m + 1) / 3} \, , \span\omit \span\omit
\\ \tag*{\tagform\aktTagnr c}
" B(2 j + 1) \; " \gets \; " B'(2 j - 1) \, + \,
1/[A(2 j) - B(2 j)] \, ,
\qquad j \le \Ent { 2 m / 3} \, .
\\ \tag*{\tagform\aktTagnr d}
\endAligntags
$$

The primed array elements $B'(2 j - 2)$ and $B'(2 j - 1)$ refer to the
occupation of $B$ after the computation of the string $\theta_j^{(2 m
-\Ent {3 j/2}-1)}$ with $0 \le j \le \Ent {(4 m - 1)/3}$. Since these
elements are overwritten during the computation of the new string,
they have to be stored in auxiliary variables.

The following FORTRAN 77 subroutine THETA performs the recursive
computation of Brezinski's $\theta$ algorithm in two 1-dimensional
arrays $A$ and $B$ using the two recursive schemes (10.2-11) and
(10.2-12). It is safeguarded against an exact or approximate vanishing
of the denominators $\Delta^2 \theta_{2 j-1}^{(n-3 j)}$ and $\Delta
\theta_{2 j - 2}^{(n-3 j+1)}$ in eqs. (10.2-4c) and (10.2-4d) by using
two variables HUGE and TINY. The elements $s_n$ with $n = 0, 1, 2,
\ldots $ of the sequence to be transformed have to be computed in a DO
loop in the calling program. Whenever a new sequence element $s_n$ is
computed in the outer DO loop, this subroutine THETA has to be called
again and a new string $\theta_{j}^{(n-\Ent{3 j/2})}$ with $0 \le j \le
\Ent{(2 n+1)/3}$ will be calculated. The new sequence element $s_n$ is
read in via the variable SOFN and the approximation to the limit is
returned via the variable ESTLIM.

It is important to note that this subroutine THETA only calculates the
approximations to the limit according to eqs. (10.2-5) -- (10.2-7). The
convergence of the whole process has to be analyzed in the calling
program.

On pp. 368 - 370 of Brezinski's second book [20] the listing of a
FORTRAN IV program, which computes Brezinski's $\theta$ algorithm using
three 1-dimensional arrays, can be found.

\bigskip

\listenvon{theta.for}

\medskip

\Abschnitt The iteration of $\theta_2^{(n)}$

\smallskip

\aktTag = 0

In section 5 it was shown how Aitken's $\Delta^2$ algorithm, which
according to eq. (5.1-5) is identical with $\epsilon_2^{(n)}$, can be
iterated to give the sequence transformation ${\cal A}_k^{(n)}$. In the
same way, it was shown in section 6.3 that $\rho_2^{(n)}$ can be
iterated to give the sequence transformation ${\cal W}_k^{(n)}$. In
this section, we want to analyze how the transform $\theta_2^{(n)}$ can
be iterated. From the recursive scheme (10.1-9) we obtain the following
expression:
$$
\theta_2^{(n)} \; = \; s_{n+1} \, - \, \frac
{[\Delta s_n] \, [\Delta s_{n+1}] \, [\Delta^2 s_{n+1}]}
{[\Delta s_{n+2}] \, [\Delta^2 s_n] -
[\Delta s_n] \, [\Delta^2 s_{n+1}]}
\, , \qquad n \in \N_0 \, .
\tag
$$

It follows from this relationship that $\theta_2^{(n)}$, which is a
kind of weighted $\Delta^3$ process, is identical with Lubkin's $W$
transformation [40]. Many other representations for $\theta_2^{(n)}$
can be derived by suitable manipulations of eq. (10.3-1). Examples are:
$$
\beginAligntags
"\theta_2^{(n)} \; " = \; " \frac
{s_{n+1} \, [\Delta s_{n+2}] \, [\Delta^2 s_n] \, - \,
s_{n+2} \, [\Delta s_n] \, [\Delta^2 s_{n+1}]}
{[\Delta s_{n+2}] \, [\Delta^2 s_n] -
[\Delta s_n] \, [\Delta^2 s_{n+1}]}
\\ \tag
" " = \; " \frac
{ \Delta^2 \, [ s_{n+1} / \Delta s_n ] }
{ \Delta^2 \, [ 1 / \Delta s_n ] } \, .
\\ \tag
\endAligntags
$$

Comparison of eq. (10.3-3) with eqs. (2.4-8), (7.3-5), and ((9.5-4)
shows that $\theta_2^{(n)}$ is also a special case of Levin's $u$
transformation or Drummond's sequence transformation with $\omega_{n+1}
= \Delta s_n$:
$$
\beginAligntags
"\theta_2^{(n)} \; " = \;
" u_2^{(n+1)} (\beta, s_{n+1})
\\ \tag
" " = \; " {\cal D}_2^{(n+1)} (s_{n+1}, \Delta s_n) \, .
\\ \tag
\endAligntags
$$

In addition, a comparison of eqs. (5.1-12) and (10.3-3) shows that
$\theta_2^{(n)}$ may also be considered to be a generalization of
Aitken's $\Delta^2$ process.

If we want to iterate eq. (10.3-1), we have to use the following
recursive scheme:
$$
\beginAligntags
" {\cal J}_0^{(n)} \; " = \; " s_n \, , \qquad n \in \N_0 \, ,
\erhoehe\aktTag \\ \tag*{\tagnr a}
" {\cal J}_{k+1}^{(n)} \; " = \;
" {\cal J}_k^{(n+1)} \, - \, \frac
{[\Delta {\cal J}_k^{(n)}] \, [\Delta {\cal J}_k^{(n+1)}]
\, [\Delta^2 {\cal J}_k^{(n+1)}]}
{[\Delta {\cal J}_k^{(n+2)}] \, [\Delta^2 {\cal J}_k^{(n)}]
\, - \,
[\Delta {\cal J}_k^{(n)}] \, [\Delta^2 {\cal J}_k^{(n+1)}]}
\, , \qquad k,n \in \N_0 \, .
\\ \tag*{\tagform\aktTagnr b}
\endAligntags
$$

As usual, the difference operator $\Delta$ acts only upon the
superscript $n$ and not upon the subscript $k$. It follows from this
recursive scheme that for the computation of ${\cal J}_k^{(n)}$ the
sequence elements $s_n, s_{n+1}, \ldots, s_{n+ 3 k}$ have to be known.
Consequently, ${\cal J}_k^{(n)}$ is a transformation of order $3 k$. In
that respect, ${\cal J}_k^{(n)}$ is equivalent to $\theta_{2 k}^{(n)}$
which needs the same set $s_n, s_{n+1}, \ldots, s_{n+ 3 k}$ of sequence
elements for its computation.

In sections 13 and 14, we shall see that ${\cal J}_k^{(n)}$ is a
powerful sequence transformation which has similar properties as
$\theta_{2 k}^{(n)}$, i.e., it is able to accelerate linear and
logarithmic convergence and is also able to sum even wildly divergent
series. This may be considered to be an indirect confirmation of the
validity of Brezinski's choice (10.1-4) which is based upon the
assumption that the exact numerical values of the transforms $\theta_{2
k+1}^{(n)}$ do not really matter as long as they diverge if the
transforms $\theta_{2 k}^{(n)}$ converge. If we would replace the
4-term recursion (10.1-9b) by the 3-term recursion
$$
\theta_{2 k + 1}^{(n)} \; = \;
1 / [\Delta \theta_{2 k}^{(n)}] \, , \qquad k,n \in \N_0 \, ,
\tag
$$
then with this modified $\theta$ algorithm we would obtain
$$
\theta_{2 k}^{(n)} \; = \; {\cal J}_k^{(n)} \, .
\tag
$$

\medskip

\Abschnitt Programming the iterated theta algorithm

\smallskip

\aktTag = 0

A program, which computes the sequence transformation ${\cal
J}_k^{(n)}$, should have the same features as the other programs in
this report. This means it should read in the sequence elements $s_0,
s_1, \ldots ,s_m , \ldots$ successively starting with $s_0$. After the
input of each new sequence element $s_m$ as many new elements ${\cal
J}_k^{(n)}$ should be computed as is permitted by the recursive scheme
(10.3-6). That new element ${\cal J}_k^{(n)}$, which has the largest
subscript $k$, should be used as the new approximation to the limit of
the sequence $\Seqn s$.

Let us arrange the elements ${\cal J}_k^{(n)}$ in a rectangular scheme
in such a way that the superscript $n$ indicates the row and the
subscript $k$ the column of the 2-dimensional array:

$$
\matrix{
{\cal J}_0^{(0)} " {\cal J}_1^{(0)} " {\cal J}_2^{(0)} " \ldots " {\cal
J}_n^{(0)} " \ldots \\
{\cal J}_0^{(1)} " {\cal J}_1^{(1)} " {\cal J}_2^{(1)} " \ldots " {\cal
J}_n^{(1)} " \ldots \\
{\cal J}_0^{(2)} " {\cal J}_1^{(2)} " {\cal J}_2^{(2)} " \ldots " {\cal
J}_n^{(2)} " \ldots \\
{\cal J}_0^{(3)} " {\cal J}_1^{(3)} " {\cal J}_2^{(3)} " \ldots " {\cal
J}_n^{(3)} " \ldots \\
\vdots    " \vdots    " \vdots    " \ddots " \vdots    " \ddots \\
{\cal J}_0^{(n)} " {\cal J}_1^{(n)} " {\cal J}_2^{(n)} " \ldots " {\cal
J}_n^{(n)} " \ldots \\
\vdots  " \vdots " \vdots " \ddots " \vdots  " \ddots }
\tag
$$

The entries in the first column of the array are the starting values
${\cal J}_0^{(n)} = s_n$ of the recursion according to eq. (10.3-6a).
The remaining elements of the ${\cal J}$ table can be computed with the
help of the recurrence formula (10.3-6b). The 5 elements, which are
connected by the nonlinear recursion (10.3-6b), form the following
pattern:

$$
\beginMatrix
\beginFormat &\Formel\links \endFormat
\+ " {\cal J}_k^{(n)} \qquad " {\cal J}_{k+1}^{(n)} "
\\ [1\jot]
\+ " {\cal J}_k^{(n+1)} \qquad " " \\ [1\jot]
\+ " {\cal J}_k^{(n+2)} \qquad " " \\ [1\jot]
\+ " {\cal J}_k^{(n+3)} \qquad " " \\
\endMatrix
\tag
$$

\medskip

It was remarked earlier that $\theta_2^{(n)}$ may be considered to be a
generalization of Aitken's $\Delta^2$ process. In the same way, the
computational algorithm for the sequence transformation ${\cal
J}_k^{(n)}$ is essentially a generalization of the computational scheme
for Aitken's iterated $\Delta^2$ process ${\cal A}_k^{(n)}$ which was
discussed in section 5.2. For that purpose, the recursive scheme
(10.3-6) is rewritten in the following way:
$$
\beginAligntags
" {\cal J}_0^{(n)} \; " = \; " s_n \, , \qquad n \, \ge \, 0 \, ,
\erhoehe\aktTag \\ \tag*{\tagnr a}
" {\cal J}_{\ell}^{(n-3 \ell)} \; " = \;
" {\cal J}_{\ell-1}^{(n-3 \ell+1)} \, - \, \frac
{[\Delta {\cal J}_{\ell-1}^{(n-3 \ell)}] \,
[\Delta {\cal J}_{\ell-1}^{(n-3 \ell+1)}]
\, [\Delta^2 {\cal J}_{\ell-1}^{(n-3 \ell+1)}]}
{[\Delta {\cal J}_{\ell-1}^{(n-3 \ell+2)}] \,
[\Delta^2 {\cal J}_{\ell-1}^{(n-3 \ell)}]
\, - \,
[\Delta {\cal J}_{\ell-1}^{(n-3 \ell)}] \,
[\Delta^2 {\cal J}_{\ell-1}^{(n-3 \ell+1)}]} \, , \, \\
" n \ge 3 \, \qquad 1 \le \ell \le \Ent {n/3} \, .
\span\omit \span\omit
\\ \tag*{\tagform\aktTagnr b}
\endAligntags
$$

As usual, $\Ent x$ denotes the integral part of $x$, i.e., the largest
integer $\nu$ satisfying $\nu \le x$. It follows either from the
geometric pattern (10.4-2) or from this recursive scheme that after
the input of a new sequence element $s_m$ the string ${\cal
J}_{\mu}^{(m- 3 \mu)}$ with $0 \le \mu \le \Ent {m/3}$ can be computed.

Again, the approximation to the limit of the sequence to be transformed
depends upon the index $m$ of the last sequence element $s_m$ which was
read in. If $m$ is a multiple of 3, $m = 3 \mu$, our approximation to
the limit will be the transformation
$$
\{ s_0, s_1, \ldots ,s_{3 \mu}\} \; \to \;
{\cal J}_{\mu}^{(0)} \, ,
\tag
$$
if we have $m = 3 \mu + 1$, our approximation to the limit will be
$$
\{ s_1, s_2, \ldots ,s_{3 \mu + 1}\} \; \to \;
{\cal J}_{\mu}^{(1)} \, ,
\tag
$$
and if we have $m = 3 \mu + 2$, our approximation to the limit will be
$$
\{ s_2, s_3, \ldots ,s_{3 \mu + 2}\} \; \to \;
{\cal J}_{\mu}^{(2)} \, .
\tag
$$

With the help of the notation $\Ent x$ for the integral part of $x$
these three relationships can be combined into a single equation
yielding
$$
\left\{ s_{m - 3 \Ent {m/3}}, s_{m - 3 \Ent {m/3} + 1}, \ldots , s_m
\right\} \; \to \;
{\cal J}_{\Ent {m/3}}^{(m - 3 \Ent {m/3})} \, .
\tag
$$

The recursive scheme (10.4-3) -- or equivalently the geometric pattern
(10.4-2) -- looks relatively complicated. But nevertheless, it is
possible to perform the computation of the string ${\cal J}_{\ell}^{(n-
3 \ell)}$ with $0 \le \ell \le \Ent {n/3}$ in a single 1-dimensional
array ${\mit J}$ if the elements of the table of this sequence
transformation are stored according to the following rule:
$$
{\cal J}_{\Ent {\nu/3}}^{(n - \nu)} \; \to \; {\mit J}
(n-\nu) \, , \qquad n \ge 0 \, , \quad 0 \le \nu \le n \, .
\tag
$$

With this convention, the recursive scheme (10.4-3) can be reformulated
in terms of the elements of the array ${\mit J}$:
$$
\beginAligntags
" {\mit J} (n) \; " \gets \; " s_n \, ,
\qquad n \, \ge \, 0 \, ,
\erhoehe\aktTag \\ \tag*{\tagnr a}
" {\mit J} (m) \; " \gets \;
" {\mit J} (m+1) \, - \, \frac
{[\Delta {\mit J} (m)] \, [\Delta {\mit J} (m+1)]
\, [\Delta^2 {\mit J} (m+1)]}
{[\Delta {\mit J} (m+2)] \, [\Delta^2 {\mit J} (m)]
\, - \,
[\Delta {\mit J} (m)] \, [\Delta^2 {\mit J} (m+1)]}
\, , \\
" m \, = \, n - 3 \ell \, , \qquad n \, \ge \, 3 \, ,
\qquad 1 \, \le \, \ell \, \le \, \Ent {n/3} \, .
\span\omit \span\omit
\\ \tag*{\tagform\aktTagnr b}
\endAligntags
$$

The following FORTRAN 77 subroutine THETIT performs the recursive
computation of the iterated $\theta_2$ algorithm in a 1-dimensional
array ${\mit J}$ using the recursive scheme (10.4-3). THETIT is
safeguarded against an exact or approximate vanishing of the
denominator in eq. (10.4-3c) by using two variables HUGE and TINY. The
elements $s_n$ with $n = 0, 1, 2, \ldots $ of the sequence to be
transformed have to be computed in a DO loop in the calling program.
Whenever a new sequence element $s_n$ is computed in the outer DO loop
this subroutine THETIT has to be called again and a new string ${\cal
J}_{\ell}^{(n- 3 \ell)}$ with $0 \le \ell \le \Ent {n/3}$ will be
calculated. The new sequence element $s_n$ is read in via the variable
SOFN and the approximation to the limit is returned via the variable
ESTLIM.

It is important to note that THETIT only calculates the approximations
to the limit according to eqs. (10.4-4) -- (10.4-6). The convergence of
the whole acceleration or summation process has to be analyzed in the
calling program.

\bigskip

\listenvon{thetit.for}

\medskip

\endAbschnittsebene

\neueSeite

\Abschnitt On the derivation of theta-type algorithms

\vskip - 2 \jot

\beginAbschnittsebene

\medskip

\Abschnitt New sequence transformations based upon Aitken's iterated
$\Delta^2$ process

\smallskip

\aktTag = 0

In section 10.1 it was discussed how Brezinski's $\theta$ algorithm,
eq. (10.1-9), can be derived by modifying Wynn's $\epsilon$ algorithm,
eq. (4.2-1). The $\theta$ algorithm is a very powerful sequence
transformation. It accelerates linear convergence and sums divergent
series approximately as efficiently as the $\epsilon$ algorithm.
However, unlike the $\epsilon$ algorithm the $\theta$ algorithm is
also able to accelerate many logarithmically convergent sequences.
Consequently, it is frequently emphasized in the literature that
Brezinski's $\theta$ algorithm combines the advantageous features of
both Wynn's $\epsilon$ and Wynn's $\rho$ algorithm.

Brezinski [27] suggested to use his approach, which led to the $\theta$
algorithm, also in the case of other sequence transformations. This
will be done in this section. However, one should not expect that
Brezinski's approach will automatically lead to new sequence
transformations that are more useful than the transformations from
which they were derived. For instance, in the case of Wynn's $\rho$
algorithm, eq. (6.2-2), which is formally almost identical with the
$\epsilon$ algorithm, Brezinski's approach led to the sequence
transformation $\Theta_k^{(n)}$, eq. (10.1-10), which is much less
efficient than Wynn's $\rho$ algorithm in the case of logarithmic
convergence and which is also not very powerful in the case of linear
convergence or divergence. Consequently, the sequence transformation
$\Theta_k^{(n)}$ is practically useless although it is certainly more
versatile than the $\rho$ algorithm from which it was derived.

But even if Brezinski's $\theta$ concept does not automatically lead to
practically useful new sequence transformations, it should nevertheless
be worthwhile to investigate in which cases new sequence
transformations can be obtained that are at least in some sense better
than the transformations from which they were derived.

It follows from eqs. (10.1-1) and (10.1-8) that the essential step in
the derivation of the $\theta$ algorithm consists in replacing a
recursion of the general type
$$
f_n \; = \; a_n \, + \, b_n
\tag
$$
by a more complicated modified recursion
$$
f_n \; = a_n \, - \,
\frac {\Delta a_n} {\Delta b_n} \> b_n \, .
\tag
$$

Thus, a new sequence transformation can be constructed by replacing a
recursion of the type of eq. (11.1-1) by a recursion of the type of eq.
(11.1-2) in the recursive scheme which defines a given sequence
transformation. The new sequence transformation will have a more
nonlinear structure than the original transformation. In addition, such
a modification will normally increase the order of the transformation
by one.

The probably closest relative of Wynn's $\epsilon$ algorithm is
Aitken's iterated $\Delta^2$ process, eq. (5.1-15). This follows from
the fact that because of eq. (5.1-5) the $\epsilon$ algorithm may also
be considered to be a generalization of Aitken's $\Delta^2$ process,
eq. (5.1-4). Aitken's iterated $\Delta^2$ process and Wynn's $\epsilon$
algorithm have similar properties since they are both able to
accelerate linear convergence and to sum many divergent series but are
unable to accelerate logarithmic convergence. Consequently, it would be
interesting to see how Brezinski's $\theta$ concept works in the case
of Aitken's iterated $\Delta^2$ process.

Let us now assume that a sequence transformation is defined by the
following recursive scheme:
$$
\beginAligntags
" T_0^{(n)} \; " = \; " s_n \, , \qquad n \in \N_0 \, ,
\erhoehe\aktTag \\ \tag*{\tagnr a}
" T_{k+1}^{(n)} \; " = \; " T_k^{(n)} \, + \,
D_k^{(n)} \, , \qquad " k,n \in \N_0 \, .
\\ \tag*{\tagform\aktTagnr b}
\endAligntags
$$

As in section 10.1 it is assumed that $D_k^{(n)}$ is a quantity which
depends upon one or several elements of the table of this
transformation. If this recursive scheme is to be modified along the
lines of Brezinski's $\theta$ algorithm, a comparison with eqs.
(11.1-1) and (11.1-2) shows that it must be changed in the following
way:
$$
\beginAligntags
" T_0^{(n)} \; " = \; " s_n \, , \qquad n \in \N_0 \, ,
\erhoehe\aktTag \\ \tag*{\tagnr a}
" T_{k+1}^{(n)} \; " = \; " T_k^{(n)} \, - \,
\frac {\Delta T_k^{(n)}} {\Delta D_k^{(n)}} \> D_k^{(n)}
\, , \qquad " k,n \in \N_0 \, .
\\ \tag*{\tagform\aktTagnr b}
\endAligntags
$$

Aitken's iterated $\Delta^2$ process, eq. (5.1-15), is of the form of
eq. (11.1-3). Thus, if Aitken's iterated $\Delta^2$ process is modified
according to eq. (11.1-4), a new sequence transformation ${\cal
B}_k^{(n)}$ results which is defined by the following nonlinear
recursive scheme:
$$
\beginAligntags
" {\cal B}_0^{(n)} \; " = \; " s_n \, , \qquad n \in \N_0 \, ,
\erhoehe\aktTag \\ \tag*{\tagnr a}
" {\cal B}_{k+1}^{(n)} \; " = \;
" {\cal B}_k^{(n)} \, + \,
\frac
{\bigl[\Delta {\cal B}_k^{(n)}\bigr]^3 \,
\bigl[\Delta^2 {\cal B}_k^{(n+1)}\bigr]}
{\bigl[\Delta {\cal B}_k^{(n)}\bigr]^2 \,
\bigl[\Delta^2 {\cal B}_k^{(n+1)}\bigr]
\, - \,
\bigl[\Delta {\cal B}_k^{(n+1)}\bigr]^2 \,
\bigl[\Delta^2 {\cal B}_k^{(n)}\bigr]}
\, , \quad k,n \in \N_0 \, .
\\ \tag*{\tagform\aktTagnr b}
\endAligntags
$$

Again, it is assumed that the difference operator $\Delta$ acts only
upon the superscript $n$ and not upon the subscript $k$.

However, this is not the only possibility of modifying Aitken's
iterated $\Delta^2$ process. It follows from eq. (5.1-6) that the
recursive scheme for Aitken's iterated $\Delta^2$ process can also be
written in the following way:
$$
\beginAligntags
" {\cal A}_0^{(n)} \, " = " \, s_n \; , \qquad n \in \N_0 \, ,
\hfill \erhoehe\aktTag \\ \tag*{\tagnr a}
" {\cal A}_{k+1}^{(n)} \, " = " \, {\cal A}_k^{(n+1)} \, - \, \frac
{\bigl[\Delta {\cal A}_k^{(n)}\bigr]
\bigl[\Delta {\cal A}_k^{(n+1)}\bigr]}
{\Delta^2 {\cal A}_k^{(n)}}
\; , \qquad k,n \in \N_0 \, . \\
\tag*{\tagform\aktTagnr b}
\endAligntags
$$

This version of Aitken's iterated $\Delta^2$ process is a recursive
scheme of the following type:
$$
\beginAligntags
" T_0^{(n)} \; " = \; " s_n \, , \qquad n \in \N_0 \, ,
\erhoehe\aktTag \\ \tag*{\tagnr a}
" T_{k+1}^{(n)} \; " = \; " T_k^{(n+1)} \, + \,
D_k^{(n)} \, , \qquad " k,n \in \N_0 \, .
\\ \tag*{\tagform\aktTagnr b}
\endAligntags
$$

If we compare this recursive scheme with eqs. (11.1-1) and (11.1-2) we
see that it has to be modified in the following way:
$$
\beginAligntags
" T_0^{(n)} \; " = \; " s_n \, , \qquad n \in \N_0 \, ,
\erhoehe\aktTag \\ \tag*{\tagnr a}
" T_{k+1}^{(n)} \; " = \; " T_k^{(n+1)} \, - \,
\frac {\Delta T_k^{(n+1)}} {\Delta D_k^{(n)}} \> D_k^{(n)}
\, , \qquad " k,n \in \N_0 \, .
\\ \tag*{\tagform\aktTagnr b}
\endAligntags
$$

However, a modification of the second version of Aitken's iterated
$\Delta^2$ process, eq. (11.1-6), according to eq. (11.1-8) does not
produce a new sequence transformation since we obtain the recursive
scheme for the sequence transformation ${\cal J}_k^{(n)}$, eq.
(10.3-6), which was derived by iterating the expression for
$\theta_2^{(n)}$, eq. (10.3-1).

It follows from eq. (5.1-7) that there is another possibility of
rewriting the recursive scheme for Aitken's $\Delta^2$ process in a way
which would be suited for our purposes:
$$
\beginAligntags
" {\cal A}_0^{(n)} \, " = " \, s_n \; , \qquad n \in \N_0 \, ,
\hfill \erhoehe\aktTag \\ \tag*{\tagnr a}
" {\cal A}_{k+1}^{(n)} \, " = " \, {\cal A}_k^{(n+2)} \, - \, \frac
{\bigl[\Delta {\cal A}_k^{(n+1)}\bigr]^2}
{\Delta^2 {\cal A}_k^{(n)}}
\; , \qquad k,n \in \N_0 \, . \\
\tag*{\tagform\aktTagnr b}
\endAligntags
$$

This version of Aitken's iterated $\Delta^2$ process is a recursive
scheme with the following general structure,
$$
\beginAligntags
" T_0^{(n)} \; " = \; " s_n \, , \qquad n \in \N_0 \, ,
\erhoehe\aktTag \\ \tag*{\tagnr a}
" T_{k+1}^{(n)} \; " = \; " T_k^{(n+2)} \, + \,
D_k^{(n)} \, , \qquad " k,n \in \N_0 \, ,
\\ \tag*{\tagform\aktTagnr b}
\endAligntags
$$
which according to eqs. (11.1-1) and (11.1-2) has to be modified in the
following way:
$$
\beginAligntags
" T_0^{(n)} \; " = \; " s_n \, , \qquad n \in \N_0 \, ,
\erhoehe\aktTag \\ \tag*{\tagnr a}
" T_{k+1}^{(n)} \; " = \; " T_k^{(n+2)} \, - \,
\frac {\Delta T_k^{(n+2)}} {\Delta D_k^{(n)}} \> D_k^{(n)}
\, , \qquad " k,n \in \N_0 \, .
\\ \tag*{\tagform\aktTagnr b}
\endAligntags
$$

If the third version of Aitken's iterated $\Delta^2$ process, eq.
(11.1-9), is modified according to eq. (11.1-11), we obtain a new
sequence transformation ${\cal C}_k^{(n)}$ which is defined by the
following nonlinear recursive scheme:
$$
\beginAligntags
" {\cal C}_0^{(n)} \; " = \; " s_n \, , \qquad n \in \N_0 \, ,
\erhoehe\aktTag \\ \tag*{\tagnr a}
" {\cal C}_{k+1}^{(n)} \; " = \;
" {\cal C}_k^{(n+2)} \, + \,
\frac
{\bigl[\Delta {\cal C}_k^{(n+1)}\bigr]^2 \,
\bigl[\Delta {\cal C}_k^{(n+2)}\bigr]
\, \bigl[\Delta^2 {\cal C}_k^{(n+1)}\bigr]}
{\bigl[\Delta {\cal C}_k^{(n+1)}\bigr]^2 \,
\bigl[\Delta^2 {\cal C}_k^{(n+1)}\bigr]
\, - \,
\bigl[\Delta {\cal C}_k^{(n+2)}\bigr]^2 \,
\bigl[\Delta^2 {\cal C}_k^{(n)}\bigr]}
\, , \quad k,n \in \N_0 \, . \qquad
\\ \tag*{\tagform\aktTagnr b}
\endAligntags
$$

The recursive schemes (11.1-5) and (11.1-12) for the two new sequence
transformations ${\cal B}_k^{(n)}$ and ${\cal C}_k^{(n)}$ have the same
structure as the recursive scheme (10.3-6) for ${\cal J}_k^{(n)}$ since
these transformations are all weighted $\Delta^3$ methods.
Consequently, the two new sequence transformations can be computed in
the same way as ${\cal J}_k^{(n)}$. In this context it is recommendable
to rewrite the recursive scheme (11.1-5) in the following way:
$$
\beginAligntags
" {\cal B}_0^{(n)} \; " = \; " s_n \, , \qquad n \ge 0 \, ,
\erhoehe\aktTag \\ \tag*{\tagnr a}
" {\cal B}_j^{(n-3 j)} \; " = \;
" {\cal B}_{j-1}^{(n-3 j)} \, + \,
\frac
{\bigl[\Delta {\cal B}_{j-1}^{(n-3 j)}\bigr]^3 \,
\bigl[\Delta^2 {\cal B}_{j-1}^{(n-3 j+1)}\bigr]}
{\bigl[\Delta {\cal B}_{j-1}^{(n-3 j)}\bigr]^2 \,
\bigl[\Delta^2 {\cal B}_{j-1}^{(n-3 j+1)}\bigr]
\, - \,
\bigl[\Delta {\cal B}_{j-1}^{(n-3 j+1)}\bigr]^2 \,
\bigl[\Delta^2 {\cal B}_{j-1}^{(n-3 j)}\bigr]} \, , \\
" n \ge 3 \, , \qquad 1 \le j \le \Ent {n/3} \, .
\span\omit\span\omit \\ \tag*{\tagform\aktTagnr b}
\endAligntags
$$

As usual, $\Ent x$ denotes the integral part of $x$, i.e., the largest
integer $\nu$ satisfying $\nu \le x$. The recursive scheme (11.1-12)
should also be rewritten in the same way:
$$
\beginAligntags
" {\cal C}_0^{(n)} \; " = \; " s_n \, , \qquad n \ge 0 \, ,
\erhoehe\aktTag \\ \tag*{\tagnr a}
" {\cal C}_j^{(n-3 j)} \; " = \;
" {\cal C}_{j-1}^{(n-3 j+2)} \, + \,
\frac
{\bigl[\Delta {\cal C}_{j-1}^{(n-3 j+1)}\bigr]^2 \,
\bigl[\Delta {\cal C}_{j-1}^{(n-3 j+2)}\bigr]
\, \bigl[\Delta^2 {\cal C}_{j-1}^{(n-3 j+1)}\bigr]}
{\bigl[\Delta {\cal C}_{j-1}^{(n-3 j+1)}\bigr]^2 \,
\bigl[\Delta^2 {\cal C}_{j-1}^{(n-3 j+1)}\bigr]
\, - \,
\bigl[\Delta {\cal C}_{j-1}^{(n-3 j+2)}\bigr]^2 \,
\bigl[\Delta^2 {\cal C}_{j-1}^{(n-3 j)}\bigr]} \, , \\
" n \ge 3 \, , \qquad 1 \le j \le \Ent {n/3} \, .
\span\omit\span\omit \\ \tag*{\tagform\aktTagnr b}
\endAligntags
$$

It follows from eqs. (11.1-13) and (11.1-14) that after the input of a
new sequence element $s_m$ the strings ${\cal B}_{\mu}^{(m- 3 \mu)}$
and ${\cal C}_{\mu}^{(m- 3 \mu)}$ with $0 \le \mu \le \Ent {m/3}$ can
be computed.

The approximations to the limit of the sequence to be transformed
depend upon the index $m$ of the last sequence element $s_m$ which was
read in. Let ${\cal X}_k^{(n)}$ stand for either ${\cal B}_k^{(n)}$ or
${\cal C}_k^{(n)}$ . Then, if $m$ is a multiple of 3, $m = 3 \mu$, our
approximation to the limit will be the transformation
$$
\{ s_0, s_1, \ldots ,s_{3 \mu}\} \; \to \;
{\cal X}_{\mu}^{(0)} \, ,
\tag
$$
if we have $m = 3 \mu + 1$, our approximation to the limit will be
$$
\{ s_1, s_2, \ldots ,s_{3 \mu + 1}\} \; \to \;
{\cal X}_{\mu}^{(1)} \, ,
\tag
$$
and if we have $m = 3 \mu + 2$, our approximation to the limit will be
$$
\{ s_2, s_3, \ldots , s_{3 \mu + 2}\} \; \to \;
{\cal X}_{\mu}^{(2)} \, .
\tag
$$

These three relationships can be combined into a single equation
yielding
$$
\left\{ s_{m - 3 \Ent {m/3}}, s_{m - 3 \Ent {m/3} + 1}, \ldots , s_m
\right\} \; \to \;
{\cal X}_{\Ent {m/3}}^{(m - 3 \Ent {m/3})} \, .
\tag
$$

As in the case of the iterated $\theta_2$ algorithm ${\cal J}_k^{(n)}$,
eq. (10.3-6), the recursive computation of the strings ${\cal
B}_{j}^{(n- 3 j)}$ and ${\cal C}_{j}^{(n- 3 j)}$ with $0 \le j \le \Ent
{n/3}$ can be done in 1-dimensional arrays ${\mit B}$ and ${\mit C}$,
respectively, if the elements of the tables of these sequence
transformations are stored according to the following rule:
$$
{\cal X}_{\Ent {\nu/3}}^{(n - \nu)} \; \to \; {\mit X}
(n-\nu) \, , \qquad n \ge 0 \, , \quad 0 \le \nu \le n \, .
\tag
$$

Again, ${\cal X}_k^{(n)}$ stands for either ${\cal B}_k^{(n)}$ or
${\cal C}_k^{(n)}$, and ${\mit X}$ stands for the corresponding
1-dimensional array ${\mit B}$ or ${\mit C}$.

With this convention, the recursive scheme (11.1-13) can be
reformulated in terms of the elements of the array ${\mit B}$:
$$
\beginAligntags
" {\mit B} (n) \; " \gets \; " s_n \, , \qquad n \ge 0 \, ,
\erhoehe\aktTag \\ \tag*{\tagnr a}
" {\mit B} (m) \; " \gets \;
" {\mit B} (m) \, + \,
\frac
{[\Delta {\mit B} (m)]^3 \,
[\Delta^2 {\mit B} (m+1)]}
{[\Delta {\mit B} (m)]^2 \,
[\Delta^2 {\mit B} (m+1)]
\, - \,
[\Delta {\mit B} (m+1)]^2 \,
[\Delta^2 {\mit B} (m)]} \, , \\
" m = n - 3 j \, , \qquad n \ge 3 \, ,
\qquad 1 \le j \le \Ent {n/3} \, .
\span\omit\span\omit \\ \tag*{\tagform\aktTagnr b}
\endAligntags
$$

Similarly, the recursive scheme (11.1-14) can be reformulated in terms
of the elements of the array ${\mit C}$:
$$
\beginAligntags
" {\mit C} (n) \; " \gets \; " s_n \, , \qquad n \ge 0 \, ,
\erhoehe\aktTag \\ \tag*{\tagnr a}
" {\mit C} (m) \; " \gets \;
" {\mit C} (m+2) \, + \,
\frac
{[\Delta {\mit C} (m+1)]^2 \,
[\Delta {\mit C} (m+2)]
\, [\Delta^2 {\mit C} (m+1)]}
{[\Delta {\mit C} (m+1)]^2 \,
[\Delta^2 {\mit C} (m+1)]
\, - \,
[\Delta {\mit C} (m+2)]^2 \,
[\Delta^2 {\mit C} (m)]} \, , \\
" m = n - 3 j \, , \qquad n \ge 3 \, ,
\qquad 1 \le j \le \Ent {n/3} \, .
\span\omit\span\omit \\ \tag*{\tagform\aktTagnr b}
\endAligntags
$$

It follows from their recursive schemes that ${\cal J}_k^{(n)}$, eq.
(10.3-6), ${\cal B}_k^{(n)}$, eq. (11.1-5), and ${\cal C}_k^{(n)}$, eq.
(11.1-12), are all transformations of order $\ell = 3 k$. In that
respect, they are equivalent to $\theta_{2 k}^{(n)}$, eq. (10.1-9),
which needs the same set $s_n, \ldots, s_{n+3 k}$ of sequence elements
for its computation.

Numerical tests showed that the sequence transformations ${\cal
J}_k^{(n)}$, eq. (10.3-6), ${\cal B}_k^{(n)}$, eq. (11.1-5), and ${\cal
C}_k^{(n)}$, eq. (11.1-12), are clearly more versatile than Aitken's
iterated $\Delta^2$ process from which they were derived since they are
also able to accelerate logarithmic convergence. The sequence
transformation ${\cal B}_k^{(n)}$ is in general less powerful than the
other two transformations, and ${\cal J}_k^{(n)}$ is normally the most
powerful transformation, being roughly comparable with Brezinski's
$\theta$ algorithm.

\medskip

\Abschnitt New nonlinear sequence transformations obtained from linear
transformations

\smallskip

\aktTag = 0

In this section we want to construct new sequence transformations by
modifying the recursive schemes of the linear transformations
$\Lambda_k^{(n)} (\beta, s_n)$, eq. (7.3-20), ${\cal F}_k^{(n)}
(\alpha, s_n)$, eq. (8.4-11), and ${\cal P}_k^{(n)} (\zeta, s_n)$, eq.
(9.4-11), according to eqs. (11.1-1) and (11.1-2).

The recursive schemes (7.3-21), (8.4-12), and (9.4-12) for the linear
sequence transformations $\Lambda_k^{(n)} (\beta, s_n)$, ${\cal
F}_k^{(n)} (\alpha, s_n)$, and ${\cal P}_k^{(n)} (\zeta, s_n)$ are all
of the form of eq. (11.1-7). Consequently, these recursive schemes will
be modified according to eq. (11.1-8). Hence, in the case of the
recursive scheme (7.3-21) for the sequence transformation
$\Lambda_k^{(n)} (\beta, s_n)$ we obtain a new nonlinear sequence
transformation $\lambda_k^{(n)}$ which is defined by the following
recursive scheme:
$$
\beginAligntags
" \lambda_0^{(n)} \, " = \, " s_n \; , \qquad n \in \N_0 \, ,
\hfill \erhoehe\aktTag \\ \tag*{\tagnr a}
" \lambda_{k+1}^{(n)} \, " = \,
" \lambda_k^{(n+1)} \, - \, \frac
{(\beta + n) \,
\bigl[\Delta \lambda_k^{(n)} \bigr] \,
\bigl[\Delta \lambda_k^{(n+1)} \bigr]}
{(\beta + n + 1) \, \bigl[\Delta \lambda_k^{(n+1)} \bigr]
\, - \, (\beta + n) \, \bigl[\Delta \lambda_k^{(n)} \bigr] }
\; , \quad k,n \in \N_0 \, .
\\ \tag*{\tagform\aktTagnr b}
\endAligntags
$$

As usual, it is assumed here that the difference operator $\Delta$ acts
only upon the superscript $n$ and not upon the subscript $k$.

The nonlinear sequence transformation $\sigma_k^{(n)}$ is obtained by
modifying the recursive scheme (8.4-12) for ${\cal F}_k^{(n)} (\alpha,
s_n)$ :
$$
\beginAligntags
" \sigma_0^{(n)} \, " = \, " s_n \; , \qquad n \in \N_0 \, ,
\hfill \erhoehe\aktTag \\ \tag*{\tagnr a}
" \sigma_{k+1}^{(n)} \, " = \,
" \sigma_k^{(n+1)} \, - \, \frac
{(\alpha + n + k) \,
\bigl[\Delta \sigma_k^{(n)} \bigr] \,
\bigl[\Delta \sigma_k^{(n+1)} \bigr]}
{(\alpha + n + k + 1) \,
\bigl[\Delta \sigma_k^{(n+1)} \bigr]
\, - \, (\alpha + n + k) \,
\bigl[\Delta \sigma_k^{(n)} \bigr]}
\; , \quad k,n \in \N_0 \, . \qquad
\\ \tag*{\tagform\aktTagnr b}
\endAligntags
$$

Finally, a modification of the recursive scheme (9.4-12) for ${\cal
P}_k^{(n)} (\zeta, s_n)$ yields the nonlinear sequence transformation
$\mu_k^{(n)}$:
$$
\beginAligntags
" \mu_0^{(n)} \, " = \, " s_n \; , \qquad n \in \N_0 \, ,
\hfill \erhoehe\aktTag \\ \tag*{\tagnr a}
" \mu_{k+1}^{(n)} \, " = \, " \mu_k^{(n+1)} \, - \, \frac
{(\zeta + n - k) \,
\bigl[\Delta \mu_k^{(n)} \bigr] \,
\bigl[\Delta \mu_k^{(n+1)} \bigr]}
{(\zeta + n - k + 1) \,
\bigl[\Delta \mu_k^{(n+1)} \bigr]
\, - \, (\zeta + n - k) \,
\bigl[\Delta \mu_k^{(n)} \bigr]}
\; , \quad k,n \in \N_0 \, . \qquad
\\ \tag*{\tagform\aktTagnr b}
\endAligntags
$$

These new sequence transformations $\lambda_k^{(n)}$, $\sigma_k^{(n)}$,
and $\mu_k^{(n)}$ are all weighted $\Delta^2$ methods, i.e.,
modifications of Aitken's iterated $\Delta^2$ process, eq. (5.1-15). In
this context it may be of interest that Aitken's iterated $\Delta^2$
process can also be derived by modifying a linear sequence
transformation along the lines of Brezinski's $\theta$ algorithm. Let
us consider the following recursive scheme:
$$
\beginAligntags
" {\cal I} _0^{(n)} \, " = \, " s_n \; , \qquad n \in \N_0 \, ,
\hfill \erhoehe\aktTag \\ \tag*{\tagnr a}
" {\cal I}_{k+1}^{(n)} \, " = \, " {\cal I}_k^{(n)} \, + \,
\Delta {\cal I}_k^{(n)} \; , \quad k,n \in \N_0 \, .
\\ \tag*{\tagform\aktTagnr b}
\endAligntags
$$

Obviously, the sequence transformation ${\cal I}_k^{(n)}$, which is
defined by a recursive scheme of the type of eq. (11.1-3), can also be
written in the following way:
$$
{\cal I}_k^{(n)} \; = \; E^k \, s_n \, .
\tag
$$

Here, $E$ denotes the shift operator which is defined by eq. (2.2-4).
If the recursive scheme (11.2-4) is changed according to eq. (11.1-4),
we obtain Aitken's iterated $\Delta^2$ process, eq. (5.1-15).

The similarity of $\lambda_k^{(n)}$, $\sigma_k^{(n)}$, and
$\mu_k^{(n)}$ with Aitken's iterated $\Delta^2$ process implies that
these sequence transformations can be computed in the same way as
Aitken's iterated $\Delta^2$ process. For that purpose, it is
recommendable to rewrite the recursive scheme (11.2-1) for
$\lambda_k^{(n)}$ in the following way:
$$
\beginAligntags
" \lambda_0^{(n)} \, " = \, " s_n \, , \qquad n \ge 0 \, ,
\hfill \erhoehe\aktTag \\ \tag*{\tagnr a}
" \lambda_j^{(n - 2 j)} \, " = \,
" \lambda_{j-1}^{(n - 2 j+1)} \, - \, \frac
{(\beta + n - 2 j) \,
\bigl[\Delta \lambda_{j-1}^{(n - 2 j)} \bigr] \,
\bigl[\Delta \lambda_{j-1}^{(n - 2 j+1)} \bigr]}
{(\beta + n - 2 j + 1) \,
\bigl[\Delta \lambda_{j-1}^{(n - 2 j+1)} \bigr]
\, - \, (\beta + n - 2 j) \,
\bigl[\Delta \lambda_{j-1}^{(n - 2 j)} \bigr] }
\; , \\
" n \ge 2 \, , \qquad 1 \le j \le \Ent {n/2} \, .
\span\omit \span\omit
\\ \tag*{\tagform\aktTagnr b}
\endAligntags
$$

The recursive schemes (11.2-2) and (11.2-3) should also be rewritten in
the same way:
$$
\beginAligntags
" \sigma_0^{(n)} \, " = \, " s_n \, , \qquad n \ge 0 \, ,
\hfill \erhoehe\aktTag \\ \tag*{\tagnr a}
" \sigma_j^{(n - 2 j)} \, " = \,
" \sigma_{j-1}^{(n - 2 j+1)} \, - \, \frac
{(\alpha + n - j) \,
\bigl[\Delta \sigma_{j-1}^{(n - 2 j)} \bigr] \,
\bigl[\Delta \sigma_{j-1}^{(n - 2 j+1)} \bigr]}
{(\alpha + n - j + 1) \,
\bigl[\Delta \sigma_{j-1}^{(n - 2 j+1)} \bigr]
\, - \, (\alpha + n - j) \,
\bigl[\Delta \sigma_{j-1}^{(n - 2 j)} \bigr]}
\, , \\
" n \ge 2 \, , \qquad 1 \le j \le \Ent {n/2} \, .
\span\omit \span\omit
\\ \tag*{\tagform\aktTagnr b}
\\
" \mu_0^{(n)} \, " = \, " s_n \, , \qquad n \ge 0 \, ,
\hfill \erhoehe\aktTag \\ \tag*{\tagnr a}
" \mu_j^{(n - 2 j)} \, " = \,
" \mu_{j-1}^{(n - 2 j+1)} \, - \, \frac
{(\zeta + n - 3 j) \,
\bigl[\Delta \mu_{j-1}^{(n - 2 j)} \bigr] \,
\bigl[\Delta \mu_{j-1}^{(n - 2 j+1)} \bigr]}
{(\zeta + n - 3 j + 1) \,
\bigl[\Delta \mu_{j-1}^{(n - 2 j+1)} \bigr]
\, - \, (\zeta + n - 3 j) \,
\bigl[\Delta \mu_{j-1}^{(n - 2 j)} \bigr]}
\, , \\
" n \ge 2 \, , \qquad 1 \le j \le \Ent {n/2} \, .
\span\omit \span\omit
\\ \tag*{\tagform\aktTagnr b}
\endAligntags
$$

As in the case of Aitken's iterated $\Delta^2$ process, the
approximations to the limit of these transformations depend upon the
index $m$ of the last sequence element $s_m$ which was used in these
recursions. Let ${\cal X}_k^{(n)}$ now stand for any of the sequence
transformations $\lambda_k^{(n)}$, $\sigma_k^{(n)}$, or $\mu_k^{(n)}$.
Then, if $m$ is even, $m = 2 \mu$, our approximations to the limit of
the sequence are the transformations
$$
\{ s_0, s_1, \ldots ,s_{2 \mu}\} \; \to \; {\cal X}_{\mu}^{(0)} \, ,
\tag
$$
and if $m$ is odd, $m = 2 \mu + 1$, our approximations to the limit
will be
$$
\{ s_1, s_2, \ldots ,s_{2 \mu + 1}\} \; \to \;
{\cal X}_{\mu}^{(1)} \, .
\tag
$$

These two relationships can be combined into a single equation yielding
$$
\left\{ s_{m - 2 \Ent {m/2}}, s_{m - 2 \Ent {m/2} + 1}, \ldots , s_{m}
\right\} \, \to \;
{\cal X}_{\Ent {m/2}}^{(m - 2 \Ent {m/2})} \, .
\tag
$$

Only 1-dimensional arrays $\tilde \lambda$, $\tilde \sigma$, and
$\tilde \mu$ will be needed for the computation of the sequence
transformations $\lambda_k^{(n)}$, $\sigma_k^{(n)}$, and $\mu_k^{(n)}$
if the elements of the tables of these transformations are stored in
the same way as the elements of the Aitken table according to eq.
(5.2-6), i.e.,
$$
{\cal X}_{\Ent {\nu / 2}}^{(n - \nu)} \; \to \;
{\mit X}(n - \nu) \, ,
\qquad n \ge 0 \, , \quad 0 \le \nu \le n \, .
\tag
$$

Here, ${\cal X}_k^{(n)}$ stands for any of the three sequence
transformations $\lambda_k^{(n)}$, $\sigma_k^{(n)}$, and $\mu_k^{(n)}$,
and ${\mit X}$ stands for the corresponding 1-dimensional array, i.e.,
either for $\tilde \lambda$, $\tilde \sigma$, or for $\tilde \mu$. With
this convention, the recursive scheme (11.2-6) for $\lambda_k^{(n)}$
can be reformulated in terms of the elements of the 1-dimensional array
$\tilde \lambda$:
$$
\beginAligntags
" {\tilde \lambda} (n) \, " = \, " s_n \; , \qquad n \ge 0 \, ,
\hfill \erhoehe\aktTag \\ \tag*{\tagnr a}
" {\tilde \lambda} (m) \, " = \,
" {\tilde \lambda} (m+1) \, - \, \frac
{(\beta + m) \,
\bigl[\Delta {\tilde \lambda} (m) \bigr] \,
\bigl[\Delta {\tilde \lambda} (m+1) \bigr]}
{(\beta + m + 1) \,
\bigl[\Delta {\tilde \lambda} (m+1) \bigr]
\, - \, (\beta + m) \,
\bigl[\Delta {\tilde \lambda} (m) \bigr] }
\; , \\
" m \, = \, n - 2 j \, , \qquad n \ge 2 \, ,
\qquad 1 \le j \le \Ent {n/2} \, . \span\omit \span\omit
\\ \tag*{\tagform\aktTagnr b}
\endAligntags
$$

With the help of convention (11.2-12) the recursive schemes (11.2-7)
for $\sigma_k^{(n)}$ and (11.2-8) for $\mu_k^{(n)}$ can also be
reformulated in terms of the elements of the arrays $\tilde \sigma$ and
$\tilde \mu$ yielding similar expressions.

It follows from their recursive schemes that $\lambda_k^{(n)}$, eq.
(11.2-1), $\sigma_k^{(n)}$, eq. (11.2-2), and $\mu_k^{(n)}$, eq.
(11.2-3), are all transformations of order $\ell = 2 k$. In that
respect, they are equivalent to ${\cal A}_k^{(n)}$, eq. (5.1-15), and
$\epsilon_{2 k}^{(n)}$, eq. (4.2-1), which need the same set $s_n,
\ldots, s_{n+2 k}$ of sequence elements for their computation.

We shall see later that the nonlinear sequence transformations
$\lambda_k^{(n)}$, $\sigma_k^{(n)}$, and $\mu_k^{(n)}$ are more
versatile than the linear sequence transformations $\Lambda_k^{(n)}
(\beta, s_n)$, eq. (7.3-20), ${\cal F}_k^{(n)} (\alpha, s_n)$, eq.
(8.4-11), and ${\cal P}_k^{(n)} (\zeta, s_n)$, eq. (9.4-12), from which
they were derived, since they are not only able to accelerate
logarithmic convergence but also linear convergence. In addition, they
can also sum many divergent series. Numerical tests showed that
$\lambda_k^{(n)}$, eq. (11.2-1), is normally a more powerful sequence
transformation than $\sigma_k^{(n)}$, eq. (11.2-2), or $\mu_k^{(n)}$,
eq. (11.2-3).

\endAbschnittsebene

\neueSeite

\Abschnitt A theoretical analysis of sequence transformations

\vskip - 2 \jot

\beginAbschnittsebene

\medskip

\Abschnitt Germain-Bonne's formal theory of convergence acceleration

\smallskip

\aktTag = 0

The properties of linear and nonlinear sequence transformations are in
some sense complementary. Theoretically, linear sequence
transformations are now very well understood (see for instance refs.
[4] and [7-11]) but their power as well as their practical usefulness
is very limited. Nonlinear sequence transformations are often able to
achieve spectacular results, but theoretically, only relatively little
is known.

Any theory of nonlinear sequence transformations has to say something
about the two fundamental questions which arise in connection with
acceleration of convergence. Firstly, is the transformation under
consideration regular, i.e., will the transformed sequence $\Seq
{s_n^{\prime}}$ converge to the same limit as the original sequence
$\Seqn s$. Secondly, will the transformed sequence converge more
rapidly than the original sequence.

The first attempt to develop a general theory of the regularity and the
acceleration properties of nonlinear sequence transformations is due to
Germain-Bonne [33] who considered sequence transformations $G_k$ with
$k \in \N_0$ that are functions defined on vectors ${\bf x} = (x_1,
x_2, \ldots, x_{k+1})$, i.e., functions of the type $G_k : \R^{k+1} \to
\R$. Germain-Bonne postulated that these sequence transformations $G_k$
possess some very general properties such as continuity, homogeneity
and translativity. On the basis of these postulates Germain-Bonne could
formulate some conditions which guarantee the regularity of such a
sequence transformation. In addition, Germain-Bonne succeeded in
formulating a general criterion which decides whether a sequence
transformation $G_k$ accelerates linear convergence or not. A good
treatment of Germain-Bonne's formal theory of convergence acceleration
[33] can also be found in Wimp's book (see pp. 101 - 105 of ref. [23]).

The applicability of Germain-Bonne's theory in its original version is
quite limited. The reason is that Germain-Bonne treats sequence
transformations which depend upon $n$ only implicitly via the $k+1$
sequence elements $s_n, s_{n+1}, \ldots, s_{n+k}$, on which they act,
but not explicitly. Consequently, Germain-Bonne's theory is limited to
sequence transformations as for instance Wynn's $\epsilon$ algorithm,
eq. (4.2-1), Aitken's iterated $\Delta^2$ process, eq. (5.1-15), or
Brezinski's $\theta$ algorithm, eq. (10.1-9), which are all defined by
recursive schemes that do not depend explicitly upon $n$. It cannot be
applied in the case of a sequence transformation such as
$\lambda_k^{(n)}$, eq. (11.2-1), although it is a close relative of
Aitken's iterated $\Delta^2$ process, since its recursive scheme
depends explicitly on $n$. Consequently, in this section
Germain-Bonne's theory will be modified in such a way that it can be
applied to the sequence transformations of this report which mostly
depend explicitly upon $n$.

Let us therefore assume that for fixed $k \in \N_0$ a sequence
transformation $G_k^{(n)}$ is a function which is defined on vectors
${\bf x} = (x_1, x_2, \ldots, x_{k+1}) \in \R^{k+1}$ and which may
depend explicitly on $n \in \N_0$. In addition, we assume that such a
sequence transformation $G_k^{(n)} : \R^{k+1} \to \R$ possesses for
fixed $k \in \N_0$ and for all $n \in \N_0$ the following properties:

\medskip

\beginBeschreibung \zu \Laenge{(H-0):} \sl

\item {(H-0):} $G_k^{(n)}$ is defined and continuous on a subset ${\bf
X}^{(n)}$ of $\R^{k+1}$.

\item {(H-1):} $G_k^{(n)}$ is a homogeneous function of degree one.
This means that $G_k^{(n)}$ satisfies for arbitrary vectors ${\bf x}
\in {\bf X}^{(n)}$ and for all $\lambda \in \R$ such that $\lambda {\bf
x}$ is still an element of ${\bf X}^{(n)}$
$$
G_k^{(n)} (\lambda x_1, \lambda x_2, \ldots, \lambda x_{k+1})
\; = \; \lambda \> G_k^{(n)} (x_1, x_2, \ldots, x_{k+1}) \, .
\tag
$$
\item {(H-2):} $G_k^{(n)}$ is invariant under translation in the sense
of eq. (3.1-4). Consequently, for arbitrary $t \in \R$ and for
arbitrary vectors ${\bf x} \in {\bf X}^{(n)}$ we have
$$
G_k^{(n)} (x_1 + t, x_2 + t, \ldots, x_{k+1} + t)
\; = \; G_k^{(n)} (x_1, x_2, \ldots, x_{k+1}) \; + \; t \, .
\tag
$$
\item {(H-3):} A subset ${\bf X}^{(\infty)}$ of $\R^{k+1}$ exists such
that for every vector ${\bf x}$ $=$ $(x_1, x_2, \ldots, x_{k+1})$
belonging to this subset the limiting transformation
$$
G_k^{(\infty)} (x_1, x_2, \ldots, x_{k+1}) \; = \;
\lim_{n \to \infty} \> G_k^{(n)} (x_1, x_2, \ldots, x_{k+1})
\tag
$$
\item {} is uniquely defined and continuous. In addition, it is assumed
that the limiting transformation $G_k^{(\infty)}$ is also homogeneous
and invariant under translation according to (H-1) and (H-2).

\endBeschreibung

\medskip

The first three postulates (H-0) -- (H-2) are essentially identical
with the analogous postulates made by Germain-Bonne [33]. The main
difference is that here a sequence transformation $G_k^{(n)}$ may
depend explicitly upon $n$ and not only implicitly via the $k+1$
sequence elements $s_n$, $s_{n+1}$, $\ldots$, $s_{n+k}$ on which it
acts. Some of Germain-Bonne's results are based on the limiting
behaviour of a sequence transformation $G_k (s_n, s_{n+1}, \ldots,
s_{n+k})$ as $n \to \infty$. If we want to formulate analogous results
for a sequence transformation $G_k^{(n)}$, which may depend explicitly
on $n$, we have to require that the limit $n \to \infty$ can be
performed in the expression defining $G_k^{(n)}$ and that a unique
limiting transformation $G_k^{(\infty)}$ exists which is defined and
continuous on a suitable subset ${\bf X}^{(\infty)}$ of $\R^{k+1}$ and
which is also homogeneous and translative according to (H-1) and (H-2).

Conditions for the existence of the limiting transformation
$G_k^{(\infty)}$ were discussed by Smith and Ford (see p. 226 of ref.
[29]). Their analysis was based upon that version of the Moore -- Smith
theorem [93] which can be found in Gleason's book (see p. 256 of ref.
[94]). However, for our purposes it is probably simpler to postulate
the validity of (H-3), since in all cases, in which we shall have to do
such a limit $n \to \infty$, the existence of a limiting transformation
$G_k^{(\infty)}$ with the required properties will always be quite
obvious.

In (H-1) the restriction, that $\lambda {\bf x}$ has to be an element
of ${\bf X}^{(n)}$, is necessary. The reason is that nonlinear sequence
transformations are frequently not defined for constant sequences. In
such a case, $\lambda = 0$ has to be excluded because $G_k^{(n)}$ would
not be defined for the vector ${\bf x} = (0, 0, \ldots, 0)$.

Concerning (H-3) it should be noted that $G_k^{(n)}$ with $n$ being
finite and its limiting transformation $G_k^{(\infty)}$ are not
necessarily defined and continuous on the same subset of $\R^{k+1}$,
i.e., in general ${\bf X}^{(n)} \ne {\bf X}^{(\infty)}$. For instance,
if in the recursive scheme (11.2-1) for $\lambda_k^{(n)}$ the limit $n
\to \infty$ is performed according to eq. (12.1-3), we obtain the
second version of Aitken's iterated $\Delta^2$ process, eq. (11.1-6).
If we compare the subsets of $\R^3$, for which the explicit expressions
for ${\cal A}_1^{(n)}$ and $\lambda_1^{(n)}$ are continuous, we find
that Aitken's $\Delta^2$ process is defined for vectors ${\bf x} =
(x_1, x_2, x_3)$ which satisfy $x_1 - 2 x_2 + x_3 \neq 0$. A different
restriction is necessary in the case of $\lambda_1^{(n)}$ as long as
$n$ is finite.

It was remarked above that a given sequence transformation $G_k^{(n)}$
will normally not be continuous for arbitrary vectors ${\bf x} \in
\R^{k+1}$. If, however, a sequence transformation $G_k^{(n)}$ is
defined and continuous for all vectors ${\bf x} \in \R^{k+1}$, then the
regularity of this transformation can be proved quite easily.

\medskip

\beginEinzug \sl \parindent = 0 pt

\Auszug {\bf Theorem 12-1:} If a sequence transformation $G_k^{(n)}$ is
a continuous function on $\R^{k+1}$ for all $n \in \N_0$ and if its
limiting transformation $G_k^{(\infty)}$ is also continuous on
$\R^{k+1}$ and satisfies (H-1) and (H-2), then $G_k^{(n)}$ is regular,
i.e., it preserves the limit $s$ of every convergent sequence $\Seqn s$.

\endEinzug

\medskip

\noindent {\it Proof: } Since $G_k^{(n)}$ is continuous on $\R^{k+1}$
for all $n \in \N_0$ and since $\Seqn s$ converges to some limit $s$,
we have by continuity
$$
\lim_{n \to \infty} \> G_k^{(n)} (s_n, s_{n+1}, \ldots, s_{n+k})
\; = \; G_k^{(\infty)} (s, s, \ldots, s) \, .
\tag
$$

Since (H-2) remains valid as $n \to \infty$, we may conclude that for
all constant sequences $s, s, s, \ldots $ we have
$$
G_k^{(\infty)} (s, s, \ldots, s) \; = \; s \, + \,
G_k^{(\infty)} (0, 0, \ldots, 0) \, .
\tag
$$

In the same way, (H-1) remains valid as $n \to \infty$. Thus, we may
conclude
$$
G_k^{(\infty)} (0, 0, \ldots, 0) \; = \; 0 \, .
\tag
$$

Hence, it follows from eqs. (12.1-4) - (12.1-6) that a sequence
transformation $G_k^{(n)}$, which satisfies the above assumptions, is
regular.

\medskip

Unfortunately, theorem 12-1 will be of little use in the case of
nonlinear sequence transformations which are in general nonregular. A
nonlinear sequence transformation $G_k^{(n)}$ is normally a rational
function of the $k+1$ sequence elements $s_n, s_{n+1}, \ldots, s_{n+k}$
which are used for its computation. Since rational functions have
poles, we cannot expect that nonlinear sequence transformations will be
continuous on $\R^{k+1}$. Consequently, in the case of a nonlinear
sequence transformation the convergence of an arbitrary sequence $\Seqn
s$ to some limit $s$ does not imply that the transformed sequence
converges at all, let alone to the same limit. In addition, eqs.
(12.1-5) and (12.1-6) need not be valid since nonlinear sequence
transformations are not necessarily defined for constant sequences.

\medskip

\beginEinzug \sl \parindent = 0 pt

\Auszug {\bf Theorem 12-2:} Let $\D^{k+1}$ be the set of vectors ${\bf
x} \in \R^{k+1}$ with distinct components. Every sequence
transformation $G_k^{(n)}$, which is a continuous function for all
vectors ${\bf x}\in \D^{k+1}$ and which also satisfies (H-1) and (H-2),
can be expressed in the following way:
$$
G_k^{(n)} (x_1, x_2, \ldots, x_{k+1}) \; = \;
x_1 \, + \, (x_2 - x_1) \>
g_k^{(n)} \biggl( \frac {x_3 - x_2}{x_2 - x_1}, \ldots ,
\frac {x_{k+1} - x_k}{x_k - x_{k-1}} \biggr) \, .
\tag
$$

The associated transformation $g_k^{(n)}$ which is defined and
continuous on a subset of $\F^{k-1}$, the set of vectors ${\bf y} \in
\R^{k-1}$ with nonzero components, is given by
$$
g_k^{(n)} \biggl( \frac {x_3 - x_2}{x_2 - x_1}, \ldots ,
\frac {x_{k+1} - x_k}{x_k - x_{k-1}} \biggr) \; = \;
G_k^{(n)} \biggl( 0, 1, 1 + \frac {x_3 - x_2}{x_2 - x_1},
\ldots , \sum_{j=0}^{k-1} \> \prod_{i=0}^{j-1}
\frac {x_{i+3} - x_{i+2}}{x_{i+2} - x_{i+1}} \biggr) \, .
\tag
$$

\endEinzug

\medskip

\noindent {\it Proof: } According to (H-2) we can subtract $x_1$ from
$G_k^{(n)}$, and according to (H-1) we can divide the $k+1$ arguments
of $G_k^{(n)}$ by $x_2 - x_1$. This yields:
$$
G_k^{(n)} (x_1, x_2, \ldots, x_{k+1}) \; = \;
x_1 \, + \, (x_2 - x_1) \>
G_k^{(n)} \biggl(0, 1, \frac {x_3 - x_1} {x_2 - x_1},
\ldots, \frac {x_{k+1} - x_1} {x_2 - x_1} \biggr) \, .
\tag
$$

We now need the following relationship which can be proved quite easily:
$$
\frac {x_{\nu+1} - x_1} {x_2 - x_1} \; = \;
\sum_{j=0}^{\nu - 1} \> \prod_{i=0}^{j-1} \>
\frac {x_{i+3} - x_{i+2}}{x_{i+2} - x_{i+1}}
\, , \qquad \nu \in \N_0 \, .
\tag
$$

If we insert eq. (12.1-10) into eq. (12.1-9) we obtain eqs. (12.1-7)
and (12.1-8). The continuity of $g_k^{(n)}$ on the subset of vectors
${\bf y} = (y_1, y_2, \ldots, y_{k-1}) \in \F^{k-1}$, which are
generated from vectors ${\bf x} = (x_1, x_2, \ldots, x_{k+1}) \in
\D^{k+1}$ via the relationship
$$
y_{\mu} \; = \; (x_{\mu+2} - x_{\mu+1}) \, / \,
(x_{\mu+1} - x_{\mu}) \, , \qquad 1 \le \mu \le k - 1 \, ,
\tag
$$
follows via eq. (12.1-7) from the continuity of $G_k^{(n)}$ on
$\D^{k+1}$. However, an associated transformation $g_k^{(n)}$ in the
sense of eqs. (12.1-7) and (12.1-8) is not defined for all vectors
${\bf y} \in \F^{k-1}$. For instance, if we choose ${\bf y} = (1, -1,
1, -1, \ldots)$, which clearly belongs to $\F^{k-1}$, we are not able
to find a vector ${\bf x} = (x_1, x_2, \ldots, x_{k+1}) \in \D^{k+1}$
which satisfies eq. (12.1-11).

\medskip

With the help of theorem 12-2 it can be analyzed for which convergent
sequences $\Seqn s$ a sequence transformation $G_k^{(n)}$ satisfying
(H-0) -- (H-3) will be regular.

\medskip

\beginEinzug \sl \parindent = 0 pt

\Auszug {\bf Theorem 12-3:} Let $G_k^{(n)}$ be a sequence
transformation which satisfies theorem 12-2. This means that according
to eqs. (12.1-7) and (12.1-8) an associated transformation $g_k^{(n)}$
exists which is defined and continuous on that subset of vectors ${\bf
y} = (y_1, y_2, \ldots, y_{k-1}) \in \F^{k-1}$ which is generated by
all vectors ${\bf x} = (x_1, x_2, \ldots, x_{k+1}) \in \D^{k+1}$
according to eq. (12.1-11). If the limiting associated transformation,
which is defined by
$$
g_k^{(\infty)} (y_1, y_2, \ldots , y_{k-1}) \; = \;
\lim_{n \to \infty} g_k^{(n)} (y_1, y_2, \ldots , y_{k-1}) \, ,
\tag
$$
is also defined and continuous on the same subset of $\F^{k-1}$ as
$g_k^{(n)}$, then $G_k^{(n)}$ will preserve the limit $s$ of every
convergent sequence $\Seqn s$ having the following properties:

\beginBeschreibung \zu \Laenge{(ii):}

\item{(i) \hfill :} For sufficiently large values of $n \in \N_0$ the
sequence elements $s_n$ are all distinct.
\item{(ii):} For sufficiently large values of $n \in \N_0$ the ratios
$\Delta s_{n+1} / \Delta s_n$ all satisfy the inequality
$$
c \le \vert \Delta s_{n+1} / \Delta s_n \vert \le c^{\prime} \, ,
\qquad 0 < c < c^{\prime} < \infty \, .
\tag
$$

\endBeschreibung

\endEinzug

\medskip

\noindent {\it Proof: } It follows from (i) that $\Delta s_n \ne 0$ for
sufficiently large values of $n$. Consequently, $G_k^{(n)}$ can be
rewritten for sufficiently large values of $n$ according to eq.
(12.1-7) yielding
$$
G_k^{(n)} (s_n, s_{n+1}, \ldots, s_{n+k}) \; = \;
s_n \, + \, \Delta s_n \>
g_k^{(n)} \biggl( \frac {\Delta s_{n+1}}{\Delta s_n}, \ldots ,
\frac {\Delta s_{n+k-1}}{\Delta s_{n+k-2}} \biggr) \, .
\tag
$$

According to eq. (12.1-8) the associated transformation is given by
$$
g_k^{(n)} \biggl( \frac {\Delta s_{n+1}}{\Delta s_n}, \ldots ,
\frac {\Delta s_{n+k-1}}{\Delta s_{n+k-2}} \biggr) \; = \;
G_k^{(n)} \biggl( 0, 1, 1 + \frac {\Delta s_{n+1}}{\Delta s_n},
\ldots , \sum_{j=0}^{k-1} \> \prod_{i=0}^{j-1}
\frac {\Delta s_{n+i+1}}{\Delta s_{n+i}} \biggr) \, .
\tag
$$

$G_k^{(n)}$ will preserve the limit of $\Seqn s$ if the second term on
the right-hand side of eq. (12.1-14) vanishes as $n \to \infty$. Since
$\Seqn s$ converges, $\Delta s_n$ vanishes as $n \to \infty$ and we
only have to show that the associated transformation $g_k^{(n)}$
remains bounded as $n \to \infty$. The limiting associated
transformation $g_k^{(\infty)}$ is by assumption defined and continuous
on that subset of $\F^{k-1}$ which is generated from vectors ${\bf x}
\in \D^{k+1}$ according to eq. (12.1-11). Continuity of a given
function on a certain subset implies the boundedness of this function
for all bounded elements belonging to this subset. It follows from (i)
that for sufficiently large values of $n$ every string $s_n, s_{n+1},
\ldots, s_{n+k}$ belongs to $\D^{k+1}$, and it follows from (ii) that
for sufficiently large values of $n$ the $k-1$ arguments of $g_k^{(n)}$
in eq. (12.1-15) are all nonzero and bounded. This implies that
$g_k^{(n)}$ will remain bounded as $n \to \infty$. This concludes the
proof of theorem 12-3.

\medskip

Next, a criterion will be formulated which decides whether a sequence
transformation $G_k^{(n)}$, which may depend explicitly upon $n$, is
able to accelerate linear convergence or not. This is probably the most
important result of our adaptation of Germain-Bonne's formal theory of
convergence acceleration [33].

\medskip

\beginEinzug \sl \parindent = 0 pt

\Auszug {\bf Theorem 12-4:} Let us assume that a sequence $\Seqn s$
converges linearly to some limit $s$, i.e.,
$$
\lim_{n \to \infty} \> \frac {s_{n+1} - s}{s_n - s} \; = \; \rho
\, , \qquad 0 < |\rho| < 1 \, .
\tag
$$
Then, a necessary and sufficient condition that a sequence
transformation $G_k^{(n)} : \R^{k+1} \to \R$ accelerates the
convergence of $\Seqn s$ is that its associated transformation
$g_k^{(n)}$ satisfies:
$$
\lim_{n \to \infty} \>
g_k^{(n)} (\rho_n, \rho_{n+1}, \ldots, \rho_{n+k-2}) \; = \;
g_k^{(\infty)} (\rho, \rho, \ldots, \rho) \; = \;
\frac {1} {1 - \rho} \, .
\tag
$$
Here, $\Seqn {\rho}$ is an arbitrary sequence which converges to
$\rho$. The above statement can also be formulated in terms of the
limiting transformation $G_k^{(\infty)}$:
$$
\beginAligntags
" \lim_{n \to \infty} G_k^{(n)} \biggl(0, 1, 1+\rho_n , \ldots,
\sum_{j=0}^{k-1} \, \prod_{i=0}^{j-1} \, \rho_{n+i} \biggr) \\
" = \; G_k^{(\infty)} \biggl(0, 1, 1+\rho, \ldots,
\sum_{j=0}^{k-1}\rho^j \biggr) \; = \; \frac {1} {1 - \rho} \, ,
\qquad 0 < \vert \rho \vert < 1 \, .
\\ \tag
\endAligntags
$$

\endEinzug

\medskip

\noindent {\it Proof: } Since $\Seqn s$ converges linearly, it follows
from eqs. (2.6-3) and (2.6-4) that $\Delta s_n \sim \rho^n$ as $n \to
\infty$. Thus, $G_k^{(n)}$ can for sufficiently large values of $n$ be
rewritten according to eq. (12.1-14). If we subtract $s$ from both
sides of eq. (12.1-14) and divide the resulting expression by $s_n - s$
we obtain:
$$
\frac
{G_k^{(n)} (s_n, s_{n+1}, \ldots, s_{n+k}) - s} {s_n - s}
\; = \; 1 \, + \, \frac {\Delta s_n} {s_n - s} \,
g_k^{(n)} \biggl( \frac {\Delta s_{n+1}}{\Delta s_n}, \ldots ,
\frac {\Delta s_{n+k-1}}{\Delta s_{n+k-2}} \biggr) \, .
\tag
$$

According to eq. (2.6-6) the convergence of a sequence $\Seqn s$ to
its limit $s$ is accelerated by a sequence transformation $G_k^{(n)}$ if
$$
\lim_{n \to \infty} \> \frac
{G_k^{(n)} (s_n, s_{n+1}, \ldots, s_{n+k}) - s} {s_n - s}
\; = \; 0 \, .
\tag
$$

Hence , if we perform the limit $n \to \infty$ in eq. (12.1-19), the
left-hand side has to vanish if $G_k^{(n)}$ accelerates convergence.
Thus, we only have to investigate under which conditions the right-hand
side of eq. (12.1-19) also vanishes as $n \to \infty$. With the help of
eq. (12.1-16) we find:
$$
\lim_{n \to \infty} \> \frac {\Delta s_n} {s_n - s} \; = \;
\lim_{n \to \infty} \> \frac {s_{n+1} - s}{s_n - s} \; - \, 1
\; = \; \rho - 1 \, .
\tag
$$

Next, we observe that because of the equivalence of eqs. (2.6-3) and
(2.6-4) for sequences satisfying eq. (12.1-16) the arguments
$$
\rho_n \; = \; \Delta s_{n+1} / \Delta s_n \, ,
\qquad n \in \N_0 \, ,
\tag
$$
of the associated transformation $g_k^{(n)}$ in eq. (12.1-19) converge
to $\rho$ as $n \to \infty$. Now, the right-hand side of eq. (12.1-19)
can only vanish as $n \to \infty$ if the limiting associated
transformation $g_k^{(\infty)}$ satisfies eq. (12.1-17). In the same
way, if eq. (12.1-17) is satisfied, then because of eq. (12.1-21) the
right-hand side of eq. (12.1-19) vanishes. Consequently, the validity
of eq. (12.1-17) is equivalent to the statement that $G_k^{(n)}$
accelerates linear convergence.

\medskip

The sequence $0, 1, 1 + \rho, \ldots$, on which the limiting
transformation $G_k^{(\infty)}$ acts in eq. (12.1-18), is apart from
its first element and apart from a shift of the indices identical with
the sequence of partial sums of the geometric series, eq. (2.6-2). This
is best seen by rewriting the elements of this sequence in the
following way:
$$
\sigma_n (\rho) \; = \; \frac {1 - \rho^n} {1 - \rho} \; = \;
\sum_{\nu=0}^{n-1} \, \rho^{\nu} \, , \qquad 0 < |\rho| < 1 \, ,
\quad n \in \N_0 \, .
\tag
$$

Obviously, this shifted sequence also converges to $1 / (1 - \rho)$ as
$n \to \infty$. Hence, we see that theorem 12-4, which decides whether
a given sequence transformation $G_k^{(n)}$ accelerates linear
convergence or not, requires that the limiting transformation
$G_k^{(\infty)}$ is exact for the shifted sequence $\Seq {\sigma_n
(\rho)}$ of partial sums of the geometric series. This fact indicates
that there is a close connection between the exactness of a sequence
transformation for the partial sums of the geometric series and its
ability of accelerating linear convergence.

\medskip

\Abschnitt Applications of Germain-Bonne's theory

\smallskip

\aktTag = 0

In this section, the properties of certain sequence transformations
will be analyzed with the help of Germain-Bonne's formal theory of
convergence acceleration. The most interesting feature of
Germain-Bonne's theory is its treatment of the acceleration of linear
convergence. Consequently, theorem 12-4 and related questions such as
the exactness of a sequence transformation for the geometric series
will be emphasized in this section.

First, we want to investigate the sequence transformations
$\lambda_k^{(n)}$, eq. (11.2-1), $\sigma_k^{(n)}$, eq. (11.2-2), and
$\mu_k^{(n)}$, eq. (11.2-3), which are close relatives of Aitken's
iterated $\Delta^2$ process, eq. (5.1-15). However, it follows from
their recursive schemes (11.2-1) - (11.2-3) that unlike Aitken's
$\Delta^2$ process ${\cal A}_1^{(n)}$, eq. (5.1-4), which is by
construction exact for the geometric series, $\lambda_1^{(n)}$,
$\sigma_1^{(n)}$, and $\mu_1^{(n)}$ are not exact for the geometric
series. Consequently, it is not obvious whether $\lambda_1^{(n)}$,
$\sigma_1^{(n)}$, and $\mu_1^{(n)}$, accelerate linear convergence.

\medskip

\beginEinzug \sl \parindent = 0 pt

\Auszug {\bf Theorem 12-5:} The sequence transformations
$\lambda_k^{(n)}$, eq. (11.2-1), $\sigma_k^{(n)}$, eq. (11.2-2), and
$\mu_k^{(n)}$, eq. (11.2-3), accelerate linear convergence if and only
if Aitken's iterated $\Delta^2$ process ${\cal A}_k^{(n)}$, eq.
(5.1-15), accelerates linear convergence.

\endEinzug

\medskip

\noindent {\it Proof: } According to theorem 12-4 a sequence
transformation $G_k^{(n)}$ accelerates linear convergence if it
possesses a limiting transformation $G_k^{(\infty)}$ which satisfies
eq. (12.1-18). If we perform the limit $n \to \infty$ in the recursive
scheme (11.2-1) for $\lambda_k^{(n)}$ according to eq. (12.1-3), we
find that the limiting transformation of $\lambda_k^{(n)}$ is defined
by a recursive scheme which is identical with the second version of
Aitken's iterated $\Delta^2$ process, eq. (11.1-6).

Hence, it follows from eq. (12.1-18) that $\lambda_k^{(n)}$ accelerates
linear convergence if the elements of the sequence $\Seq {\sigma_n
(\rho)}$, eq. (12.1-23) are transformed into $1 / (1 - \rho)$ by the
recursive scheme (11.1-6).

However, since Aitken's iterated $\Delta^2$ process does not depend
explicitly upon $n$, this is at the same time the condition which
determines whether Aitken's iterated $\Delta^2$ process is able to
accelerate linear convergence.

If we perform the limit $n \to \infty$ in the recursive schemes
(11.2-2) for $\sigma_k^{(n)}$ and (11.2-3) for $\mu_k^{(n)}$ in the
sense of eq. (12.1-3), we find that they have the same limiting
transformation which is again defined by the recursive scheme for
Aitken's iterated $\Delta^2$ process, eq. (11.1-6). This concludes the
proof of theorem 12-5.

\medskip

Aitken's iterated $\Delta^2$ process is defined by a recursive scheme
and no explicit expression is known. Consequently, a general proof,
that ${\cal A}_k^{(n)}$ accelerates linear convergence for all $k \ge
1$, does not seem possible. Instead, one can only construct an explicit
rational expression for some special ${\cal A}_k^{(n)}$ with $k,n$
being fixed integers. It can then be checked whether this explicit
expression produces $1 / (1 - \rho)$ if it is applied to the elements
of the sequence $\Seq {\sigma_n (\rho)}$, eq. (12.1-23). Unfortunately,
these rational expressions become very complicated for larger values of
the subscript $k$. Therefore, it is recommendable to consider only the
simplest case.

\medskip

\beginEinzug \sl \parindent = 0 pt

\Auszug {\bf Theorem 12-6:} Aitken's $\Delta^2$ process ${\cal
A}_1^{(n)}$, eq. (5.1-4), accelerates linear convergence for all $n \in
\N_0$.

\endEinzug

\medskip

\noindent {\it Proof: } According to theorem 12-4 we have to show that
Aitken's $\Delta^2$ process ${\cal A}_1^{(0)}$, eq. (5.1-4), produces
$1 / ( 1 - \rho )$ if it acts upon the sequence elements $\sigma_0
(\rho)$, $\sigma_1 (\rho)$, and $\sigma_2 (\rho)$ which are defined by
eq. (12.1-23). Straightforward computation shows:
$$
{\cal A}_1^{(0)} \; = \; \sigma_0 (\rho) - \frac
{[\sigma_1 (\rho) - \sigma_0 (\rho)]^2}
{\sigma_2 (\rho) - 2 \sigma_1 (\rho) + \sigma_0 (\rho)}
\; = \; \frac {1} {1 - \rho} \, .
\tag
$$

\medskip

It follows from theorems 12-5 and 12-6 that $\lambda_1^{(n)}$,
$\sigma_1^{(n)}$, and $\mu_1^{(n)}$ also accelerate linear convergence
for all $n \in \N_0$.

Here, it must be emphasized that theorem 12-6 does not imply that
${\cal A}_k^{(n)}$ with $k > 1$ will also accelerate linear
convergence. This has to be checked separately and independently for
every $k > 1$. For instance, it follows from eq. (5.1-5) and theorem
12-6 that $\epsilon_2^{(n)}$ also accelerates linear convergence.
However, in Wimp's book it is shown that in the case $\epsilon_4^{(n)}$
the limited associated transformation according to eqs. (12.1-7) and
(12.1-8) cannot satisfy eq. (12.1-17) since it is unbounded in the
vicinity of any $(\rho, \rho, \rho) \in \R^3$ (see pp. 127 - 128 of ref.
[23]).

Next, the sequence transformations ${\cal J}_k^{(n)}$, eq. (10.3-6),
${\cal B}_k^{(n)}$, eq. (11.1-5), and ${\cal C}_k^{(n)}$, eq.
(11.1-12), will be analyzed. Again, no explicit expressions are known
for these transformations which are all defined by recursive schemes.
Hence, as in the case of Aitken's iterated $\Delta^2$ process, eq.
(5.1-15), only special cases can be considered and it is again
recommendable to consider only the simplest case.

A simple calculation shows that ${\cal J}_1^{(n)}$, ${\cal B}_1^{(n)}$,
and ${\cal C}_1^{(n)}$ are exact for the partial sums (2.6-2) of the
geometric series. Consequently, it is not surprising that these
transformations accelerate linear convergence.

\medskip

\beginEinzug \sl \parindent = 0 pt

\Auszug {\bf Theorem 12-7:} The sequence transformations ${\cal
J}_1^{(n)}$, eq. (10.3-6), ${\cal B}_1^{(n)}$, eq. (11.1-5), and ${\cal
C}_1^{(n)}$, eq. (11.1-12), accelerate linear convergence for all $n
\in \N_0$.

\endEinzug

\medskip

\noindent {\it Proof: } The recursive schemes, which define these
sequence transformations, do not depend explicitly on $n$.
Consequently, one only has to show that the explicit expressions for
${\cal J}_1^{(0)}$, ${\cal B}_1^{(0)}$, and ${\cal B}_1^{(0)}$ produce
$1 / (1 - \rho)$ if they are applied to the first four elements of the
sequence $\Seq {\sigma_n (\rho)}$, eq. (12.1-23). Straightforward
computation shows that this and consequently theorem 12-7 is indeed
true.

\medskip

Since ${\cal J}_1^{(n)}$ is identical with $\theta_2^{(n)}$, theorem
12-7 also implies that $\theta_2^{(n)}$ accelerates linear convergence.

Next, we shall analyze those variants of the sequence transformations
${\cal L}_k^{(n)} (\beta, s_n, \omega_n)$, eq. (7.1-7), ${\cal
S}_k^{(n)} (\beta, s_n, \omega_n)$, eq. (8.2-7), ${\cal M}_k^{(n)}
(\gamma, s_n, \omega_n)$, eq. (9.2-6), and ${\cal D}_k^{(n)} (s_n,
\omega_n)$, eq. (9.5-4), which are based upon Levin's [28] explicit
remainder estimates (7.3-4), (7.3-6), and (7.3-10), as well as Smith
and Ford's [29] modified remainder estimate (7.3-8). First, it will be
shown that these sequence transformations are exact for the geometric
series. For that purpose it is recommendable to modify a theorem, which
was originally used by Smith and Ford (see p. 227 of ref. [29]) to
prove that $u_k^{(n)} (\beta, s_n)$, eq. (7.3-5), $d_k^{(n)} (\beta,
s_n)$, eq. (7.3-9), and $v_k^{(n)} (\beta, s_n)$, eq. (7.3-11), are
exact for the geometric series, in such a way that it can also be
applied in the case of the analogous variants of ${\cal S}_k^{(n)}
(\beta, s_n, \omega_n)$, ${\cal M}_k^{(n)} (\gamma, s_n, \omega_n)$,
and ${\cal D}_k^{(n)} (s_n, \omega_n)$.

\medskip

\beginEinzug \sl \parindent = 0 pt

\Auszug {\bf Theorem 12-8:} Assume that a sequence transformation
$T_k^{(n)} (s_n, \omega_n)$ is defined in the following way:
$$
T_k^{(n)} (s_n, \omega_n) \; = \; \frac
{ \Delta^k \> [ P_{k-1} (n) \, s_n / \omega_n ] }
{ \Delta^k \> [ P_{k-1} (n) \, / \omega_n ] }
\, , \qquad k,n \in \N_0 \, .
\tag
$$
Here, $P_{k-1} (n)$ is a polynomial of degree $\le k - 1$ in $n$. The
sequence transformation $T_k^{(n)}$ is defined if the elements of
$\Seqn s$ are bounded in magnitude and if the sequence $\Seqn \omega$
of remainder estimates is chosen in such a way that the denominator in
eq. (12.2-2) does not vanish.

Let us assume that the sequence $\Seqn s$ converges to some limit $s$.
Then, for $k \ge 1$ and for $n \ge 0$ this sequence transformation
$T_k^{(n)} (s_n, \omega_n)$ is exact for the sequence $\Seqn s$ if the
sequence $\Seqn {\omega}$ of remainder estimate can be chosen in such a
way that the denominator in eq. (12.2-2) does not vanish and that
$\omega_n$ is proportional to $s_n - s$, i.e.,
$$
\omega_n \; = \; c \, (s_n - s) \, ,
\qquad c \neq 0 \, , \quad n \in \N_0 \, .
\tag
$$

\endEinzug

\medskip

\noindent {\it Proof: } Since this sequence transformation $T_k^{(n)}
(s_n, \omega_n)$ is obviously invariant under translation in the sense
of eq. (3.1-4), we can write
$$
T_k^{(n)} (s_n, \omega_n) \; = \; s \, + \, \frac
{ \Delta^k \> [ P_{k-1} (n) \, (s_n - s) / \omega_n ] }
{ \Delta^k \> [ P_{k-1} (n) \, / \omega_n ] }
\, , \qquad k,n \in \N_0 \, .
\tag
$$

If $\omega_n$ is proportional to $s_n - s$ according to eq. (12.2-3),
the difference operator $\Delta^k$ in the numerator on the right-hand
side acts only on $P_{k-1} (n)$ which is a annihilated because it is a
polynomial of degree $\le k -1$ in $n$. Since the denominator on the
right-hand side of eq. (12.2-4) does not vanish by assumption, we see
that $T_k^{(n)}$ is exact for $\Seqn s$.

\medskip

It is immediately obvious that the sequence transformations ${\cal
L}_k^{(n)} (\beta, s_n, \omega_n)$, eq. (7.1-7), ${\cal S}_k^{(n)}
(\beta, s_n, \omega_n)$, eq. (8.2-7), ${\cal M}_k^{(n)} (\gamma, s_n,
\omega_n)$, eq. (9.2-6), and ${\cal D}_k^{(n)} (s_n, \omega_n)$, eq.
(9.5-4), are all of the form of eq. (12.2-2). Hence, according to
theorem 12-8 these sequence transformations are exact for the partial
sums (2.6-2) of the geometric series if the remainder estimates
$\omega_n$ satisfy
$$
\omega_n \; = \; c \, \frac{z^{n+1}} {1 - z} \ ,
\qquad c \neq 0 \, , \quad n \in \N_0 \, .
\tag
$$

Since $z^{\alpha} / (1 - z)$ is for fixed $z$ and $\alpha$ also a
constant, an equivalent condition for the exactness would be:
$$
\omega_n \; = \; c^{\prime} \, z^{n - \alpha + 1} \, ,
\qquad c^{\prime} \neq 0 \, , \quad n \in \N_0 \, ,
\quad \alpha \in \R \, .
\tag
$$

\medskip

\beginEinzug \sl \parindent = 0 pt

\Auszug {\bf Theorem 12-9:} The sequence transformations $u_k^{(n)}
(\beta, s_n)$, eq. (7.3-5), $y_k^{(n)} (\beta, s_n)$, eq. (8.4-2), and
$Y_k^{(n)} (\gamma, s_n)$, eq. (9.4-2), are all exact for the geometric
series for $k \ge 2$ and $n \ge 0$, whereas $t_k^{(n)} (\beta, s_n)$,
eq. (7.3-7), $d_k^{(n)} (\beta, s_n)$, eq. (7.3-9), $v_k^{(n)} (\beta,
s_n)$, eq. (7.3-11), $\tau_k^{(n)} (\beta, s_n)$, eq. (8.4-3),
$\delta_k^{(n)} (\beta, s_n)$, eq. (8.4-4), $\phi_k^{(n)} (\beta,
s_n)$, eq. (8.4-5), $T_k^{(n)} (\gamma, s_n)$, eq. (9.4-3),
$\Delta_k^{(n)} (\gamma, s_n)$, eq. (9.4-4), and $\Phi_k^{(n)} (\gamma,
s_n)$, eq. (9.4-5), are all exact for the geometric series for $k \ge
1$ and $n \ge 0$. Drummond's sequence transformation ${\cal D}_k^{(n)}
(s_n, \omega_n)$, eq. (9.5-4), is also exact for $k \ge 1$ and $n \ge
0$ if the remainder estimates (7.3-6), (7.3-8), and (7.3-10) are used.

\endEinzug

\medskip

\noindent {\it Proof: } We only have to prove that in the case of the
partial sums (2.6-2) of the geometric series the remainder estimates,
which define these variants of the sequence transformations ${\cal
L}_k^{(n)} (\beta, s_n, \omega_n)$, eq. (7.1-7), ${\cal S}_k^{(n)}
(\beta, s_n, \omega_n)$, eq. (8.2-7), ${\cal M}_k^{(n)} (\gamma, s_n,
\omega_n)$, eq. (9.2-6), and ${\cal D}_k^{(n)} (s_n, \omega_n)$, eq.
(9.5-4), lead to sequence transformations of the type of eq. (12.2-2)
with remainder estimates $\omega_n$ that satisfy either eq. (12.2-5) or
(12.2-6).

The remainder estimate (7.3-4) leads to $\omega_n = (\beta + n) z^n$.
In the case of $u_k^{(n)} (\beta, s_n)$, eq. (7.3-5), and $y_k^{(n)}
(\beta, s_n)$, eq. (8.4-2), the factor $\beta + n$ can be absorbed in
eqs. (7.1-6) and (8.2-6), respectively, leading for $k \ge 2$ to new
sequence transformations which are of the type of eq. (12.2-2). This
proves the exactness of $u_k^{(n)} (\beta, s_n)$ and $y_k^{(n)} (\beta,
s_n)$. The exactness of $Y_k^{(n)} (\gamma, s_n)$, eq. (9.4-2), can be
proved in the same way because the remainder estimate (9.4-1) yields
$\omega_n = - (\gamma + n) z^n$ and because $- \gamma - n$ can for $k
\ge 2$ be absorbed in eq. (9.2-5) yielding a new sequence
transformation satisfying eq. (12.2-2). However, in the case of
Drummond's sequence transformation ${\cal D}_k^{(n)} (s_n, \omega_n)$,
eq. (9.5-4), we would not obtain a sequence transformation of the type
of eq. (12.2-2) if we absorb either $\beta + n$ or $- \gamma - n$.
Thus, with neither of the two remainder estimates (7.3-4) or (9.4-1)
Drummond's sequence transformation is exact for the geometric series.

The remainder estimate (7.3-6) leads to $ \omega_n = z^n$. Because of
eq. (12.2-6) this proves the exactness of $t_k^{(n)} (\beta, s_n)$, eq.
(7.3-7), $\tau_k^{(n)} (\beta, s_n)$, eq. (8.4-3), and $T_k^{(n)}
(\gamma, s_n)$, eq. (9.4-3), for $k \ge 1$. With this remainder
estimate ${\cal D}_k^{(n)} (s_n, \omega_n)$ is also exact for $k \ge 1$.

The remainder estimate (7.3-8) leads to $ \omega_n = z^{n+1}$. Because
of eq. (12.2-6) this proves the exactness of $d_k^{(n)} (\beta, s_n)$,
eq. (7.3-9), $\delta_k^{(n)} (\beta, s_n)$, eq. (8.4-4), and
$\Delta_k^{(n)} (\gamma, s_n)$, eq. (9.4-4), for $k \ge 1$. With this
remainder estimate ${\cal D}_k^{(n)} (s_n, \omega_n)$ is also exact for
$k \ge 1$.

The remainder estimate (7.3-10) leads to $ \omega_n = z^{n+1} / (1-z)$.
Because of eq. (12.2-5) this proves the exactness of $v_k^{(n)} (\beta,
s_n)$, eq. (7.3-11), $\phi_k^{(n)} (\beta, s_n)$, eq. (8.4-5), and
$\Phi_k^{(n)} (\gamma, s_n)$, eq. (9.4-5), for $k \ge 1$. With this
remainder estimate ${\cal D}_k^{(n)} (s_n, \omega_n)$ is also exact for
$k \ge 1$.

\medskip

Since the sequence transformations, which are listed in theorem 12-9,
are all exact for the geometric series, it is not surprising that they
are also able to accelerate linear convergence.

\medskip

\beginEinzug \sl \parindent = 0 pt

\Auszug {\bf Theorem 12-10:} The sequence transformations $u_k^{(n)}
(\beta, s_n)$, eq. (7.3-5), $t_k^{(n)} (\beta, s_n)$, eq. (7.3-7),
$y_k^{(n)} (\beta, s_n)$, eq. (8.4-2), $\tau_k^{(n)} (\beta, s_n)$, eq.
(8.4-3), $Y_k^{(n)} (\gamma, s_n)$, eq. (9.4-2), and $T_k^{(n)}
(\gamma, s_n)$, eq. (9.4-3), accelerate linear convergence if and only
if Drummond's sequence transformation ${\cal D}_k^{(n)} (s_n,
\omega_n)$, eq. (9.5-4), with $\omega_n = a_n$ accelerates linear
convergence.

The sequence transformations $d_k^{(n)} (\beta, s_n)$, eq. (7.3-9),
$\delta_k^{(n)} (\beta, s_n)$, eq. (8.4-4), and $\Delta_k^{(n)}
(\gamma, s_n)$, eq. (9.4-4), accelerate linear convergence if and only
if Drummond's sequence transformation ${\cal D}_k^{(n)} (s_n,
\omega_n)$, eq. (9.5-4), with $\omega_n = a_{n+1}$ accelerates linear
convergence.

The sequence transformations $v_k^{(n)} (\beta, s_n)$, eq. (7.3-11),
$\phi_k^{(n)} (\beta, s_n)$, eq. (8.4-5), and $\Phi_k^{(n)} (\gamma,
s_n)$, eq. (9.4-5), accelerate linear convergence if and only if
Drummond's sequence transformation ${\cal D}_k^{(n)} (s_n, \omega_n)$,
eq. (9.5-4), with $\omega_n = a_n a_{n+1} / (a_n - a_{n+1})$
accelerates linear convergence.

\endEinzug

\medskip

\noindent {\it Proof: } According to theorem 12-4 the sequence
transformations mentioned above accelerate linear convergence if they
possess limiting transformations which satisfy eq. (12.1-18). If we
perform the limit $n \to \infty$ in the explicit expressions for these
sequence transformations according to eq. (12.1-3), we find that their
limiting transformations $G_k^{(\infty)}$ are Drummond's sequence
transformation with $\omega_n = a_n$, $\omega_n = a_{n+1}$, or
$\omega_n = a_n a_{n+1} / (a_n - a_{n+1})$, respectively.

Since Drummond's sequence transformation does not explicitly depend on
$n$, this is at the same time the condition which determines whether
${\cal D}_k^{(n)} (s_n, \omega_n)$ with either $\omega_n = a_n$,
$\omega_n = a_{n+1}$, or $\omega_n = a_n a_{n+1} / (a_n - a_{n+1})$
accelerates linear convergence or not.

\medskip

\beginEinzug \sl \parindent = 0 pt

\Auszug {\bf Theorem 12-11:} Drummond's sequence transformation ${\cal
D}_k^{(n)} (s_n, \omega_n)$, eq. (9.5-4), with either $\omega_n = a_n$,
$\omega_n = a_{n+1}$, or $\omega_n = a_n a_{n+1} / (a_n - a_{n+1})$
accelerates linear convergence for $k \ge 1$ and $n \ge 0$.

\endEinzug

\medskip

\noindent {\it Proof: } According to theorem 12-4 Drummond's sequence
transformation ${\cal D}_k^{(n)} (s_n, \omega_n)$ accelerates linear
convergence if ${\cal D}_k (s_0, \omega_0)$ is exact for the sequence
$\Seq {\sigma_n (\rho)}$, eq. (12.1-23). The remainder estimate
$\omega_n = a_n$ leads to $\omega_n = \rho^{n-1}$, the remainder
estimate $\omega_n = a_{n+1}$ leads to $\omega_n = \rho^n$, and the
remainder estimate $\omega_n = a_n a_{n+1} / (a_n - a_{n+1})$ leads to
$\omega_n = \rho^n / (1 - \rho)$. Since these remainder estimates are
all of the form of either eq. (12.2-5) or (12.2-6), it follows from
theorem 12-8 that for $k \ge 1$ ${\cal D}_k^{(0)} (s_0, \omega_0)$ with
either $\omega_n = a_n$, $\omega_n = a_{n+1}$, or $\omega_n = a_n
a_{n+1} / (a_n - a_{n+1})$ is exact for the sequence $\Seq {\sigma_n
(\rho)}$, eq. (12.1-23). This completes the proof of theorem 12-11.

\medskip

Similar results as in theorems 12-9 and 12-10 can also be derived for
the analogous variants of the generalized transformations ${\cal L}_{k,
\ell}^{(n)} (\beta, s_n, \omega_n)$, eq. (7.1-8), ${\cal S}_{k,
\ell}^{(n)} (\beta, s_n, \omega_n)$, eq. (8.2-8), and ${\cal M}_{k,
\ell}^{(n)} (\gamma, s_n, \omega_n)$, eq. (9.2-7), with $\ell \ge 1$.

\medskip

\Abschnitt A modification of Germain-Bonne's theory for sequence
transformations involving remainder estimates.

\smallskip

\aktTag = 0

In the last section Germain-Bonne's formal theory of convergence
acceleration was applied to those variants of the sequence
transformations ${\cal L}_k^{(n)} (\beta, s_n, \omega_n)$, eq. (7.1-7),
${\cal S}_k^{(n)} (\beta, s_n, \omega_n)$, eq. (8.2-7), ${\cal
M}_k^{(n)} (\gamma, s_n, \omega_n)$, eq. (9.2-6), and ${\cal D}_k^{(n)}
(s_n, \omega_n)$, eq. (9.5-4), which are based upon Levin's [28]
explicit remainder estimates (7.3-4), (7.3-6), and (7.3-10), as well as
Smith and Ford's [29] modified remainder estimate (7.3-8).

The sequence $\Seqn {\omega}$ of remainder estimates plays a central
r\^ole in the sequence transformations mentioned above since its choice
will ultimately determine success or failure. Experience shows that the
simple remainder estimates (7.3-4), (7.3-6), (7.3-8), and (7.3-10)
often work remarkably well in a variety of situations. However, one
cannot expect that these simple remainder estimates will always lead to
satisfactory results and in some cases it may well be much more
efficient to use other remainder estimates $\Seq {\omega_n}$.

It is often possible to obtain explicit expressions for the remainders
$\Seqn r$ of a sequence $\Seqn s$. Unfortunately, expressions of that
kind are in most cases practically useless since they are normally too
complicated. In some cases, however, it may be possible to derive with
the help of simplifying assumptions, which are valid in the limit of
large indices $n$, simple explicit expressions which can be used as
remainder estimates $\Seqn {\omega}$. If such an explicit remainder
estimate $\omega_n$ does not depend explicitly upon one or several
elements of $\Seqn s$, the sequence transformations ${\cal L}_k^{(n)}
(\beta, s_n, \omega_n)$, ${\cal S}_k^{(n)} (\beta, s_n, \omega_n)$,
${\cal M}_k^{(n)} (\gamma, s_n, \omega_n)$, and ${\cal D}_k^{(n)} (s_n,
\omega_n)$ are linear sequence transformations. In addition, they are
also defined and exact for constant sequences.

In such a case, one would of course like to know how the two sequences
$\Seqn s$ and $\Seqn {\omega}$ have to be related in order to guarantee
at least the regularity of the transformation, and under which
circumstances the convergence of $\Seqn s$ will be accelerated. These
questions can at least partially be answered by a suitable modification
of Germain-Bonne's theory.

Our approach is inspired by a modification of Germain-Bonne's theory
which can be found in Brezinski's first book (see pp. 126 - 132 of ref.
[19]). Brezinski considered sequence transformations which
simultaneously act upon $k+1$ consecutive elements of the sequence
$\Seqn s$ to be transformed and on $k+1$ consecutive elements of an
auxiliary sequence $\Seqn x$. Brezinski's modification of
Germain-Bonne's theory is suited for algorithms which remain
well-defined if some or all elements of the auxiliary sequence $\Seqn
x$ are zero. This cannot be assumed here.

It follows from eqs. (7.1-6), (8.2-6), (9.2-5), and (9.5-3) that the
sequence transformations ${\cal L}_k^{(n)} (\beta, s_n, \omega_n)$,
${\cal S}_k^{(n)} (\beta, s_n, \omega_n)$, ${\cal M}_k^{(n)} (\gamma,
s_n, \omega_n)$, and ${\cal D}_k^{(n)} (s_n, \omega_n)$ are all of the
form of eq. (12.2-2). If we assume that the remainder estimates $\Seqn
\omega$ do not explicitly depend upon the elements of $\Seqn s$, then
it is a typical feature of the sequence transformations mentioned above
that they are linear functions of their first $k+1$ variables.
Consequently, these transformations are all continuous in their first
$k+1$ variables provided that the elements of $\Seqn s$ are bounded.
Much more critical is the continuity of these sequence transformations
with respect to their second $k + 1$ variables, the remainder estimates
$\omega_n, \omega_{n+1}, \ldots, \omega_{n+k}$. Since the remainder
estimates always occur in denominators, the elements of $\Seqn
{\omega}$ must not be zero for all finite values of $n$. In addition,
the remainder estimates have to be chosen in such a way that the
denominator of such a sequence transformation, which is the $k$-th
difference of $P_{k-1} (n)$, a polynomial of degree $\le k - 1$ in $n$,
divided by $\omega_n$, will not vanish. Hence, a necessary but
unfortunately not sufficient condition, which an admissible sequence
$\Seqn {\omega}$ of remainder estimates would have to satisfy, is that
its elements are nonzero and distinct for all finite values of $n$.
This implies that for every finite value of $n$ a substring $\omega_n,
\omega_{n+1}, \ldots, \omega_{n+k}$ has to belong to a suitable subset
of $\H^{k+1}$, the intersection of $\F^{k+1}$ and $\D^{k+1}$.

In this section $\Gamma_k^{(n)}$ stands for a sequence transformation
which acts upon $k+1$ consecutive elements of a convergent sequence
$\Seqn s$ and which also requires $k+1$ consecutive elements of a
sequence $\Seqn {\omega}$ of remainder estimates. The superscript $n$
indicates that $\Gamma_k^{(n)}$ may depend explicitly upon $n$.

Hence, for fixed $k \in \N_0$ a sequence transformation
$\Gamma_k^{(n)}$ is a function which may depend explicitly on $n \in
\N_0$ and which is defined on vectors ${\bf x} = (x_1, x_2, \ldots,
x_{k+1}) \in \R^{k+1}$ and ${\bf z} = (z_1, z_2, \ldots, z_{k+1}) \in
\H^{k+1}$. In addition, we assume that such a sequence transformation
$\Gamma_k^{(n)} : \R^{k+1} \times \H^{k+1} \to \R$ possesses for fixed
$k \in \N_0$ and for all $n \in \N_0$ the following properties:

\medskip

\beginBeschreibung \zu \Laenge{(A-0):} \sl

\item {(A-0):} $\Gamma_k^{(n)}$ is defined and continuous on a subset
of $\R^{k+1} \times \H^{k+1}$.
\item {(A-1):} $\Gamma_k^{(n)}$ is a homogeneous function of degree one
in its first $k+1$ variables and a homogeneous function of degree zero
in its second $k+1$ variables. This means that for all vectors ${\bf x}
\in \R^{k+1}$ and ${\bf z} \in \H^{k+1}$, for which $\Gamma_k^{(n)}$ is
defined and continuous, and for all $\lambda, \mu \in \R$ with $\mu
\neq 0$ we have
$$
\beginAligntags
" \Gamma_k^{(n)}
(\lambda x_1, \lambda x_2, \ldots, \lambda x_{k+1} \Bar
z_1, z_2, \ldots, z_{k+1}) \\
" = \; \lambda \> \Gamma_k^{(n)} (x_1, x_2, \ldots, x_{k+1} \Bar
z_1, z_2, \ldots, z_{k+1}) \, ,
\erhoehe\aktTag \\ \tag*{\tagnr a}
" \Gamma_k^{(n)}
(x_1, x_2, \ldots, x_{k+1} \Bar
\mu z_1, \mu z_2, \ldots, \mu z_{k+1}) \\
" = \; \Gamma_k^{(n)} (x_1, x_2, \ldots, x_{k+1} \Bar
z_1, z_2, \ldots, z_{k+1}) \, .
\\ \tag*{\tagform\aktTagnr b}
\endAligntags
$$
\item {(A-2):} $\Gamma_k^{(n)}$ is linear in its first $k+1$ variables.
Consequently, for all vectors ${\bf x}, {\bf y} \in \R^{k+1}$ and ${\bf
z} \in \H^{k+1}$, for which is $\Gamma_k^{(n)}$ defined and continuous,
we have
$$
\beginAligntags
" \Gamma_k^{(n)} (x_1 + y_1, x_2 + y_2, \ldots, x_{k+1} + y_{k+1}
\Bar z_1, z_2, \ldots, z_{k+1}) \\
" = \; \Gamma_k^{(n)} (x_1, x_2, \ldots, x_{k+1} \Bar
z_1, z_2, \ldots, z_{k+1})
\; + \;
\Gamma_k^{(n)} (y_1, y_2, \ldots, y_{k+1} \Bar
z_1, z_2, \ldots, z_{k+1}) \, . \qquad
\\ \tag
\endAligntags
$$
\item {(A-3):} Let ${\bf c} = (c, c, \ldots, c) \in \R^{k+1}$ be a
vector with constant components and let ${\bf z}$ belong to the subset
of $\H^{k+1}$ for which $\Gamma_k^{(n)}$ is defined and continuous.
Then, $\Gamma_k^{(n)}$ is exact, i.e.,
$$
\Gamma_k^{(n)} (c, c, \ldots, c \Bar z_1, z_2, \ldots, z_{k+1})
\; = \; c \, .
\tag
$$
\item {(A-4):} There exists a subset of $\H^{k+1}$ such that for all
bounded vectors ${\bf x} \in \R^{k+1}$ and for all vectors ${\bf z}$
belonging to this subset the limiting transformation
$$
\beginAligntags
" \Gamma_k^{(\infty)} (x_1, x_2, \ldots, x_{k+1} \Bar
z_1, z_2, \ldots, z_{k+1}) \\
" = \; \lim_{n \to \infty} \>
\Gamma_k^{(n)} (x_1, x_2, \ldots, x_{k+1} \Bar
z_1, z_2, \ldots, z_{k+1})
\\ \tag
\endAligntags
$$
\item {} is uniquely defined and continuous on this subset of $\R^{k+1}
\times \H^{k+1}$. In addition, it is assumed that the limiting
transformation $\Gamma_k^{(\infty)}$ is also homogeneous and linear
according to (A-1) and (A-2).

\endBeschreibung

\medskip

Similarly as in the case of the sequence transformations $G_k^{(n)}$ it
cannot be assumed that a sequence transformation $\Gamma_k^{(n)}$ and
its limiting transformation $\Gamma_k^{(\infty)}$ will be defined and
continuous on the same subset of $\R^{k+1} \times \H^{k+1}$.

It follows from their explicit expressions that the sequence
transformations ${\cal L}_k^{(n)} (\beta, s_n, \omega_n)$, eq. (7.1-7),
${\cal S}_k^{(n)} (\beta, s_n, \omega_n)$, eq. (8.2-7), ${\cal
M}_k^{(n)} (\gamma, s_n, \omega_n)$, eq. (9.2-6), and ${\cal D}_k^{(n)}
(s_n, \omega_n)$, eq. (9.5-4), satisfy (A-0) -- (A-4) if suitable
sequences $\Seqn {\omega}$ of remainder estimates are used. Since these
transformations are linear functions of the $k+1$ sequence elements
$s_n, s_{n+1}, \ldots, s_{n+k}$ if the elements of $\Seqn {\omega}$ do
not depend explicitly upon the elements of $\Seqn s$, they are defined
and continuous for arbitrary sequences $\Seqn s$ if the elements of
$\Seqn {\omega}$ are nonzero for all finite values of $n$ and if the
remainder estimates are chosen in such a way that the denominators of
these transformations do not vanish. The denominators of ${\cal
L}_k^{(n)} (\beta, s_n, \omega_n)$, ${\cal S}_k^{(n)} (\beta, s_n,
\omega_n)$, ${\cal M}_k^{(n)} (\gamma, s_n, \omega_n)$, and ${\cal
D}_k^{(n)} (s_n, \omega_n)$, which are all of the type of eq. (12.2-2),
will be nonzero for all $k,n \in \N_0$ if the remainder estimates
satisfy $\Delta^k (P_{k-1} (n) / \omega_n) \neq 0$. In the case of
${\cal L}_k^{(n)} (\beta, s_n, \omega_n)$, we have according to eq.
(7.1-6) $P_{k-1} (n) = (n + \beta)^{k-1}$, in the case of ${\cal
S}_k^{(n)} (\beta, s_n, \omega_n)$ we have according to eq. (8.2-6)
$P_{k-1} (n) = (n + \beta)_{k-1}$, in the case of ${\cal M}_k^{(n)}
(\gamma, s_n, \omega_n)$ we have according to eq. (9.2-5) $P_{k-1} (n)
= (- n - \gamma)_{k-1}$, and in the case of ${\cal D}_k^{(n)} (s_n,
\omega_n)$ we have according to eq. (9.5-3) $P_{k-1} (n) = 1$.

We are now in a position to formulate an analogue of theorem 12-2 for
sequence transformations of the type $\Gamma_k^{(n)} : \R^{k+1} \times
\H^{k+1} \to \R$.

\medskip

\beginEinzug \sl \parindent = 0 pt

\Auszug {\bf Theorem 12-12:} Every sequence transformation
$\Gamma_k^{(n)}$, which is defined and continuous for all vectors ${\bf
x}\in \R^{k+1}$ and for all ${\bf z}$ belonging to a suitable subset of
$\H^{k+1}$ and which also satisfies (A-0) - (A-4), can be expressed in
the following way:
$$
\beginAligntags
" \Gamma_k^{(n)} (x_1, x_2, \ldots, x_{k+1} \Bar
z_1, z_2, \ldots, z_{k+1}) \\
" = \; x_1 \, + \, z_1 \>
\gamma_k^{(n)} \biggl(\frac {x_2 - x_1} {z_1}, \ldots ,
\frac {x_{k+1} - x_k} {z_k} \biggBar \frac {z_2} {z_1},
\ldots, \frac {z_{k+1}} {z_k} \biggr) \, .
\\ \tag
\endAligntags
$$
The associated transformation $\gamma_k^{(n)}$, which is defined and
continuous on a suitable subset of $\R^k \times \F^k$, is given by
$$
\beginAligntags
" \gamma_k^{(n)} \biggl(\frac {x_2 - x_1} {z_1}, \ldots ,
\frac {x_{k+1} - x_k} {z_k} \biggBar \frac {z_2} {z_1},
\ldots, \frac {z_{k+1}} {z_k} \biggr) \\
" = \; \Gamma_k^{(n)} \biggl( 0, \frac {x_2 - x_1}{z_1}, \ldots ,
\sum_{j=0}^{k-1} \, \frac {x_{j+2} - x_{j+1}} {z_{j+1}} \,
\prod_{i=0}^{j-1} \, \frac {z_{i+2}}{z_{i+1}} \biggBar
1, \frac {z_2}{z_1}, \ldots,
\prod_{j=0}^{k-1} \, \frac {z_{j+2}}{z_{j+1}} \biggr) \, .
\\ \tag
\endAligntags
$$

\endEinzug

\medskip

\noindent {\it Proof: } It follows from (A-2) and (A-3) that
$\Gamma_k^{(n)}$ is invariant under translation in the sense of eq.
(3.1-4). Consequently, we can subtract $x_1$ from $\Gamma_k^{(n)}$.
Since by assumption $z_1 \neq 0$, it follows from (A-1) that we can
divide the $2 k + 2$ arguments of $\Gamma_k^{(n)}$ by $z_1$. This
yields:
$$
\beginAligntags
" \Gamma_k^{(n)} (x_1, x_2, \ldots, x_{k+1} \Bar
z_1, z_2, \ldots, z_{k+1}) \\
" = \; x_1 \, + \, z_1 \>
\Gamma_k^{(n)} \biggl(0, \frac {x_2 - x_1}{z_1}, \ldots ,
\frac {x_{k+1} - x_1}{z_1} \biggBar
1, \frac{z_2}{z_1}, \ldots, \frac {z_{k+1}}{z_1} \biggr) \, .
\\ \tag
\endAligntags
$$

We now need the following two relationships which can be proved quite
easily:
$$
\beginAligntags
" \frac {z_{\nu+1}} {z_1} \; " = \;
" \prod_{i=0}^{\nu-1} \, \frac {z_{i+2}} {z_{i+1}}
\, , \qquad \nu \in \N_0 \, ,
\\ \tag
" \frac {x_{\nu+1} - x_1} {z_1} \; " = \;
" \sum_{j=0}^{\nu-1} \, \frac {x_{j+2} - x_{j+1}} {z_{j+1}} \,
\prod_{i=0}^{j-1} \, \frac {z_{i+2}} {z_{i+1}}
\, , \qquad \nu \in \N_0 \, .
\\ \tag
\endAligntags
$$

If we insert eqs. (12.3-8) and (12.3-9) into eq. (12.3-7) we obtain
eqs. (12.3-5) and (12.3-6). The continuity of $\gamma_k^{(n)}$ on a
suitable subset of $\R^k \times \F^k$ follows from the continuity of
$\Gamma_k^{(n)}$ on a subset of $\R^{k+1} \times \H^{k+1}$ via eq.
(12.3-5).

\medskip

With the help of theorem 12-12 it can be analyzed for which convergent
sequences $\Seqn s$ and for which sequences $\Seqn \omega$ of remainder
estimates a sequence transformation $\Gamma_k^{(n)}$ satisfying (A-0)
-- (A-4) will be regular.

\medskip

\beginEinzug \sl \parindent = 0 pt

\Auszug {\bf Theorem 12-13:} Let $\Gamma_k^{(n)}$ be a sequence
transformation which satisfies theorem 12-12. This means that according
to eqs. (12.3-5) and (12.3-6) an associated transformation
$\gamma_k^{(n)}$ exists which is continuous on a suitable subset of
$\R^k \times \F^k$. Let us assume that a sequence $\Seqn s$ converges
to some limit $s$, and that the elements of a sequence $\Seqn {\omega}$
of remainder estimates -- although they are different from zero for all
finite values of $n$ -- approach zero as $n \to \infty$. Then,
$\Gamma_k^{(n)}$ is regular if the elements of $\Seqn s$ and $\Seqn
{\omega}$ satisfy:

\beginBeschreibung \zu \Laenge{(ii):}

\item{(i) \hfill :} For sufficiently large values of $n \in \N_0$ the
ratios $\Delta s_n / \omega_n$ are all bounded, i.e.,
$$
\vert \Delta s_n / \omega_n \vert \le c \, ,
\qquad 0 \le c < \infty \, .
\tag
$$
\item{(ii):} For all bounded vectors ${\bf y} = (y_1, y_2, \ldots, y_k)
\in \R^k$ the associated transformation $\gamma_k^{(n)}$ remains
bounded as $n \to \infty$:
$$
\lim_{n \to \infty} \> \biggl| \gamma_k^{(n)}
\bigl( y_1, \ldots, y_k \bigBar \omega_{n+1} / \omega_n,
\ldots, \omega_{n+k} / \omega_{n+k-1} \bigr) \biggr|
\; \le \; M \, , \qquad 0 < M < \infty \, .
\tag
$$

\endBeschreibung

\endEinzug

\medskip

\noindent {\it Proof: } It follows from eq. (12.3-5) that
$\Gamma_k^{(n)}$ can be written in the following way:
$$
\beginAligntags
" \Gamma_k^{(n)} (s_n, s_{n+1}, \ldots, s_{n+k} \Bar
\omega_n, \omega_{n+1}, \ldots, \omega_{n+k}) \\
" = \; s_n \, + \, \omega_n \>
\gamma_k^{(n)} \biggl( \frac {\Delta s_n} {\omega_n}, \ldots ,
\frac {\Delta s_{n+k-1}}{\omega_{n+k-1}} \biggBar
\frac {\omega_{n+1}}{\omega_n}, \ldots,
\frac {\omega_{n+k}}{\omega_{n+k-1}} \biggr) \, . \\
\tag
\endAligntags
$$

According to eq. (12.3-6) the associated transformation is given by
$$
\beginAligntags
" \gamma_k^{(n)} \biggl( \frac {\Delta s_n} {\omega_n}, \ldots ,
\frac {\Delta s_{n+k-1}}{\omega_{n+k-1}} \biggBar
\frac {\omega_{n+1}}{\omega_n}, \ldots,
\frac {\omega_{n+k}}{\omega_{n+k-1}} \biggr) \\
" = \; \Gamma_k^{(n)} \biggl( 0, \frac{\Delta s_{n+1}}{\omega_n},
\ldots , \sum_{j=0}^{k-1} \frac {\Delta s_{n+j}}{\omega_{n+j}}
\> \prod_{i=0}^{j-1} \frac {\omega_{n+i+1}}{\omega_{n+i}}
\biggBar 1, \frac {\omega_{n+1}}{\omega_n}, \ldots,
\prod_{j=0}^{k-1} \frac {\omega_{n+j+1}}{\omega_{n+j}} \biggr)
\, . \\
\tag
\endAligntags
$$

$\Gamma_k^{(n)}$ preserves the convergence of $\Seqn s$ to its limit
$s$ if the second term on the right-hand side of eq. (12.3-12) vanishes
as $n \to \infty$. Since $\Seqn {\omega}$ approaches zero as $n \to
\infty$, we only have to show that the associated transformation
$\gamma_k^{(n)}$ remains bounded as $n \to \infty$. Since
$\Gamma_k^{(n)}$ is according to (A-2) linear in its first $k+1$
components, we may conclude from eq. (12.3-13) that
$\gamma_k^{(\infty)}$ is bounded for all $n \in \N_0$ if its $2 k$
arguments remain bounded as $n \to \infty$, and if it remains
continuous in its second $k$ variables as $n \to \infty$. It follows
from (i) and (ii) that this is indeed the case which proves theorem
12-13.

\medskip

If we compare theorem 12-13 with the analogous theorem 12-3, which
formulates criteria for the regularity of sequence transformations
$G_k^{(n)} : \H^{k+1} \to \R$, we see that theorem 12-13 is quite
liberal with respect to the set of admissible sequences $\Seqn s$ since
only convergence to some limit is assumed. However, given a convergent
sequence $\Seqn s$, the criteria, which have to be satisfied by an
admissible sequence $\Seqn \omega$ of remainder estimates, are quite
restrictive.

The next theorem deals with the acceleration of linear convergence by
sequence transformations $\Gamma_k^{(n)}$. The following theorem is
virtually identical with the analogous theorem 12-4 which deals with
sequence transformations $G_k^{(n)}$. In both cases the decisive
criterion is that the limiting transformations $\Gamma_k^{(\infty)}$
and $G_k^{(\infty)}$ have to be exact for a shifted sequence of partial
sums of the geometric series. This again emphasizes the importance of
the geometric series for a theoretical analysis of the acceleration of
linear convergence.

\medskip

\beginEinzug \sl \parindent = 0 pt

\Auszug {\bf Theorem 12-14:} Let us assume that the elements of the
sequences $\Seqn s$ and $\Seqn {\omega}$ satisfy:
$$
\beginAligntags
" \text{(i)} ": \qquad
" \lim_{n \to \infty} \> s_n \; = \; s \, , \hfill
\\ \tag
" \text{(ii)} ": \qquad " \lim_{n \to \infty} \>
\frac {s_n - s} {\omega_n} \; = \; c \, ,
\qquad c \ne 0 \, , \hfill
\\ \tag
" \text{(iii)} ": \qquad " \lim_{n \to \infty} \>
\frac {\omega_{n+1}} {\omega_n} \; = \; \rho \, ,
\qquad 0 < \vert \rho \vert < 1 \, . \hfill
\\ \tag
\endAligntags
$$

Then, a necessary and sufficient condition that a sequence
transformation $\Gamma_k^{(n)}$ accelerates the convergence of the
sequence $\Seqn s$ is that its associated transformation
$\gamma_k^{(n)}$ satisfies:
$$
\beginAligntags
" \lim_{n \to \infty} \>
\gamma_k^{(n)} (y_n, y_{n+1}, \ldots, y_{n+k-1} \Bar
z_n, z_{n+1}, \ldots, z_{n+k-1}) \\
" = \; \gamma_k^{(\infty)} (y, y, \ldots, y \Bar
z, z, \ldots, z) \; = \;
\frac {y} {1 - z} \, .
\\ \tag
\endAligntags
$$
Here, $\Seqn y$ and $\Seqn z$ are essentially arbitrary sequences which
converge to $y$ and $z$, respectively. The above statement can also be
formulated in terms of the limiting sequence transformation
$\Gamma_k^{(\infty)}$:
$$
\beginAligntags
" \lim_{n \to \infty} \> \Gamma_k^{(n)} \biggl(0, y_n, \ldots,
\, \sum_{j=0}^{k-1} y_{n+j} \, \prod_{i=0}^{j-1} \, z_{n+i} \biggBar
1, z_n, \ldots, \prod_{j=0}^{k-1} \, z_{n+j} \biggr) \\
" = \; \Gamma_k^{(\infty)} \biggl(0, y, \ldots,
y \, \sum_{j=0}^{k-1}\, z^j \biggBar 1, z, \ldots, z^k \biggr)
\; = \; \frac {y} {1 - z} \, ,
\qquad 0 < \vert z \vert < 1 \, .
\\ \tag
\endAligntags
$$

\endEinzug

\medskip

\noindent {\it Proof: } If we subtract $s$ from both sides of eq.
(12.3-12) and divide the resulting expression by $s_n - s$ we obtain:
$$
\beginAligntags
" \frac
{\Gamma_k^{(n)} (s_n, s_{n+1}, \ldots, s_{n+k} \Bar
\omega_n, \omega_{n+1}, \ldots, \omega_{n+k}) - s}
{s_n - s} \\
" = \; 1 \, + \, \frac {\omega_{n}} {s_n - s} \,
\gamma_k^{(n)} \biggl( \frac {\Delta s_n}{\omega_n}, \ldots ,
\frac {\Delta s_{n+k-1}}{\omega_{n+k-1}} \biggBar
\frac {\omega_{n+1}}{\omega_n}, \ldots,
\frac {\omega_{n+k}}{\omega_{n+k-1}} \biggr) \, .
\\ \tag
\endAligntags
$$

According to eq. (2.6-6) the convergence of a sequence $\Seqn s$ to
its limit $s$ is accelerated by a sequence transformation
$\Gamma_k^{(n)}$ if
$$
\lim_{n \to \infty} \> \frac
{\Gamma_k^{(n)} (s_n, s_{n+1}, \ldots, s_{n+k} \Bar
\omega_n, \omega_{n+1}, \ldots, \omega_{n+k}) - s}
{s_n - s} \; = \; 0 \, .
\tag
$$

Hence , if we perform the limit $n \to \infty$ in eq. (12.3-19), the
left-hand side has to vanish if $\Gamma_k^{(n)}$ accelerates
convergence. Thus, we only have to investigate under which conditions
the right-hand side of eq. (12.3-19) also vanishes as $n \to \infty$.
It follows from eq. (12.3-15) that $\omega_n/(s_n - s) \to 1/c$ as $n
\to \infty$. In addition, it follows from eq. (12.3-16) that the second
$k$ arguments of the associated transformation $\gamma_k^{(n)}$ all
approach $\rho$ as $n \to \infty$, and with the help of eqs. (12.3-15)
and (12.3-16) we find that the first $k$ arguments of $\gamma_k^{(n)}$
all satisfy:
$$
\lim_{n \to \infty} \> \frac {\Delta s_n} {\omega_n} \; = \;
\lim_{n \to \infty} \> \biggl\{ \frac {\omega_{n+1}}{\omega_n} \,
\frac {s_{n+1} - s} {\omega_{n+1}} \; - \;
\frac {s_n - s} {\omega_n} \biggr\}
\; = \; c \, (\rho - 1) \, .
\tag
$$

Hence, if we perform the limit $n \to \infty$ in eq. (12.3-19) we find:
$$
1 \, + \, \frac {1}{c} \, \gamma_k^{(\infty)}
\bigl( c (\rho - 1), c (\rho - 1), \ldots, c (\rho - 1) \bigBar
\rho, \rho, \ldots, \rho \bigr) \; = \; 0 \, .
\tag
$$

Now we only have to set $c (\rho - 1) = y$ and $\rho = z$ in order to
see that if a sequence transformation $\Gamma_k^{(n)}$ accelerates the
convergence of $\Seqn s$, then its associated transformation has to
satisfy eq. (12.3-17). In the same way, if eq. (12.3-17) is satisfied
by the associated transformation of a sequence transformation
$\Gamma_k^{(n)}$, then it follows from eqs. (12.3-15), (12.3-16), and
(12.3-19) that $\Gamma_k^{(n)}$ accelerates the convergence of $\Seqn
s$ according to eq. (12.3-20). This proves theorem 12-14.

\medskip

First, it should be remarked that condition (ii), eq. (12.3-15), is
identical with eq. (7.3-1). This is another confirmation that the
elements of a sequence of remainder estimates should be chosen in such
a way that $\omega_n$ is proportional to the leading term of an
asymptotic expansion of $s_n - s$ as $n \to \infty$.

The first $k+1$ arguments of $\Gamma_k^{(\infty)}$ in eq. (12.3-18) are
apart from the factor $y$ identical with the elements of the sequence
$\Seq {\sigma_n (z)}$ which are defined by eq. (12.1-23). Consequently,
the right-hand side of eq. (12.3-18) can be rewritten in the following
way:
$$
\Gamma_k^{(\infty)} \biggl( 0, y, \ldots, y \frac{1-z^k} {1-z}
\biggBar 1,z, \ldots, z^k \biggr) \; = \;
\frac {y}{1 - z} \, , \qquad 0 < | z | < 1 \, .
\tag
$$

Next, theorem 12-14 will be used to prove that the sequence
transformations ${\cal L}_k^{(n)} (\beta, s_n, \omega_n)$, eq. (7.1-7),
${\cal S}_k^{(n)} (\beta, s_n, \omega_n)$, eq. (8.2-7), ${\cal
M}_k^{(n)} (\gamma, s_n, \omega_n)$, eq. (9.2-6), and ${\cal D}_k^{(n)}
(s_n, \omega_n)$, eq. (9.5-4), are able to accelerate linear
convergence if the remainder estimates are chosen in such a way that
$\omega_n$ is proportional to the leading term on an asymptotic
expansion of $s_n - s$ as $n \to \infty$.

\medskip

\beginEinzug \sl \parindent = 0 pt

\Auszug {\bf Theorem 12-15:} We assume that a sequence transformation
$T_k^{(n)} (s_n, \omega_n)$ can be written in the following way:
$$
T_k^{(n)} (s_n, \omega_n) \; = \; \frac
{\displaystyle
\sum_{j=0}^{k} \; ( - 1)^{j} \; \binom {k} {j} \;
f_k (n+j) \; \frac {s_{n+j}} {\omega_{n+j}} }
{\displaystyle
\sum_{j=0}^{k} \; ( - 1)^{j} \; \binom {k} {j} \;
f_k (n+j) \; \frac {1} {\omega_{n+j}} }
\; , \qquad k,n \in \N_0 \, .
\tag
$$
If the sequences $\Seqn s$ and $\Seqn \omega$ are as in theorem 12-14
and if the coefficients $f_k (n)$ satisfy
$$
\lim_{n \to \infty} \; f_k (n) \; = \; 1 \, ,
\qquad k \in \N_0 \, ,
\tag
$$
then $T_k^{(n)} (s_n, \omega_n)$ accelerates the convergence of $\Seqn
s$ for $k \ge 1$.

\endEinzug

\medskip

\noindent {\it Proof: } According to theorem 12-14 we have to show that
the limiting transformation $T_k^{(\infty)}$ satisfies either eq.
(12.3-18) or (12.3-23). However, the limiting transformation of all
sequence transformations $T_k^{(n)} (s_n, \omega_n)$ satisfying eqs.
(12.3-24) and (12.3-25) is Drummond's sequence transformation ${\cal
D}_k^{(n)} (s_n, \omega_n)$, eq. (9.5-4). Consequently, it is
sufficient to show that Drummond's sequence transformation is exact for
the sequence $y \Seq { \sigma_n (z)}$, with $\Seq {\sigma_n (z)}$
defined by eq. (12.1-23), if the remainder estimates $\omega_n = z^n$
are used, i.e.,
$$
\frac
{\displaystyle
\sum_{j=0}^{k} \; ( - 1)^{j} \; \binom {k} {j} \;
\frac {y} {z^j} \, \frac {1 - z^j} {1-z} }
{\displaystyle
\sum_{j=0}^{k} \; ( - 1)^{j} \; \binom {k} {j} \;
\frac {1} {z^j} }
\; = \; \frac {y} {1-z} \, .
\tag
$$

The numerator sum in eq. (12.3-26) may be rewritten in the following
way:
$$
\sum_{j=0}^{k} \; ( - 1)^{j} \; \binom {k} {j} \;
\frac {y} {z^j} \, \frac {1 - z^j} {1-z} \; = \;
\frac {y} {1-z} \sum_{j=0}^{k} \; ( - 1)^{j} \; \binom {k} {j} \;
\frac {1} {z^j} \; - \;
\frac {y} {1-z} \sum_{j=0}^{k} \; ( - 1)^{j} \; \binom {k} {j} \, .
\tag
$$

It follows from eq. (2.4-8) that the second sum on the right-hand side
is zero for $k \ge 1$. This shows that eq. (12.3-26) is correct and
proves theorem 12-15.

\medskip

It can be deduced directly from their explicit representations that
${\cal L}_k^{(n)} (\beta, s_n, \omega_n)$, eq. (7.1-7), ${\cal
S}_k^{(n)} (\beta, s_n, \omega_n)$, eq. (8.2-7), ${\cal M}_k^{(n)}
(\gamma, s_n, \omega_n)$, eq. (9.2-6), and ${\cal D}_k^{(n)} (s_n,
\omega_n)$, eq. (9.5-4), satisfy eqs. (12.3-24) and (12.3-25). In
addition, it can be shown that the generalized transformations ${\cal
L}_{k, \ell}^{(n)} (\beta, s_n, \omega_n)$, eq. (7.1-8), ${\cal S}_{k,
\ell}^{(n)} (\beta, s_n, \omega_n)$, eq. (8.2-8), and ${\cal M}_{k,
\ell}^{(n)} (\gamma, s_n, \omega_n)$, eq. (9.2-7), with $\ell \ge 1$
also satisfy eqs. (12.3-24) and (12.3-25) for sufficiently large values
of $k$. Consequently, these sequence transformations accelerate the
convergence of a linearly convergent sequence $\Seqn s$ if the
remainder estimates $\Seqn \omega$ are chosen in such a way that
conditions (i) -- (iii) of theorem 12-14 are fulfilled.

It is in fact by no means trivial that the generalized sequence
transformations ${\cal L}_{k, \ell}^{(n)} (\beta, s_n, \omega_n)$,
${\cal S}_{k, \ell}^{(n)} (\beta, s_n, \omega_n)$, ${\cal M}_{k,
\ell}^{(n)} (\gamma, s_n, \omega_n)$ with $\ell \ge 1$ also accelerate
convergence if the sequences $\Seqn s$ and $\Seqn \omega$ satisfy
conditions (i) -- (iii) of theorem 12-14. It follows from the model
sequences (7.1-10), (8.2-10), and (9.2-9) that these transformations
were derived assuming that $s_n - s \sim P_{\ell} (n) \omega_n$ as $n
\to \infty$ with $P_{\ell} (n)$ being a polynomial of degree $\ell$ in
$n$. However, in theorem 12-14 it is assumed that $s_n - s \sim
\omega_n$ as $n \to \infty$. Hence, we see that sequence
transformations $T_k^{(n)} (s_n, \omega)$ satisfying eqs. (12.3-24) and
(12.3-25) accelerate linear convergence even if instead of the
``right'' sequence $\Seqn \omega$ of remainder estimates a ``wrong''
sequence $\Seq {\omega_n^{\prime}}$ of remainder estimates with
$\omega_n^{\prime} = P_{\ell} (n) \omega_n$ with $\ell \in \N_0$ is
used.

This behaviour is quite typical of all sequence transformations of the
type of eq. (12.2-2) which are defined as ratios of finite differences.

With the help of the following theorem it can also be shown that the
sequence transformations ${\cal L}_k^{(n)} (\beta, s_n, \omega_n)$, eq.
(7.1-7), ${\cal S}_k^{(n)} (\beta, s_n, \omega_n)$, eq. (8.2-7), ${\cal
M}_k^{(n)} (\gamma, s_n, \omega_n)$, eq. (9.2-6), and ${\cal D}_k^{(n)}
(s_n, \omega_n)$, eq. (9.5-4), as well as the generalized
transformations ${\cal L}_{k, \ell}^{(n)} (\beta, s_n, \omega_n)$, eq.
(7.1-8), ${\cal S}_{k, \ell}^{(n)} (\beta, s_n, \omega_n)$, eq.
(8.2-8), and ${\cal M}_{k, \ell}^{(n)} (\gamma, s_n, \omega_n)$, eq.
(9.2-7), with $\ell \ge 1$ are all exact for the geometric series.

\medskip

\beginEinzug \sl \parindent = 0 pt

\Auszug {\bf Theorem 12-16:} We assume that a sequence transformation
$T_k^{(n)} (s_n, \omega_n)$ can be written in the following way:
$$
T_k^{(n)} (s_n, \omega_n) \; = \; \frac
{\displaystyle
\sum_{j=0}^{k} \; ( - 1)^{j} \; \binom {k} {j} \;
\phi_{k-1} (n+j) \; \frac {s_{n+j}} {\omega_{n+j}} }
{\displaystyle
\sum_{j=0}^{k} \; ( - 1)^{j} \; \binom {k} {j} \;
\phi_{k-1} (n+j) \; \frac {1} {\omega_{n+j}} }
\; , \qquad k,n \in \N_0 \, .
\tag
$$

If for sufficiently large values of $k$ the coefficients $\phi_{k-1}
(n)$ are polynomials of degree $\le k - 1$ in $n$, then for
sufficiently large value of $k$ such a sequence transformation
$T_k^{(n)} (s_n, \omega_n)$ is exact for the partial sums of the
geometric series, eq. (2.6-2), if the remainder estimates are chosen
according to $\omega_n = z^{n + \alpha}$ with $\alpha \in \R$.

\endEinzug

\medskip

\noindent {\it Proof: } The numerator polynomial in eq. (12.3-28) can
then be rewritten in the following way:
$$
\beginAligntags
" \frac {1} {1 - z} \sum_{j=0}^{k} \, ( - 1)^{j} \,
\binom {k} {j} \, \phi_{k-1} (n+j) \,
\frac {1 - z^{n+j}} {z^{n+j+\alpha}} \\
" = \; \frac {1} {1-z} \sum_{j=0}^{k} \, ( - 1)^{j} \,
\binom {k} {j} \, \phi_{k-1} (n+j) \, \frac {1} {z^{n+j+\alpha}} \\
" - \; \frac {1} {z^{\alpha} \, (1-z)} \,
\sum_{j=0}^{k} \, ( - 1)^{j} \, \binom {k} {j}
\phi_{k-1} (n+j) \, .
\\ \tag
\endAligntags
$$
Let us now assume that $k$ is large enough such that $\phi_{k-1} (n)$
is a polynomial of degree $\le k - 1$ in $n$. Then it follows from eq.
(2.4-8) that the second sum in eq. (12.3-29) vanishes. This proves
theorem 12-16.

\medskip

Obviously, all sequence transformations mentioned above satisfy the
prerequisites of theorem 12-16.

\medskip

\Abschnitt A critical assessment of Germain-Bonne's theory

\smallskip

\aktTag = 0

With the help of either the original version of Germain-Bonne's formal
theory of convergence acceleration or its modifications it can be
decided whether a sequence transformation is regular, i.e., whether the
transformed sequence $\Seq {s^{\prime}_n}$ converges to the same limit
as the original sequence $\Seqn s$. In addition, a necessary and
sufficient condition could be formulated by means of which it can be
decided whether a sequence transformation is able to accelerate linear
convergence or not. Theoretically, these results are certainly
remarkable achievements, in particular since for their derivation only
some very general properties of the sequence transformation such as
continuity, homogeneity, and translativity had to be assumed. Also,
concerning the sequences $\Seqn s$, which are to be transformed, only
relatively little has to be assumed. In most cases it is sufficient
that the sequences converge or -- if the acceleration of linear
convergence is analyzed -- that they converge linearly.

However, it must not be overlooked that Germain-Bonne's formal theory
of convergence acceleration has some serious shortcomings which
definitely limit its practical usefulness, although it certainly is a
beautiful mathematical theory. For instance, the generality of
Germain-Bonne's theory and its modifications -- although highly
desirable from a theoretical point of view -- is at the same time a
major weakness since it implies that the results of this theory are
quite general and cannot be as specific as one would like them to be.

Germain-Bonne's theory is only able to make statements as for instance
that a sequence transformation is regular or that it is able to
accelerate linear convergence. However, from a practical point of view
the statement that a given sequence transformation is able to
accelerate linear convergence is just as useful -- or as useless -- as
the statement that a given series converges.

Germain-Bonne's theory is essentially asymptotic in nature because only
the sequence elements $s_n$ and the transforms $G_k^{(n)}$ or
$\Gamma_k^{(n)}$ with large values of $n$ matter. This asymptotic
nature is essential because it makes a theoretical analysis possible.
However, it refers to a situation -- sequence elements $s_n$ or
transforms $s^{\prime}_n$ with large indices $n$ -- which one would
like to avoid by using sequence transformations. In addition, the
predictive value of an asymptotic theory is often quite limited. A
given numerical technique need not be particularly powerful for
moderately large values of $n$, let alone for small values of $n$ even
if it is guaranteed that this technique will work well in the limit $n
\to \infty$.

In actual computations only a relatively small number of sequence
elements will normally be known, say $s_n, s_{n+1}, \ldots, s_{n+k}$,
and one would like to know how the information, which is contained in
these sequences elements, can be extracted and utilized in an optimal
way. Here, Germain-Bonne's theory or its variants cannot help at all
since it does not discriminate among sequence transformations which all
have the same properties in the limit $n \to \infty$. For instance,
according to theorem 12-10 certain variants of the sequence
transformations ${\cal L}_k^{(n)} (\beta, s_n, \omega_n)$, eq. (7.1-7),
${\cal S}_k^{(n)} (\beta, s_n, \omega_n)$, eq. (8.2-7), and ${\cal
M}_k^{(n)} (\gamma, s_n, \omega_n)$, eq. (9.2-6), are able to
accelerate linear convergence if and only if the analogous variants of
Drummond's sequence transformation ${\cal D}_k^{(n)} (s_n, \omega_n)$,
eq. (9.5-4), are able to accelerate linear convergence. Consequently,
Drummond's sequence transformation plays a very important r{\^o}le in
our modification of Germain-Bonne's formal theory of convergence
acceleration since it allows a unified treatment of a large class of
sequence transformations. But it is wrong to assume that Drummond's
sequence transformation will be particularly useful in actual
computations. In fact, we shall see later that Drummond's sequence
transformation is normally significantly less powerful than the other
sequence transformations mentioned above.

Germain-Bonne's theory is essentially a successful theory of the
acceleration of linear convergence and it does not help at all if for
instance logarithmic convergence is to be accelerated. This is quite
deplorable since the acceleration of logarithmic convergence is a much
more annoying problem than the acceleration of linear convergence --
both theoretically and computationally.

Another problem of great practical relevance is the determination of
the antilimit $s$ of a divergent sequence $\Seqn s$ by employing a
suitable sequence transformation $G_k^{(n)}$ or $\Gamma_k^{(n)}$. In
such a case a formal theory, which involves a limit $n \to \infty$,
makes no sense. Instead, any theoretical analysis of such a summation
process would have to say something about the convergence of a sequence
of transforms $G_k^{(n)}$ or $\Gamma_k^{(n)}$ to the antilimit $s$ and
how this convergence is affected if the subscript $k$ is increased
while the superscript $n$ is held fixed. Again, Germain-Bonne's theory
and its modifications cannot contribute anything.

Hence, Germain-Bonne's theory of convergence acceleration is not able
to treat several problems of great practical relevance and has to be
supplemented by other theoretical approaches. However, it is not likely
that Germain-Bonne's theory can be improved significantly without
making much more detailed assumptions about both the sequence
transformations, which are to be analyzed, and the sequences, which are
to be accelerated or summed.

\endAbschnittsebene

\neueSeite

\Abschnitt Summation and convergence acceleration of Stieltjes series

\vskip - 2 \jot

\beginAbschnittsebene

\medskip

\Abschnitt Stieltjes series and Stieltjes functions

\smallskip

\aktTag = 0

This section deals with the summation of divergent asymptotic series,
as they for instance occur in the theory of special functions or in
quantum mechanical perturbation theory. It is well known that a given
function $f(z)$ admits at most one asymptotic power series. The
converse, however, is not true, i.e., it may happen that there are
several different functions which all admit the same asymptotic power
series. In a theoretical analysis of summation processes, the possible
nonuniqueness of asymptotic expansions is certainly quite annoying and
the set of admissible asymptotic series should be suitably restricted
in order to avoid these complications.

These nonuniqueness problems can be avoided in the case of Stieltjes
series which assume an exceptional position among divergent series. For
Stieltjes series there exists a highly developed convergence and
representation theory (see refs. [18,22,79,95,96]). For instance, it
can be shown that Pad\'e approximants are able to sum even wildly
divergent Stieltjes series if the terms $a_n$ of these series do not
grow faster in magnitude than $c^{n+1} (2n)!$ as $n \to \infty$ with
$c$ being a suitable positive constant (see theorems 1.2 and 1.3 of
ref. [87]). This implies that Pad\'e approximants are able to sum the
divergent Euler series, eq. (1.1-7), which is, as we shall see later, a
Stieltjes series.

Stieltjes series are also of considerable physical interest since many
quantum mechanical perturbation expansions are Stieltjes series. For
instance, it could be proved rigorously that the Rayleigh-Schr\"odinger
perturbation expansions for the energy eigenvalues of the quartic
anharmonic oscillator are Stieltjes series [86,97]. It also follows
from the asymptotic behaviour (1.1-5) of the perturbation series
coefficients that Pad\'e approximants are able to sum the divergent
perturbation series (1.1-4) for the ground state energy of the quartic
anharmonic oscillator.

In this section, the summation of divergent Stieltjes series by means
of nonlinear sequence transformations will be investigated both
theoretically and numerically. It will be one of the main results of
this section that the sequence transformations ${\cal L}_k^{(n)}
(\beta, s_n, \omega_n)$, eq. (7.1-7), ${\cal S}_k^{(n)} (\beta, s_n,
\omega_n)$, eq. (8.2-7), and ${\cal M}_k^{(n)} (\gamma, s_n,
\omega_n)$, eq. (9.2-6), which all require a sequence $\Seqn \omega$ of
remainder estimates for their computation, sum divergent Stieltjes
series much more efficiently than for instance Pad\'e approximants, which
may be computed with the help of Wynn's $\epsilon$ algorithm, eq.
(4.2-1), or also Brezinski's $\theta$ algorithm, eq. (10.1-9). This is
probably due to the fact that in the case of a Stieltjes series it is
comparatively easy to find simple remainder estimates $\Seqn \omega$,
which may be used in the sequence transformations mentioned above and
which in spite of their simplicity yield rigorous and tight upper
bounds for the remainders of truncated Stieltjes series.

In order to make this section more self-contained, first those
properties of Stieltjes series and Stieltjes functions, which are of
particular importance for our purposes, will be reviewed.

A function $f(z)$ with $z \in \C$ will be called {\it Stieltjes
function} if it can be expressed as a {\it Stieltjes integral},
$$
f(z) \; = \; \int\nolimits_{0}^{\infty} \,
\frac {\d \psi (t)} {1 + z t} \, ,
\qquad \vert \arg (z) \vert < \pi \, .
\tag
$$

Here, $\psi (t)$ is a positive measure on $0 \le t < \infty$ which has
for all $m \in \N_0$ finite and positive moments $\mu_m$ defined by
$$
\mu_m \; = \; \int\nolimits_{0}^{\infty} \, t^m \d \psi (t)
\, , \qquad m \in \N_0 \, .
\tag
$$

A formal series expansion of the following type, which need not be
convergent,
$$
f(z) \; = \; \sum_{m=0}^{\infty} \; (-1)^m \, \mu_m \, z^m \, ,
\tag
$$
is called a {\it Stieltjes series} if its coefficients $\mu_m$ are
moments of a positive measure $\psi (t)$ on $0 \le t < \infty$
according to eq. (13.1-2), i.e.,
$$
f(z) \; = \; \sum_{m=0}^{\infty} \; (-1)^m \, z^m \,
\int\nolimits_{0}^{\infty} \, t^m \d \psi (t) \, .
\tag
$$

A good example for a Stieltjes function with a wildly divergent
Stieltjes series is the so-called Euler integral, eq. (1.1-6), and its
associated asymptotic series, the so-called Euler series, eq. (1.1-7).

\medskip

\beginEinzug \sl \parindent = 0 pt

\Auszug {\bf Theorem 13-1:} Every Stieltjes function $f (z)$ can be
written in the following way:
$$
f (z) \; = \; \sum_{m=0}^{n} \, (-1)^m \, \mu_m \, z^m \; + \;
(- z)^{n+1} \,
\int\nolimits_{0}^{\infty} \,
\frac {t^{n+1} \, \d \psi (t)} {1 + z t}
\, , \qquad \vert \arg (z) \vert < \pi \, .
\tag
$$

\endEinzug

\medskip

\noindent {\it Proof: } We only have to insert the relationship
$$
\sum_{m=0}^{n} \, x^m \; = \; \frac {1 - x^{n+1}} {1 - x}
\tag
$$
with $x = - z t$ into eq. (13.1-1) and do the moment integrals
according to eq. (13.1-2).

\medskip

Whether a Stieltjes series converges or diverges depends upon the
behaviour of the remainder integral on the right-hand side of eq.
(13.1-5). The next theorem shows that this remainder integral is
bounded by the first term of the power series (13.1-3) which was not
included in the partial sum in eq. (13.1-5). This bound will also help
us to find a simple and -- as we shall see later -- efficient sequence
of remainder estimates for the sequence transformations ${\cal
L}_k^{(n)} (\beta, s_n, \omega_n)$, eq. (7.1-7), ${\cal S}_k^{(n)}
(\beta, s_n, \omega_n)$, eq. (8.2-7), ${\cal M}_k^{(n)} (\gamma, s_n,
\omega_n)$, eq. (9.2-6), and ${\cal D}_k^{(n)} (s_n, \omega_n)$, eq.
(9.5-4).

\medskip

\beginEinzug \sl \parindent = 0 pt

\Auszug {\bf Theorem 13-2:} The remainder term in theorem 13-1,
$$
R_n (z) \; = \; (-z)^{n+1} \, \int\nolimits_{0}^{\infty} \,
\frac {t^{n+1} \, \d \psi (t)} {1 + z t} \, ,
\tag
$$
satisfies depending upon the value of $\theta = \arg (z)$ the following
inequalities:
$$
\beginAligntags
" \vert R_n (z) \vert \; \le \;
\mu_{n+1} \, \vert z^{n+1} \vert \, ,
\qquad " \vert \theta \vert \le \pi /2 \, ,
\erhoehe\aktTag \\ \tag*{\tagnr a}
" \vert R_n (z) \vert \; \le \;
\mu_{n+1} \, \vert z^{n+1} \Funk{cosec} \, \theta \vert \, ,
\qquad " \pi / 2 < \vert \theta \vert < \pi \, .
\\ \tag*{\tagform\aktTagnr b}
\endAligntags
$$

\endEinzug

\medskip

\noindent {\it Proof: } Setting $z = r \e^{i \theta}$ gives
$$
\vert 1 + z t \vert =
\bigl[ 1 + 2 r t \cos \theta + r^2 t^2 \bigr]^{1/2} \, .
\tag
$$

Next, one has to look for the value of $t$ with $0 \le t < \infty$ for
which $\vert 1 + z t \vert$ assumes its minimal value. Differentiation
with respect to $t$ gives the extremal condition
$$
t \; = \; - \; \frac {\cos \theta} {r} \, .
\tag
$$
Now, two different cases have to be distinguished:

\smallskip

\beginBeschreibung \zu \Laenge{(ii):}

\item{(i) \hfill :} $\vert \theta \vert \le \pi / 2$. Then, $\cos
\theta \ge 0$. Consequently, $\vert 1 + z t \vert$ assumes its minimal
value for $t = 0$ and we obtain the estimate (13.1-8a).
\item{(ii):} $\pi / 2 < \vert \theta \vert < \pi$. Then, $\cos \theta <
0$. Combination of eqs. (13.1-9) and (13.1-10) yields
$$
\vert \sin \theta \vert \le \vert 1 + z t \vert \, .
\tag
$$
If this inequality is used in eq. (13.1-7), estimate (13.1-8b) follows.
This shows that theorem 13-2 is correct.

\endBeschreibung

\medskip

It also follows from theorem 13-2 that every Stieltjes function
possesses an asymptotic series valid uniformly in every sector $\vert
\arg (z) \vert < \theta$ for any $\theta < \pi$ and that this
asymptotic series is a Stieltjes series (see p. 398 of ref. [86]). It
can also be proved that for every Stieltjes series there exists at
least one associated Stieltjes function. Since this possible
nonuniqueness is very inconvenient in summation processes, a criterion
would be needed which makes it possible to prove that there is a
one-to-one correspondence at least between certain divergent Stieltjes
series and certain Stieltjes functions. Thus, a condition would be
needed which is stronger than the existence of an asymptotic power
series of the type of eq. (13.1-3) but weaker than the existence of a
convergent Stieltjes series.

On the basis of Carleman's theorem (see p. 39 of ref. [98]) a
sufficient condition can be formulated which guarantees that there
exists a one-to-one correspondence between a Stieltjes function and its
associated asymptotic series.

A Stieltjes function $f (z)$, which is analytic in a sectorial region
of the complex plane, is said to satisfy a {\it strong asymptotic
condition} and its associated Stieltjes series is called a {\it strong
asymptotic series} if suitable positive constants $A$ and $\xi$ can be
found such that
$$
\biggl\vert f (z) \; - \;
\sum_{m=0}^n \, (-1)^m \mu_m z^m \biggr\vert \; \le \;
A \xi^{n+1} (n+1)! \, {\vert z \vert}^{n+1}
\tag
$$
holds for all $n \in \N_0$ and for all $z$ in this sectorial region.

The validity of such a strong asymptotic condition implies that a
Stieltjes function $f (z)$ is uniquely determined by its asymptotic
series (see p. 40 of ref. [98]). Such a strong asymptotic condition can
only be valid if the Stieltjes moments $\mu_n$, which are defined by
eqs. (13.1-3) and (13.1-4), satisfy for all $n \in \N_0$ (see p. 43 of
ref. [98])
$$
\mu_n \; \le \; A \, \xi^n \, n! \, .
\tag
$$

The moments of the Euler series, eq. (1.1-7), satisfy this inequality.
Hence, we may conclude that the Euler integral, eq. (1.1-6), is
uniquely determined by its asymptotic series (1.1-7). In the same way,
it follows from the asymptotic behaviour (1.1-5) of the series
coefficients that the perturbation series (1.1-4) for the ground state
energy of the quartic anharmonic oscillator is a strong asymptotic
series. Consequently, the ground state energy of the quartic anharmonic
oscillator is uniquely determined by its divergent perturbation series
(see also p. 41 of ref. [98]).

However, there are Stieltjes series of considerable physical interest
which have moments $\mu_n$ that behave like $(k n)!$ with $k > 1$ as $n
\to \infty$. For instance, in the case of the perturbation expansions
for the energy eigenvalues of the sextic or octic anharmonic oscillator
there is ample numerical evidence that the coefficients of these series
grow as $(2 n)!$ or $(3 n)!$, respectively, as $n \to \infty$ (see p.
43 of ref. [98]). Obviously, a strong asymptotic condition cannot be
valid in such a case. However, it can be shown (see p. 43 of ref. [98])
that a function $f (z)$, which is analytic within a sectorial region of
the complex plane, is also uniquely determined by its asymptotic series
if $f (z)$ satisfies a {\it modified strong asymptotic condition of
order $k$} and if its asymptotic series is a {\it modified strong
asymptotic series of order $k$}. This means that suitable positive
constants $A$ and $\xi$ can be found such that
$$
\biggl\vert f (z) \; - \;
\sum_{m=0}^n \, (-1)^m \mu_m z^m \biggr\vert \; \le \;
A \xi^{n+1} [k (n+1)]! \, {\vert z \vert}^{n+1}
\tag
$$
holds for all $n \in \N_0$ and for all $z$ in this sectorial region.

Again, such a modified strong asymptotic condition of order $k$ can
only be valid if the Stieltjes moments $\mu_n$ satisfy for all $n \in
\N_0$ (see p. 406 of ref. [86])
$$
\mu_n \; \le \; A \, \xi^n \, (k n)! \, .
\tag
$$

The bounds for the remainders $R_n (z)$ in theorem 13-2 are also of
considerable importance for convergence acceleration and summation
processes because it helps us to find simple and manageable remainder
estimates $\Seqn \omega$. Because of the specific structure of the
sequence transformations ${\cal L}_k^{(n)} (\beta, s_n, \omega_n)$,
${\cal S}_k^{(n)} (\beta, s_n, \omega_n)$, ${\cal M}_k^{(n)} (\gamma,
s_n, \omega_n)$, and ${\cal D}_k^{(n)} (s_n, \omega_n)$, all those
quantities which are independent of $n$, do not have to be included in
the remainder estimate $\omega_n$. Consequently, it is not necessary to
distinguish the two different cases $\vert \theta \vert \le \pi$ and
$\pi / 2 < \vert \theta \vert < \pi$ in theorem 13-2 and for every
sector $\vert \arg (z) \vert < \theta$ with $\theta < \pi$ a suitable
estimate $\omega_n$ for the remainder $R_n (z)$ of a Stieltjes series
would be
$$
\omega_n \; = \; (-1)^{n+1} \, \mu_{n+1} \, z^{n+1} \, ,
\qquad n \in \N_0 \, .
\tag
$$

This choice is identical with the remainder estimate (7.3-8) of Smith
and Ford [29]. Hence, for the summation of divergent Stieltjes series
the most natural choices among the numerous variants of the sequence
transformations ${\cal L}_k^{(n)} (\beta, s_n, \omega_n)$, ${\cal
S}_k^{(n)} (\beta, s_n, \omega_n)$, and ${\cal M}_k^{(n)} (\gamma, s_n,
\omega_n)$ would be the transformations $d_k^{(n)} (\beta, s_n)$, eq.
(7.3-9), $\delta_k^{(n)} (\beta, s_n)$, eq. (8.4-4), and
$\Delta_k^{(n)} (\gamma, s_n)$, eq. (9.4-4).

\medskip

\Abschnitt Theoretical error estimates

\smallskip

\aktTag = 0

In this section theoretical error estimates for the summation of a
divergent Stieltjes series by means of nonlinear sequence
transformations will be derived. However, the error estimates of this
section can also be applied if the convergence of sequences with
strictly alternating remainders is accelerated.

There are only relatively few references in the literature in which the
summation of divergent Stieltjes series by means of nonlinear sequence
transformations is analyzed. In articles by Wynn [99], Common [100],
Allen, Chui, Madych, Narcowich, and Smith [101], and Karlsson and Sydow
[102] the Pad\'e summation of Stieltjes series was analyzed. Then, there
is an article by Sidi [103] on the summation of certain wildly
divergent series by Levin's $u$ and $t$ transformations, eqs. (7.3-5)
and (7.3-7). Sidi could show that if the divergent series satisfies
certain conditions, Levin's $u$ and $t$ transformation produce
sequences of approximants which converge to the Borel sum of the
divergent series. Other sequence transformations were apparently not
yet treated in the literature. This is not too surprising since many
nonlinear sequence transformations as for instance Brezinski's $\theta$
algorithm, eq. (10.1-9), are defined by relatively complicated
recursive schemes and otherwise only very little is known about these
transformations. Currently, a detailed theoretical analysis of the
efficiency of such a sequence transformation in convergence
acceleration and summation processes seems to be more or less
impossible.

However, sequence transformations of the type of eq. (12.2-2) can be
analyzed relatively easily if suitable assumptions concerning the
sequences $\Seqn s$ and $\Seqn \omega$ are made. Consequently, in this
section the emphasis will be on the sequence transformations ${\cal
L}_k^{(n)} (\beta, s_n, \omega_n)$, eq. (7.1-7), ${\cal S}_k^{(n)}
(\beta, s_n, \omega_n)$, eq. (8.2-7), ${\cal M}_k^{(n)} (\gamma, s_n,
\omega_n)$, eq. (9.2-6), and ${\cal D}_k^{(n)} (s_n, \omega_n)$, eq.
(9.5-4), as well as on the mild generalizations ${\cal L}_{k,
\ell}^{(n)} (\beta, s_n, \omega_n)$, eq. (7.1-8), ${\cal S}_{k,
\ell}^{(n)} (\beta, s_n, \omega_n)$, eq. (8.2-8), and ${\cal M}_{k,
\ell}^{(n)} (\gamma, s_n, \omega_n)$, eq. (9.2-7). The following
theorem will be the basis of our analysis.

\medskip

\beginEinzug \sl \parindent = 0 pt

\Auszug {\bf Theorem 13-3:} Let us assume that a sequence
transformation $G_k^{(n)} (s_n, s_{n+1}, \ldots, s_{n+k})$ with $k,n
\in \N_0$ is invariant under translation according to eq. (3.1-4). Then
a necessary and sufficient condition that this sequence transformation
is able to sum a divergent sequence $\Seqn s$ to its antilimit $s$ on a
path ${\cal P} = \Seq {(n_j, k_j)}$ with $j \in \N_0$ is that
$$
\lim_{j \to \infty} \; G_{k_j}^{({n_j})}
(s_{n_j} - s, s_{n_j + 1} - s, \ldots, s_{n_j + k_j} - s)
\; = \; 0 \, .
\tag
$$

\endEinzug

\medskip

\noindent {\it Proof: } Since $G_k^{(n)}$ is by assumption invariant
under translation according to eq. (3.1-4), we have for arbitrary
integers $n_j$ and $k_j$:
$$
G_{k_j}^{({n_j})} (s_{n_j}, s_{n_j + 1}, \ldots, s_{n_j + k_j})
\; = \; s \, + \, G_{k_j}^{({n_j})}
(s_{n_j} - s, s_{n_j + 1} - s, \ldots, s_{n_j + k_j} - s) \, .
\tag
$$
Performing the limit $j \to \infty$ shows that theorem 13-3 is correct.

\medskip

It can be proved quite easily by a typical $2 \varepsilon$ proof that
if the antilimit $s$ of a divergent series $\Seqn s$ exists on a given
path ${\cal P}$, then it is uniquely determined on this path. For
different paths, however, no general statement concerning the
uniqueness of the antilimit $s$ can be made. In summation processes,
one is normally only interested in horizontal paths, i.e., in paths in
which $n_j$ is ultimately constant and in which only $k_j$ is
increased. Of course, theorem 13-3 can be reformulated in such a way
that it applies to convergence acceleration processes.

\medskip

\beginEinzug \sl \parindent = 0 pt

\Auszug {\bf Theorem 13-4:} Let us assume that a sequence
transformation $G_k^{(n)} (s_n, s_{n+1}, \ldots, s_{n+k})$ with $k,n
\in \N_0$ is invariant under translation according to eq. (3.1-4). Then
a necessary and sufficient condition that this sequence transformation
preserves the limit $s$ of a convergent sequence $\Seqn s$ on a path
${\cal P} = \Seq {(n_j, k_j)}$ with $j \in \N_0$ is that
$$
\lim_{j \to \infty} \; G_{k_j}^{({n_j})}
(s_{n_j} - s, s_{n_j + 1} - s, \ldots, s_{n_j + k_j} - s)
\; = \; 0 \, .
\tag
$$

\endEinzug

\medskip

What is gained if summation and convergence acceleration processes are
analyzed with the help of theorems 13-3 and 13-4. Since the limit or
antilimit $s$ of a sequence $\Seqn s$ is normally not known, it is in
most cases very hard or even impossible to estimate, how close
$G_k^{(n)} (s_n, s_{n+1}, \ldots, s_{n+k})$ and $s$ are. However, it
will become clear in the sequel that it is frequently comparatively
easy to obtain a theoretical estimate for the magnitude of the error
term $G_k^{(n)} (s_n - s, s_{n+1} - s, \ldots, s_{n+k} - s)$ and its
dependence upon $k$ and $n$.

Theorems 13-3 and 13-4 remain of course valid if the sequence
transformation $G_k^{(n)}$, which only depends upon $k+1$ sequence
elements $s_n, s_{n+1}, \ldots, s_{n+k}$, is replaced by a sequence
transformation $\Gamma_k^{(n)}$ which in addition to the $k+1$ sequence
elements $s_n, s_{n+1}, \ldots, s_{n+k}$ also depends upon $k+1$
remainder estimates $\omega_n, \omega_{n+1}, \ldots, \omega_{n+k}$.

In this section we shall try to make some quantitative predictions
about the magnitude of the summation error if a given sequence
transformation, which is of the type of eq. (12.2-2), is used for the
summation of a divergent Stieltjes series. Unfortunately, it seems that
such an error analysis cannot be done in the case of a completely
arbitrary Stieltjes series. However, if we apply sequence
transformations of the type of eq. (12.2-2) to some suitably chosen
model sequences, valuable insight into the mechanism as well as the
power of these sequence transformations can be gained.

Our error analysis will be based upon theorems 13-3 and 13-4, i.e., we
shall try to estimate the magnitude of the error term $G_k^{(n)} (s_n -
s, s_{n+1} - s, \ldots, s_{n+k} - s)$ and its dependence upon $k$ and
$n$. In addition, we assume that the sequences $\Seqn s$ and $\Seqn
\omega$ possess the following properties:

\medskip

\beginBeschreibung \zu \Laenge{(S-0):} \sl

\item{(S-0):} The elements of $\Seqn s$ are the partial sums of an
infinite series which either converges to some limit $s$ or in the case
of divergence can be summed to give $s$.
\item{(S-1):} The elements of the sequence $\Seqn \omega$ of remainder
estimates for $\Seqn s$ are strictly alternating in sign.
\item{(S-2):} For all $n \in \N_0$ the ratio $(s_n - s) / \omega_n$ can
be expressed as a factorial series, i.e.,
$$
\frac{s_n - s} {\omega_n} \; = \; \sum_{j=0}^{\infty}
\frac {c_j} {(\beta + n)_j} \, ,
\qquad \beta \in \R_{+} \, , \quad n \in \N_0 \, .
\tag
$$

\endBeschreibung

\medskip

On the basis of these assumptions the summation of divergent Stieltjes
series as well as the acceleration of the convergence of certain
alternating series can be analyzed.

Concerning (S-1) it should be remarked that if we chose the remainder
estimates $\Seqn \omega$ according to eq. (13.1-16) then the positivity
of the Stieltjes moments $\mu_n$ according to eq. (13.1-2) implies that
we are restricted to power series with positive arguments $z$. If $z$
would be an arbitrary complex number, it could not be guaranteed that
our remainder estimates $\omega_n$ will have strictly alternating signs
if they are chosen according to eq. (13.1-16).

The requirement that $(s_n - s)/\omega_n$ can be expressed as a
factorial series according to eq. (13.2-4) may appear to be somewhat
restrictive. However, this is not necessarily more restrictive than the
analogous requirement that $(s_n - s)/\omega_n$ can be expressed as a
series in inverse powers of $\beta+n$,
$$
\frac{s_n - s} {\omega_n} \; = \; \sum_{j=0}^{\infty}
\frac {c_j^{\prime}} {(\beta + n)^j} \, ,
\qquad \beta \in \R_{+} \, , \quad n \in \N_0 \, .
\tag
$$
In Nielsen's book it is described how inverse power series and
factorial series can be transformed into each other (see pp. 272 - 282
of ref. [77]). Assumptions (S-0) - (S-2) will now be used to obtain
quantitative error estimates in summation and convergence acceleration
processes.

\medskip

\beginEinzug \sl \parindent = 0 pt

\Auszug {\bf Theorem 13-5:} Let us assume that the sequences $\Seqn s$
and $\Seqn \omega$ satisfy (S-0) - (S-2) and that the sequence
transformation ${\cal S}_k^{(n)} (\beta, s_n, \omega_n)$, eq. (8.2-7),
is used for the transformation of $\Seqn s$. Then we obtain for fixed
$k \in \N$ and for all $n \in \N_0$ the following estimate for the
error term:
$$
\left\vert
{\cal S}_k^{(n)} (\beta, s_n, \omega_n) \, - \, s
\right\vert
\; \le \;
\left\vert
\frac {\omega_n} {(\beta + n)_{2 k} } \, \sum_{j=0}^{\infty} \,
\frac {c_{k + j} \, (j+1)_k } {(\beta + n + 2 k)_j }
\right\vert \, .
\tag
$$
This implies for fixed $k \in \N$ and for large values of $n$ the
following order estimate:
$$
\frac
{{\cal S}_k^{(n)} (\beta, s_n, \omega_n) \, - \, s} {s_n - s}
\; = \; O (n^{- 2 k}) \, , \qquad n \to \infty \, .
\tag
$$

\endEinzug

\medskip

\noindent {\it Proof: } First, we observe that ${\cal S}_k^{(n)}
(\beta, s_n, \omega_n)$, eq. (8.2-7), is invariant under translation
according to eq. (3.1-4). This implies that the magnitude of the
transformation error and its dependence upon $k$ and $n$ can be
analyzed by estimating the magnitude of ${\cal S}_k^{(n)} (\beta, s_n -
s, \omega_n)$. The starting point of our analysis is eq. (8.2-6) which
is rewritten in the following way:
$$
{\cal S}_k^{(n)} (\beta, s_n - s, \omega_n) \; = \;
\frac
{\Delta^k \{ (\beta + n)_{k-1} \, (s_n - s) \, / \, \omega_n \} }
{\Delta^k \{ (\beta + n)_{k-1} \, / \, \omega_n \} } \, .
\tag
$$

In the numerator in eq. (13.2-8) $(s_n - s) / \omega_n$ is replaced by
the factorial series (13.2-4) yielding
$$
\beginAligntags
" \Delta^k \, \frac {(\beta + n)_{k-1} \, (s_n - s)} {\omega_n}
\; = \; \Delta^k \, (\beta + n)_{k-1} \, \sum_{j=0}^{\infty} \,
\frac {c_j} {(\beta + n)_j}
\\ \tag
" = \; \Delta^k \, \sum_{j=0}^{k-1} \,
c_j \, (\beta + n + j)_{k-j-1} \, + \,
\Delta^k \, \sum_{j=0}^{\infty}
\frac {c_{k+j} } {(\beta + n + k - 1)_{j+1} } \, .
\\ \tag
\endAligntags
$$

Next, we derive from eq. (8.4-10) the following two relationships:
$$
\beginAligntags
" \Delta^k \, (a + n)_m \; " = \;
" (-1)^k \, (-m)_k \, (a + n + k)_{m-k} \, ,
\\ \tag
" \Delta^k \, [1 / (b + n)_m] \; " = \;
" (-1)^k \, \, (m)_k / (b + n)_{k + m} \, .
\\ \tag
\endAligntags
$$

If these two relationships are used in eq. (13.2-10), we find that the
first sum vanishes since it is a polynomial of degree $k-1$ in $n$ and
we obtain for the numerator in eq. (13.2-8):
$$
\Delta^k \, \frac {(\beta + n)_{k-1} \, (s_n - s)} {\omega_n}
\; = \; \frac
{\Gamma (\beta + n + k - 1)} {\Gamma (\beta + n + 2 k)} \,
\sum_{j=0}^{\infty} \, \frac
{c_{k+j} \, (j+1)_k} {(\beta + n + 2 k)_j} \, .
\tag
$$

With the help of eq. (2.4-8) we obtain for the denominator in eq.
(13.2-8):
$$
\Delta^k \, \frac {(\beta + n)_{k-1}} {\omega_n} \; = \;
(-1)^k \, \sum_{j=0}^{k} \, (-1)^j \binom {k} {j}
\frac {(\beta + n + j)_{k-1}} {\omega_{n + j} } \, .
\tag
$$

According to (S-1) the elements of $\Seqn \omega$ strictly alternate in
sign. This gives us immediately the following estimate:
$$
\vert (\beta + n)_{k-1} \, / \, \omega_n \vert \; \le \;
\vert \Delta^k \{ (\beta + n)_{k-1} \, / \, \omega_n \} \vert
\, .
\tag
$$
Combination of eqs. (13.3-8), (13.2-13) and (13.2-15) gives eq.
(13.2-6). The order estimate (13.2-7) follows from the fact that
according to eq. (13.2-4)
$$
(s_n - s) / \omega_n \; = \; c_0 \, [1 + O (n^{-1})] \, ,
\qquad n \to \infty \, ,
\tag
$$
and that $(\beta + n)_{2 k} \; = \; O (n^{2 k})$ as $n \to \infty$.
This proves theorem 13-5.

\medskip

Since the denominator sum (13.2-14), which consists of $k+1$ terms, is
estimated by a single term according to eq. (13.2-15), the error
estimate (13.2-6) is quite conservative.

It is a typical feature of the error estimate (13.2-6) and also of some
analogous error estimates for other sequence transformations, which
will be derived later in this section, that the error estimate is
directly proportional to $\omega_n$. Consequently, no distinction
between convergent and divergent sequences $\Seqn s$ of partial sums
has to be made. It also follows from the error estimate (13.2-6) that
${\cal S}_k^{(n)} (\beta, s_n, \omega_n)$, eq. (8.2-7), is able to sum
a divergent series satisfying (S-0) - (S-2) on a horizontal path if the
coefficients $c_j$ of the factorial series (13.2-4) do not grow too
fast in magnitude as $j \to \infty$.

The next theorem, which can be proved in essentially the same way as
theorem 13-5, shows that our error analysis is able to distinguish
between the sequence transformation ${\cal S}_k^{(n)} (\beta, s_n,
\omega_n)$, eq. (8.2-7), and its mild generalization ${\cal S}_{k,
\ell}^{(n)} (\beta, s_n, \omega_n)$, eq. (8.2-8), with $\ell \ge 1$.

\medskip

\beginEinzug \sl \parindent = 0 pt

\Auszug {\bf Theorem 13-6:} Let us assume that the sequences $\Seqn s$
and $\Seqn \omega$ satisfy (S-0), (S-1), and
$$
\frac{s_n - s} {\omega_n} \; = \; \sum_{j=0}^{\infty}
\frac {d_j} {(\beta + n + \ell)_j} \, , \qquad \beta \in \R_{+}
\, , \quad \ell \in \N \, , \quad n \in \N_0 \, ,
\tag
$$
and that the sequence transformation ${\cal S}_{k, \ell}^{(n)} (\beta,
s_n, \omega_n)$, eq. (8.2-8), with $\ell \ge 1$ is used for the
transformation of $\Seqn s$. Then we obtain for fixed $k, \ell \in \N$
with $k \ge \ell + 1$ and for all $n \in \N_0$ the following estimate
for the error term:
$$
\left\vert
{\cal S}_{k, \ell}^{(n)} (\beta, s_n, \omega_n) \, - \, s
\right\vert
\; \le \;
\left\vert
\frac {\omega_n} {(\beta + n + \ell)_{2 k - \ell} } \,
\sum_{j=0}^{\infty} \,
\frac {d_{k - \ell + j} \, (j+1)_k } {(\beta + n + 2 k)_j }
\right\vert \, .
\tag
$$
This implies for fixed $k, \ell \in \N$ with $k \ge \ell + 1$ and for
large values of $n$ the following order estimate:

$$
\frac
{{\cal S}_{k, \ell}^{(n)} (\beta, s_n, \omega_n) \, - \, s}
{s_n - s}
\; = \; O (n^{\ell - 2 k}) \, , \qquad n \to \infty \, .
\tag
$$

\endEinzug

\medskip

If we compare eqs. (13.2-4) and (13.2-17), it seems that in theorems
13-5 and 13-6 the existence of two different factorial series
expansions for the ratio $(s_n - s)/\omega_n$ are assumed. However, the
factorial series (13.2-4) and (13.2-17) are not independent. In
Nielsen's book it is shown how a factorial series of the type of eq.
(13.2-4) can be transformed into a factorial series of the type of eq.
(13.2-17) (see pp. 252 - 253 of ref. [77]).

Theorems 13-5 and 13-6 indicate that ${\cal S}_k^{(n)} (\beta, s_n,
\omega_n)$, eq. (8.2-7), should normally be more efficient than its
generalization ${\cal S}_{k, \ell}^{(n)} (\beta, s_n, \omega_n)$, eq.
(8.2-8), with $\ell \ge 1$.

In the following theorem, the efficiency of Drummond's sequence
transformation ${\cal D}_k^{(n)} (s_n, \omega_n)$, eq. (9.5-4), is
analyzed.

\medskip

\beginEinzug \sl \parindent = 0 pt

\Auszug {\bf Theorem 13-7:} Let us assume that the sequences $\Seqn s$
and $\Seqn \omega$ satisfy (S-0) - (S-2) and that Drummond's sequence
transformation ${\cal D}_k^{(n)} (s_n, \omega_n)$, eq. (9.5-4), is used
for the transformation of $\Seqn s$. Then we obtain for fixed $k \in
\N$ and for all $n \in \N_0$ the following estimate for the error term:
$$
\left\vert
{\cal D}_k^{(n)} (\beta, s_n, \omega_n) \, - \, s
\right\vert
\; \le \;
\left\vert
\frac {\omega_n} {(\beta + n)_{k + 1} } \, \sum_{j=0}^{\infty} \,
\frac {c_{j + 1} \, (j+1)_k } {(\beta + n + k + 1)_j }
\right\vert \, .
\tag
$$
This implies for fixed $k \in \N$ and for large values of $n$ the
following order estimate:
$$
\frac
{{\cal D}_k^{(n)} (s_n, \omega_n) \, - \, s} {s_n - s}
\; = \; O (n^{- k - 1}) \, , \qquad n \to \infty \, .
\tag
$$

\endEinzug

\medskip

This theorem, which can be proved in exactly the same way as theorem
13-5, indicates that in particular for larger values of $k$ Drummond's
sequence transformation ${\cal D}_k^{(n)} (s_n, \omega_n)$, eq.
(9.5-4), should be significantly less powerful than ${\cal S}_k^{(n)}
(\beta, s_n, \omega_n)$, eq. (8.2-7), or its mild generalization ${\cal
S}_{k, \ell}^{(n)} (\beta, s_n, \omega_n)$, eq. (8.2-8), with $\ell \ge
1$.

It would be interesting to do the same kind of error analysis also for
the sequence transformation ${\cal M}_k^{(n)} (\gamma, s_n, \omega_n)$,
eq. (9.2-6), and its mild generalization ${\cal M}_{k, \ell}^{(n)}
(\gamma, s_n, \omega_n)$, eq. (9.2-7). However, if we would try to
estimate the error term of this transformation for arbitrary sequences
$\Seqn s$ and $\Seqn \omega$ satisfying (S-0) - (S-2) we would in
general end up with very complicated formulas which would contribute
little to our understanding. This is due to the fact that for arbitrary
$\beta$ and $\gamma$ we would have to use Leibniz' theorem for finite
differences (see p. 35 of ref. [72]),
$$
\Delta^k \, [ f (n) g (n) ] \; = \; \sum_{j=0}^k \,
\binom {k} {j} \, [\Delta^j \, f (n)] \, [\Delta^{k-j} g(n)] \, .
\tag
$$

Much more revealing and enlightening is, however, the following
observation:

\medskip

\beginEinzug \sl \parindent = 0 pt

\Auszug {\bf Theorem 13-8:} Assume that $\gamma = \beta + k - 2$ holds.
Then,
$$
{\cal M}_k^{(n)} (\gamma, s_n, \omega_n) \; = \;
{\cal S}_k^{(n)} (\beta, s_n, \omega_n) \, .
\tag
$$

\endEinzug

\medskip

\noindent {\it Proof: } If we use the following relationship for
Pochhammer symbols (see eq. (I.5) on p. 239 of ref. [104]),
$$
(a - m)_m \; = \, (-1)^m \, (1 - a)_m \, ,
\tag
$$
we obtain
$$
(- \gamma - n)_{k-1} \; = \;
(-1)^{k-1} \, (n + \gamma - k + 2)_{k-1} \, .
\tag
$$

If we insert this relationship into eq. (9.2-6) and use $\gamma = \beta
+ k - 2$, we obtain eq. (8.2-7).

\medskip

Theorem 13-8 does not imply that the two strings ${\cal M}_j^{(n)}
(\gamma, s_n, \omega_n)$ and ${\cal S}_j^{(n)} (\beta, s_n, \omega_n)$
with $0 \le j \le k$ are identical if $\gamma = \beta + k -2$ holds.
Only the last elements of the two strings are guaranteed to be
identical, but not the others.

In the case of Levin's sequence transformation ${\cal L}_k^{(n)}
(\beta, s_n, \omega_n)$, eq. (7.1-7), an analysis of the magnitude of
the error term as in theorem 13-5 would again be very complicated and
would lead to lengthy and messy expressions. This is due to the fact
that Levin's sequence transformation is constructed on the basis of the
model sequence (7.1-1), which is merely a truncation of the inverse
power series (13.2-5) after $k$ terms, and that in the calculus of
finite differences Pochhammer symbols and not powers are the most
simple functions. However, at least some order estimates can be
obtained relatively easily in the case of the Levin transformation.

\medskip

\beginEinzug \sl \parindent = 0 pt

\Auszug {\bf Theorem 13-9:} Let us assume that the sequences $\Seqn s$
and $\Seqn \omega$ satisfy (S-0), (S-1), and eq. (13.2-5) and that
Levin's sequence transformation ${\cal L}_k^{(n)} (\beta, s_n,
\omega_n)$, eq. (7.1-7), is used for the transformation of $\Seqn s$.
Then we obtain for large values of $n$ and for fixed $k \in \N$ the
following order estimate:
$$
\frac
{{\cal L}_k^{(n)} (\beta, s_n, \omega_n) \, - \, s} {s_n - s}
\; = \; O (n^{- 2 k}) \, , \qquad n \to \infty \, .
\tag
$$

\endEinzug

\medskip

\noindent {\it Proof: } Obviously, ${\cal L}_k^{(n)} (\beta, s_n,
\omega_n)$ is invariant under translation according to eq. (3.1-4).
This implies that the magnitude of the error term and its dependence
upon $k$ and $n$ can be analyzed by estimating the magnitude of ${\cal
L}_k^{(n)} (\beta, s_n - s, \omega_n)$. The starting point of our
analysis is eq. (7.1-6) which is rewritten in the following way:
$$
{\cal L}_k^{(n)} (\beta, s_n - s, \omega_n) \; = \;
\frac
{\Delta^k \{ (\beta + n)^{k-1} \, (s_n - s) \, / \, \omega_n \} }
{\Delta^k \{ (\beta + n)^{k-1} \, / \, \omega_n \} } \, .
\tag
$$

In the numerator in eq. (13.2-27) $(s_n - s) / \omega_n$ is replaced by
the power series (13.2-5) yielding
$$
\beginAligntags
" \Delta^k \frac {(\beta + n)^{k-1} \, (s_n - s)} {\omega_n}
\; = \; \Delta^k \, (\beta + n)^{k-1} \, \sum_{j=0}^{\infty} \,
\frac {c_j^{\prime} } {(\beta + n)^j}
\\ \tag
" = \; \Delta^k \, \sum_{j=0}^{k-1} \,
c_j^{\prime} (\beta + n)^{k-j-1} \, + \,
\Delta^k \, \sum_{j=0}^{\infty}
\frac {c_{k+j}^{\prime} } {(\beta + n)^{j+1} } \, .
\\ \tag
\endAligntags
$$

The first sum on the right-hand side of eq. (13.2-29) is annihilated by
$\Delta^k$ since it is a polynomial of degree $k - 1$ in $n$. The large
$n$ behaviour of the second sum can be estimated with the help of the
relationship
$$
\Delta^m \, n^{ - \alpha} \; = \; O (n^{ - \alpha - m}) \, ,
\qquad \alpha > 0 \, ,
\tag
$$
to give
$$
\Delta^k \frac {(\beta + n)^{k-1} \, (s_n - s)} {\omega_n}
\; = \; O (n^{- k - 1}) \, , \qquad n \to \infty \, .
\tag
$$

Since the remainder estimates $\Seqn \omega$ are strictly alternating
in sign according to (S-1), we obtain the following estimate for the
denominator in eq. (13.2-27),
$$
\vert (\beta + n)^{k-1} / \omega_n \vert \le
\vert \Delta^k \{(\beta + n)^{k-1} / \omega_n\} \vert \, .
\tag
$$

If we combine eqs. (13.2-27), (13.2-31), and (13.2-32) and take into
account that $(\beta + n)^{k-1} = O (n^{k-1})$ as $n \to \infty$, we
see that theorem 13-9 is correct.

\medskip

In the following theorem, which can be proved in essentially the same
way as theorem 13-9, an order estimate for Levin's generalized sequence
transformation ${\cal L}_{k, \ell}^{(n)} (\beta, s_n, \omega_n)$, eq.
(7.1-8), with $\ell \ge 1$ is derived.

\medskip

\beginEinzug \sl \parindent = 0 pt

\Auszug {\bf Theorem 13-10:} Let us assume that the sequences $\Seqn s$
and $\Seqn \omega$ satisfy (S-0), (S-1), and eq. (13.2-5) and that
Levin's generalized sequence transformation ${\cal L}_{k, \ell}^{(n)}
(\beta, s_n, \omega_n)$, eq. (7.1-8), with $\ell \ge 1$ is used for the
transformation of $\Seqn s$. Then we obtain for fixed $k, \ell \in \N$
with $k \ge \ell + 1$ and for large values of $n$ the following order
estimate:
$$
\frac
{{\cal L}_{k, \ell}^{(n)} (\beta, s_n, \omega_n)
\, - \, s} {s_n - s} \; = \; O (n^{\ell - 2 k})
\, , \qquad n \to \infty \, .
\tag
$$

\endEinzug

\medskip

A comparison of the order estimates (13.2-7), (13.2-19), (13.2-26), and
(13.2-33) shows that Levin's sequence transformation ${\cal L}_k^{(n)}
(\beta, s_n, \omega_n)$, eq. (7.1-7), should be roughly comparable with
${\cal S}_k^{(n)} (\beta, s_n, \omega_n)$, eq. (8.2-7), and that for
fixed $\ell \ge 1$ Levin's generalized transformation ${\cal L}_{k,
\ell}^{(n)} (\beta, s_n, \omega_n)$, eq. (7.1-8), should be roughly
comparable with ${\cal S}_{k, \ell}^{(n)} (\beta, s_n, \omega_n)$, eq.
(8.2-8). In addition, a comparison with the order estimate (13.2-21)
shows that these sequence transformations should all be significantly
more powerful than Drummond's sequence transformation ${\cal D}_k^{(n)}
(s_n, \omega_n)$, eq. (9.5-4). A more detailed comparison cannot be
made here since this would require additional knowledge about the
sequence $\Seqn s$ and the remainder estimates $\Seqn \omega$.

The error analysis of this section is restricted to convergent or
divergent sequences $\Seqn s$ with strictly alternating remainder
estimates $\Seqn \omega$. This restriction is essential because
otherwise the denominator sums of the pertaining sequence
transformations cannot be estimated by a single term as it was for
instance done in eqs. (13.2-15) and (13.2-32). If we want to analyze
the transformation of sequences with nonalternating remainders,
additional assumptions about the behaviour of the remainders have to be
made. For instance, in the case of logarithmic convergence we could
assume something like
$$
\frac {(\beta+n)_{k-1}}{\omega_n} \; = \;
\sum_{j=0}^{\infty} \, c_j \,
\frac {\Gamma (\beta + n + k - 1)} {\Gamma (\delta + n + j)}
\, , \qquad \beta, \delta \in \R_{+} \, ,
\quad n \in \N_0 \, ,
\tag
$$
because then the denominator of eq. (13.2-8) could be computed with the
help of eq. (8.4-10) and we would obtain an explicit expression for the
transformation error ${\cal S}_k^{(n)} (\beta, s_n - s, \omega_n)$.
With the help of similar assumptions Sidi [56,105] could derive various
error estimates for Levin's sequence transformation ${\cal L}_k^{(n)}
(\beta, s_n, \omega_n)$, eq. (7.1-7), in convergence acceleration and
summation processes.

\medskip

\Abschnitt Summation of the Euler series

\smallskip

\aktTag = 0

Pad\'e approximants are generally accepted to be valuable numerical tools
for the treatment of scientific problems. Therefore, it is certainly
interesting to compare Pad\'e approximants with the other sequence
transformations of this report.

Unfortunately, the theoretical error estimates for the Pad\'e summation
of a Stieltjes series, which can be found in the literature, are not
directly comparable with the error estimates of section 13.2, in which
the error is always directly proportional to the remainder estimate
$\omega_n$. For instance, in the articles by Allen, Chui, Madych,
Narcowich, and Smith [101], and by Karlsson and Sydow [102] the
summation error is expressed in terms of polynomials which are
orthogonal with respect to the measure $\psi (t)$ in the Stieltjes
integral (13.1-1).

Consequently, we first would have to derive something like theorem 13-5
for Pad\'e approximants before we could compare Pad\'e approximants and the
sequence transformations ${\cal L}_k^{(n)} (\beta, s_n, \omega_n)$, eq.
(7.1-7), ${\cal S}_k^{(n)} (\beta, s_n, \omega_n)$, eq. (8.2-7), ${\cal
M}_k^{(n)} (\gamma, s_n, \omega_n)$, eq. (9.2-6), and ${\cal D}_k^{(n)}
(s_n, \omega_n)$, eq. (9.5-4), and their mild generalizations ${\cal
L}_{k, \ell}^{(n)} (\beta, s_n, \omega_n)$, eq. (7.1-8), ${\cal S}_{k,
\ell}^{(n)} (\beta, s_n, \omega_n)$, eq. (8.2-8), and ${\cal M}_{k,
\ell}^{(n)} (\gamma, s_n, \omega_n)$, eq. (9.2-7), with respect to
their ability of summing divergent Stieltjes series or accelerating the
convergence of some alternating series. Unfortunately, no such theorem
could be derived which treats the Pad\'e summation of an arbitrary
Stieltjes series.

However, there is a notable exception. In the case of the Euler
integral, eq. (1.1-6), and its associated asymptotic series, the
so-called Euler series, eq. (1.1-7), Sidi [84] could show that their
Pad\'e approximants can be expressed in closed form via Drummond's
sequence transformation ${\cal D}_k^{(n)} (s_n, \omega_n)$, eq.
(9.5-4), with $\omega_n = a_{n+1}$. Sidi's proof is based upon the
well-known fact that Wynn's $\epsilon$ algorithm, eq. (4.2-1), which
according to eq. (4.2-10) is able to compute the Pad\'e approximants $[n
+ k / k]$, is exact for the model sequence (4.1-3). In the case of the
Euler series, eq. (1.1-7), we have
$$
\beginAligntags
" s_n \; " = \;
" \sum_{\nu=0}^{n} \, (-1)^{\nu} \, {\nu}! \, z^{\nu} \, ,\\ \tag
" a_n \; " = \; " (-1)^n \, n! \, z^n \, .
\\ \tag
\endAligntags
$$

If we insert these relationships into the model sequence (4.1-3), we
obtain
$$
s_n \; = \; s \, + \, \sum_{j=0}^{k-1} \,
c_j \, (-1)^{n+j+1} \, (n+j+1)! \, z^{n+j+1} \, .
\tag
$$

This model sequence for the Pad\'e approximants $[n + k / k]$ of the
Euler series can be rewritten in the following way:
$$
\frac {s_n - s} {(-1)^{n+1} (n+1)! z^{n+1} } \; = \;
\sum_{j=0}^{k-1} \, c_j \, (-1)^j \, z^j \, (n + 2)_j \, .
\tag
$$

The sum on the right-hand side of eq. (13.3-4) is a polynomial of
degree $k-1$ in $n$. Consequently, it can be annihilated by the
difference operator $\Delta^k$. Hence, it follows from eq. (9.5-2) that
the Pad\'e approximants $[n + k / k]$ for the Euler series can be
expressed in closed form in terms of Drummond's sequence
transformation, eq. (9.5-4), $$
[n + k / k] \; = \; \epsilon_{2 k}^{(n)}
\; = \; {\cal D}_k^{(n)} (s_n, a_{n+1}) \, ,
\qquad k,n \in \N_0 \, .
\tag
$$

From this relationship we may conclude that in the case of the Euler
series Drummond's sequence transformation is much more efficient than
Wynn's $\epsilon$ algorithm. If the diagonal Pad\'e approximant $[n /
n]$, which according to eq. (4.1-8) is a ratio of two polynomials $p_n
(z)$ and $q_n (z)$ of degree $n$ in $z$, is computed with the help of
Wynn's $\epsilon$ algorithm as $\epsilon_{2 n}^{(0)}$, the partial sums
$s_0, s_1, \ldots, s_{2 n}$ of the Euler series will be needed. If the
same diagonal Pad\'e approximant $[n / n]$ is computed as ${\cal
D}_n^{(0)} (s_0, a_1)$, then according to eq. (9.5-4) only the partial
sums $s_0, s_1, \ldots, s_{n + 1}$ will be needed. Consequently, in the
case of the Euler series, Drummond's sequence transformation is
approximately twice as efficient as Wynn's $\epsilon$ algorithm.

How can this behaviour be explained? If Drummond's sequence
transformation is applied to a sequence of partial sums of the Euler
series, then we find that ${\cal D}_n^{(0)} (s_0, a_1)$ is the ratio of
two polynomials $p_n (z)$ and $q_n (z)$ of degree $n$ in $z$. However,
it follows from eq. (9.5-4) that the $2 n + 2$ coefficients of the two
polynomials are not all independent. In fact, these two polynomials
$p_n (z)$ and $q_n (z)$ are completely determined by the $n + 2$ terms
$a_0, a_1, \ldots, a_{n+1}$. If the same ratio $p_n (z) / q_n (z)$ is
computed via Wynn's $\epsilon$ algorithm, it is implicitly assumed that
the $2 n + 2$ coefficients of the two polynomials are independent apart
from a common normalization condition. This implies that Wynn's
$\epsilon$ algorithm needs $2 n + 1$ independent conditions -- in this
case the $2 n + 1$ partial sums $s_0, s_1, \ldots, s_{2 n}$ -- for the
construction of the ratio $p_n (z) / q_n (z)$.

Actually, it is a typical feature of all sequence transformation
$T_k^{(n)} (s_n, \omega_n)$ of the type of eq. (12.2-2) that the
coefficients of the numerator and denominator sums are not independent.

It is a natural idea to try to apply theorem 13-7, which gives an error
estimate for Drummond's sequence transformation, also for of the Pad\'e
summation of the Euler series. Assumptions (S-0) and (S-1) are
obviously satisfied. However, it is not clear whether and how a
sequence $\Seqn \omega$ of remainder estimates can be found such that
assumption (S-2), which requires that $(s_n - s)/\omega_n$ can be
represented as a factorial series according to eq. (13.2-4), is valid.
No explicit proof for the existence of such a factorial series could be
found in the case of the Euler series, if the remainder estimates were
chosen according to
$$
\omega_n \; = \; (-1)^{n+1} \, (n+1)! \, z^{n+1} \, ,
\qquad n \in \N_0 \, .
\tag
$$

Consequently, it can only be investigated numerically whether the error
analysis of theorem 13-7 provides an adequate description of the Pad\'e
summation of the Euler series.

A close relative of the Euler integral, eq. (1.1-6), is the so-called
exponential integral
$$
E_1 (z) \; = \; \int\nolimits_{z}^{\infty} \,
\frac {\e^{- x} } {x} \d x \, .
\tag
$$

By means of some elementary operations we find:
$$
z \, \e^z \, E_1 (z) \; = \;
\int\nolimits_{0}^{\infty} \,
\frac {\e^{- t} \d t} {1 + t / z} \, .
\tag
$$

If we compare this relationship with eq. (1.1-6) and also use eq.
(1.1-7), we see that the associated Stieltjes series of the integral in
eq. (13.3-8) is the Euler series with argument $1/z$,
$$
z \, \e^z \, E_1 (z) \; \sim \;
\sum_{m=0}^{\infty} \, (-1)^m \, m! \, z^{- m} \; = \;
{}_2 F_0 (1,1; - 1/z) \, , \qquad z \to \infty \, .
\tag
$$

The radius of convergence of the hypergeometric series ${}_2 F_0$ on
the right-hand side of eq. (13.3-9) is zero, i.e., the series diverges
quite rapidly for all finite values of $z$. Since reliable programs for
the exponential integral $E_1 (z)$ with $z \in \R_{+}$ are available,
eq. (13.3-9) is well suited to test the ability of a sequence
transformation of summing even wildly divergent series. In this report,
the exponential integral $E_1 (z)$ will be computed with the help of
the routine S13AAF of the NAG FORTRAN library [106]. This function
computes an approximation for the exponential integral in DOUBLE
PRECISION (15 -- 16 decimal digits) using appropriate Chebyshev
expansions.

\beginFloat

\medskip

\beginTabelle [to \kolumnenbreite]
\beginFormat \rechts " \rechts " \mitte " \mitte " \mitte
\endFormat
\+ " \links {\bf Table 13-1} \@ \@ \@ \@ " \\
\+ " \links {Summation of the asymptotic series ${}_2 F_0 (1,1; -1/z)
\; = \; z \, \e^z \, E_1 (z)$ for $z \; = \; 3$}
\@ \@ \@ \@ " \\
\- " \- " \- " \- " \- " \- " \\ \sstrut {} {1.5 \jot} {1.5 \jot}
\+ " \rechts {$n$} " \mitte {partial sum $s_n$}
" ${\cal A}_{\Ent {n/2}}^{(n - 2 \Ent {n/2})}$
" ${\cal D}_n^{(0)} (s_0, a_1)$
" $\epsilon_{2 \Ent {n/2}}^{(n - 2 \Ent {n/2})}$ " \\
\+ " " \mitte {eq. (13.1-10)} " eq. (5.1-15) " eq. (9.5-4) " eq.
(4.2-1) " \\
\- " \- " \- " \- " \- " \- " \\ \sstrut {} {1 \jot} {1 \jot}
\+ " 10 " $ 0.4831550069 \times 10^{02}$ " 0.78625130019479 "
0.78625125348502 " 0.78626367674141 " \\
\+ " 11 " $ -0.1770160037 \times 10^{03}$ "  0.78625114835779 "
0.78625123263883 "  0.78624220653206 " \\
\+ "  12 "  $ 0.7243100137 \times 10^{03}$  "  0.78625122394910 "
0.78625122525386 "  0.78625447790898 " \\
\+ "  13 " $-0.3181436062 \times 10^{04}$  "  0.78625121766831 "
0.78625122252501 "  0.78624881508686 " \\
\+ "  14 " $ 0.1504537896 \times 10^{05}$  "  0.78625122089403 "
0.78625122147819 "  0.78625215335611 " \\
\+ "  15 " $-0.7608869613 \times 10^{05}$  "  0.78625122063943 "
0.78625122106292 "  0.78625052018310 " \\
\+ "  16 " $ 0.4099597043 \times 10^{06}$  "  0.78625122077179 "
0.78625122089311 "  0.78625150842397 " \\
\+ "  17 " $-0.2344314565 \times 10^{07}$  "  0.78625122076057 "
0.78625122082175 "  0.78625100153477 " \\
\+ "  18 " $ 0.1418133105 \times 10^{08}$  "  0.78625122076626 "
0.78625122079099 "  0.78625131522011 " \\
\+ "  19 " $-0.9048109119 \times 10^{08}$  "  0.78625122076568 "
0.78625122077742 "  0.78625114787954 " \\
\+ "  20 " $ 0.6072683904 \times 10^{09}$  "  0.78625122076597 "
0.78625122077131 "  0.78625125348502 " \\
\+ "  21 " $-0.4276977981 \times 10^{10}$  "  0.78625122076594 "
0.78625122076850 "  0.78625119524201 " \\
\+ "  22 " $ 0.3154082874 \times 10^{11}$  "  0.78625122076595 "
0.78625122076718 "  0.78625123263883 " \\
\+ "  23 " $-0.2430623561 \times 10^{12}$  "  0.78625122076596 "
0.78625122076656 "  0.78625121141456 " \\
\+ "  24 " $ 0.1953763123 \times 10^{13}$  "  0.78625122076595 "
0.78625122076626 "  0.78625122525386 " \\
\+ "  25 " $-0.1635311587 \times 10^{14}$  "  0.78625122076596 "
0.78625122076611 "  0.78625121720071 " \\
\+ "  26 " $ 0.1423065021 \times 10^{15}$  "  0.78625122076596 "
0.78625122076603 "  0.78625122252501 " \\
\+ "  27 " $-0.1285630059 \times 10^{16}$  "  0.78625122076596 "
0.78625122076600 "  0.78625121935772 " \\
\+ "  28 " $ 0.1204177785 \times 10^{17}$  "  0.78625122076596 "
0.78625122076598 "  0.78625122147819 " \\
\+ "  29 " $-0.1167898319 \times 10^{18}$  "  0.78625122076596 "
0.78625122076597 "  0.78625122019177 " \\
\+ "  30 " $ 0.1171526266 \times 10^{19}$  "  0.78625122076596 "
0.78625122076596 "  0.78625122106292 " \\
\- " \- " \- " \- " \- " \- " \\ \sstrut {} {1 \jot} {1 \jot}
\+ " \links {NAG function S13AAF} \@        "  0.78625122076594 "
0.78625122076594 "  0.78625122076594 " \\
\- " \- " \- " \- " \- " \- " \\ \sstrut {} {1 \jot} {1 \jot}
\endTabelle

\medskip

\endFloat

In table 13-1 the effect of Aitken's iterated $\Delta^2$ process, eq.
(5.1-15), of Drummond's sequence transformation, eq. (9.5-4), and of
Wynn's $\epsilon$ algorithm, eq. (4.2-1), on the partial sums
$$
s_n \; = \; \sum_{m=0}^n \, (-1)^m \, m! \; z^{- m} \, ,
\qquad n \in \N_0
\tag
$$
of the divergent series ${}_2 F_0$ in eq. (13.3-9) with $z = 3$ is
compared. In Drummond's sequence transformation, eq. (9.5-4), the
remainder estimates are chosen according to eq. (13.1-16) which in this
case means
$$
\omega_n \; = \; a_{n+1} \; = \;
(-1)^{n+1} \, (n+1)! \, z^{- n - 1} \, , \qquad n \in \N_0 \, .
\tag
$$

The partial sums and the three different transforms in table 13-1 were
computed in QUADRUPLE PRECISION (31 - 32 decimal digits). When these
computations were repeated in DOUBLE PRECISION (15 - 16 digits) in
order to study the numerical stability of the pertaining numerical
processes, it turned out that the two computations agreed at least up
to 12 decimal digits.

In all cases, the approximants were chosen in such a way that the
information, which is contained in the finite string $s_0, s_1, \ldots,
s_n$ of partial sums, is exploited optimally. This means that in the
case of Aitken's iterated $\Delta^2$ process, eq. (5.1-15), the
approximants were chosen according to eq. (5.2-6), and in the case of
Wynn's $\epsilon$ algorithm, eq. (4.2-1), they were chosen according to
eq. (4.3-6).

A comparison of these three sequence transformations is quite
interesting. Aitken's iterated $\Delta^2$ process and Wynn's $\epsilon$
algorithm are closely related since they are both generalizations of
Aitken's $\Delta^2$ process, eq. (5.1-4), and one would like to know
which one of these two generalizations fares better. In addition, since
the series (13.3-9) is the Euler series with argument $1/z$, the
validity of eq. (13.3-5) can be checked numerically by comparing the
results for Wynn's $\epsilon$ algorithm and Drummond's sequence
transformation.

The clear winner in table 13-1 is Aitken's iterated $\Delta^2$ process
which produces 14 decimal digits after $n = 23$ (there is strong
independent evidence that the last digit produced by the NAG function
S13AAF in table 13-1 is incorrect and that Aitken's iterated $\Delta^2$
process and Drummond's sequence transformation produce the correct
result). It is followed by Drummond's sequence transformation, which
reaches an accuracy of 14 decimal digits after $n = 30$, and the clear
loser is Wynn's $\epsilon$ algorithm.

The results in table 13-1 show that eq. (13.3-5) is obviously valid in
the case of the divergent series ${}_2 F_0 (1,1; -1/z)$ because we find
$$
{\cal D}_n^{(0)} (s_0, a_1) \; = \; \epsilon_{2 n}^{(0)} \, .
\tag
$$

Since ${\cal D}_{30}^{(0)} (s_0, a_1)$ is able to produce an accuracy
of 14 decimal digits, it follows from eq. (13.3-12) that Wynn's
$\epsilon$ algorithm will need the partial sums $s_0, s_1, \ldots,
s_{60}$ of the asymptotic series in eq.(13.3-9) to produce the same
accuracy.

\beginFloat

\medskip

\beginTabelle [to \kolumnenbreite]
\beginFormat \rechts " \rechts " \mitte " \mitte " \mitte
\endFormat
\+ " \links {\bf Table 13-2} \@ \@ \@ \@ " \\
\+ " \links {Summation of the asymptotic series ${}_2 F_0 (1,1; -1/z)
\; = \; z \, \e^z \, E_1 (z)$ for $z \; = \; 3$}
\@ \@ \@ \@ " \\
\- " \- " \- " \- " \- " \- " \\ \sstrut {} {1 \jot} {1 \jot}
\+ " \rechts {$n$} " \mitte {partial sum $s_n$} " $d_n^{(0)} (1, s_0)$
" ${\delta}_n^{(0)} (1, s_0)$ " ${\Delta}_n^{(0)} (17, s_0)$ " \\
\+ " " \mitte {eq. (13.1-10)} " eq. (7.3-9) " eq. (8.4-4) " eq.
(9.4-4)" \\
\- " \- " \- " \- " \- " \- " \\ \sstrut {} {1 \jot} {1 \jot}
\+ "    3 " $  0.6666666667 \times 10^{00}$  "  0.78709677419355 "
0.78672985781991 "  0.78633660627852 " \\
\+ "    4 " $  0.9629629630 \times 10^{00}$  "  0.78607714016933 "
0.78622197922362 "  0.78625813355638 " \\
\+ "    5 " $  0.4691358025 \times 10^{00}$  "  0.78628225839245 "
0.78625036724446 "  0.78625167667778 " \\
\+ "    6 " $  0.1456790123 \times 10^{01}$  "  0.78624675493384 "
0.78625141640628 "  0.78625123654802 " \\
\+ "    7 " $ -0.8477366255 \times 10^{00}$  "  0.78625162955159 "
0.78625123162756 "  0.78625121997903 " \\
\+ "    8 " $  0.5297668038 \times 10^{01}$  "  0.78625123599599 "
0.78625121903376 "  0.78625122068020 " \\
\+ "    9 " $ -0.1313854595 \times 10^{02}$  "  0.78625120523222 "
0.78625122051031 "  0.78625122077447 " \\
\+ "   10 " $  0.4831550069 \times 10^{02}$  "  0.78625122396512 "
0.78625122077239 "  0.78625122076641 " \\
\+ "   11 " $ -0.1770160037 \times 10^{03}$  "  0.78625122056582 "
0.78625122077131 "  0.78625122076576 " \\
\+ "   12 " $  0.7243100137 \times 10^{03}$  "  0.78625122068924 "
0.78625122076646 "  0.78625122076598 " \\
\+ "   13 " $ -0.3181436062 \times 10^{04}$  "  0.78625122079175 "
0.78625122076590 "  0.78625122076596 " \\
\+ "   14 " $  0.1504537896 \times 10^{05}$  "  0.78625122076354 "
0.78625122076593 "  0.78625122076595 " \\
\+ "   15 " $ -0.7608869613 \times 10^{05}$  "  0.78625122076528 "
0.78625122076595 "  0.78625122076596 " \\
\+ "   16 " $  0.4099597043 \times 10^{06}$  "  0.78625122076622 "
0.78625122076596 "  0.78625122076596 " \\
\+ "   17 " $ -0.2344314565 \times 10^{07}$  "  0.78625122076593 "
0.78625122076596 "  0.78625122076596 " \\
\+ "   18 " $  0.1418133105 \times 10^{08}$  "  0.78625122076595 "
0.78625122076596 "  0.78625122076596 " \\
\- " \- " \- " \- " \- " \- " \\ \sstrut {} {1 \jot} {1 \jot}
\+ " \links {NAG function S13AAF} \@   "  0.78625122076594 "
0.78625122076594 "  0.78625122076594 " \\
\- " \- " \- " \- " \- " \- " \\ \sstrut {} {1 \jot} {1 \jot}
\endTabelle

\medskip

\endFloat

In table 13-2 the same divergent series ${}_2 F_0$ in eq. (13.3-9) with
$z = 3$ is summed by the sequence transformations $d_n^{(0)} (\beta,
s_0)$, eq. (7.3-9), and ${\delta}_n^{(0)} (\beta, s_0)$, eq. (8.4-4),
with $\beta = 1$ and ${\Delta}_n^{(0)} (\gamma, s_0)$, eq. (9.4-4),
with $\gamma = 17$. These three sequence transformations use the same
remainder estimate (13.3-11) as ${\cal D}_n^{(0)} (s_0, a_1)$ in table
13-1.

Table 13-2 was also produced in QUADRUPLE PRECISION. When this
computation was repeated in DOUBLE PRECISION, it turned out that in the
last two columns all 14 digits agreed. Only in the case of the Levin
transformation $d_n^{(0)} (\beta, s_0)$, eq. (7.3-9), it occasionally
happened that the last digit disagreed. Thus, numerical instabilities
are no problem here.

A comparison of the results in tables 13-1 and 13-2 confirms the error
analysis in section 13.2, which indicates that Drummond's sequence
transformation ${\cal D}_k^{(n)} (s_n, \omega_n)$, eq. (9.5-4), should
be significantly less powerful than the sequence transformations ${\cal
L}_k^{(n)} (\beta, s_n, \omega_n)$, eq. (7.1-7), ${\cal S}_k^{(n)}
(\beta, s_n, \omega_n)$, eq. (8.2-7), and ${\cal M}_k^{(n)} (\gamma,
s_n, \omega_n)$, eq. (9.2-6), if the same remainder estimates $\Seqn
\omega$ are used. Even $d_n^{(0)} (\beta, s_0)$, which is somewhat
weaker than the other two transformations in table 13-2, is clearly
more powerful than the transformations in table 13-1, and both
$\Delta_{n}^{(0)} (\gamma, s_0)$ and $\delta_{n}^{(0)} (\beta, s_0)$
are approximately twice as efficient as ${\cal D}_n^{(0)} (s_0, a_1)$.
This observation is at least qualitatively in agreement with the order
estimates (13.2-7) and (13.2-21).

In view of its slow convergence a Pad\'e summation of the divergent
series ${}_2 F_0$ in eq. (13.3-9) does not seem feasible if its
argument is significantly smaller than $z = 3$ as it was chosen in
table 13-1. If, however, variants of the sequence transformations
${\cal L}_k^{(n)} (\beta, s_n, \omega_n)$, eq. (7.1-7), ${\cal
S}_k^{(n)} (\beta, s_n, \omega_n)$, eq. (8.2-7), and ${\cal M}_k^{(n)}
(\gamma, s_n, \omega_n)$, eq. (9.2-6), are used, the summation of the
divergent series in eq. (13.3-9) can be done even for relatively small
arguments. Table 13-3 shows that the sequence transformations
$d_n^{(0)} (\beta, s_0)$, eq. (7.3-9), and ${\delta}_n^{(0)} (\beta,
s_0)$, eq. (8.4-4), with $\beta = 1$ and ${\Delta}_n^{(0)} (\gamma,
s_0)$, eq. (9.4-4), with $\gamma = 29$ are able to sum the divergent
series ${}_2 F_0$ in eq. (13.3-9) with an accuracy of 14 decimal digits
even if the argument of the series is as small as $z = 1/2$.

\beginFloat

\medskip

\beginTabelle [to \kolumnenbreite]
\beginFormat \rechts " \rechts " \mitte " \mitte " \mitte
\endFormat
\+ " \links {\bf Table 13-3} \@ \@ \@ \@ " \\
\+ " \links {Summation of the asymptotic series ${}_2 F_0 (1,1; -1/z)
\; = \; z \, \e^z \, E_1 (z)$ for $z \; = \; 1/2$}
\@ \@ \@ \@ " \\
\- " \- " \- " \- " \- " \- " \\ \sstrut {} {1 \jot} {1 \jot}
\+ " \rechts {$n$} " \mitte {partial sum $s_n$} " $d_n^{(0)} (1, s_0)$
" ${\delta}_n^{(0)} (1, s_0)$ " ${\Delta}_n^{(0)} (29, s_0)$ " \\
\+ " " \mitte {eq. (13.1-10)} " eq. (7.3-9) " eq. (8.4-4) " eq.
(9.4-4)" \\
\- " \- " \- " \- " \- " \- " \\ \sstrut {} {1 \jot} {1 \jot}
\+ "   15 " $ -0.4147067254 \times 10^{17}$  "  0.46145531715043 "
0.46145531958535 "  0.46145595366489 " \\
\+ "   16 " $  0.1329725286 \times 10^{19}$  "  0.46145530923846 "
0.46145531701552 "  0.46145551453546 " \\
\+ "   17 " $ -0.4529093729 \times 10^{20}$  "  0.46145531613431 "
0.46145531625982 "  0.46145536941468 " \\
\+ "   18 " $  0.1633052915 \times 10^{22}$  "  0.46145531735759 "
0.46145531613493 "  0.46145532757622 " \\
\+ "   19 " $ -0.6214401349 \times 10^{23}$  "  0.46145531627646 "
0.46145531616450 "  0.46145531778365 " \\
\+ "   20 " $  0.2488938643 \times 10^{25}$  "  0.46145531605612 "
0.46145531620445 "  0.46145531622965 " \\
\+ "   21 " $ -0.1046565329 \times 10^{27}$  "  0.46145531622971 "
0.46145531622787 "  0.46145531618769 " \\
\+ "   22 " $  0.4609744216 \times 10^{28}$  "  0.46145531627375 "
0.46145531623807 "  0.46145531623838 " \\
\+ "   23 " $ -0.2122526902 \times 10^{30}$  "  0.46145531624564 "
0.46145531624153 "  0.46145531624494 " \\
\+ "   24 " $  0.1019714416 \times 10^{32}$  "  0.46145531623631 "
0.46145531624231 "  0.46145531624191 " \\
\+ "   25 " $ -0.5102726985 \times 10^{33}$  "  0.46145531624080 "
0.46145531624227 "  0.46145531624156 " \\
\+ "   26 " $  0.2655415912 \times 10^{35}$  "  0.46145531624283 "
0.46145531624210 "  0.46145531624194 " \\
\+ "   27 " $ -0.1434925159 \times 10^{37}$  "  0.46145531624214 "
0.46145531624197 "  0.46145531624189 " \\
\+ "   28 " $  0.8040791666 \times 10^{38}$  "  0.46145531624170 "
0.46145531624191 "  0.46145531624184 " \\
\+ "   29 " $ -0.4666476909 \times 10^{40}$  "  0.46145531624180 "
0.46145531624188 "  0.46145531624188 " \\
\+ "   30 " $  0.2801466126 \times 10^{42}$  "  0.46145531624189 "
0.46145531624187 "  0.46145531624187 " \\
\- " \- " \- " \- " \- " \- " \\ \sstrut {} {1 \jot} {1 \jot}
\+ " \links {NAG function S13AAF} \@  "  0.46145531624187 "
0.46145531624187 "  0.46145531624187 " \\
\- " \- " \- " \- " \- " \- " \\ \sstrut {} {1 \jot} {1 \jot}
\endTabelle

\medskip

\endFloat

In the case of table 13-3 it is essential to use QUADRUPLE PRECISION.
In DOUBLE PRECISION, a heavy loss of significant digits occurs. The
best results in DOUBLE PRECISION are obtained by ${\Delta}_n^{(0)}
(\gamma, s_0)$ for $n = 20$ (10 decimal digits). For larger values of
$n$, the accuracy of the results deteriorates rapidly due to numerical
instabilities, leading to nonsensical results for the Levin
transformation $d_{30}^{(0)} (\beta, s_0)$ and to only 3 digits
accuracy for ${\delta}_{30}^{(0)} (\beta, s_0)$ and
${\Delta}_{30}^{(0)} (\gamma, s_0)$.

If the other variants of the sequence transformations ${\cal L}_k^{(n)}
(\beta, s_n, \omega_n)$, eq. (7.1-7), ${\cal S}_k^{(n)} (\beta, s_n,
\omega_n)$, eq. (8.2-7), and ${\cal M}_k^{(n)} (\gamma, s_n,
\omega_n)$, eq. (9.2-6), which are based upon the remainder estimates
(7.3-4), (7.3-6), (7.3-10), and (9.4-1), are used for the summation of
the divergent series ${}_2 F_0$ in eq. (13.3-9), it turns out that
these transformations are roughly comparable with $d_n^{(0)} (\beta,
s_0)$, ${\delta}_n^{(0)} (\beta, s_0)$, or ${\Delta}_n^{(0)} (\gamma,
s_0)$, which were used in tables 13-2 and 13-3.

The other sequence transformations of this report do not sum the
divergent series ${}_2 F_0$ in eq. (13.3-9) as efficiently as the
transformations mentioned above. For instance, Brezinski's $\theta$
algorithm, eq. (10.1-9), or other transformations, which are based upon
the $\theta$ algorithm, as for instance ${\cal J}_k^{(n)}$, eq.
(10.3-6), ${\cal B}_k^{(n)}$, eq. (11.1-5), and ${\cal C}_k^{(n)}$, eq.
(11.1-12), all rank between Aitken's iterated $\Delta^2$ algorithm, eq.
(5.1-15), and Drummond's sequence transformation, eq. (9.5-4), with
respect to their ability of summing the divergent series in eq.
(13.3-9). The sequence transformations $\lambda_k^{(n)}$, eq. (11.2-1),
$\sigma _k^{(n)}$, eq. (11.2-2), and $\mu_k^{(n)}$, eq. (11.2-3), sum
the divergent series ${}_2 F_0$ in eq. (13.3-9) slightly less efficient
than Drummond's sequence transformation.

It is a remarkable fact that compared with $d_n^{(0)} (\beta, s_0)$,
eq. (7.3-9), ${\delta}_n^{(0)} (\beta, s_0)$, eq. (8.4-4), or
${\Delta}_n^{(0)} (\gamma, s_0)$, eq. (9.4-4), which were used in
tables 13-2 and 13-3, the Pad\'e summation of the divergent series ${}_2
F_0$ in eq. (13.3-9) is hopelessly inefficient, even if the Pad\'e
approximants are computed via Drummond's sequence transformation
according to eq. (13.3-5) and not via Wynn's $\epsilon$ algorithm.

How can this inferiority of Pad\'e approximants be explained? It was
remarked earlier, that the Pad\'e approximants $[n +k / k]$ for the Euler
series can be constructed on the basis of the model sequence (13.3-5).
The remainder $r_n$ of this model sequence is of order $O (z^{n+k}
n^{n+k})$ as $n \to \infty$. However, it follows from theorem 13-2 that
the remainder integral $R_n (z)$ of the Euler series (1.1-7) with $z
\in \R_{+}$ is rigorously bounded by $(n+1)! z^{n+1}$ which is only of
order $O (z^{n+1} n^{n+1})$ as $n \to \infty$. Hence, we see that the
remainder of the model sequence (13.3-3) for Pad\'e approximants $[n + k/
k]$ yields unrealistically large estimates in the case of a wildly
divergent series such as the Euler series. Consequently, it is to be
expected that in the case of such a wildly divergent series Pad\'e
approximants will be less efficient than sequence transformations which
use tighter remainder estimates as for instance suitable variants of
${\cal L}_k^{(n)} (\beta, s_n, \omega_n)$, eq. (7.1-7), ${\cal
S}_k^{(n)} (\beta, s_n, \omega_n)$, eq. (8.2-7), and ${\cal M}_k^{(n)}
(\gamma, s_n, \omega_n)$, eq. (9.2-6).

There is considerable numerical evidence that the inferiority of Pad\'e
approximants in summation processes is not restricted to the Euler
series. For instance, the numerical tests performed by Smith and Ford
[29,30], who also considered the summation of several alternating
divergent series, showed that Levin's $u$ transformation, eq. (7.3-5),
is significantly more powerful than Wynn's $\epsilon$ algorithm.

The inferiority of Pad\'e approximants in summation processes also
becomes quite obvious in the case of the following class of auxiliary
functions,
$$
F_m (z) \; = \; \int\nolimits_{0}^{1} u^{2 m} \e^{- z u^2} \d u
\, , \qquad m \in \N_0 \, , \quad z \in \R_{+} \, .
\tag
$$

These auxiliary functions $F_m (z)$ are of considerable importance in
molecular {\it ab initio} calculations with Gaussian-type basis
functions since the nuclear attraction and interelectronic repulsion
integrals are ultimately expressed in terms of these functions. In
molecular calculations these auxiliary functions have to be computed
over a wide range of parameters $m$ and arguments $z$ so frequently
that it amounts to a significant part of the whole integral evaluation
time. In the case of larger arguments $z$, it is recommendable to
compute this auxiliary function via its asymptotic expansion,
$$
F_m (z) \; \sim \; \frac{\Gamma (m + 1/2)} {2 z^{m + 1/2} }
\, - \, \frac {\e^{- z}} {2 z} \,
{}_2 F_0 (1, 1/2 - m; -1/z) \, , \qquad z \to \infty \, .
\tag
$$

In ref. [63] it was shown that Levin's $d$ transformation, eq. (7.3-9),
sums this divergent series much more efficiently than Wynn's $\epsilon$
algorithm. Later, in ref. [107] the effect of the sequence
transformations $u_k^{(n)} (\beta, s_n)$, eq. (7.3-5), $y_k^{(n)}
(\beta, s_n)$, eq. (8.4-2), and $Y_k^{(n)} (\gamma, s_n)$, eq. (9.4-2),
on the divergent series ${}_2 F_0$ in eq. (13.3-14) was compared.
Similarly as in the case of the divergent series ${}_2 F_0$ in eq.
(13.3-9) it was found that Levin's $u$ transformation is slightly less
efficient than the analogous new sequence transformations $y_k^{(n)}
(\beta, s_n)$ and $Y_k^{(n)} (\gamma, s_n)$.

\medskip

\Abschnitt A Stieltjes series with a finite radius of convergence

\smallskip

\aktTag = 0

Let us consider the following integral representation for the logarithm
which is defined for all $z$ belonging to the cut complex plane which
is cut along $- \infty < z \le - 1$,
$$
\frac {1} {z} \ln (1+z) \; = \; \int\nolimits_{0}^{1}
\frac {\d t}{1 + z t} \, .
\tag
$$

The integral in eq. (13.4-1) is a Stieltjes integral as the one in eq.
(13.1-1). To see this we only have to set $\psi (t) = t$ for $0 \le t
\le 1$ and $\psi (t) = 1$ for $1 < t < \infty$ in eq. (13.1-1). The
moments $\mu_m$ of this positive measure $\psi (t)$ are given by
$$
\int\nolimits_{0}^{\infty} t^m \d \psi (t) \; = \;
\int\nolimits_{0}^{1} t^m \d t \; = \; \frac {1} {m + 1} \, ,
\qquad m \in \N_0 \, .
\tag
$$

If we use these moments $\mu_m$ in eq. (13.1-3), we obtain the
following power series for the logarithm which is by construction a
Stieltjes series:
$$
\ln (1+z) \; = \;
\sum_{m=0}^{\infty} \, \frac {(-1)^m z^{m+1}} {m+1}
\; = \; z \, {}_2 F_1 (1,1;2, -z) \, .
\tag
$$

The power series in eq. (13.4-3) converges absolutely for all $z \in
\C$ with $|z| < 1$, for $z = 1$ the series converges conditionally, and
all for $z \in \C$ with $|z| > 1$ the series diverges. However, as long
as the argument $z \in \C$ does not lie on the cut, the divergent
series can at least in principle be summed.

It may be interesting to note that the infinite series (13.4-3) for
$\ln (2)$ occurs also in solid state physics since it gives the
Madelung constant of a 1-dimensional lattice of oppositely charged ions
(see pp. 74 - 75 of ref. [108]). According to Killingbeck the infinite
series (13.4-3) occurs also if correlation effects in atoms are treated
via perturbation theory (see p. 969 of ref. [109]).

In the last section, it was demonstrated both theoretically and
numerically that the sequence transformations ${\cal L}_k^{(n)} (\beta,
s_n, \omega_n)$, eq. (7.1-7), ${\cal S}_k^{(n)} (\beta, s_n,
\omega_n)$, eq. (8.2-7), and ${\cal M}_k^{(n)} (\gamma, s_n,
\omega_n)$, eq. (9.2-6), sum the wildly divergent series ${}_2 F_0$ in
eq. (13.3-9), which is essentially the Euler series, significantly more
efficiently than Pad\'e approximants. In addition, some arguments were
presented which indicate that this inferiority of Pad\'e approximants is
not restricted to the Euler series (1.1-7) and will occur also in the
case of other wildly divergent series. For $1 < z < \infty$, the
sequence of partial sums of the Stieltjes series in eq. (13.4-3),
$$
s_n \; = \; \sum_{m=0}^n \, \frac {(-1)^m z^{m+1}} {m+1}
\, , \qquad n \in \N_0 \, ,
\tag
$$
obviously diverges but not as wildly as the partial sums of the
divergent series ${}_2 F_0$ in eq. (13.3-9). Consequently, it should be
interesting to investigate whether the striking superiority of the
sequence transformations ${\cal L}_k^{(n)} (\beta, s_n, \omega_n)$,
${\cal S}_k^{(n)} (\beta, s_n, \omega_n)$, and ${\cal M}_k^{(n)}
(\gamma, s_n, \omega_n)$, over Pad\'e approximants is also observed in
the case of the Stieltjes series (13.4-3).

In this context it would of course be helpful to have some theoretical
summation error estimates. In the case of the sequence transformations
${\cal L}_k^{(n)} (\beta, s_n, \omega_n)$, eq. (7.1-7), ${\cal
S}_k^{(n)} (\beta, s_n, \omega_n)$, eq. (8.2-7), and ${\cal M}_k^{(n)}
(\gamma, s_n, \omega_n)$, eq. (9.2-6), this poses no problems. If the
argument $z$ of the Stieltjes series in eq. (13.4-3) is positive, the
error analysis of section 13-2 can be used since the remainders of the
power series are then strictly alternating.

In the case of Pad\'e approximants the error analysis of section 13-2
cannot be applied because only the Pad\'e approximants of the Euler
series, eq. (1.1-7), can be computed via Drummond's sequence
transformation. In Wimp's book [23] the effect of Wynn's $\epsilon$
algorithm, eq. (4.2-2), on the following model sequence is studied. The
elements of this model sequence are defined by Poincar\'e-type asymptotic
expansions in inverse powers of $n$,
$$
s_n \; \sim \; s \, + \, \lambda^n \, n^{\theta}
\sum_{j=0}^{\infty} \, c_j \, / \, n^j \, ,
\qquad c_0 \ne 0 \, , \quad n \to \infty \, .
\tag
$$

A sequence of the type of eq. (13.4-5) should be a reasonably good
model for the behaviour of the partial sums $s_n$ of the series
(13.4-3) as $n \to \infty$. The sequence (13.4-5) obviously converges
linearly if $| \lambda | < 1$ and it diverges if $| \lambda | > 1$.
Assuming $\lambda \ne 1$ and $\theta \ne 0, 1, \ldots , k - 1$ in eq.
(13.4-5), Wimp obtained for fixed $k \in \N$ the following order
estimate (see p. 127 of ref. [23]):
$$
\epsilon_{2 k}^{(n)} \, - \, s \; = \; \frac
{c_0 \, \lambda^{n + 2 k} \, n^{\theta - 2 k}
\, k! \, (- \theta)_k }
{ (\lambda - 1)^{2 k} }
\left[ 1 + O \left( \frac 1 n \right) \right]
\, , \qquad n \to \infty \, .
\tag
$$

Since $s_n -s \sim \lambda^n n^{\theta}$ as $n \to \infty$, we obtain
from eq. (13.4-6) the following order estimate:
$$
\frac {\epsilon_{2 k}^{(n)} \, - \, s} {s_n - s} \; \sim \;
O (n^{- 2 k}) \, , \qquad n \to \infty \, .
\tag
$$

The error estimate (13.4-6) shows quite clearly that for $| \lambda | <
1$ Wynn's $\epsilon$ algorithm accelerates the convergence of the
linearly convergent sequence. It also follows from the error estimate
(13.4-6) that the limit $s$ can be determined more easily if $\lambda$
is negative which is well in agreement with experience.

\beginFloat

\medskip

\beginTabelle [to \kolumnenbreite]
\beginFormat \rechts " \rechts " \mitte " \mitte " \mitte
\endFormat
\+ " \links {\bf Table 13-4} \@ \@ \@ \@ " \\
\+ " \links {Summation of the divergent series $z {}_2 F_1 (1,1;2; - z)
\; = \; \ln (1+z)$ for $z \; = \; 5$}
\@ \@ \@ \@ " \\
\- " \- " \- " \- " \- " \- " \\ \sstrut {} {1.5 \jot} {1.5 \jot}
\+ " \rechts {$n$} " \mitte {partial sum $s_n$}
" $\epsilon_{2 \Ent {n/2}}^{(n - 2 \Ent {n/2})}$
" $t_n^{(0)} (1, s_0)$ " ${\tau}_n^{(0)} (1, s_0)$ " \\
\+ " " \mitte {eq. (13.4-4)} " eq. (4.2-1) " eq. (7.3-7) " eq. (8.4-3)
" \\
\- " \- " \- " \- " \- " \- " \\ \sstrut {} {1 \jot} {1 \jot}
\+ "   10 " $   0.3639603183 \times 10^{07}$  "  1.79198007997771 "
1.79175951159974 "  1.79175959220168 " \\
\+ "   11 " $  -0.1670544890 \times 10^{08}$  "  1.79159768463775 "
1.79175946864530 "  1.79175949178480 " \\
\+ "   12 " $   0.7719479148 \times 10^{08}$  "  1.79179758764032 "
1.79175946794412 "  1.79175947333854 " \\
\+ "   13 " $  -0.3587706103 \times 10^{09}$  "  1.79173348919423 "
1.79175946933559 "  1.79175946997338 " \\
\+ "   14 " $   0.1675734598 \times 10^{10}$  "  1.79176609278102 "
1.79175946926071 "  1.79175946936268 " \\
\+ "   15 " $  -0.7861008566 \times 10^{10}$  "  1.79175520193427 "
1.79175946922241 "  1.79175946925230 " \\
\+ "   16 " $   0.3701778279 \times 10^{11}$  "  1.79176062438322 "
1.79175946922743 "  1.79175946923241 " \\
\+ "   17 " $  -0.1749098431 \times 10^{12}$  "  1.79175875744767 "
1.79175946922828 "  1.79175946922884 " \\
\+ "   18 " $   0.8289578584 \times 10^{12}$  "  1.79175967119854 "
1.79175946922806 "  1.79175946922819 " \\
\+ "   19 " $  -0.3939413724 \times 10^{13}$  "  1.79175934919749 "
1.79175946922805 "  1.79175946922808 " \\
\+ "   20 " $   0.1876711762 \times 10^{14}$  "  1.79175950460547 "
1.79175946922806 "  1.79175946922806 " \\
\+ "   21 " $  -0.8960496379 \times 10^{14}$  "  1.79175944882296 "
1.79175946922806 "  1.79175946922806 " \\
\+ "   22 " $   0.4286962951 \times 10^{15}$  "  1.79175947543322 "
1.79175946922805 "  1.79175946922806 " \\
\+ "   23 " $  -0.2054830571 \times 10^{16}$  "  1.79175946573795 "
1.79175946922805 "  1.79175946922806 " \\
\+ "   24 " $   0.9866098385 \times 10^{16}$  "  1.79175947031756 "
1.79175946922806 "  1.79175946922806 " \\
\+ "   25 " $  -0.4744606005 \times 10^{17}$  "  1.79175946862827 "
1.79175946922806 "  1.79175946922806 " \\
\- " \- " \- " \- " \- " \- " \\ \sstrut {} {1 \jot} {1 \jot}
\+ " \links {FORTRAN function QLOG} \@ "  1.79175946922806 "
1.79175946922806 "  1.79175946922806 " \\
\- " \- " \- " \- " \- " \- " \\ \sstrut {} {1 \jot} {1 \jot}
\endTabelle

\medskip

\endFloat

However, if we compare the order estimate (13.4-7) with the order
estimates (13.2-7) and (13.2-26) for ${\cal S}_k^{(n)} (\beta, s_n,
\omega_n)$, eq. (8.2-7), and ${\cal L}_k^{(n)} (\beta, s_n, \omega_n)$,
eq. (7.1-7), respectively, we find that Wynn's $\epsilon$ algorithm
should be significantly less efficient than the other two sequence
transformations mentioned above. This follows from the fact that for
the computation of $\epsilon_{2 k}^{(n)}$, which according to eq.
(13.4-7) gives an order estimate of order $O (n^{- 2 k})$, $2 k + 1$
sequence elements $s_n, s_{n+1}, \ldots, s_{n+2 k}$ will be needed,
whereas for the computation of ${\cal S}_k^{(n)} (\beta, s_n,
\omega_n)$, and ${\cal L}_k^{(n)} (\beta, s_n, \omega_n)$, which also
give an error estimate of order $O (n^{- 2 k})$, only $k + 1$ sequence
elements $s_n, s_{n+1}, \ldots, s_{n+k}$ will be needed.

In table 13-4 the effect of the sequence transformations $t_k^{(n)}
(\beta, s_n)$, eq. (7.3-7), $\tau_k^{(n)} (\beta, s_n)$, eq. (8.4-3),
with $\beta = 1$ and Wynn's $\epsilon$ algorithm, eq. (4.2-1), on the
partial sums of the divergent Stieltjes series in eq. (13.4-3) with $z
= 5$ are compared.

The results in table 13-4 are another striking example for the
inferiority of Pad\'e approximants in summation processes. The results
also indicate that our conclusions concerning the efficiency of Pad\'e
approximants in summation processes, which were based upon a comparison
of the order estimates in theorems 13-5 and 13-9 and in eq. (13.4-6),
should at least be qualitatively correct.

Table 13-4 was produced in QUADRUPLE PRECISION. When the same
computation was repeated in DOUBLE PRECISION, some loss of accuracy was
observed. The best results were produced by $t_n^{(0)} (1, s_0)$ and
$\tau_n^{(0)} (1, s_0)$ for $n$ between 15 and 18 (approximately 11
decimal digits). For larger values of $n$ the accuracy deteriorates.

\beginFloat

\medskip

\beginTabelle [to \kolumnenbreite]
\beginFormat \rechts " \mitte " \mitte " \mitte " \mitte
\endFormat
\+ " \links {\bf Table 13-5} \@ \@ \@ \@ " \\
\+ " \links {Acceleration of the conditionally convergent series $z
{}_2 F_1 (1,1;2; - z) \; = \; \ln (1+z)$ for $z \; = \; 1$}
\@ \@ \@ \@ " \\
\- " \- " \- " \- " \- " \- " \\ \sstrut {} {1.5 \jot} {1.5 \jot}
\+ " \rechts {$n$} " \mitte {partial sum $s_n$} " $\epsilon_{2 \Ent
{n/2}}^{(n - 2 \Ent {n/2})}$ " $t_n^{(0)} (1, s_0)$
" ${\tau}_n^{(0)} (1, s_0)$ " \\
\+ " " \mitte {eq. (13.4-4)} " eq. (4.2-1) " eq. (7.3-7) " eq. (8.4-3)
" \\
\- " \- " \- " \- " \- " \- " \\ \sstrut {} {1 \jot} {1 \jot}
\+ "    3 "   0.58333333333333  "  0.69047619047619 "  0.69313725490196
"  0.69321533923304 " \\
\+ "    4 "   0.78333333333333  "  0.69333333333333 "  0.69314393939394
"  0.69314971751412 " \\
\+ "    5 "   0.61666666666667  "  0.69308943089431 "  0.69314740192831
"  0.69314726571364 " \\
\+ "    6 "   0.75952380952381  "  0.69315245478036 "  0.69314717779003
"  0.69314718328808 " \\
\+ "    7 "   0.63452380952381  "  0.69314574314574 "  0.69314718001500
"  0.69314718064517 " \\
\+ "    8 "   0.74563492063492  "  0.69314733235438 "  0.69314718060123
"  0.69314718056257 " \\
\+ "    9 "   0.64563492063492  "  0.69314714248772 "  0.69314718055924
"  0.69314718056003 " \\
\+ "   10 "   0.73654401154401  "  0.69314718496213 "  0.69314718055985
"  0.69314718055995 " \\
\+ "   11 "   0.65321067821068  "  0.69314717951778 "  0.69314718055995
"  0.69314718055995 " \\
\+ "   12 "   0.73013375513376  "  0.69314718068816 "  0.69314718055995
"  0.69314718055995 " \\
\+ "   13 "   0.65870518370518  "  0.69314718053085 "  0.69314718055995
"  0.69314718055995 " \\
\+ "   14 "   0.72537185037185  "  0.69314718056369 "  0.69314718055995
"  0.69314718055995 " \\
\+ "   15 "   0.66287185037185  "  0.69314718055912 "  0.69314718055995
"  0.69314718055995 " \\
\+ "   16 "   0.72169537978362  "  0.69314718056005 "  0.69314718055995
"  0.69314718055995 " \\
\+ "   17 "   0.66613982422806  "  0.69314718055992 "  0.69314718055995
"  0.69314718055995 " \\
\+ "   18 "   0.71877140317543  "  0.69314718055995 "  0.69314718055995
"  0.69314718055995 " \\
\- " \- " \- " \- " \- " \- " \\ \sstrut {} {1 \jot} {1 \jot}
\+ " \links {FORTRAN function QLOG} \@ "  0.69314718055995 "
0.69314718055995 "  0.69314718055995 " \\
\- " \- " \- " \- " \- " \- " \\ \sstrut {} {1 \jot} {1 \jot}
\endTabelle

\medskip

\endFloat

The Stieltjes series (13.4-3) is not only suited to test the efficiency
of a sequence transformation in summation processes. A very popular
test case, which is frequently found in the literature, is the
conditionally convergent series (13.4-3) for $\ln (2)$ which converges
quite slowly. According to Bender and Orszag (see p. 372 of ref. [2])
about $7000$ terms of the series in eq. (13.4-3) with $z = 1$ will be
needed to compute $\ln (2)$ with a relative accuracy of $0.01$ percent.
The same sequence transformations as in table 13-4 accelerate the
convergence of the series for $\ln (2)$ also in table 13-5. This time,
the $\epsilon$ algorithm is comparatively successful since it only
needs the partial sums (13.3-4) up to $n = 18$ to produce an accuracy
of 14 decimal digits. However, the other two sequence transformations
in table 13-5 are still significantly more powerful.

Table 13-5 was again produced in QUADRUPLE PRECISION. When this
computation was repeated in DOUBLE PRECISION, no loss of accuracy was
observed.

The potential of the Stieltjes series (13.4-3) for $\ln (1+z)$ to test
the performance of sequence transformations is not yet exhausted. If
the argument $z$ of the power series satisfies $- 1 < z < 0$, all of
its terms have the same sign. The convergence of this series will
become quite bad if $z$ approaches $- 1$ because for $z = - 1$ it
becomes the series (1.1-2) for $\zeta (1)$ which diverges.
Consequently, it should be interesting to find out whether and how well
the convergence of the Stieltjes series (13.4-3) can be accelerated if
its argument $z$ is close to $- 1$.

In this context, it would again be helpful to have some theoretical
error estimates. In the case of Pad\'e approximants this poses no
problems. If we assume that the elements of the model sequence (13.4-5)
is still a good model for the behaviour of the partial sums (13.4-4) as
$n \to \infty$, we may conclude from eq. (13.4-6) that Wynn's
$\epsilon$ algorithm will accelerate the convergence of the sequence
(13.4-5).

In the case of the sequence transformations ${\cal L}_k^{(n)} (\beta,
s_n, \omega_n)$, eq. (7.1-7), ${\cal S}_k^{(n)} (\beta, s_n,
\omega_n)$, eq. (8.2-7), and ${\cal M}_k^{(n)} (\gamma, s_n,
\omega_n)$, eq. (9.2-6), the situation is more complicated since the
error analysis of section 13-2, which rests upon the assumption that
the remainder estimates strictly alternate in sign, cannot be applied
here. Consequently, we have to find estimates of the type of eq.
(13.4-6) for the other sequence transformations mentioned above.

\medskip

\beginEinzug \sl \parindent = 0 pt

\Auszug {\bf Theorem 13-11:} Let us assume that the elements of the
sequence $\Seqn s$ satisfy
$$
s_n \; = \; s \, + \, \lambda^n \, n^{\theta} \,
[ c_0 + O (n^{- 1}) ] \, , \qquad c_0 \ne 0 \, ,
\quad \lambda \ne 0, 1 \, , \quad n \to \infty \, ,
\tag
$$
that the elements of the sequence of remainder estimates $\Seqn \omega$
can be chosen in such a way that
$$
\omega_n \; = \; \lambda^n \, n^{\theta} \,
[ d_0 + O (n^{- 1}) ] \, ,
\qquad d_0 \ne 0 \, , \quad n \to \infty \, ,
\tag
$$
and that the ratio $(s_n - s) / \omega_n$ can for all $n \in \N_0$ be
expressed as a factorial series,
$$
\frac {s_n - s} {\omega_n} \; = \;
\sum_{j=0}^{\infty} \; \frac {\gamma_j} {(\beta+n)_j} \, ,
\qquad \beta \in \R_{+} \, .
\tag
$$

If the sequence transformation ${\cal S}_k^{(n)} (\beta, s_n,
\omega_n)$, eq. (8.2-7), is used for the transformation of $\Seqn s$,
we obtain for fixed $k \in \N$ the following order estimate:
$$
\frac
{{\cal S}_k^{(n)} (\beta, s_n, \omega_n) \, - \, s} {s_n - s}
\; = \; O (n^{- 2 k}) \, , \qquad n \to \infty \, .
\tag
$$

\endEinzug

\medskip

\noindent {\it Proof: } We can proceed as in theorem 13-5, i.e., the
starting point for the proof of theorem 13-11 is eq. (13.2-8). Since
eqs. (13.2-4) and (13.4-10) are identical we find that the numerator of
this expression is also given by eq. (13.2-13) which is obviously of
order $O (n^{- k - 1})$ as $n \to \infty$.

In order to obtain an estimate for the denominator $\Delta^k
[(\beta+n)_{k-1} / \omega_n ]$ we use (see eq. (41) on p. 21 of ref.
[23])
$$
\Delta^k [ z^n n^{\alpha} ] \; \sim \; z^n (z - 1)^k n^{\alpha}
\, , \qquad z \ne 1 \, , \quad n \to \infty \, .
\tag
$$

This relationship gives us the following asymptotic estimate for the
denominator in eq. (13.2-8):
$$
\Delta^k [(\beta+n)_{k-1} / \omega_n ] \; \sim \;
[1 - \lambda]^k \, \lambda^{- n - k} \, n^{k - \theta - 1} \, ,
\qquad n \to \infty \, .
\tag
$$

If we combine the expressions for the numerator and the denominator and
take into account that $s_n - s \sim \lambda^n n^{\theta}$ as $n \to
\infty$, we obtain eq. (13.4-11).

\medskip

In the next theorem, which can be proved in essentially the same way as
theorem 13-11, it is shown that Levin's sequence transformation ${\cal
L}_k^{(n)} (\beta, s_n, \omega_n)$, eq. (7.1-7), also leads to an error
estimate of order $O (n^{- 2 k})$.

\medskip

\beginEinzug \sl \parindent = 0 pt

\Auszug {\bf Theorem 13-12:} Let us assume that the elements of the
sequence $\Seqn s$ satisfy
$$
s_n \; = \; s \, + \, \lambda^n \, n^{\theta} \,
[ c_0 + O (n^{- 1}) ] \, , \qquad c_0 \ne 0 \, ,
\quad \lambda \ne 0, 1 \, , \quad n \to \infty \, ,
\tag
$$
that the elements of the sequence of remainder estimates $\Seqn \omega$
can be chosen in such a way that
$$
\omega_n \; = \; \lambda^n \, n^{\theta} \,
[ d_0 + O (n^{- 1}) ] \, ,
\qquad d_0 \ne 0 \, , \quad n \to \infty \, ,
\tag
$$
and that the ratio $(s_n - s)/\omega_n$ can for all $n \in \N_0$ be
expressed as a power series of the following type,
$$
\frac {s_n - s} {\omega_n} \; = \;
\sum_{j=0}^{\infty} \; \frac {\gamma'_j} {(\beta+n)^j} \, ,
\qquad \beta \in \R_{+} \, .
\tag
$$

If the sequence transformation ${\cal L}_k^{(n)} (\beta, s_n,
\omega_n)$, eq. (7.1-7), is used for the transformation of $\Seqn s$,
we obtain for fixed $k \in \N$ the following order estimate:
$$
\frac
{{\cal L}_k^{(n)} (\beta, s_n, \omega_n) \, - \, s} {s_n - s}
\; = \; O (n^{- 2 k}) \, , \qquad n \to \infty \, .
\tag
$$

\endEinzug

\medskip

The error estimate (13.4-17) for Levin's sequence transformation has in
principle already been derived by Sidi (see eq. (3.14) on p. 840 of
ref. [105]).

In the same way, it can be proved that for fixed $k \in \N$ and for
large values of $n$ the sequence transformations ${\cal L}_{k,
\ell}^{(n)} (\beta, s_n, \omega_n)$, eq. (7.1-8), and ${\cal S}_{k,
\ell}^{(n)} (\beta, s_n, \omega_n)$, eq. (8.2-8), with $\ell \in \N$
lead to error estimates of order $O (n^{\ell - 2 k})$ and that
Drummond's sequence transformation ${\cal D}_k^{(n)} (s_n, \omega_n)$,
eq. (9.5-4), leads to an error estimate of order $O (n^{- k - 1})$.

On the basis of these order estimates it is to be expected that in the
case of negative arguments $z$ the different variants of ${\cal
L}_k^{(n)} (\beta, s_n, \omega_n)$, eq. (7.1-7), and ${\cal S}_k^{(n)}
(\beta, s_n, \omega_n)$, eq. (8.2-7), should accelerate the convergence
of the Stieltjes series in eq. (13.4-3) more efficiently than the
analogous variants of the generalized transformations ${\cal L}_{k,
\ell}^{(n)} (\beta, s_n, \omega_n)$, eq. (7.1-8), and ${\cal S}_{k,
\ell}^{(n)} (\beta, s_n, \omega_n)$, eq. (8.2-8), with $\ell \ge 1$.
Also, Drummond's sequence transformation, eq. (9.5-4), should be
roughly as efficient as Wynn's $\epsilon$ algorithm, eq. (4.2-1).

In table 13-6 the convergence of the absolutely convergent Stieltjes
series (13.4-3) with $z = - 0.9$ is accelerated by the same sequence
transformations as in tables 13-4 and 13-5. The inferiority of Wynn's
$\epsilon$ algorithm is again evident.

\beginFloat

\medskip

\beginTabelle [to \kolumnenbreite]
\beginFormat \rechts " \mitte " \mitte " \mitte " \mitte
\endFormat
\+ " \links {\bf Table 13-6} \@ \@ \@ \@ " \\
\+ " \links {Acceleration of the absolutely convergent series $z {}_2
F_1 (1,1;2; - z) \; = \; \ln (1+z)$ for $z \; = \; - 0.9$}
\@ \@ \@ \@ " \\
\- " \- " \- " \- " \- " \- " \\ \sstrut {} {1.5 \jot} {1.5 \jot}
\+ " \rechts {$n$} " \mitte {partial sum $s_n$} " $\epsilon_{2 \Ent
{n/2}}^{(n - 2 \Ent {n/2})}$ " $t_n^{(0)} (1, s_0)$
" ${\tau}_n^{(0)} (1, s_0)$ " \\
\+ " " \mitte {eq. (13.4-4)} " eq. (4.2-1) " eq. (7.3-7) " eq. (8.4-3)
" \\
\- " \- " \- " \- " \- " \- " \\ \sstrut {} {1 \jot} {1 \jot}
\+ "   15 " $-2.23201245730299$ " $-2.30249119271252$ "
$-2.30258308878949$ " $-2.30258507564758$ " \\
\+ "   16 " $-2.24182256418514$ " $-2.30253886165435$ "
$-2.30258429564850$ " $-2.30258508829481$ " \\
\+ "   17 " $-2.25016115503498$ " $-2.30255980521704$ "
$-2.30258477584394$ " $-2.30258509172157$ " \\
\+ "   18 " $-2.25727090091747$ " $-2.30257257511841$ "
$-2.30258496686622$ " $-2.30258509264961$ " \\
\+ "   19 " $-2.26334973364700$ " $-2.30257828055989$ "
$-2.30258504284117$ " $-2.30258509290084$ " \\
\+ "   20 " $-2.26856016170088$ " $-2.30258170574611$ "
$-2.30258507305389$ " $-2.30258509296883$ " \\
\+ "   21 " $-2.27303639307444$ " $-2.30258325724815$ "
$-2.30258508506691$ " $-2.30258509298723$ " \\
\+ "   22 " $-2.27688984443081$ " $-2.30258417685642$ "
$-2.30258508984292$ " $-2.30258509299220$ " \\
\+ "   23 " $-2.28021344622568$ " $-2.30258459820761$ "
$-2.30258509174154$ " $-2.30258509299355$ " \\
\+ "   24 " $-2.28308503817645$ " $-2.30258484529651$ "
$-2.30258509249623$ " $-2.30258509299391$ " \\
\+ "   25 " $-2.28557006967230$ " $-2.30258495961133$ "
$-2.30258509279620$ " $-2.30258509299401$ " \\
\+ "   26 " $-2.28772376363538$ " $-2.30258502604192$ "
$-2.30258509291542$ " $-2.30258509299404$ " \\
\+ "   27 " $-2.28959286232476$ " $-2.30258505703207$ "
$-2.30258509296280$ " $-2.30258509299404$ " \\
\+ "   28 " $-2.29121704463416$ " $-2.30258507490093$ "
$-2.30258509298163$ " $-2.30258509299404$ " \\
\+ "   29 " $-2.29263008324333$ " $-2.30258508329701$ "
$-2.30258509298911$ " $-2.30258509299405$ " \\
\+ "   30 " $-2.29386079429003$ " $-2.30258508810542$ "
$-2.30258509299209$ " $-2.30258509299405$ " \\
\+ " \- " \- " \- " \- " \- " \\ \sstrut {} {1 \jot} {1 \jot}
\+ " \links {FORTRAN function QLOG} \@ " $-2.30258509299405$ "
$-2.30258509299405$ " $-2.30258509299405$ " \\
\- " \- " \- " \- " \- " \- " \\ \sstrut {} {1 \jot} {1 \jot}

\endTabelle

\medskip

\endFloat

Table 13-6 was produced in QUADRUPLE PRECISION. When the same
computation was repeated in DOUBLE PRECISION, a heavy loss of accuracy
was observed. No transformation was able to produce an accuracy of more
than 8 decimal digits. These 8 digits were produced by $\tau_n^{(0)}
(1, s_0)$ for $n = 17$, by $t_n^{(0)} (1, s_0)$ for $n = 21$, and by
Wynn's $\epsilon$ algorithm for $n = 30$. For larger values of $n$ the
accuracy deteriorated again. However, it seems that Wynn's $\epsilon$
algorithm is not as much affected by numerical instabilities as the
other two sequence transformations in table 13-6.

More extensive numerical tests showed that in the case of the
convergent or divergent Stieltjes series in eq. (13.4-3) ${\cal
S}_k^{(n)} (\beta, s_n, \omega_n)$, eq. (8.2-7), had a slight plus over
${\cal L}_k^{(n)} (\beta, s_n, \omega_n)$, eq. (7.1-7). A comparison of
the numerous variants of the sequence transformations ${\cal L}_k^{(n)}
(\beta, s_n, \omega_n)$ and ${\cal S}_k^{(n)} (\beta, s_n, \omega_n)$
showed that they differ in their ability of summing or accelerating the
Stieltjes series in eq. (13.4-3). Those variants, which are based upon
the remainder estimates (7.3-4) and (7.3-8), are approximately as
efficient as $t_k^{(n)} (\beta, s_n)$, eq. (7.3-7), and $\tau_k^{(n)}
(\beta, s_n)$, eq. (7.3-7). However, $v_k^{(n)} (\beta, s_n)$, eq.
(7.3-11), and $\phi_k^{(n)} (\beta, s_n)$, eq. (8.4-5), which are both
based upon the remainder estimate (7.3-10), were significantly less
efficient. Also, ${\cal M}_k^{(n)} (\gamma, s_n, \omega_n)$, eq.
(9.2-6), and its variants were somewhat less efficient than the
analogous variants of ${\cal L}_k^{(n)} (\beta, s_n, \omega_n)$ and
${\cal S}_k^{(n)} (\beta, s_n, \omega_n)$.

Almost as efficient as the sequence transformations mentioned above
were ${\cal J}_k^{(n)}$, eq. (10.3-6), and ${\cal C}_k^{(n)}$, eq.
(11.1-12), followed by Brezinski's $\theta$ algorithm, eq. (10.1-9).
Next came Aitken's iterated $\Delta^2$ process, eq. (5.1-15), which was
again more efficient than Wynn's $\epsilon$ algorithm, eq. (4.2-1). The
transformations ${\cal B}_k^{(n)}$, eq. (11.1-5), $\lambda_k^{(n)}$,
eq. (11.2-1), $\sigma_k^{(n)}$, eq. (11.2-2), and $\mu_k^{(n)}$, eq.
(11.2-3), were weaker than the $\epsilon$ algorithm.

Of all transformations tested Drummond's sequence transformation, eq.
(9.5-4), was least efficient. It was able to sum or accelerate the
Stieltjes series (13.4-3) moderately well for positive arguments $z$,
i.e., as long as the terms of the power series for $\ln (1+z)$ had
alternating signs. However, it failed completely if the argument $z$
approached $- 1$. For instance, for $z = - 0.9$ Drummond's sequence
transformation produced a sequence of transforms ${\cal D}_n^{(0)}
(s_0, a_1)$ which rapidly diverged with increasing $n$. This example
shows that asymptotic order estimates -- although undeniably quite
helpful for the classification of sequence transformations -- do not
necessarily tell the whole truth about the capability of a sequence
transformation.

\endAbschnittsebene

\neueSeite

\Abschnitt The acceleration of logarithmic convergence

\vskip - 2 \jot

\beginAbschnittsebene

\medskip

\Abschnitt Properties of logarithmically convergent sequences and
series

\smallskip

\aktTag = 0

It is tempting to believe that the most formidable task for a nonlinear
sequence transformation is the summation of a wildly divergent series
such as the Euler series, eq. (1.1-7), and that convergence
acceleration should not be overly troublesome. In the case of
alternating series or sequences with strictly alternating remainders,
this is indeed normally true. However, it will become clear later in
this section that the acceleration of the convergence of a monotonic
sequence or a series with terms, that all have the same sign, can be a
more formidable computational problem than the summation of an
alternating divergent series.

The numerical examples presented in sections 13.3 and 13.4 showed that
several sequence transformations are able to sum efficiently divergent
series with alternating terms. In addition, it was shown that it is
frequently possible to sum alternating divergent series with an
accuracy that is close to or identical with machine accuracy.
Particularly efficient and also remarkably reliable were variants of
Levin's sequence transformation ${\cal L}_k^{(n)}(\beta, s_n,
\omega_n)$, eq. (7.1-7), and partly even more so variants of the new
sequence transformations ${\cal S}_k^{(n)}(\beta, s_n, \omega_n)$, eq.
(8.2-7), and ${\cal M}_k^{(n)}(\gamma, s_n, \omega_n)$, eq. (9.2-6).

The situation is much less satisfactory if logarithmic convergence has
to be accelerated. Many series with positive terms are known which
converge so slowly that the evaluation of such a series by successively
adding up the terms would be hopeless. In such a case, the use of a
convergence acceleration method is indispensable. However, there is a
considerable amount of theoretical and numerical evidence which
indicates that convergence acceleration methods are generally less
efficient and also more susceptible to rounding errors in the case of
series with positive terms than in the case of alternating series.
Consequently, it is often easier to sum an alternating divergent
series, even if it diverges quite wildly, than to accelerate the
convergence of a slowly convergent series with terms that all have the
same sign.

A good example of a very slowly convergent series with positive terms
is the series (1.1-2) for the Riemann zeta function. It is well known
that this series converges for all $z \in \C$ with $Re (z) > 1$.
However, it follows from eqs. (7.3-12) and (7.3-14) that the remainder
$r_n$ of the series (1.1-2), which is defined by
$$
r_n \; = \; \sum_{m=n+1}^{\infty} \, (m + 1)^{- z} \, ,
\tag
$$
is of order $O (n^{1-z})$ as $n \to \infty$. Consequently, the
computation of $\zeta (z)$ with the help of the series (1.1-2) would
only be feasible if $Re (z)$ is relatively large. But even then, the
use of convergence acceleration techniques would be recommendable.

In order to make this section, in which the acceleration of logarithmic
convergence by means of nonlinear sequence transformations will be
treated, more selfcontained, first some properties of logarithmically
convergent sequences and series are reviewed. A sequence $\Seqn s$,
which converges to some limit $s$, is said to converge logarithmically
if
$$
\lim_{n \to \infty} \frac {s_{n+1} - s} {s_n - s} \; = \; 1 \, .
\tag
$$

This definition of logarithmic convergence is inconvenient since it
involves the limit $s$ of the sequence $\Seqn s$ which is normally not
known. Thus, it would be advantageous to have an alternative criterion
for logarithmic convergence which only involves the differences $\Delta
s_n$. Such an equivalent criterion can be formulated if it is assumed
that the elements of the sequence $\Seqn s$ are partial sums of an
infinite series with real terms $a_m$ that all have the same sign.
Then, Clark, Gray, and Adams could show (see theorem 2 on p. 26 of ref.
[35]) that the sequence of partial sums converges logarithmically
according to eq. (14.1-2) if the following condition holds:
$$
\lim_{n \to \infty} \frac {a_{n+1}} {a_n} \; = \;
\lim_{n \to \infty} \frac {\Delta s_n} {\Delta s_{n-1}}
\; = \; 1 \; .
\tag
$$

Eq. (14.1-3) implies that for larger indices $n$ the terms $a_n$ of a
logarithmically convergent series differ only slightly. This fact is
not only responsible for the often prohibitively slow convergence of
logarithmically convergent series but also affects convergence
acceleration processes in a very unpleasant way.

A sequence transformation can only beat the conventional process of
successively adding up the terms of a series if it does not only use
the numerical values of the terms. It also has to extract from the
terms of the series some additional information about the behaviour of
the partial sums $s_n$ as $n$ increases.

All sequence transformations of this report retrieve and utilize this
additional information by computing rational expressions of weighted
differences of partial sums. In the case of alternating series the
computation of these weighted differences normally does not lead to a
serious loss of significant digits. If, however, the terms $a_n$ of a
series all have the same sign and do not differ much in magnitude, this
additional information, which has to be retrieved, is hidden somewhere
in the later digits. Consequently, in the case of logarithmic
convergence the computation of weighted differences is likely to lead
to a cancellation of significant digits and ultimately, i.e., in the
case of large transformation orders, completely nonsensical results are
to be expected. This explains why rounding errors are more or less
inevitable if logarithmic convergence is accelerated and why the
acceleration or summation of alternating series is frequently
remarkably stable. A good discussion of these stability problems can be
found in an article by Longman [58].

A theoretical analysis of the acceleration of logarithmic convergence
is also far from simple. For instance, it would certainly be quite
helpful if an analogue of Germain-Bonne's theory of the acceleration of
linear convergence [33] could be developed because then a decision
based on some theoretical criteria could be made whether a given
sequence transformation is able to accelerate logarithmic convergence
or not. In the case of linear convergence this question can be decided
on the basis of theorems 12-4 and 12-14. It only has to be shown that
the sequence transformation under consideration is exact for a sequence
which apart from a shift of indices consists of the partial sums of the
geometric series, and it is guaranteed that linear convergence will be
accelerated.

Smith and Ford had speculated whether it might be possible to develop
an analogue of Germain-Bonne's theory also in the case of logarithmic
convergence, i.e., whether some special sequence could be found such
that the exactness of a sequence transformation for this sequence would
imply that all logarithmically convergent sequences will be
accelerated. They also had presented some potential candidates which in
their opinion might possibly be suited to serve as this special
sequence (see p. 238 of ref. [29]). In the meantime, this question has
been answered by Delahaye and Germain-Bonne [110], but unfortunately
the answer is negative. Delahaye and Germain-Bonne [110] showed that no
algorithm exists which would be able to accelerate the convergence of
every logarithmically convergent sequence. Consequently, a general
theory in the spirit of Germain-Bonne's theory, which would cover the
acceleration of all logarithmically convergent sequences, cannot exist.
Such an analogue of Germain-Bonn's theory can exist at most for
suitably restricted subsets of the set of logarithmically convergent
sequences. But it seems that even this has not yet been accomplished so
far.

\medskip

\Abschnitt Exactness results and error estimates

\smallskip

\aktTag = 0

As mentioned in the last section, a theoretical analysis of the
acceleration of logarithmic convergence is far from simple.
Particularly hard is the analysis of the acceleration properties of
those sequence transformations which are defined by a complicated
nonlinear recursive scheme as for instance Brezinski's $\theta$
algorithm, eq. (10.1-9). In those cases, apart from the defining
recursive scheme only very little else is normally known. However, at
least for Levin's sequence transformation ${\cal L}_k^{(n)}(\beta, s_n,
\omega_n)$, eq. (7.1-7), some exactness results and asymptotic error
estimates can be derived quite easily if suitable model sequences are
considered.

Similar exactness results and error estimates as for Levin's sequence
transformation can also be derived for the sequence transformations
${\cal S}_k^{(n)}(\beta, s_n, \omega_n)$, eq. (8.2-7), and ${\cal
M}_k^{(n)}(\gamma, s_n, \omega_n)$, eq. (9.2-6). However, they will not
be considered here. The reason is that numerical tests showed that
those variants of these sequence transformations, which should be able
to accelerate logarithmic convergence, are significantly less efficient
than the analogous variants of Levin's sequence transformation. For the
moment, no completely satisfactory explanation can be given why these
otherwise very powerful sequence transformations perform so weakly in
the case of logarithmic convergence. It is at least conceivable that in
the case of logarithmically convergent sequences inverse powers of $n$
are better suited for a description of the $n$-dependence of the ratios
$(s_n - s)/\omega_n$ than Pochhammer symbols, which are the basis of
the sequence transformations ${\cal S}_k^{(n)}(\beta, s_n, \omega_n)$,
eq. (8.2-7), and ${\cal M}_k^{(n)}(\gamma, s_n, \omega_n)$, eq.
(9.2-6). However, this is only speculation.

The explicit expressions of Levin's $u$ transformation, eq. (7.3-5),
and $t$ transformation, eq. (7.3-7), are very similar. Also, with
respect to the acceleration of linear convergence or the summation of
alternating divergent series these two sequence transformations have
virtually identical properties. However, in the literature on
convergence acceleration it is always emphasized that Levin's $u$
transformation is one of the best accelerators of logarithmic
convergence, whereas Levin's $t$ transformation completely fails to
accelerate logarithmic convergence. In view of the otherwise close
similarity of these two sequence transformations this different
behaviour with respect to the acceleration of logarithmic convergence
is certainly quite puzzling. It will now be shown that the different
properties of Levin's $u$ and $t$ transformation can be understood on
the basis of the different exactness properties of these two sequence
transformations.

Quite common in practical applications are logarithmically convergent
sequences $\Seqn s$ with remainders $\Seqn r$ that are of order $O
(n^{- \alpha})$ with $\alpha \in \R_{+}$ as $n \to \infty$. Hence, for
large values of $n$ the elements of such a logarithmically convergent
sequence can be characterized in the following way: $$
s_n \; = \; s \, + \, n^{- \alpha} \, [ c + O (n^{-1})] \, ,
\qquad c \ne 0 \, , \quad \alpha \in \R_{+} \, ,
\quad n \to \infty \, .
\tag
$$

The explicit expressions for $u_k^{(n)} (\beta, s_n)$, eq. (7.3-5), and
$t_k^{(n)} (\beta, s_n)$, eq. (7.3-7), contain the terms $a_n$ of the
series which is to be transformed. If these two transformations are to
be applied to sequences of the type of eq. (14.2-1), the terms $a_n$ in
the explicit expressions have to be replaced by the differences $\Delta
s_{n-1}$. If we compute these differences and apply some simplifying
assumptions, which are permitted if $n$ is large, we see that sequences
of the type of eq. (14.2-1) satisfy:
$$
(s_n - s) \, / \, \Delta s_{n-1} \; = \; O (n) \, ,
\qquad n \to \infty \, .
\tag
$$

This relationship is quite typical of logarithmically convergent
sequences of the type of eq. (14.2-1). Also, eq. (14.2-2) is
essentially identical with the remainder estimate (7.3-4), which is the
basis of Levin's $u$ transformation.

In the following theorem the exactness properties of Levin's $u$ and
$t$ transformation for a special class of logarithmically convergent
model sequences are analyzed. Model sequences belonging to this special
class have the same behaviour as $n \to \infty$ as the dominant term of
sequences of the type of eq. (14.2-1). Consequently, the following
theorem makes it plausible why Levin's $u$ transformation accelerates
the convergence of sequences of the type of eq. (14.2-1), and why
Levin's $t$ transformation fails to accelerate the convergence of these
sequences.

\medskip

\beginEinzug \sl \parindent = 0 pt

\Auszug {\bf Theorem 14-1:} Assume that a sequence transformation
${\cal T}_k^{(n)} (s_n)$ is defined in the following way:
$$
{\cal T}_k^{(n)} (s_n) \; = \; \frac
{ \Delta^k \> [ P_{k-1} (n) \, s_n / \Delta s_{n-1} ] }
{ \Delta^k \> [ P_{k-1} (n) \, / \Delta s_{n-1} ] }
\, , \qquad k,n \in \N_0 \, .
\tag
$$

$P_{k-1} (n)$ is for sufficiently large values of $k$ a polynomial of
degree $\le k - 1$ in $n$. Obviously, ${\cal T}_k^{(n)} (s_n)$ is
defined as long as $\Delta^k$ does not annihilate $P_{k-1} (n) / \Delta
s_{n-1}$, i.e., as long as $P_{k-1} (n) / \Delta s_{n-1}$ is not a
polynomial of degree $\le k - 1$ in $n$.

Let us assume that the sequence $\Seqn s$, which converges to some
limit $s$, belongs to the domain of the sequence transformation ${\cal
T}_k^{(n)}$ and that its elements satisfy for all $n \in \N_0$
$$
(s_n - s) \, / \, \Delta s_{n-1} \; = \; \gamma n \, + \, \delta \, ,
\qquad \gamma, \delta \in \R \, , \quad \gamma \neq 0 \, .
\tag
$$

If for sufficiently large values of $k$ the degree of $P_{k-1} (n)$ is
exactly $k - 1$, ${\cal T}_k^{(n)} (s_n)$ does not accelerate the
convergence of $\Seqn s$, and if $P_{k-1} (n)$ is a polynomial of
degree $\le k - 2$ in $n$, ${\cal T}_k^{(n)} (s_n)$ is exact for $\Seqn
s$.

\endEinzug

\medskip

\noindent {\it Proof: } Since ${\cal T}_k^{(n)} (s_n)$ is obviously
invariant under translation in the sense of eq. (3.1-4), we can write
$$
{\cal T}_k^{(n)} (s_n) \; = \; s \, + \, \frac
{ \Delta^k \> [ P_{k-1} (n) \, (s_n - s) / \Delta s_{n-1} ] }
{ \Delta^k \> [ P_{k-1} (n) \, / \Delta s_{n-1} ] }
\, , \qquad k,n \in \N_0 \, .
\tag
$$

Next, the ratio $(s_n - s) / \Delta s_{n-1}$ in eq. (14.2-5) is
replaced by $\gamma n + \delta$ according to eq. (14.2-4). This yields
$$
{\cal T}_k^{(n)} (s_n) \; = \; s \, + \, \frac
{ \Delta^k \> [ P_{k-1} (n) \, (\gamma n + \delta) ] }
{ \Delta^k \> [ P_{k-1} (n) \, / \Delta s_{n-1} ] }
\, , \qquad k,n \in \N_0 \, , \quad \gamma \neq 0 \, .
\tag
$$

Let us now assume that $k$ is large enough such that $P_{k-1} (n)$ is a
polynomial of degree $\le k - 1$ in $n$. If the degree of $P_{k-1} (n)$
in eq. (14.2-6) is exactly $k-1$, the product $ (\gamma n + \delta)
P_{k-1} (n)$ is a polynomial of degree $k$ in $n$. Consequently, this
product will not be annihilated by $\Delta^k$ and ${\cal T}_k^{(n)}
(s_n)$ will not accelerate $\Seqn s$. If, however, $P_{k-1} (n)$ is of
degree $\le k - 2$ in $n$, the product $ (\gamma n + \delta) P_{k-1}
(n)$ is a polynomial of degree $\le k - 1$ in $n$. Consequently, this
product will be annihilated by $\Delta^k$ and ${\cal T}_k^{(n)} (s_n)$
is exact for $\Seqn s$. This proves theorem 14-1.

\medskip

In the case of Levin's $u$ transformation, eq. (7.3-5), the polynomial
$P_{k-1} (n)$ is given by $(\beta + n)^{k-2}$, i.e., it is a polynomial
of degree $k - 2$ in $n$. Consequently, $u_k^{(n)} (\beta, s_n)$ will
be exact for sequences $\Seqn s$ satisfying eq. (14.2-4). In the case
of Levin's $t$ transformation, eq. (7.3-7), the polynomial $P_{k-1}
(n)$ is given by $(\beta + n)^{k-1}$, i.e., it is a polynomial of
degree $k-1$ in $n$. This implies that $t_k^{(n)} (\beta, s_n)$ will
not accelerate the convergence of a sequence satisfying eq. (14.2-4).
It also follows from theorem 14-1 that for sufficiently large values of
$k$ Levin's generalized sequence transformation ${\cal L}_{k,
\ell}^{(n)}(\beta, s_n, \omega_n)$, eq. (7.1-8), with $\ell \ge 2$ and
$\omega_n = \Delta s_{n-1}$ is also exact for every sequence $\Seqn s$
satisfying eq. (14.2-4).

A simple example of a logarithmically convergent sequence, which for
all $n \in \N_0$ satisfies eq. (14.2-4), would be
$$
s_n \; = \; s \, + \, \frac {( a )_{n+1}} {( b )_{n+1}} \, ,
\qquad a, b \in \R_{+} \, , \quad a < b \, .
\tag
$$

In sections 13.2 and 13.4 error estimates for the summation of
divergent Stieltjes series and the acceleration of the convergence of
Stieltjes series by means of sequence transformations as for instance
${\cal L}_k^{(n)}(\beta, s_n, \omega_n)$, eq. (7.1-7), or ${\cal
S}_k^{(n)}(\beta, s_n, \omega_n)$, eq. (8.2-7), were derived. It could
be shown that the application of these sequence transformations to the
partial sums of convergent or divergent Stieltjes series lead to
asymptotic error estimates which were of order $O (n^{- 2 k})$ as $n
\to \infty$. In the next theorem a similar asymptotic error analysis is
done for a large class of logarithmically convergent sequences.

\medskip

\beginEinzug \sl \parindent = 0 pt

\Auszug {\bf Theorem 14-2:} Let us assume that the elements of the
sequence $\Seqn s$, which converges logarithmically to some limit $s$,
satisfy
$$
s_n \; = \; s \, + \, n^{- \alpha} \, [ b_0 + O (n^{- 1}) ]
\, , \qquad b_0 \ne 0 \, , \quad \alpha \in \R_{+}
\, , \quad n \to \infty \, .
\tag
$$

Let us also assume that the elements of a sequence of remainder
estimates $\Seqn \omega$ can be chosen in such a way that
$$
\omega_n \; = \; n^{- \alpha} \,
[ d_0 + O (n^{- 1}) ] \, ,
\qquad d_0 \ne 0 \, , \quad n \to \infty \, ,
\tag
$$
and that the ratio $(s_n - s)/(\omega_n)$ can for all $n \in \N_0$ be
expanded in a power series of the following type,
$$
\frac {s_n - s} {\omega_n} \; = \;
\sum_{j=0}^{\infty} \; \frac {c_j} {(\beta+n)^j} \, ,
\qquad \beta \in \R_{+} \, .
\tag
$$

If the sequence transformation ${\cal L}_k^{(n)} (\beta, s_n,
\omega_n)$, eq. (7.1-7), is used for the acceleration of the
convergence of $\Seqn s$, we obtain for fixed $k \in \N$ and for $n \to
\infty$ the following order estimate:
$$
\frac
{{\cal L}_k^{(n)} (\beta, s_n, \omega_n) \, - \, s} {s_n - s}
\; = \; O (n^{- k}) \, , \qquad n \to \infty \, .
\tag
$$

\endEinzug

\medskip

\noindent {\it Proof: } We can proceed as in theorem 13-9, i.e., the
starting point for the proof of theorem 14-2 is the representation of
the transformation error ${\cal L}_k^{(n)} (\beta, s_n - s, \omega_n)$
as in eq. (13.2-27). Since the series expansions (13.2-4) and
(14.2-10) are structurally identical, we can conclude that according to
eq. (13.2-31) the numerator of the transformation error is also of
order $O (n^{- k - 1})$ as $n \to \infty$.

In order to obtain an estimate for the denominator $\Delta^k [(\beta +
n)^{k-1} / \omega_n]$ we take into account that according to eq.
(14.2-9) the remainder estimate $\omega_n$ is of order $O (n^{-
\alpha})$ as $n \to \infty$. This implies that $(\beta + n)^{k-1} /
\omega_n$ is of order $O (n^{k + \alpha - 1})$ as $n \to \infty$.
Hence, with the help of eq. (13.2-30) we obtain the following order
estimate for the denominator of the transformation error:
$$
\Delta^k \, [(\beta + n)^{k-1} / \omega_n ] \; = \;
O (n^{\alpha - 1}) \, , \qquad n \to \infty \, .
\tag
$$

If we combine this relationship with the expression for the numerator
of the transformation error according to eq. (13.2-21), which is of
order $O (n^{- k - 1})$ as $n \to \infty$, we find that the
transformation error ${\cal L}_k^{(n)} (\beta, s_n - s, \omega_n)$ is
of order $O (n^{- \alpha - k})$ as $n \to \infty$. If we next divide
the transformation error by $s_n - s$ and use eq. (14.2-8), we obtain
eq. (14.2-11) which proves theorem 14-2.

\medskip

An essentially identical result as theorem 14-2 was derived previously
by Sidi (see theorem 4.2 on p. 320 of ref. [56]). A comparison of the
asymptotic order estimate (14.2-11) with the analogous order estimates
(13.2-26) and (13.4-17) in theorems 13-9 and 13-12, respectively, which
are all of order $O (n^{- 2 k})$ as $n \to \infty$, shows that the
acceleration of logarithmic convergence is indeed a much more
formidable task than the acceleration of linear convergence, and it can
be even harder than the summation of wildly divergent series.

In extensive numerical studies performed by Smith and Ford [29,30] it
was demonstrated that Levin's $u$ and $v$ transformations, eqs. (7.3-5)
and (7.3-11), respectively, are among the best accelerators for
logarithmic convergence. Hence, we have to conclude that the relative
inefficiency of Levin's sequence transformation according to theorem
14-2 is entirely due to the complicated nature of logarithmically
convergent sequences and that it cannot be attributed to an intrinsic
weakness of Levin's sequence transformation.

The next theorem, which can be proved in the same way as theorem 14-2,
shows that Levin's generalized sequence transformation ${\cal L}_{k,
\ell}^{(n)} (\beta, s_n, \omega_n)$, eq. (7.1-8), with $\ell \ge 1$, is
also able to accelerate logarithmically convergent sequences of the
type of eq. (14.2-8). However, with increasing $\ell \in \N$ the
efficiency of the acceleration process deteriorates.

\medskip

\beginEinzug \sl \parindent = 0 pt

\Auszug {\bf Theorem 14-3:} Let us assume that the sequences $\Seqn s$
and $\Seqn \omega$ are chosen as in theorem 14-2. If the sequence
transformation ${\cal L}_{k, \ell}^{(n)} (\beta, s_n, \omega_n)$, eq.
(7.1-8), is used for the acceleration of the convergence of $\Seqn s$,
we obtain for fixed $k, \ell \in \N$ with $k \ge \ell + 1$ and for $n
\to \infty$ the following order estimate:
$$
\frac
{{\cal L}_{k, \ell}^{(n)} (\beta, s_n, \omega_n) \, - \, s}
{s_n - s} \; = \; O (n^{\ell - k}) \, ,
\qquad n \to \infty \, .
\tag
$$

\endEinzug

\medskip

\Abschnitt Some numerical test series

\smallskip

\aktTag = 0

Due to the lack of theoretical criteria, by means of which it can be
decided whether a given sequence transformation is able to accelerate
logarithmic convergence or not, numerical testing will be of particular
importance. In the literature on convergence acceleration the partial
sums of the series expansion (1.1-2) for $\zeta (2)$,
$$
\zeta (2) \; = \; \sum_{m=0}^{\infty} \, (m + 1)^{-2}
\; = \; \pi^2 / 6 \, ,
\tag
$$
are frequently used to test the ability of a sequence transformation of
accelerating logarithmic convergence. It follows from eqs. (7.3-12) and
(7.3-14) that the remainders $r_n$ of this series are of order $O (n^{-
1})$ as $n \to \infty$, i.e., the convergence of the sequence of
partial sums of this series is prohibitively slow and the series for
$\zeta (2)$ indeed appears to be a good test problem for logarithmic
convergence. However, the next theorem shows that it can happen that a
sequence transformation is able to accelerate the convergence of a
sequence if the remainders behave like an integral power of $1/n$ as $n
\to \infty$ but fails to accelerate convergence if the remainders
behave like a nonintegral power of $1/n$. Consequently, the partial
sums of the series (14.3-1) for $\zeta (2)$ are not suited to test the
ability of a sequence transformation of accelerating the convergence of
a large class of logarithmically convergent sequences.

\medskip

\beginEinzug \sl \parindent = 0 pt

\Auszug {\bf Theorem 14-4:} Let us assume that the linear sequence
transformation $\Lambda_k^{(n)} (\beta, s_n)$, eq. (7.3-20), which
corresponds to the special case $x_n = 1/(\beta + n)$ of the Richardson
extrapolation scheme, eq. (6.1-5), is applied to the following
logarithmically convergent model sequence:
$$
s_n \; = \; s \, + \, \sum_{j=0}^{\infty} \,
\frac {c_j} {(\beta + n)^{\alpha + j} } \, , \qquad n \in \N_0
\, , \quad \alpha, \beta \in \R_{+} \, , \quad c_0 \neq 0 \, .
\tag
$$

If $\alpha$ is a positive integer, i.e., $\alpha = 1, 2, \ldots$, we
obtain for fixed $k \in \N$ and $k \ge \alpha$ and for large values of
$n$ the following asymptotic order estimate
$$
\Lambda_k^{(n)} (\beta, s_n) - s \; = \; O (n^{- k - 1}) \, ,
\tag
$$
which shows that $\Lambda_k^{(n)} (\beta, s_n)$ accelerates the
convergence of this sequence according to eq. (2.7-7) for sufficiently
large values of $k$.

If, however, $\alpha$ is not a positive integer, $\Lambda_k^{(n)}
(\beta, s_n)$ does not accelerate the convergence of the sequence
(14.3-2).

\endEinzug

\medskip

\noindent {\it Proof: } Obviously, $\Lambda_k^{(n)} (\beta, s_n)$, eq.
(7.3-20), is invariant with respect to translation according to eq.
(3.1-4). Hence, with the help of eqs. (7.3-17) and (7.3-19) we can
write:
$$
\Lambda_k^{(n)} (\beta, s_n) \, - \, s \; = \; \frac
{\Delta^k \, [(\beta + n)^k (s_n - s)]} {k!} \, .
\tag
$$

Let us now assume that $\alpha$ is a positive integer, i.e., $\alpha =
m$ with $m \in \N$, and that $k \ge m$. Then, with the help of eq.
(14.3-2) we obtain for the numerator of the right-hand side of eq.
(14.3-4):
$$
\Delta^k \, [(\beta + n)^k (s_n - s)] \; = \;
\Delta^k \, \sum_{j=0}^{k - m} \, c_j (\beta + n)^{k - m - j}
\, + \, \Delta^k \, \sum_{j=0}^{\infty} \, c_{k - m + j + 1}
(\beta + n)^{ - j - 1} \, .
\tag
$$

The first sum on the right-hand side is a polynomial of degree $k - m$
in $n$, i.e., it is annihilated by $\Delta^k$, and according to eq.
(13.2-30) the second sum will produce a result which is of order $O
(n^{- k - 1})$ as $n \to \infty$. This proves eq. (14.3-3).

Let us now assume that $\alpha$ is not a positive integer. Then, with
the help of eqs. (13.2-30) and (14.3-2) we obtain for the numerator of
the right-hand side of eq. (14.3-4) the following asymptotic estimate
$$
\Delta^k \, [(\beta + n)^k (s_n - s)] \; = \; \Delta^k \,
\sum_{j=0}^{\infty} \, c_j \, (\beta + n)^{k - \alpha - j}
\; = \; O (n^{- \alpha}) \, , \qquad n \to \infty \, ,
\tag
$$
which proves the second part of theorem 14-4.

\medskip

Hence, if we want to use the Richardson extrapolation scheme, eq.
(6.1-5), for the acceleration of the logarithmically convergent
sequence (14.3-2) if $\alpha$ is not a positive integer, we cannot use
the interpolation points $x_n = 1/(\beta + n)$ and we would at least
have to find a different set of interpolation points $\Seqn x$. There
is considerable numerical evidence that the situation is quite
analogous in the case of Wynn's $\rho$ algorithm, eq. (6.2-2). The
standard form of Wynn's $\rho$ algorithm, eq. (6.2-4), corresponds to
the choice $x_n = \beta + n$ for the interpolation points. It is
together with its iteration ${\cal W}_k^{(n)}$, eq. (6.3-4), probably
the best accelerator for the partial sums of the series (14.3-1) for
$\zeta (2)$, but is apparently not able to accelerate the convergence
of a sequence with remainders that behave like $n^{- 1/2}$ as $n \to
\infty$. However, it will be shown later that the general form of
Wynn's $\rho$ algorithm, eq. (6.2-2), is able to accelerate the
convergence of sequences with remainders $r_n \sim n^{- 1/2}$ as $n \to
\infty$ if the interpolation points $\Seqn x$ are chosen according to
$x_n = (\beta + n)^{1/2}$ with $\beta > 0$.

We shall see later that for instance Brezinski's $\theta$ algorithm,
eq. (10.1-9), or its iteration ${\cal J}_k^{(n)}$, eq. (10.3-6), are
able to accelerate the convergence of sequences, whose remainders
behave like $n^{- 1/2}$ as $n \to \infty$. Consequently,
$\theta_k^{(n)}$ and ${\cal J}_k^{(n)}$ are more flexible and versatile
than the Richardson extrapolation scheme or Wynn's $\rho$ algorithm
since these sequence transformations only work if appropriate
interpolation points $\Seqn x$ are used.

This implies that because of theorem 14-4 the infinite series (14.3-1)
for $\zeta (2)$ is not suited to test the ability of a sequence
transformation of accelerating logarithmic convergence. Instead, one
should use test problems with remainders that behave like a nonintegral
power of $1/n$ as $n \to \infty$. A simple idea would be to use the
infinite series (1.1-2) for $\zeta (z)$ for nonintegral arguments and
not for $z = 2$. However, only if $z$ is an even positive integer, a
simple explicit expression for $\zeta (z)$ is known (see p. 19 of ref.
[34]). Therefore, the use of the infinite series (1.1-2) for $\zeta
(z)$ with $z \notin \N$ as a test problem for logarithmic convergence
would be somewhat inconvenient and the emphasis in this report will be
on other test problems.

Well suited for our purposes is the following series expansion (see p.
14 of ref. [111]):
$$
A \; = \; \sum_{m=0}^{\infty} \, \frac {(2 m - 1)!!}{(2 m)!!}
\, \frac {1} {4 m + 1} \, .
\tag
$$

Here, $A$ stands for the so-called lemniscate constant which can be
expressed in terms of the complete elliptic integral $K$ (see pp. 358 -
359 of ref. [34]),
$$
A \; = \; \int\nolimits_{0}^{1} \, \frac {\d t} {[1 - t^4]^{1/2}}
\; = \; 2^{- 1/2} K (1/2) \; = \;
\frac {[ \Gamma (1/4)]^2} {4 \, (2 \pi)^{1/2}} \, .
\tag
$$

If we use the following expression for the ratio of two gamma functions
which holds for $\vert z \vert \to \infty$ (see p. 12 of ref. [34]),
$$
\frac {\Gamma (z + \alpha)} {\Gamma (z + \beta)} \; = \;
z^{\alpha - \beta} \, [ 1 + O (z^{-1}) ] \, , \qquad \vert \arg (z)
\vert < \pi \, ,
\tag
$$
we find that the terms of the series (14.3-7) behave like $n^{-3/2}$ as
$n \to \infty$. Consequently, it follows from eqs. (7.3-12) and
(7.3-14) that the remainders of this series behave like $n^{-1/2}$ as
$n \to \infty$.

When Smith and Ford [29] investigated numerically the performance of
sequence transformations in convergence acceleration and summation
processes, they used the series (14.3-7) for the lemniscate constant
$A$ to test the ability of a sequence transformation of accelerating
logarithmic convergence. Smith and Ford observed that the standard
version of Wynn's $\rho$ algorithm, eq. (6.2-4), was not able to
accelerate the convergence of this series (see p. 235 of ref. [29]).

Another test problem, which is well suited for our purposes, is the
following series expansion for $1/z$ in terms of the so-called reduced
Bessel functions,
$$
1 / z \; = \; \sum_{m=0}^{\infty} \,
{\widehat k}_{m-1/2} (z) \, / \, [2^m \, m!] \, ,
\qquad z \in \R_{+} \, .
\tag
$$

This series expansion was derived and used in connection with explicit
expressions for certain molecular multicenter integrals of
exponentially declining basis functions (see eq. (6.5) of ref. [112]).
In table I of ref. [64] it was shown that this series converges
extremely slowly. For $z = 1$ the infinite series (14.3-10) produced an
accuracy of only three decimal digits after $1\ 000\ 000$ terms.

The so-called reduced Bessel function ${\widehat k}_{\nu} (z)$ of real
or complex order $\nu$, which was introduced by Steinborn and Filter
(see eqs. (3.1) and (3.2) of ref. [113]) as an exponentially declining
basis function in electronic structure calculations, is defined by
$$
{\widehat k}_{\nu} (z) \; = \;
(2/\pi)^{1/2} \, z^{\nu} \, K_{\nu} (z) \, .
\tag
$$
Here, $K_{\nu} (z)$ is a modified Bessel function of the second kind
(see p. 66 of ref. [34]). If the order $\nu$ of the reduced Bessel
function is half-integral and positive, $\nu = n + 1/2$ with $n \in
\N_0$, a reduced Bessel function can be represented as an exponential
multiplied by a terminating confluent hypergeometric series ${}_1 F_1$
(see eq. (3.7) of ref. [114]),
$$
{\widehat k}_{n+1/2} (z) \; = \; 2^n \, (1/2)_n \, \e^{- z} \,
{}_1 F_1 (- n; - 2 n;2 z) \, .
\tag
$$

The polynomial part of these reduced Bessel functions with
half-integral orders $\nu = n + 1/2$ with $n \in \N_0$ has also been
investigated independently in the mathematical literature. There, the
following notation is used (see p. 34 of ref. [115]):
$$
\theta_n (z) \; = \; \e^z \, {\widehat k}_{n+1/2} (z) \; = \;
2^n \, (1/2)_n \, {}_1 F_1 (- n; - 2 n;2 z) \, ,
\qquad n \in \N_0 \, .
\tag
$$

Together with some other, closely related polynomials, which are
denoted by $y_n (z)$, the polynomials $\theta_n (z)$ are called Bessel
polynomials. In Grosswald's book on Bessel polynomials [115] numerous
applications of these polynomials in vastly differing fields are
described. For instance, they are applied in number theory, in
statistics, or for the analysis of complex electrical networks.

In the context of convergence acceleration and summation it may be
interesting to note that Bessel polynomials occur also in the theory of
Pad\'e approximants. In the book by Baker and Graves-Morris it is shown
that the Pad\'e approximants $[ \, \ell \, / \, m \, ]$ for $\e^z$ are
given by (see eq. (2.12) of part I of ref. [22])
$$
[ \, \ell \, / \, m \, ] \; = \; \frac
{ {}_1 F_1 (- \ell; - \ell - m; z)}
{ {}_1 F_1 (- m; - \ell - m;- z)} \, ,
\qquad \ell, m \in \N_0 \, .
\tag
$$
Comparison of eqs. (14.3-13) and (14.3-14) shows that the diagonal
elements $[ \, n \, / \, n \, ]$ of the Pad\'e table for $\e^z$ can be
expressed as ratios of Bessel polynomials,
$$
[ \, n \, / \, n \, ] \; = \; \frac
{ \theta_n (z/2) } { \theta_n (- z/2) } \, ,
\qquad n \in \N_0 \, .
\tag
$$

With the help of some well known monotonicity properties of the
modified Bessel function of the second kind, $K_{\nu} (z)$, it can be
shown that the reduced Bessel functions ${\widehat k}_{\nu} (z)$ are
positive and bounded by their values at the origin provided that $\nu >
0$ and $z \ge 0$ (see eq. (3.1) of ref. [114]). In the case of reduced
Bessel functions with half-integral orders this implies:
$$
0 \; < \; {\widehat k}_{n+1/2} (z) \; \le \;
{\widehat k}_{n+1/2} (0) \; = \; 2^n \, (1/2)_n \, ,
\qquad 0 \le z < \infty \, , \quad n \in \N_0 \, .
\tag
$$

Grosswald's book [115] also contains a chapter on the asymptotic
properties of Bessel polynomials. There, it is shown that for fixed and
finite argument $z$ the Bessel polynomials $\theta_n (z)$ satisfy (see
p. 125 of ref. [115])
$$
\theta_n (z) \; \sim \; \frac {(2 n)!} {2^n n!} \, \e^z
\, , \qquad n \to \infty \, .
\tag
$$

If we combine eqs. (14.3-13) and (14.3-17) we find that the dominant
term of the Poincar\'e-type asymptotic expansion in inverse powers of $n$
of a reduced Bessel function ${\widehat k}_{n+1/2} (z)$ with fixed and
finite argument $z$ corresponds to its value at the origin,
$$
{\widehat k}_{n+1/2} (z) \; = \;
2^n \, (1/2)_n \, [1 + O (n^{- 1})] \; = \;
{\widehat k}_{n+1/2} (0) \, [1 + O (n^{- 1})]
\, , \qquad n \to \infty \, .
\tag
$$

Higher terms of the asymptotic expansion of a reduced Bessel function
${\widehat k}_{n+1/2} (z)$ in inverse powers of $n$ can in principle be
obtained from related expansions for Bessel polynomials $y_n (z)$. In
Grosswald's book the coefficients for terms up to an order $O (n^{-
3})$ can be found (see p. 130 of ref. [115]) and in an article by
Salzer [116] the coefficients for the terms up to an order $O (n^{-
4})$.

Starting from eq. (14.3-18) it can be proved quite easily with the help
of eq. (14.3-9) that the terms of the infinite series (14.3-10) behave
like $n^{- 3/2}$ as $n \to \infty$ (see p. 3709 of ref. [64]). In view
of eqs. (7.3-12) and (7.3-14) this implies that the remainders of the
partial sums of the series (14.3-10) behave like $n^{- 1/2}$ as $n \to
\infty$.

From the series expansion (14.3-10) in terms of reduced Bessel
functions another series of the same convergence type can be derived.
The new series is essentially the $z$-independent part of the infinite
series (14.3-10). If we take into account that ${\widehat k}_{- 1/2}
(z) = \e^{- z}/z$, we can conclude from eq. (14.3-12) that only the
first term of the infinite series (14.3-10) is singular at the origin.
Consequently, the following limit exists,
$$
\lim_{z \to 0} \, [ \, 1/z \, - \, {\widehat k}_{- 1/2} (z) \, ]
\; = \; \lim_{z \to 0} \, \sum_{m=1}^{\infty} \,
{\widehat k}_{m - 1/2} (z) \, / \, [2^m \, m!] \; = \; 1 \, ,
\tag
$$
and we obtain with the help of eq. (14.3-16):
$$
\sum_{m=0}^{\infty} \, \frac {(2 m - 1)!!} {( 2 m + 2)!!} \; = \;
\frac {1} {2} \, \sum_{m=0}^{\infty} \, \frac {(1/2)_m} {(m+1)!}
\; = \; 1 \, .
\tag
$$

Again, it follows from eq. (14.3-9) that the terms of this series
behave like $n^{- 3/2}$ as $n \to \infty$, which according to eqs.
(7.3-12) and (7.3-14) implies that the remainders of this series behave
like $n^{- 1/2}$.

Several sequence transformations are exact for the infinite series
(14.3-20) since its partial sums satisfy the prerequisites of theorem
14-1. With the help of a summation theorem by Gauss for a
hypergeometric series ${}_2 F_1$ with unit argument (see p. 40 of ref.
[34]) it can be proved quite easily that the remainders of the infinite
series (14.3-20) satisfy
$$
\sum_{m=n+1}^{\infty} \,
\frac {(2 m - 1)!!} {( 2 m + 2)!!} \; = \;
\frac {(1/2)_{n+1}} {(n+1)!} \, ,
\tag
$$
which shows that the partial sums of the infinite series (14.3-20) are
of the type of the sequence (14.2-7). Hence, it follows from theorem
14-1 that Levin's $u$ transformation, eq. (7.3-5), will only need the
partial sums $s_{n-1}, s_n, s_{n+1}$, and $s_{n+2}$ to be exact for the
infinite series (14.3-20). Also, from eq. (10.3-4) it follows that
Brezinski's $\theta$ algorithm, eq. (10.1-9), and its iteration ${\cal
J}_k^{(n)}$, eq. (10.3-6), only need the partial sums $s_n, s_{n+1},
s_{n+2}$, and $s_{n+3}$ to be exact for the infinite series (14.3-20).

Levin's $v$ transformation, eq. (7.3-11), is even more efficient than
the sequence transformations mentioned above because only the partial
sums $s_{n-1}$, $s_n$, and $s_{n+1}$ are needed to sum the infinite
series (14.3-20) exactly. This follows from the fact that in the case
of the infinite series (14.3-20) the remainder estimate (7.3-10), which
is the basis for the $v$ transformation,
$$
\frac {a_{n+1} a_n} {a_{n+1} - a_n} \; = \;
- \, \frac {1} {3} \, \frac {(1/2)_{n+1}} {(n+1)!} \, ,
\tag
$$
is proportional to the remainder (14.3-21) of the infinite series
(14.3-20). Consequently, in this case the ratio $(s_n - s)/\omega_n$ is
independent of $n$ which implies that $(\beta + n)^{k-1} (s_n -
s)/\omega_n$ will be annihilated by $\Delta^k$ for $k \ge 1$.

The terms of the three test series (14.3-7), (14.3-10), and (14.3-20)
all behave like $n^{- 3/2}$ as $n \to \infty$. Consequently, we may
expect that these three series should be roughly comparable with
respect to their rates of convergence as well as in convergence
acceleration processes. However, at least in convergence acceleration
processes these three test series are definitely not equivalent. The
acceleration of the convergence of the series expansion (14.3-10) in
terms of reduced Bessel functions is a much more formidable task, in
particular for larger values of $z$, than the acceleration of the
convergence of the other two series (14.3-7) and (14.3-20). In fact,
for sufficiently large values of $z$ it is virtually impossible to
accelerate the convergence of the infinite series (14.3-10). This is
probably a consequence of the exponential decline of the terms and also
of the partial sums of this series. Numerical tests showed that for
larger values of $z \in \R_{+}$ a reduced Bessel function is
approximated by its value at the origin,
$$
{\widehat k}_{n+1/2} (z) \; \approx \; {\widehat k}_{n+1/2} (0)
\; = \; 2^n \, (1/2)_n \, , \qquad n \in \N_0 \, ,
\tag
$$
with reasonable accuracy only if $n$ is very large. For instance, if we
require that eq. (14.3-23) should be accurate to one percent for $z =
8$ then we would need $n \ge 1400$, and for $z = 4$ we would still need
$n \ge 400$. Consequently, in particular for larger values of $z$ a
partial sum $s_n$ of the series (14.3-10) is essentially a linear
combination of some quantities which decline exponentially, and only
for relatively large values of $n$ it can actually be observed that
$s_n$ behaves like $n^{- 1/2}$ as $n \to \infty$. Thus, it is likely
that for a sequence transformation, which tries to extract and utilize
some regularity in the behaviour of the partial sums, the partial sums
of the infinite series (14.3-10) appear for larger values of $z$ to be
much more irregular than the partial sums of the other two infinite
series (14.3-7) and (14.3-20).

\medskip

\Abschnitt Numerical examples

\smallskip

\aktTag = 0

In this section, the acceleration of logarithmically convergent
sequences will be studied numerically. The emphasis will be on the test
series (14.3-7), (14.3-10), and (14.3-20), which should converge
approximately as slowly as the infinite series (1.1-2) for $\zeta
(3/2)$. But since the infinite series (14.3-1) for $\zeta (2)$ is the
most popular test problem for logarithmic convergence in the
literature, it is of interest to see how in particular the new sequence
transformations fare if they are applied to the partial sums
$$
s_n \; = \; \sum_{m=0}^n \, (m+1)^{- 2}
\, , \qquad n \in \N_0 \, ,
\tag
$$
of the infinite series (14.3-1) for $\zeta (2)$.

In table 14-1 the partial sums (14.4-1) are transformed by the standard
form of ${\cal W}_k^{(n)}$, eq. (6.3-4), by $\Lambda_k^{(n)} (\beta,
s_n)$, eq. (7.3-20), with $\beta = 1$, which corresponds to the special
case $x_n = 1/(\beta + n)$ of the Richardson extrapolation scheme, eq.
(6.1-5), and by ${\cal J}_k^{(n)}$, eq. (10.3-6). In all cases, the
approximants were chosen in such a way that the information, which is
contained in the finite string $s_0, s_1, \ldots, s_n$ of partial sums,
is exploited in an optimal way. This means that in the case of ${\cal
W}_k^{(n)}$, eq. (6.3-4), the approximations to $\zeta (2)$ were chosen
according to eq. (6.3-9), and in the case of ${\cal J}_k^{(n)}$, they
were chosen according to eq. (10.4-7).

\beginFloat

\medskip

\beginTabelle [to \kolumnenbreite]
\beginFormat \rechts " \mitte " \mitte " \mitte " \mitte
\endFormat
\+ " \links {\bf Table 14-1} \@ \@ \@ \@ " \\
\+ " \links {Acceleration of the series (14.3-1) for $\zeta (2) = \pi^2
/ 6$} \@ \@ \@ \@ " \\
\- " \- " \- " \- " \- " \- " \\ \sstrut {} {1.5 \jot} {1.5 \jot}
\+ " \rechts {$n$} " \mitte {partial sum $s_n$} " ${\cal W}_{\Ent
{n/2}}^{(n - 2 \Ent {n/2})}$ " $\Lambda_n^{(0)} (1, s_0)$ " ${\cal
J}_{\Ent {n/3}}^{(n - 3 \Ent {n/3})}$ " \\
\+ " " \mitte {eq. (14.4-1)} " eq. (6.3-4) " eq. (7.3-20) " eq.
(10.3-6) " \\
\- " \- " \- " \- " \- " \- " \\ \sstrut {} {1 \jot} {1 \jot}
\+ "    0 "    1.0000000000000   "   1.0000000000000 "
1.0000000000000 "   1.0000000000000 " \\
\+ "    1 "    1.2500000000000   "   1.2500000000000 "
1.5000000000000 "   1.2500000000000 " \\
\+ "    2 "    1.3611111111111   "   1.6500000000000 "
1.6250000000000 "   1.3611111111111 " \\
\+ "    3 "    1.4236111111111   "   1.6468253968254 "
1.6435185185185 "   1.6388888888889 " \\
\+ "    4 "    1.4636111111111   "   1.6449013949014 "
1.6449652777778 "   1.6423611111111 " \\
\+ "    5 "    1.4913888888889   "   1.6449244489889 "
1.6449513888889 "   1.6436111111111 " \\
\+ "    6 "    1.5117970521542   "   1.6449342449874 "
1.6449351851852 "   1.6449225865209 " \\
\+ "    7 "    1.5274220521542   "   1.6449341126465 "
1.6449339434186 "   1.6449297924298 " \\
\+ "    8 "    1.5397677311665   "   1.6449340660297 "
1.6449340411698 "   1.6449321959755 " \\
\+ "    9 "    1.5497677311665   "   1.6449340666548 "
1.6449340662475 "   1.6449340557022 " \\
\+ "   10 "    1.5580321939765   "   1.6449340668515 "
1.6449340671488 "   1.6449340629267 " \\
\+ "   11 "    1.5649766384209   "   1.6449340668489 "
1.6449340668835 "   1.6449340652730 " \\
\+ "   12 "    1.5708937981842   "   1.6449340668482 "
1.6449340668472 "   1.6449340668410 " \\
\+ "   13 "    1.5759958390005   "   1.6449340668482 "
1.6449340668476 "   1.6449340668458 " \\
\+ "   14 "    1.5804402834450   "   1.6449340668482 "
1.6449340668482 "   1.6449340668473 " \\
\+ "   15 "    1.5843465334450   "   1.6449340668482 "
1.6449340668482 "   1.6449340668482 " \\
\- " \- " \- " \- " \- " \- " \\ \sstrut {} {1 \jot} {1 \jot}
\+ " \links {$\pi^2 /6$} \@    "   1.6449340668482 "   1.6449340668482
"   1.6449340668482 " \\
\- " \- " \- " \- " \- " \- " \\ \sstrut {} {1 \jot} {1 \jot}

\endTabelle

\medskip

\endFloat

All sequence transformations in table 14-1 accelerate the convergence
of the infinite series (14.3-1) for $\zeta (2)$ quite efficiently. The
winner in table 14-1 is ${\cal W}_k^{(n)}$, eq. (6.3-4), which together
with the standard form of Wynn's $\rho$ algorithm, eq. (6.2-4), is the
best accelerator for the series for $\zeta (2)$. Somewhat less
efficient are $\Lambda_k^{(n)} (\beta, s_n)$, eq. (7.3-20), which in
the case of the partial sums (14.4-1) is identical with Levin's $u$
transformation, eq. (7.3-5), and ${\cal J}_k^{(n)}$, eq. (10.3-6).

Other good accelerators for the partial sums (14.4-1) are Levin's $v$
transformation, eq. (7.3-11), which is as efficient as Levin's $u$
transformation, eq. (7.3-5), and Brezinski's $\theta$ algorithm, eq.
(10.1-9), and $\sigma_k^{(n)}$, eq. (11.2-2), which are almost as
efficient as ${\cal J}_k^{(n)}$, eq. (10.3-6).

The partial sums and the three transforms in table 14-1 were computed
in QUADRUPLE PRECISION (31 - 32 decimal digits). When these
computations were repeated in DOUBLE PRECISION (15 - 16 decimal digits)
the loss of some significant digits was observed. This is not
surprising since the acceleration of logarithmic convergence is -- as
emphasized previously -- an inherently unstable process. Of the three
sequence transformations in table 14-1 it was again ${\cal W}_k^{(n)}$,
eq. (6.3-4), which turned out to be the numerically most stable
transformation since it lost at most 3 significant digits in DOUBLE
PRECISION. The other two transformations, which lost up to 5 decimal
digits, are apparently more sensitive to rounding errors.

Next, we want to see how the convergence of the test series (14.3-7),
(14.3-10), and (14.3-20) can be accelerated. Since the remainders $r_n$
of these series behave like $n^{- 1/2}$ as $n \to \infty$, we expect
that these series will converge significantly more slowly than the
series (14.3-1) for $\zeta (2)$, whose remainders are of order $O (n^{-
1})$. Here, it must be emphasized that it is not clear how and to what
extent the slower convergence of the test series (14.3-7), (14.3-10),
and (14.3-20) will affect convergence acceleration processes.

\beginFloat

\medskip

\beginTabelle [to \kolumnenbreite]
\beginFormat \rechts " \mitte " \mitte " \mitte " \mitte
\endFormat
\+ " \links {\bf Table 14-2} \@ \@ \@ \@ " \\
\+ " \links {Acceleration of the series (14.3-7) for the lemniscate
constant $A$} \@ \@ \@ \@ " \\
\- " \- " \- " \- " \- " \- " \\ \sstrut {} {1.5 \jot} {1.5 \jot}
\+ " \rechts {$n$} " \mitte {partial sum $s_n$} " ${\cal L}_{n,
2}^{(0)} (1, s_0, a_0)$ " $u_n^{(0)} (1, s_0)$ " ${\cal J}_{\Ent
{n/3}}^{(n - 3 \Ent {n/3})}$ " \\
\+ " " \mitte {eq. (14.4-2)} " eq. (7.1-8) " eq. (7.3-5) " eq. (10.3-6)
" \\
\- " \- " \- " \- " \- " \- " \\ \sstrut {} {1 \jot} {1 \jot}
\+ "    3  "    1.1657051282051    "   1.2190476190476  "
1.3163120567376  "   1.3037037037037  " \\
\+ "    4  "    1.1817896870287    "   1.3343421605717  "
1.3108727079053  "   1.3080867850099  " \\
\+ "    5  "    1.1935084370287    "   1.3103293923028  "
1.3109952008776  "   1.3095200070979  " \\
\+ "    6  "    1.2025318745287    "   1.3108082123295  "
1.3110289627926  "   1.3110119624014  " \\
\+ "    7  "    1.2097550695718    "   1.3110277257117  "
1.3110291499078  "   1.3110229739315  " \\
\+ "    8  "    1.2157059973061    "   1.3110318476640  "
1.3110287979182  "   1.3110263202535  " \\
\+ "    9  "    1.2207187157130    "   1.3110290080304  "
1.3110287737803  "   1.3110287611468  " \\
\+ "   10  "    1.2250162047862    "   1.3110287461269  "
1.3110287766205  "   1.3110287718416  " \\
\+ "   11  "    1.2287537180105    "   1.3110287708982  "
1.3110287771522  "   1.3110287750875  " \\
\+ "   12  "    1.2320431110268    "   1.3110287771312  "
1.3110287771540  "   1.3110287771349  " \\
\+ "   13  "    1.2349672811610    "   1.3110287772480  "
1.3110287771466  "   1.3110287771425  " \\
\+ "   14  "    1.2375891404731    "   1.3110287771553  "
1.3110287771460  "   1.3110287771447  " \\
\+ "   15  "    1.2399574101139    "   1.3110287771452  "
1.3110287771460  "   1.3110287771461  " \\
\+ "   16  "    1.2421104860230    "   1.3110287771458  "
1.3110287771461  "   1.3110287771461  " \\
\+ "   17  "    1.2440790912340    "   1.3110287771461  "
1.3110287771461  "   1.3110287771461  " \\
\+ "   18  "    1.2458881405432    "   1.3110287771461  "
1.3110287771461  "   1.3110287771461  " \\
\- " \- " \- " \- " \- " \- " \\ \sstrut {} {1 \jot} {1 \jot}
\+ " \links {$[\Gamma (1/4)]^2 / [4 (2 \pi)^{1/2}]$} \@ "
1.3110287771461 "   1.3110287771461 "   1.3110287771461 " \\ \- " \- "
\- " \- " \- " \- " \\ \sstrut {} {1 \jot} {1 \jot}

\endTabelle

\medskip

\endFloat

First, we shall accelerate the convergence of the sequence of partial
sums of the infinite series (14.3-7) for the lemniscate constant $A$,
$$
s_n \; = \; \sum_{m=0}^n \, \frac {(2 m - 1)!!} {(2 m)!!} \,
\frac {1} {4 m + 1} \, .
\tag
$$

In table 14-2 we see the effect of Levin's generalized sequence
transformation ${\cal L}_{k, \ell}^{(n)} (\beta, s_n, \omega_n)$, eq.
(7.1-8), with $\ell = 2$, $\omega_n = a_n$, and $\beta = 1$, of Levin's
$u$ transformation, eq. (7.3-5), with $\beta = 1$, and of ${\cal
J}_k^{(n)}$ on the partial sums (14.4-2). As usual, the approximants
were chosen in such a way that the information, which is contained in
the finite string $s_0, s_1, \ldots, s_n$ of partial sums, is exploited
optimally. This means that in the case of the Levin transformations the
approximations to the lemniscate constant $A$ were chosen according to
eq. (7.5-4).

If we compare the results in tables 14-1 and 14-2, we see that the
significantly slower convergence of the infinite series (14.3-7) does
not affect the efficiency of convergence acceleration too much. The
convergence of the transforms in table 14-2 is almost as fast as in
table 14-1. However, it seems that the slower convergence of the test
series (14.3-7) has a detrimental effect on the numerical stability of
the transformations. As usual, table 14-2 was produced in QUADRUPLE
PRECISION. When the same computations were repeated in DOUBLE
PRECISION, a larger number of significant digits were lost than in
table 14-1. The best results were obtained by Levin's $u$
transformation which achieved a relative accuracy of 11 decimal digits
after $n = 11$. For larger values of $n$, the accuracy deteriorated
again. For instance, for $n = 18$ the best results were obtained by
${\cal J}_{12}^{(0)}$ which achieved an accuracy of 8 decimal digits.

If we compare eqs. (7.1-8) and (7.3-5), we find that Levin's $u$
transformation may also be considered to be a special case of Levin's
generalized sequence transformation ${\cal L}_{k, \ell}^{(n)} (\beta,
s_n, \omega_n)$,
$$
u_k^{(n)} (\beta, s_n) \; = \;
{\cal L}_{k, 1}^{(n)} (\beta, s_n, a_n) \, .
\tag
$$

This relationship shows that the numerical data in table 14-2 are in
agreement with theorem 14-3 which predicts that the efficiency of
Levin's generalized sequence transformation ${\cal L}_{k, \ell}^{(n)}
(\beta, s_n, \omega_n)$, eq. (7.1-8), for the acceleration of the
convergence of sequences of the type of eq. (14.2-8) should decrease
with increasing $\ell \in \N$. This behaviour is apparently quite
typical since it was also observed when the partial sums of the test
series (14.3-10) were accelerated by $u_k^{(n)} (\beta, s_n)$, eq.
(7.3-5), and ${\cal L}_{k, \ell}^{(n)} (\beta, s_n, a_n)$, eq. (7.1-8),
with $\ell \ge 2$.

Other good sequence transformations for the test series (14.3-7) are
Levin's $v$ transformation, eq. (7.3-11), which is as good or even
slightly better than Levin's $u$ transformation, eq. (7.3-5), and
Brezinski's $\theta$ algorithm, eq. (10.1-9), which is as good as
${\cal L}_{k, 2}^{(n)} (\beta, s_n, a_n)$.

Theorem 14-4 predicts that $\Lambda_k^{(n)} (\beta, s_n)$, eq.
(7.3-20), which corresponds to the special case $x_n = 1/(\beta + n)$
of the Richardson extrapolation scheme, eq. (6.1-5), is not able to
accelerate the convergence of the series (14.3-7), (14.3-10), and
(14.3-20). Numerical tests confirmed this prediction. In addition, the
standard forms of Wynn's $\rho$ algorithm, eq. (6.2-4), and of its
iteration ${\cal W}_k^{(n)}$, eq. (6.3-4), also fail to accelerate the
convergence of these series. Since the remainders of the test series
mentioned above all behave like $n^{- 1/2}$ as $n \to \infty$, and
since the Richardson extrapolation scheme ${\cal N}_k^{(n)} (s_n,
x_n)$, eq. (6.1-5), is by construction exact for the model sequence
(6.1-6), it is an obvious idea to choose an alternative set of
interpolation points $\Seqn x$ according to
$$
x_n \; = \; (\beta + n)^{- 1/2} \, , \qquad n \in \N_0 \, ,
\quad \beta \in \R_{+} \, .
\tag
$$

Practical experience has shown that if the Richardson extrapolation
scheme, eq. (6.1-5), can successfully handle a certain problem if the
interpolation points $\Seqn x$ are used, then the general forms of
Wynn's $\rho$ algorithm, eq. (6.2-2), and of its iteration ${\cal
W}_k^{(n)}$, eq. (6.3-3), are usually able to handle the same problem
if the interpolation points $\Seqn \xi$ with $\xi_n = 1/x_n$ are used.
Hence, if we want to use these nonlinear sequence transformations for
the acceleration of the convergence of the series (14.3-7), (14.3-10),
and (14.3-20), we should choose the elements of the set $\Seqn \xi$ of
interpolation points according to
$$
\xi_n \; = \; (\beta + n)^{1/2} \, , \qquad n \in \N_0 \, ,
\quad \beta \in \R_{+} \, .\tag
$$

\beginFloat

\medskip

\beginTabelle [to \kolumnenbreite]
\beginFormat \rechts " \mitte " \mitte " \mitte " \mitte
\endFormat
\+ " \links {\bf Table 14-3} \@ \@ \@ \@ " \\
\+ " \links {Acceleration of the series (14.3-7) for the lemniscate
constant $A$} \@ \@ \@ \@ " \\
\- " \- " \- " \- " \- " \- " \\ \sstrut {} {1.5 \jot} {1.5 \jot}
\+ " \rechts {$n$} " \mitte {partial sum $s_n$} " ${\cal N}_n^{(0)}
(s_0, x_0)$ " $\rho_{2 \Ent {n/2}}^{(n - 2 \Ent {n/2})}$ " ${\cal
W}_{\Ent {n/2}}^{(n - 2 \Ent {n/2})}$ " \\
\+ " " \mitte {eq. (14.4-2)} " eq. (6.1-5) " eq. (6.2-2) " eq. (6.3-3)
" \\
\+  " " " $x_n = (n+1)^{- 1/2}$ " $\xi_n = (n+1)^{1/2}$ " $\xi_n =
(n+1)^{1/2}$ " \\
\- " \- " \- " \- " \- " \- " \\ \sstrut {} {1 \jot} {1 \jot}
\+ "   10  "    1.2250162047862     "   1.3110281470344  "
1.3110287489950  "   1.3110289097428 " \\
\+ "   11  "    1.2287537180105     "   1.3110291844571  "
1.3110287903217  "   1.3110290181070 " \\
\+ "   12  "    1.2320431110268     "   1.3110287720479  "
1.3110287927338  "   1.3110289480394 " \\
\+ "   13  "    1.2349672811610     "   1.3110287371888  "
1.3110287908112  "   1.3110286302938 " \\
\+ "   14  "    1.2375891404731     "   1.3110287841474  "
1.3110287783445  "   1.3110287340891 " \\
\+ "   15  "    1.2399574101139     "   1.3110287800479  "
1.3110287774952  "   1.3110287694697 " \\
\+ "   16  "    1.2421104860230     "   1.3110287759641  "
1.3110287770923  "   1.3110287666673 " \\
\+ "   17  "    1.2440790912340     "   1.3110287770396  "
1.3110287771290  "   1.3110287675098 " \\
\+ "   18  "    1.2458881405432     "   1.3110287772843  "
1.3110287771523  "   1.3110287670755 " \\
\+ "   19  "    1.2475580797723     "   1.3110287771349  "
1.3110287771467  "   1.3110287674918 " \\
\+ "   20  "    1.2491058660392     "   1.3110287771336  "
1.3110287771459  "   1.3110287670448 " \\
\+ "   21  "    1.2505456974656     "   1.3110287771492  "
1.3110287771460  "   1.3110287770066 " \\
\+ "   22  "    1.2518895646423     "   1.3110287771468  "
1.3110287771461  "   1.3110287770883 " \\
\+ "   23  "    1.2531476731141     "   1.3110287771456  "
1.3110287771461  "   1.3110287771050 " \\
\+ "   24  "    1.2543287710802     "   1.3110287771461  "
1.3110287771461  "   1.3110287770977 " \\
\+ "   25  "    1.2554404064530     "   1.3110287771461  "
1.3110287771461  "   1.3110287770960 " \\
\- " \- " \- " \- " \- " \- " \\ \sstrut {} {1 \jot} {1 \jot}
\+ " \links {$[\Gamma (1/4)]^2 / [4 (2 \pi)^{1/2}]$} \@ "
1.3110287771461 "   1.3110287771461 "   1.3110287771461 " \\ \- " \- "
\- " \- " \- " \- " \\ \sstrut {} {1 \jot} {1 \jot}

\endTabelle

\medskip

\endFloat

In table 14-3 the partial sums (14.4-2) are accelerated by the
Richardson extrapolation scheme, eq. (6.1-5), and by the general forms
of Wynn's $\rho$ algorithm, eq. (6.2-2), and of ${\cal W}_k^{(n)}$, eq.
(6.3-3). In the case of the Richardson extrapolation scheme the
interpolation points $\Seqn x$ were chosen according to eq. (14.4-4)
with $\beta = 1$, and in the case of $\rho_k^{(n)}$ and ${\cal
W}_k^{(n)}$ the interpolation points $\Seqn \xi$ were chosen according
to eq. (14.4-5) with $\beta = 1$.

The most efficient transformation in table 14-3 is the general form of
Wynn's $\rho$ algorithm, eq. (6.2-2), followed by the Richardson
extrapolation scheme, eq. (6.1-5), and the general form of ${\cal
W}_k^{(n)}$, eq. (6.3-3). However, a comparison of tables 14-2 and 14-3
shows that even if we choose the interpolation points according to eqs.
(14.4-4) and (14.4-5), the transformations in table 14-3 are clearly
less efficient than Levin's $u$ transformation or ${\cal J}_k^{(n)}$.

As usual, table 14-3 was produced in QUADRUPLE PRECISION. When the same
computations were repeated in DOUBLE PRECISION, it was observed that
the Richardson extrapolation scheme is much more unstable than the
other two transformations. The Richardson extrapolation scheme achieved
a relative accuracy of 8 decimal digits for $n = 12$. For larger values
of $n$ the accuracy of the transforms deteriorated rapidly, yielding
totally nonsensical results for $n \ge 22$. The other two
transformations also did not accomplish more than a relative accuracy
of 7 or 8 decimal digits. However, they maintained this relative
accuracy throughout the whole range of $n$ between $n = 10$ and $n =
25$.

The results in table 14-3 and similar results for the other two test
series (14.3-10) and (14.3-20) indicate that the Richardson
extrapolation scheme, eq. (6.1-5), is able to accelerate the
convergence of a sequence, whose remainders $r_n$ behave like $n^{-
1/2}$ as $n \to \infty$, if the interpolation points $\Seqn x$ are
chosen according to eq. (14.4-4). Similarly, the general forms of Wynn's
$\rho$ algorithm, eq. (6.2-2), or of ${\cal W}_k^{(n)}$, eq. (6.3-3),
should be able to accelerate the convergence of such a sequence if the
interpolation points $\Seqn \xi$ are chosen according to eq. (14.4-5).

Unfortunately, in practical applications these observations are not
necessarily very helpful. Let us assume that only the numerical values
of a few sequence elements are known but nothing about the behaviour of
the remainders. In such a situation, it will be very hard or even
impossible to find out whether the remainders of this sequence behave
like an integral or like a nonintegral power of $1/n$. If one wants to
use the Richardson extrapolation scheme, eq. (6.1-5), in such an
unfavourable situation, it may be a good idea to follow a
recommendation by Beleznay [117]. He suggested to choose the
interpolation points according to
$$
x_n \; = \; (n + \beta)^{- \alpha} \, , \qquad n \in \N_0 \, ,
\quad \alpha, \beta \in \R_{+} \, ,
\tag
$$
and to optimize the exponent $\alpha$ of the interpolation points in
such a way that the error term $\vert {\cal N}_{n-1}^{(1)} - {\cal
N}_{n-1}^{(0)} \vert$ becomes minimal. This technique was later used by
Liegener, Beleznay, and Ladik [118] to extrapolate the results of
Hartree-Fock calculations on periodic chains. A similar approach is of
course also possible in the case of the general forms of Wynn's $\rho$
algorithm, eq. (6.2-2), or of its iteration ${\cal W}_k^{(n)}$, eq.
(6.3-3). In that case, the interpolation points would have to be chosen
according to
$$
\xi_n \; = \; (n + \beta)^{\alpha} \, , \qquad n \in \N_0 \, ,
\quad \alpha, \beta \in \R_{+} \, .
\tag
$$

There is still another possibility of modifying either the Richardson
extrapolation scheme, eq. (6.1-5), or the general forms of the rational
transformations $\rho_k^{(n)}$, eq. (6.2-2), and ${\cal W}_k^{(n)}$,
eq. (6.3-3), in such a way that it will be unimportant whether the
remainder of the sequence to be transformed behaves like an integral or
a nonintegral power of $1/n$. This approach is inspired by a nonlinear
variant of the Richardson extrapolation scheme, which in Wimp's book is
called GBW (Germain-Bonne Wimp) transformation (see p. 106 of ref.
[23]). This GBW transformation is obtained from the Richardson
extrapolation scheme by choosing the interpolation points according to
$$
x_n \; = \; \Delta s_n \; = \; a_{n+1} \, ,
\qquad n \in \N_0 \, .
\tag
$$

If the interpolation points $\Seqn x$ are chosen according to eq.
(14.4-8), then it follows from eq. (6.1-6) that the Richardson
extrapolation scheme, eq. (6.1-5), is exact for the following model
sequence:
$$
s_n \; = \; s \, + \, \sum_{j=0}^{k-1} \, c_j \,
[ \Delta s_n ]^{j+1} \, , \qquad k, n \in \N_0
\, , \quad \beta \in \R_{+} \, .
\tag
$$

According to Wimp this GBW transformation works quite well in the case
of iteration sequences (see pp. 106 - 108 of ref. [23]). However, the
GBW transformation is apparently not able to accelerate logarithmic
convergence.

Let us now assume that the remainders $r_n$ of a sequence behave like
$n^{- \alpha}$ as $n \to \infty$. Then it follows from eq. (14.2-2)
that the product $[ n \, \Delta s_{n-1} ]$ also behaves like $n^{-
\alpha}$ as $n \to \infty$. Hence, if we choose the interpolation
points $\Seqn x$ for the Richardson extrapolation scheme according to
$$
x_n \; = \; (\beta + n) \, \Delta s_{n-1} \, ,
\qquad n \in \N_0 \, , \quad \beta \in \R_{+} \, ,
\tag
$$
it is at least guaranteed that the interpolation points $x_n$ behave
like the $r_n$ as $n \to \infty$.

If the interpolation points $\Seqn x$ are chosen according to eq.
(14.4-10), then obviously the Richardson extrapolation scheme, eq.
(6.1-5), is exact for the following model sequence:
$$
s_n \; = \; s \, + \, \sum_{j=0}^{k-1} \, c_j \,
[(\beta + n) \, \Delta s_{n-1}]^{j+1} \, , \qquad k, n \in \N_0
\, , \quad \beta \in \R_{+} \, .
\tag
$$

It was remarked previously that if the Richardson extrapolation scheme,
eq. (6.1-5), is able to handle a certain problem using the
interpolation points $\Seqn x$, then the extrapolation points $\Seqn
\xi$ with $\xi_n = 1/x_n$ should be used if the same problem is to be
treated by the general forms of Wynn's $\rho$ algorithm, eq. (6.2-2),
and of its iteration ${\cal W}_k^{(n)}$, eq. (6.3-3). Hence, the
appropriate interpolation points $\Seqn \xi$ for $\rho_k^{(n)}$ and
${\cal W}_k^{(n)}$ would be
$$
\xi_n \; = \; 1 \, / \, [(\beta + n) \, \Delta s_{n-1}] \, ,
\qquad n \in \N_0 \, , \quad \beta \in \R_{+} \, .
\tag
$$
\beginFloat

\medskip

\beginTabelle [to \kolumnenbreite]
\beginFormat \rechts " \mitte " \mitte " \mitte " \mitte
\endFormat
\+ " \links {\bf Table 14-4} \@ \@ \@ \@ " \\
\+ " \links {Acceleration of the series (14.3-7) for the lemniscate
constant $A$} \@ \@ \@ \@ " \\
\- " \- " \- " \- " \- " \- " \\ \sstrut {} {1.5 \jot} {1.5 \jot}
\+ " \rechts {$n$} " \mitte {partial sum $s_n$} " ${\cal N}_n^{(0)}
(s_0, x_0)$ " $\rho_{2 \Ent {n/2}}^{(n - 2 \Ent {n/2})}$ " ${\cal
W}_{\Ent {n/2}}^{(n - 2 \Ent {n/2})}$ " \\
\+ " " \mitte {eq. (14.4-2)} " eq. (6.1-5) " eq. (6.2-2) " eq. (6.3-3)
" \\
\+  " " " $x_n = (n+1) a_n$ " $\xi_n = 1/[(n+1) a_n]$ " $\xi_n =
1/[(n+1) a_n]$ " \\
\- " \- " \- " \- " \- " \- " \\ \sstrut {} {1 \jot} {1 \jot}
\+ "   10  "    1.2250162047862     "   1.3106718871541  "
1.3108757199219  "   1.3110586974716 " \\
\+ "   11  "    1.2287537180105     "   1.3109702421830  "
1.3110068620892  "   1.3110406715649 " \\
\+ "   12  "    1.2320431110268     "   1.3110807598150  "
1.3110296350804  "   1.3110524519639 " \\
\+ "   13  "    1.2349672811610     "   1.3110705535246  "
1.3110300492139  "   1.3110491470884 " \\
\+ "   14  "    1.2375891404731     "   1.3110388651127  "
1.3110297665932  "   1.3110531590519 " \\
\+ "   15  "    1.2399574101139     "   1.3110247934048  "
1.3110284817234  "   1.3110289425634 " \\
\+ "   16  "    1.2421104860230     "   1.3110247919628  "
1.3110287561737  "   1.3110288767162 " \\
\+ "   17  "    1.2440790912340     "   1.3110277629872  "
1.3110287810924  "   1.3110291518189 " \\
\+ "   18  "    1.2458881405432     "   1.3110291063652  "
1.3110287803056  "   1.3110289945504 " \\
\+ "   19  "    1.2475580797723     "   1.3110291052833  "
1.3110287815108  "   1.3110288722980 " \\
\+ "   20  "    1.2491058660392     "   1.3110288500150  "
1.3110287770120  "   1.3110288696319 " \\
\+ "   21  "    1.2505456974656     "   1.3110287479173  "
1.3110287771445  "   1.3110288718241 " \\
\+ "   22  "    1.2518895646423     "   1.3110287537771  "
1.3110287771529  "   1.3110288775601 " \\
\+ "   23  "    1.2531476731141     "   1.3110287734917  "
1.3110287771494  "   1.3110288714005 " \\
\+ "   24  "    1.2543287710802     "   1.3110287796054  "
1.3110287771460  "   1.3110288733834 " \\
\+ "   25  "    1.2554404064530     "   1.3110287785488  "
1.3110287771461  "   1.3110288715363 " \\
\- " \- " \- " \- " \- " \- " \\ \sstrut {} {1 \jot} {1 \jot}
\+ " \links {$[\Gamma (1/4)]^2 / [4 (2 \pi)^{1/2}]$} \@ "
1.3110287771461 "   1.3110287771461 "   1.3110287771461 " \\
\- " \- " \- " \- " \- " \- " \\ \sstrut {} {1 \jot} {1 \jot}

\endTabelle

\medskip

\endFloat

In table 14-4 the partial sums (14.4-2) are accelerated by the
Richardson extrapolation scheme, eq. (6.1-5), and by the general forms
of Wynn's $\rho$ algorithm, eq. (6.2-2), and of its iteration ${\cal
W}_k^{(n)}$, eq. (6.3-3). In the case of the Richardson extrapolation
scheme the interpolation points $\Seqn x$ were chosen according to eq.
(14.4-10) with $\beta = 1$, and in the case of $\rho_k^{(n)}$ and
${\cal W}_k^{(n)}$ the interpolation points $\Seqn \xi$ were chosen
according to eq. (14.4-12) with $\beta = 1$. If we compare tables 14-3
and 14-4, we find that the rate of convergence of the transforms is
slower in table 14-4, but otherwise, the results are quite similar.

As usual, table 14-4 was produced in QUADRUPLE PRECISION. When the same
computation was repeated in DOUBLE PRECISION, it was again observed
that the Richardson extrapolation scheme is much more sensitive to
rounding errors than the other two transformations. For $n = 15$ the
Richardson extrapolation scheme achieved a relative accuracy of 6
decimal digits, and for larger values of $n$ the accuracy deteriorated
rapidly yielding nonsensical results for $n \ge 22$. Of the other two
transformations in table 14-4, ${\cal W}_k^{(n)}$, eq. (6.3-3), was
this time the numerically more stable transformation. Wynn's $\rho$
algorithm, eq. (6.2-2), achieved for $n \ge 19$ a relative accuracy for
7 decimal digits, whereas ${\cal W}_k^{(n)}$, eq. (6.3-3), achieved for
$n \ge 18$ a relative accuracy of 8 decimal digits.

Next, the acceleration of the convergence of the other two test series
(14.3-10) and (14.3-20) will be considered. The infinite series
(14.3-20) may be considered to be a special case of the infinite series
(14.3-10) since it was derived from it by performing the limit $z \to
0$. Because of eq. (14.3-18) we expect that these two series (14.3-10)
and (14.3-20) should have roughly the same convergence properties.
However, in convergence acceleration processes these two series differ
significantly. As remarked previously, Levin's $u$ and $v$
transformation, eqs. (7.3-5) and (7.3-11), respectively, Brezinski's
$\theta$ algorithm, eq. (10.1-9), and ${\cal J}_k^{(n)}$, eq. (10.3-6),
are all exact for the partial sums of the infinite series (14.3-20),
whereas no sequence transformation is known which is exact for the
sequence of partial sums
$$
s_n \; = \; \sum_{m=0}^n \,
{\widehat k}_{m-1/2} (z) \, / \, [2^m \, m!] \, ,
\qquad n \in \N_0 \, , \quad z \in \R_{+} \, ,
\tag
$$
of the infinite series (14.3-10). In the case of those sequence
transformations, which are not exact for the series (14.3-20), it was
observed quite consistently that the series expansion (14.3-20) can be
accelerated more easily than the series expansion (14.3-10).
Consequently, we shall not consider explicitly the acceleration of the
convergence of the series (14.3-20). Instead, we shall focus our
attention on the acceleration of the convergence of the infinite series
(14.3-10) which is much more interesting in this context. The
acceleration of the convergence of the series expansion (14.3-10) is
particularly hard for larger values of $z$. In fact, for sufficiently
large values of $z$, every sequence transformation has so far been
brought down to its knees. This strong dependence of the success of a
convergence acceleration process on the the magnitude of the argument
$z$ make the series expansion (14.3-10) of $1/z$ in terms of reduced
Bessel functions a very interesting test problem.

\beginFloat

\medskip

\beginTabelle [to \kolumnenbreite]
\beginFormat \rechts " \mitte " \mitte " \mitte " \mitte
\endFormat
\+ " \links {\bf Table 14-5} \@ \@ \@ \@ " \\
\+ " \links {Acceleration of the series expansion (14.3-10) for $z =
4/5$} \@ \@ \@ \@ " \\
\- " \- " \- " \- " \- " \- " \\ \sstrut {} {1.5 \jot} {1.5 \jot}
\+ " \rechts {$n$} " \mitte {partial sum $s_n$}
" $\theta_{2 \Ent {n/3}}^{(n - 3 \Ent {n/3})})$
" ${\cal J}_{\Ent {n/3}}^{(n - 3 \Ent {n/3})}$
" $\lambda_{\Ent {n/2}}^{(n - 2 \Ent {n/2})}$ " \\
\+ " " \mitte {eq. (14.4-13)} " eq. (10.1-9) " eq. (10.3-6) " eq.
(11.2-1) " \\
\- " \- " \- " \- " \- " \- " \\ \sstrut {} {1 \jot} {1 \jot}
\+ "     7  "    1.0422312196170     "   1.2497381860187  "
1.2479931939358  "   1.2173888687023 " \\
\+ "     8  "    1.0550056275790     "   1.2499224049805  "
1.2494113420204  "   1.2172621135266 " \\
\+ "     9  "    1.0656857514131     "   1.2497936062022  "
1.2499952677367  "   1.2015462070872 " \\
\+ "    10  "    1.0747865667307     "   1.2500104486053  "
1.2499875622609  "   1.2083962702181 " \\
\+ "    11  "    1.0826618965033     "   1.2500117160875  "
1.2499952780892  "   1.2511499315352 " \\
\+ "    12  "    1.0895638413456     "   1.2500026110013  "
1.2499738050830  "   1.2501603786137 " \\
\+ "    13  "    1.0956774851981     "   1.2500122711095  "
1.2499996707633  "   1.2500076445011 " \\
\+ "    14  "    1.1011421634246     "   1.2500120099841  "
1.2499999749148  "   1.2500005196325 " \\
\+ "    15  "    1.1060650318428     "   1.2500122009821  "
1.2499999979586  "   1.2500007343758 " \\
\+ "    16  "    1.1105300244656     "   1.2500168769473  "
1.2500000006228  "   1.2500006080295 " \\
\+ "    17  "    1.1146039429560     "   1.2500002080018  "
1.2499999999894  "   1.2500010863749 " \\
\+ "    18  "    1.1183407028756     "   1.2500032915285  "
1.2499999999664  "   1.2500008892374 " \\
\+ "    19  "    1.1217843613599     "   1.2500000030342  "
1.2499999999776  "   1.2499999946191 " \\
\+ "    20  "    1.1249713188304     "   1.2500000009966  "
1.2499999999625  "   1.2499999323380 " \\
\+ "    21  "    1.1279319483516     "   1.2500000107858  "
1.2499999999720  "   1.2499999648641 " \\
\+ "    22  "    1.1306918204772     "   1.2499999999866  "
1.2500000000029  "   1.2499999366280 " \\
\- " \- " \- " \- " \- " \- " \\ \sstrut {} {1 \jot} {1 \jot}
\+ " \links {exact} \@ "   1.2500000000000  "   1.2500000000000  "
1.2500000000000 " \\
\- " \- " \- " \- " \- " \- " \\ \sstrut {} {1 \jot} {1 \jot}

\endTabelle

\medskip

\endFloat

In table 14-5 the convergence of the partial sums (14.4-13) is
accelerated by Brezinski's $\theta$ algorithm, eq. (10.1-9), by its
iteration ${\cal J}_k^{(n)}$, eq. (10.3-6), and by $\lambda_k^{(n)}$,
eq. (11.2-1), with $\beta = 1$. As in tables 14-1 and 14-2 ${\cal
J}_k^{(n)}$ has a slight plus over Brezinski's $\theta$ algorithm. The
third transformation in table 14-5, $\lambda_k^{(n)}$, is clearly less
efficient than the other two.

As usual, table 14-5 was produced in QUADRUPLE PRECISION. When the same
computations were repeated in DOUBLE PRECISION, the loss of some
significant digits was again observed. Relatively insensitive to
rounding errors was ${\cal J}_k^{(n)}$, which for $n = 22$ reproduced
10 decimal digits. For $n = 22$ Brezinski's $\theta$ algorithm
reproduced 8 decimal digits, whereas $\lambda_k^{(n)}$ reproduced 7
decimal digits.

If we compare the results produced by ${\cal J}_k^{(n)}$, eq. (10.3-6),
in tables 14-2 and 14-5 we see that the convergence of the series
(14.3-7) for the lemniscate constant $A$ can apparently be accelerated
much more easily than the convergence of the series expansion (14.3-10)
for $1/z$ in terms of reduced Bessel functions. Similar results were
observed also in the case of other sequence transformations.

The third sequence transformation in table 14-5, $\lambda_k^{(n)}$, eq.
(11.2-1), was derived by modifying the recursive scheme (7.3-21) for
$\Lambda_k^{(n)} (\beta, s_n)$ along the lines of Brezinski's $\theta$
algorithm. There are some interesting differences between
$\lambda_k^{(n)}$, which may be considered to be an iterated weighted
$\Delta^2$ process, and $\Lambda_k^{(n)} (\beta, s_n)$, which
corresponds to the special case $x_n = 1/(\beta + n)$ of the Richardson
extrapolation scheme, with respect to their ability of accelerating
logarithmic convergence. The linear sequence transformation
$\Lambda_k^{(n)} (\beta, s_n)$ is one of the best accelerators for the
series (14.3-1) for $\zeta (2)$ but according to theorem 14-4 is not
able to accelerate the convergence of the test series (14.3-7),
(14.3-10), and (14.3-20). The nonlinear sequence transformation
$\lambda_k^{(n)}$ is clearly less efficient than $\Lambda_k^{(n)}
(\beta, s_n)$ in the case of the series for $\zeta (2)$, but is at
least moderately powerful in the case of the test series (14.3-7),
(14.3-10), and (14.3-20). This example shows once more that the
modification of the recursive scheme of a sequence transformation along
the lines of Brezinski's $\theta$ algorithm does not automatically lead
to a sequence transformation which is able to outperform the
transformation, from which it was derived, in all respects. However, it
is quite likely that the new transformation will be more versatile than
the transformation from which it was derived.

The greater flexibility of those sequence transformations, which are
derived along the lines of Brezinski's $\theta$ algorithm, is probably
responsible for their ability of accelerating the convergence of the
test series (14.3-7), (14.3-10), and (14.3-20), whose remainders all
behave like $n^{- 1/2}$ as $n \to \infty$. A sequence transformation
like the standard form of Wynn's $\rho$ algorithm, eq. (6.2-4), is
often able to achieve really spectacular results if the remainders of
the sequence to be transformed behave like an integral power of $1/n$,
but it fails completely if the remainders behave like a nonintegral
powers of $1/n$. In such a case, the general form of Wynn's $\rho$
algorithm, eq. (6.2-2), together with an appropriate set of
interpolation points $\Seqn x$ has to be used. However, one should not
expect that it will always be an easy task to find an appropriate set
of interpolation points $\Seqn x$.

\beginFloat

\medskip

\beginTabelle [to \kolumnenbreite]
\beginFormat \rechts " \mitte " \mitte " \mitte " \mitte
\endFormat
\+ " \links {\bf Table 14-6} \@ \@ \@ \@ " \\
\+ " \links {Acceleration of the series expansion (14.3-10) for $z =
4/5$} \@ \@ \@ \@ " \\
\- " \- " \- " \- " \- " \- " \\ \sstrut {} {1.5 \jot} {1.5 \jot}
\+ " \rechts {$n$} " \mitte {partial sum $s_n$}
" $u_n^{(0)} (1/2, s_0)$
" ${\cal L}_n^{(0)} (1/2, s_0, \omega_0)$
" ${\cal L}_n^{(0)} (1/2, s_0, \omega_0)$
\+ " " \mitte {eq. (14.4-13)} " eq. (7.3-5) " eq. (7.1-7) " eq. (7.1-7)
" \\
\+  " " " " $\omega_n = (n+1)^{- 1/2}$ " $\omega_n = (2 n - 1)!! / (2
n)!!$ " \\
\- " \- " \- " \- " \- " \- " \\ \sstrut {} {1 \jot} {1 \jot}
\+ "    7  "    1.0422312196170     "   1.2472807413200  "
1.2500112144531  "   1.2500305329006 " \\
\+ "    8  "    1.0550056275790     "   1.2519888543148  "
1.2499669957586  "   1.2499700063760 " \\
\+ "    9  "    1.0656857514131     "   1.2498098463228  "
1.2500047257200  "   1.2500032427111 " \\
\+ "   10  "    1.0747865667307     "   1.2498185973113  "
1.2500012644944  "   1.2500013433652 " \\
\+ "   11  "    1.0826618965033     "   1.2500735225247  "
1.2499994685087  "   1.2499995333305 " \\
\+ "   12  "    1.0895638413456     "   1.2499966179905  "
1.2500000284158  "   1.2500000125919 " \\
\+ "   13  "    1.0956774851981     "   1.2499942161928  "
1.2500000253626  "   1.2500000250149 " \\
\+ "   14  "    1.1011421634246     "   1.2500019544458  "
1.2499999930747  "   1.2499999939416 " \\
\+ "   15  "    1.1060650318428     "   1.2499999104639  "
1.2500000000686  "   1.2499999999124 " \\
\+ "   16  "    1.1105300244656     "   1.2499998754938  "
1.2500000003615  "   1.2500000003511 " \\
\+ "   17  "    1.1146039429560     "   1.2500000414919  "
1.2499999999164  "   1.2499999999261 " \\
\+ "   18  "    1.1183407028756     "   1.2499999971441  "
1.2499999999999  "   1.2499999999984 " \\
\+ "   19  "    1.1217843613599     "   1.2499999980000  "
1.2500000000043  "   1.2500000000042 " \\
\+ "   20  "    1.1249713188304     "   1.2500000007366  "
1.2499999999990  "   1.2499999999991 " \\
\+ "   21  "    1.1279319483516     "   1.2499999999235  "
1.2500000000000  "   1.2500000000000 " \\
\+ "   22  "    1.1306918204772     "   1.2499999999760  "
1.2500000000000  "   1.2500000000000 " \\
\- " \- " \- " \- " \- " \- " \\ \sstrut {} {1 \jot} {1 \jot}
\+ " \links {exact} \@ "   1.2500000000000  "   1.2500000000000  "
1.2500000000000 " \\
\- " \- " \- " \- " \- " \- " \\ \sstrut {} {1 \jot} {1 \jot}

\endTabelle

\medskip

\endFloat

If the convergence of a given sequence $\Seqn s$ is to be accelerated
by Levin's sequence transformation ${\cal L}_k^{(n)} (\beta, s_n,
\omega_n)$, eq. (7.1-7), together with one of the explicit remainder
estimates (7.3-4), (7.3-6), (7.3-8), and (7.3-10), then the elements of
the sequence $\Seqn s$ also supply the remainder estimates $\Seqn
\omega$. We can expect that Levin's sequence transformation will
produce good results if the the remainder estimates $\Seqn \omega$ can
be chosen in such a way that the ratio $(s_n - s)/\omega_n$ will depend
upon $n$ only quite weakly, i.e., if it is a constant apart from terms,
that are at least of order $O (n^{- 1})$ as $n \to \infty$ or smaller.
However, in practical applications it may happen that unless $n$ is
very large, the explicit remainder estimates (7.3-4), (7.3-6), (7.3-8),
and (7.3-10) yield only relatively bad approximations for the actual
remainders $\Seqn r$ of the sequence to be transformed. In such a case,
it is to be expected that Levin's sequence transformation ${\cal
L}_k^{(n)} (\beta, s_n, \omega_n)$, eq. (7.1-7), will be a relatively
weak sequence transformation if it uses one the the explicit remainder
estimates (7.3-4), (7.3-6), (7.3-8), and (7.3-10).

According to eq. (14.3-16) the reduced Bessel functions ${\widehat
k}_{n+1/2} (z)$ with $z \in \R_{+}$ and $n \in \N_0$ are positive and
bounded by their values at the origin. However, it was already remarked
in section 14.3 that due to the exponential decline of the reduced
Bessel functions ${\widehat k}_{n+1/2} (0)$ is a good approximation for
${\widehat k}_{n+1/2} (z)$ with a larger argument $z$ only if $n$ is
relatively large. This has some unpleasant consequences if for instance
Levin's $u$ transformation, eq. (7.3-5), is to be used for an
acceleration of the convergence of the series expansion (14.3-10) of
$1/z$ in terms of reduced Bessel functions. For larger values of $z$
the product $(\beta + n) \, {\widehat k}_{n-1/2} (z)$ with $\beta \in
\R_{+}$  will be a good approximation for the remainder
$$
r_n \; = \; \sum_{m=n+1}^{\infty} \,
{\widehat k}_{m-1/2} (z) \, / \, [2^m \, m!] \, ,
\qquad n \in \N_0 \, \quad z \in \R_{+} \, ,
\tag
$$
of the infinite series (14.3-10) only if $n$ is very large. In such a
situation, it should be worthwhile to look for other sets of remainder
estimates $\Seqn \omega$ even if Levin's sequence transformation ${\cal
L}_k^{(n)} (\beta, s_n, \omega_n)$, eq. (7.1-7), would then be a linear
sequence transformation. One simple possibility would be to choose
$$
\omega_n \; = \; (n + 1)^{- 1/2} \, ,
\qquad n \in \N_0 \, .
\tag
$$

Another possibility, which would also lead to remainder estimates that
behave like $n^{- 1/2}$ as $n \to \infty$, would be to choose the
remainder estimates according to eq. (14.3-21),
$$
\omega_n \; = \; (2 n - 1)!! \, / \, (2 n)!! \, ,
\qquad n \in \N_0 \, .
\tag
$$

In table 14-6 the partial sums (14.4-13) are accelerated by $u_k^{(n)}
(\beta, s_n)$, eq. (7.3-5), and by ${\cal L}_k^{(n)} (\beta, s_n,
\omega_n)$, eq. (7.1-7), with either $\omega_n = (n + 1)^{- 1/2}$ or
$\omega_n = (2 n - 1)!! / (2 n)!!$. In all cases $\beta = 1/2$ was used
which gives slightly better results than $\beta = 1$.

The results in table 14-6 show quite clearly that the remainder
estimates (14.4-15) and (14.4-16) produce significantly better results
than the remainder estimate (7.3-4) which is the basis of Levin's $u$
transformation. There is indirect evidence that this improved
convergence of the transforms is indeed due to the better approximation
of the remainders (14.4-14) by the remainder estimates (14.4-15) and
(14.4-16). If the convergence of the series (14.3-7) for the lemniscate
constant $A$ is accelerated by Levin's sequence transformation, eq.
(7.1-7), with the remainder estimates being chosen according to eq.
(14.4-15), then the results obtained in this way are as good or only
marginally better than the results obtained by Levin's $u$
transformation. Hence, in the case of the series (14.3-7) for the
lemniscate constant $A$ the remainder estimates (14.4-15) do not lead
to a spectacular improvement.

As usual, table 14-6 was produced in QUADRUPLE PRECISION. When the same
computation was repeated in DOUBLE PRECISION, Levin's $u$
transformation produced for $n = 15$ a relative accuracy of 8 decimal
digits. The other two transformations produced also for $n = 15$ a
relative accuracy of 10 decimal digits. For larger values of $n$, the
accuracy deteriorated again.

\endAbschnittsebene

\neueSeite

\Abschnitt Synopsis

\vskip - 2 \jot

\beginAbschnittsebene

\medskip

\Abschnitt General considerations

\smallskip

\aktTag = 0

In this report a large number of mainly nonlinear sequence
transformations for the acceleration of convergence and the summation
of divergent series were discussed. Some of those sequence
transformations as for instance Wynn's $\epsilon$ algorithm are well
established in the literature, while many others are new. The
properties of these sequence transformations were analyzed and
efficient algorithms for their evaluation were derived. In sections 13
and 14 the performance of these sequences transformations was tested by
applying them to certain slowly convergent and divergent series, which
are hopefully realistic models for a large part of the slowly
convergent or divergent series that can occur in scientific problems
and in applied mathematics.

It still has to be discussed how one should actually proceed if the
convergence of a slowly convergent sequence or series has to be
accelerated or if a divergent series has to be summed. In view of the
numerous different types of sequences and series, which can occur in
practical problems, and because of the large number of sequence
transformations, which are known, the selection of an appropriate
sequence transformation is certainly a nontrivial problem.

If the terms of the series, which is to be transformed, are known
analytically or if it is known how the elements of the sequence $\Seqn
s$ of partial sums behave as $n \to \infty$, it is normally
comparatively easy to find a suitable sequence transformation.
Unfortunately, it can happen that only a few elements of a slowly
convergent or divergent sequence $\Seqn s$ are available and that the
behaviour of the sequence elements $s_n$ as $n \to \infty$ is not
known. In such an unfavourable case, in which it is often not easy to
decide whether $\Seqn s$ converges at all, and if it does, whether it
converges linearly or logarithmically, the choice of an appropriate
sequence transformation is by no means simple and also of decisive
importance for the success of the whole approach.

It is well known that the performance of a sequence transformation
depends in most cases quite strongly upon the type of convergence of
the sequence to which it is applied. Apparently, there is no sequence
transformation which excells in every respect. Of all the sequence
transformations in this report, only Levin's $u$ transformation, eq.
(7.3-5), to a somewhat lesser extent also Levin's $v$ transformation,
eq. (7.3-11), Brezinski's $\theta$ algorithm, eq. (10.1-9), and its
iteration ${\cal J}_k^{(n)}$, eq. (10.3-6), are powerful accelerators
for both linear and logarithmic convergence and are also able to sum
efficiently even wildly divergent series. In all test cases considered
they were among the better sequence transformations.

It is tempting to believe that it would be sufficient to use only the
four sequence transformations mentioned above in situations, in which
apart from the numerical values of a few sequence elements little is
known. Since these sequence transformations are known to work well in a
variety of different situations, it seems reasonable to expect that
they will accomplish at least something. However, in many of the test
problems of this report other, less versatile sequence transformations
were actually more efficient. Hence, even if the four transformations
mentioned above are successful, they do not necessarily give the best
results, and it may well be worthwhile to look for other sequence
transformations which are possibly more efficient, in particular if
only relatively few sequence elements are available.

There are also some other aspects which should be taken into
consideration. In this report, only the most common types of sequences
and series were considered, i.e., either linearly and logarithmically
convergent sequences and series or alternating divergent series. This
does not exhaust all possibilities. Therefore, it is not certain
whether the four sequence transformations mentioned above will also be
able to handle successfully other types of convergence. For instance,
Smith and Ford report that in the case of some slowly convergent series
with terms having irregular sign patterns Wynn's $\epsilon$ algorithm,
eq. (4.2-1), clearly outperformed Levin's $u$ transformation, eq.
(7.3-5), and Brezinski's $\theta$ algorithm, eq. (10.1-9), which both
did not accomplish much (see table 5 on p. 484 of ref. [30]). Also, in
a situation, in which apart from the numerical values of only a few
sequence elements very little is known, it is often not clear whether
and how well the whole process has already converged. Even if a
sequence transformation produces a sequence of transforms which
apparently converges to some limit it cannot be excluded that this
convergence is an artifact. If two different sequence transformations
converge to the same limit, an artifact still cannot be ruled out but
it is much less likely. Consequently, in such a situation it should be
worthwhile to use more than a single sequence transformation.

In the opinion of the author the approximate determination of the limit
or antilimit $s$ of a slowly convergent or divergent sequence $\Seqn s$
is essentially an experimental problem which should be handled with
utmost care. The numerical evidence supplied by a single sequence
transformation is not necessarily sufficient, and it is usually a good
idea to compare the results produced by several sequence
transformations. In order to facilitate the task of selecting
appropriate sequence transformations, short r\'esum\'es of the properties
of all sequence transformations, which occur in this report, will now
be be given.

\medskip

\Abschnitt Wynn's epsilon algorithm and related transformations

\smallskip

\aktTag = 0

Wynn's $\epsilon$ algorithm, eq. (4.2-1), and its close relative,
Aitken's iterated $\Delta^2$ process, eq. (5.1-15), are both able to
accelerate linear convergence and to sum alternating divergent series
even if they diverge as wildly as the Euler series, eq. (1.1-7), but
they are not able to accelerate logarithmic convergence.

Practical experience and also some theoretical estimates indicate that
Wynn's $\epsilon$ algorithm is only a moderately powerful sequence
transformation, in particular if wildly divergent series must be
summed. In the case of the two Stieltjes series (13.3-9) for the
exponential integral and (13.4-3) for the logarithm, and also in other
tests not discussed in this report, Wynn's $\epsilon$ algorithm was not
only clearly outperformed by Levin's sequence transformation ${\cal
L}_k^{(n)} (\beta, s_n, \omega_n)$, eq. (7.1-7), and the new sequence
transformations ${\cal S}_k^{(n)} (\beta, s_n, \omega_n)$, eq. (8.2-7),
and ${\cal M}_k^{(n)} (\gamma, s_n, \omega_n)$, eq. (9.2-6), but also
frequently by Aitken's iterated $\Delta^2$ algorithm, albeit to a
lesser extent.

On the basis of these results it looks as if Wynn's $\epsilon$
algorithm should be dismissed. However, the real strength of Wynn's
$\epsilon$ algorithm is not its efficiency but its robustness. The
$\epsilon$ algorithm is remarkably insensitive to rounding errors and
can also tolerate input data which either have a low relative accuracy
or which behave in a comparatively irregular way. Due to its
robustness, Wynn's $\epsilon$ algorithm is often able to produce
meaningful and reliable results in situations in which other sequence
transformations, which are in principle much more powerful, fail. For
instance, in ref. [64] the convergence of some infinite series with
very complicated terms was accelerated by Wynn's $\epsilon$ algorithm
and by Levin's $u$ transformation, eq. (7.3-5). Since these infinite
series converge linearly, it was to be expected that Levin's $u$
transformation would do better than Wynn's $\epsilon$ algorithm.
However, it was found that the $\epsilon$ algorithm converged more
rapidly than the $u$ transformation. In addition, the $\epsilon$
algorithm was apparently not affected by numerical instabilities
whereas in the case of the $u$ transformation a dangerous accumulation
of rounding errors was observed (see pp. 3716 - 3717 of ref. [64]).

Superficially, Aitken's iterated $\Delta^2$ process appears to be a
better sequence transformation than Wynn's $\epsilon$ algorithm.
However, to a certain extent Aitken's iterated $\Delta^2$ process
combines the disadvantageous features of both Wynn's $\epsilon$
algorithm, which is only moderately powerful, and of the
transformations ${\cal L}_k^{(n)} (\beta, s_n, \omega_n)$, eq. (7.1-7),
${\cal S}_k^{(n)} (\beta, s_n, \omega_n)$, eq. (8.2-7), and ${\cal
M}_k^{(n)} (\gamma, s_n, \omega_n)$, eq. (9.2-6), which are not very
robust since they are powerful sequence transformations only if the
remainder estimates $\Seqn \omega$ are good approximations of the
actual remainders $\Seqn r$. Experience indicates that Aitken's
iterated $\Delta^2$ process is in general less efficient than ${\cal
L}_k^{(n)} (\beta, s_n, \omega_n)$, ${\cal S}_k^{(n)} (\beta, s_n,
\omega_n)$, and ${\cal M}_k^{(n)} (\gamma, s_n, \omega_n)$, and that it
is less robust and more susceptible to rounding errors than Wynn's
$\epsilon$ algorithm.

\medskip

\Abschnitt Wynn's rho algorithm and related transformations

\smallskip

\aktTag = 0

The properties of Wynn's $\epsilon$ algorithm, eq. (4.2-1), and of
Wynn's $\rho$ algorithm, eq. (6.2-2), are complementary. The $\rho$
algorithm is often a good or even very good accelerator for logarithmic
convergence but is unable to accelerate linear convergence or to sum
divergent series. But with respect to robustness, the $\epsilon$ and
the $\rho$ algorithm are very similar. Experience indicates that the
$\rho$ algorithm is in general less sensitive to rounding errors than
other sequence transformations which are also able to accelerate
logarithmic convergence. This is certainly no mean accomplishment, in
particular since the acceleration of logarithmic convergence is an
inherently unstable process.

The power of Wynn's $\rho$ algorithm, which is essentially an
intelligent way of computing and extrapolating to infinity an
interpolating rational function of the type of eq. (6.2-1), depends
decisively upon an appropriate choice of the interpolation points
$\Seqn x$. In this respect, the $\rho$ algorithm closely resembles the
Richardson extrapolation scheme, eq. (6.1-5), which is essentially an
efficient way of computing and extrapolating to zero an interpolating
polynomial of the type of eq. (6.1-3). According to theorem 14-4 the
linear sequence transformation $\Lambda_k^{(n)} (\beta, s_n)$, eq.
(7.3-20), which corresponds to the special case $x_n = 1/(\beta + n)$
of the Richardson extrapolation scheme, eq. (6.1-5), is only able to
accelerate the convergence of logarithmically convergent sequences if
the remainders of these sequences behave like integral powers of $1/n$
as $n \to \infty$.

The standard form of Wynn's $\rho$ algorithm, eq. (6.2-4), which
corresponds to the choice $x_n = \beta + n$ for the interpolation
points in eq. (6.2-2), is a very powerful accelerator for
logarithmically convergent sequences with remainders that behave like
integral powers of $1/n$ as $n \to \infty$. However, experience
indicates that the standard form of the $\rho$ algorithm is unable to
accelerate convergence if the remainders of the sequence to be
transformed behave like a nonintegral power of $1/n$. In section 14.4
it was shown that the general form of Wynn's $\rho$ algorithm, eq.
(6.2-2), is apparently able to accelerate the convergence of sequences
with remainders that behave like $n^{- 1/2}$ as $n \to \infty$ if the
interpolation points $\Seqn x$ are chosen in such a way that $x_n \sim
n^{1/2}$ as $n \to \infty$.

The iterated sequence transformation ${\cal W}_k^{(n)}$, eq. (6.3-3),
is also only able to accelerate logarithmic convergence. The power of
${\cal W}_k^{(n)}$ depends as in the case of Wynn's $\rho$ algorithm,
from which it was derived, strongly upon an appropriate choice of the
interpolation points $\Seqn x$. The numerical results presented in
section 14.4 indicate that ${\cal W}_k^{(n)}$ has similar properties as
Wynn's $\rho$ algorithm. The standard form of ${\cal W}_k^{(n)}$, eq.
(6.3-4), is apparently not able to accelerate the convergence of
sequences whose remainders behave like nonintegral powers of $1/n$.
However, the results in section 14.4 also show that the general form of
${\cal W}_k^{(n)}$, eq. (6.3-3), is apparently able to accelerate the
convergence of sequences with remainders that behave like $n^{- 1/2}$
as $n \to \infty$ if the interpolation points $\Seqn x$ are chosen in
such a way that $x_n \sim n^{1/2}$ as $n \to \infty$. It also seems
that ${\cal W}_k^{(n)}$ is relatively insensitive to rounding errors.
However, because of the limited experience with this transformation it
seems to be too early for a definite assessment of its merits as well
as its weaknesses. Further tests of this sequence transformation should
therefore be of interest.

The fact, that Wynn's $\rho$ algorithm and its iteration ${\cal
W}_k^{(n)}$ are only successful if an appropriate set of interpolation
points is used, severely limits the practical usefulness of these
transformations in situations in which only the numerical values of a
few sequence elements are known. In such a situation, it may be a good
idea to use a modification of a suggestion by Beleznay [117]. In this
approach, the interpolation points are chosen according to eq. (14.4-7)
and the free parameter $\alpha$ is optimized after the input of every
new sequence $s_n$. Another possibility would be to choose the
interpolation points according to eq. (14.4-12). However, these two
suggestions are not yet sufficiently tested and it seems to be too
early for a definite assessment of their practical usefulness.

\medskip

\Abschnitt Levin's sequence transformation and related transformations

\smallskip

\aktTag = 0

It is a typical feature of Levin's sequence transformations ${\cal
L}_k^{(n)} (\beta, s_n, \omega_n)$, eq. (7.1-7), and of the related
transformations ${\cal S}_k^{(n)} (\beta, s_n, \omega_n)$, eq. (8.2-7),
and ${\cal M}_k^{(n)} (\gamma, s_n, \omega_n)$, eq. (9.2-6), that they
not only require the sequence elements $s_n, s_{n+1}, \ldots, s_{n+k}$,
but also the remainder estimates $\omega_n, \omega_{n+1}, \ldots,
\omega_{n+k}$. This explicit incorporation of remainder estimates is
both the strength as well as the weakness of these sequence
transformations. If it is possible to find a sequence of remainder
estimates $\Seqn \omega$ that are good approximations of the remainders
$\Seqn r$ of the sequence to be transformed, then experience as well as
some theoretical estimates indicate that such a sequence transformation
is extremely powerful. If, however, a good sequence of remainder
estimates cannot be found, such a sequence transformation will probably
perform quite poorly.

Theoretical estimates as well as practical experience indicate that the
remainder estimates $\Seqn \omega$ should be chosen in such a way that
the ratios $(s_n - s)/\omega_n$ depend on $n$ only quite weakly and
approach a constant as $n \to \infty$. In practical applications
Levin's sequence transformation ${\cal L}_k^{(n)} (\beta, s_n,
\omega_n)$ has so far exclusively been used in connection with the
simple remainder estimates (7.3-4), (7.3-6), (7.3-8), and (7.3-10),
which can also be used in the case of the new transformations ${\cal
S}_k^{(n)} (\beta, s_n, \omega_n)$ and ${\cal M}_k^{(n)} (\gamma, s_n,
\omega_n)$. The remainder estimate (7.3-4) gives Levin's $u$
transformation, eq. (7.3-5), which is certainly one of the most
powerful and most versatile sequence transformations. It is a powerful
accelerator for both linear and logarithmic convergence and is able to
sum efficiently divergent alternating series. The remainder estimates
(7.3-6) and (7.3-8) give Levin's $t$ and $d$ transformations, eqs.
(7.3-7) and (7.3-9), respectively, which are powerful accelerators for
linear convergence and are able to sum efficiently divergent
alternating series. However, they are unable to accelerate logarithmic
convergence. The remainder estimate (7.3-10) gives Levin's $v$
transformation, eq. (7.3-11), which has similar properties as Levin's
$u$ transformation. Finally, the remainder estimate $\omega_n =
1/(\beta + n)$ gives the linear sequence transformation
$\Lambda_k^{(n)} (\beta, s_n)$, eq. (7.3-20), which can also be
obtained from the Richardson extrapolation scheme, eq. (6.1-5), by
choosing $x_n = 1/(\beta + n)$. It is able to accelerate logarithmic
convergence if the remainders of the sequence to be transformed behave
like an integral power of $1/n$ as $n \to \infty$.

The simple remainder estimates (7.3-4), (7.3-6), (7.3-8), and (7.3-10)
are essentially asymptotic in nature because they were derived using
some simplifications which are valid for large values of $n$. However,
in convergence acceleration or summation processes it is tried to
approximate the limit or antilimit of a sequence $\Seqn s$ using only
the information stored in the first few sequence elements, say $s_0,
s_1, \ldots, s_m$, with $m$ being a relatively small number. Therefore,
it is by no means clear whether the simple remainder estimates (7.3-4),
(7.3-6), (7.3-8), and (7.3-10) lead to good approximations of the
actual remainders for only moderately large or even small indices.

There is some evidence that the simple remainder estimates (7.3-4),
(7.3-6), (7.3-8), and (7.3-10) lead to efficient sequence
transformations if the terms $a_n$ of the series, which is to be
accelerated or summed, approach their asymptotic limits relatively
fast. Let us for instance assume that $a_n$ is the term of a power
series in $z$ and that $a_n$ behaves like $n^\alpha z^n$ as $n \to
\infty$. Then, the sequence transformations ${\cal L}_k^{(n)} (\beta,
s_n, \omega_n)$, ${\cal S}_k^{(n)} (\beta, s_n, \omega_n)$, and ${\cal
M}_k^{(n)} (\gamma, s_n, \omega_n)$ should work well in combination
with one of the simple remainder estimates (7.3-4), (7.3-6), (7.3-8),
and (7.3-10) if the leading term $n^\alpha z^n$ is a good approximation
for $a_n$ already for moderately large or even small values of $n$. If
this is not the case, the simple remainder estimates (7.3-4), (7.3-6),
(7.3-8), and (7.3-10) will probably not work particularly well.

The infinite series (14.3-10) is a good example for the complications
which can occur in this context. It was remarked earlier that the
reduced Bessel functions, eq. (14.3-11), approach their asymptotic
limits according to eq. (14.3-18) quite slowly. This slow approach
decreases the efficiency of Levin's $u$ transformation, eq. (7.3-5),
considerably. The results presented in table 14-6 show that in this
case it is advantageous to use other, explicit remainder estimates,
which are not obtained from the elements of the sequence $\Seqn s$ to
be transformed, even if Levin's sequence transformation is then a
linear sequence transformation. Unfortunately, such an approach is only
possible if the remainders are known analytically and if simple and yet
good approximations for the remainders can be derived.

Also, under unfavourable circumstances the simple remainder estimates
(7.3-4), (7.3-6), (7.3-8), and (7.3-10) may have a detrimental effect
on the robustness of the sequence transformations ${\cal L}_k^{(n)}
(\beta, s_n, \omega_n)$, ${\cal S}_k^{(n)} (\beta, s_n, \omega_n)$, and
${\cal M}_k^{(n)} (\gamma, s_n, \omega_n)$. The elements of the
sequence $\Seqn s$ are not only input data, but they are also used to
compute the remainder estimates $\Seqn \omega$. Consequently, the
elements of $\Seqn s$ induce two fundamentally different kinds of
errors. More or less inevitable are the errors due to the limited
accuracy of the input data. However, the elements of the sequence
$\Seqn s$ induce also potentially large errors among the remainder
estimates $\Seqn \omega$, either because they are not accurate enough
or because they deviate too much from their asymptotic limits and
therefore produce bad remainder estimates. The worst scenario, which
can be imagined in this context, would be that the terms $a_n$ of a
series are not very accurate and that the terms approach their limiting
expressions only quite slowly and in an irregular fashion. In such a
situation the sequence transformations ${\cal L}_k^{(n)} (\beta, s_n,
\omega_n)$, ${\cal S}_k^{(n)} (\beta, s_n, \omega_n)$, and ${\cal
M}_k^{(n)} (\gamma, s_n, \omega_n)$ will be in trouble and it is likely
that Wynn's $\epsilon$ algorithm, although in principle only moderately
powerful, will produce better results.

Levin's sequence transformation ${\cal L}_k^{(n)} (\beta, s_n,
\omega_n)$, eq. (7.1-7), is based upon the assumption that the ratio
$(s_n - s)/\omega_n$ can be approximated by a polynomial in $1/(\beta +
n)$, whereas the new sequence transformations ${\cal S}_k^{(n)} (\beta,
s_n, \omega_n)$, eq. (8.2-7), and ${\cal M}_k^{(n)} (\gamma, s_n,
\omega_n)$, eq. (9.2-6), were derived assuming that the ratio $(s_n -
s)/\omega_n$ can be approximated by truncated factorial series or
related expressions. Since power series and factorial series have
different properties, it is not surprising that the new sequence
transformations and Levin's transformation behave differently in
convergence acceleration and summation processes. With respect to the
acceleration of linear convergence the new sequence transformations
${\cal S}_k^{(n)} (\beta, s_n, \omega_n)$ and ${\cal M}_k^{(n)}
(\gamma, s_n, \omega_n)$ are approximately as efficient as Levin's
sequence transformation. The new sequence transformations are
particularly well suited to sum wildly divergent alternating series such
as the Euler series, eq. (1.1-7). In that respect, they are usually at
least as good as Levin's sequence transformation and often they are
even clearly better. However, the new transformations perform quite
poorly if logarithmic convergence is to be accelerated. Also, the
linear transformations ${\cal F}_k^{(n)} (\alpha, s_n)$, eq. (8.4-11),
and ${\cal P}_k^{(n)} (\gamma, s_n)$, eq.(9.4-11), are much less
efficient than their analogue ${\Lambda}_k^{(n)} (\gamma, s_n)$, eq.
(7.3-20).

The stability properties of the sequence transformations ${\cal
L}_k^{(n)} (\beta, s_n, \omega_n)$, ${\cal S}_k^{(n)} (\beta, s_n,
\omega_n)$, and ${\cal M}_k^{(n)} (\gamma, s_n, \omega_n)$ depend very
much upon the sequence which is to be transformed and upon the
remainder estimates being used. However, at least some statements of a
more general nature, which are based upon experience, can be made. It
seems that the transformation of both convergent and divergent
alternating series is in general remarkably stable. Also, the
acceleration of linear convergence usually poses no particular
stability problems. The acceleration of logarithmic convergence is
always a problem which may easily lead to a serious loss of accuracy.
But it cannot be said that Levin's sequence transformation is more
sensitive to rounding errors than most other sequence transformations.
Also, it is probably safe to say that the sequence transformations
${\cal L}_k^{(n)} (\beta, s_n, \omega_n)$, ${\cal S}_k^{(n)} (\beta,
s_n, \omega_n)$, and ${\cal M}_k^{(n)} (\gamma, s_n, \omega_n)$ are in
general more efficient and at the same time less robust than Wynn's
$\epsilon$ algorithm if linear convergence is accelerated or if
divergent alternating series are summed. In the same way, Wynn's $\rho$
algorithm is apparently more robust than the $u$ and $v$ transformation
or the linear sequence transformation $\Lambda_k^{(n)} (\beta, s_n)$ in
the case of logarithmic convergence.

Finally, there is Drummond's sequence transformation ${\cal D}_k^{(n)}
(s_n, \omega_n)$, eq. (9.5-4), which is another relative of Levin's
sequence transformation since it also uses a sequence of remainder
estimates $\Seqn \omega$. Drummond's sequence transformation is very
important theoretically, in particular in connection with
Germain-Bonne's formal theory of convergence acceleration [33] and the
explicit construction of Pad\'e approximants for the Euler series, eq.
(1.1-7). However, in practical applications Drummond's sequence
transformations is at most moderately powerful. It is significantly
less powerful than the sequence transformations ${\cal L}_k^{(n)}
(\beta, s_n, \omega_n)$, eq. (7.1-7), ${\cal S}_k^{(n)} (\beta, s_n,
\omega_n)$, eq. (8.2-7), and ${\cal M}_k^{(n)} (\gamma, s_n,
\omega_n)$, eq. (9.2-6), but has the same weaknesses as these
transformations.

\medskip

\Abschnitt Brezinski's theta algorithm and related transformations

\smallskip

\aktTag = 0

Brezinski's $\theta$ algorithm, eq. (10.1-9), and its iteration ${\cal
J}_k^{(n)}$, eq. (10.3-6), combine many of the advantageous features
of Wynn's $\epsilon$ algorithm, eq. (4.2-1), and of Wynn's $\rho$
algorithm, eq. (6.2-2). They are able to accelerate linear convergence
and to sum even wildly divergent alternating series, and they are also
able to accelerate logarithmic convergence.

In those tests, in which linear convergence had to be accelerated or
divergent alternating series had to be summed, Brezinski's $\theta$
algorithm and its iteration ${\cal J}_k^{(n)}$ were usually better than
the $\epsilon$ algorithm, but less powerful than the sequence
transformations ${\cal L}_k^{(n)} (\beta, s_n, \omega_n)$, eq. (7.1-7),
${\cal S}_k^{(n)} (\beta, s_n, \omega_n)$, eq. (8.2-7), and ${\cal
M}_k^{(n)} (\gamma, s_n, \omega_n)$, eq. (9.2-6).

With respect to the acceleration of logarithmic convergence Brezinski's
$\theta$ algorithm and its iteration ${\cal J}_k^{(n)}$ are more
reliable than the standard form of Wynn's $\rho$ algorithm, eq.
(6.2-4), since they are not restricted to sequences with remainders
that behave like an integral power of $1/n$ as $n \to \infty$, and they
are easier to use than the general form of Wynn's $\rho$ algorithm, eq.
(6.2-2), since no interpolation points are needed. In those tests, in
which logarithmic convergence had to be accelerated, they were
approximately as powerful as Levin's $u$ and $v$ transformations.

It also seems that with respect to numerical stability and robustness
Brezinski's $\theta$ algorithm and its iteration ${\cal J}_k^{(n)}$ are
less robust than the $\epsilon$ or the $\rho$ algorithm and also more
susceptible to rounding errors.

The other sequence transformations, which were also derived along the
lines of Brezinski's $\theta$ algorithm are ${\cal B}_k^{(n)}$, eq.
(11.1-5), ${\cal C}_k^{(n)}$, eq. (11.1-12), ${\lambda}_k^{(n)}$, eq.
(11.2-1), ${\sigma}_k^{(n)}$, eq. (11.2-2), and ${\mu}_k^{(n)}$, eq.
(11.2-3). It is a typical feature of these transformations that they
are much more versatile than the transformations from which they were
derived. This means they are all able to accelerate linear and
logarithmic convergence and are also able to sum even wildly divergent
alternating series.

Unfortunately, it is also a typical feature of these sequence
transformations that their performance in the numerical tests described
in sections 13 and 14 was quite inconsistent and more or less
unpredictable. For instance, ${\sigma}_k^{(n)}$ turned out to be a very
powerful accelerator for the infinite series (14.3-1) for $\zeta (2)$,
a powerful accelerator for the infinite series (14.3-20), but a
relatively weak accelerator for the infinite series (14.3-7) and
(14.3-10). No explanation for this inconsistent behaviour can be given.
At best, the sequence transformations listed above were as good as
Brezinski's $\theta$ algorithm or its iteration ${\cal J}_k^{(n)}$, but
in most cases they were significantly weaker. Also, it seems that the
sequence transformations mentioned above are not more robust and less
susceptible to rounding errors than the $\theta$ algorithm or ${\cal
J}_k^{(n)}$. Hence, it seems that the most promising choices among all
sequence transformations, which were derived along the lines of
Brezinski's $\theta$ algorithm, are the $\theta$ algorithm, eq.
(10.1-9), and ${\cal J}_k^{(n)}$, eq. (10.3-6).

\bigskip

\Ueberschrift \large Acknowledgement

\medskip

I would like to thank Professor E. O. Steinborn for stimulating
discussions, for his encouragement, for his constant support, and for
the excellent working conditions at the Institut f\"ur Physikalische und
Theoretische Chemie der Universit\"at Regensburg.

The research, which ultimately led to this report, was begun during a
stay at the Faculty of Mathematics of the University of Waterloo,
Ontario, Canada. I would like to thank Professor J. {\v C}{\' \i}{\v
z}ek and Professor J. Paldus for their invitation to work with them as
a postdoctoral fellow in the Quantum Theory Group of the Department of
Applied Mathematics. Their hospitality, their generosity, and the
inspiring atmosphere which they provided is highly appreciated. Special
thanks to Professor J. {\v C}{\' \i}{\v z}ek who aroused my interest in
asymptotic techniques, the summation of divergent series, and the
construction of rational approximants.

Many thanks also to H. Homeier, who patiently read and discussed the
numerous preliminary versions of this manuscript, and to M. Middleton,
the \TeX \ expert at the Universit\"at Regensburg and the author of many
useful macros, who helped me to typeset this manuscript in \RzTeX, the
local variant of \PCTeX.

\endAbschnittsebene

\neueSeite

\Ueberschrift \large References

\Inhaltspez {\protect\bigbreak}

\Inhaltszeile {\large References}

\bigskip

\item {[1]} S. E. Haywood and J. D. Morgan III, {\it Discrete basis-set
approach for calculating Bethe logarithms}, Phys. Rev. A {\bf 32}
(1985), 3179 -- 3186.

\item {[2]} C. M. Bender and S. A. Orszag, {\it Advanced mathematical
methods for scientists and engineers} (McGraw-Hill, New York, 1978).

\item {[3]} C. M. Bender and T. T. Wu, {\it Anharmonic oscillator. II.
A study in perturbation theory in large order}, Phys. Rev. D {\bf 7}
(1973), 1620 -- 1636.

\item {[4]} K. Knopp, {\it Theorie and Anwendung der unendlichen
Reihen} (Springer-Verlag, Berlin, 1964).

\item {[5]} J. Stirling, {\it Methodus differentialis sive tractatus
de summatione et interpolatione serium infinitarum} (London, 1730).
English translation by F. Holliday, {\it The differential method, or, a
treatise concerning the summation and interpolation of infinite series}
(London, 1749).

\item {[6]} L. Euler, {\it Institutiones calculi differentialis cum
eius usu in analysi finitorum ac doctrina serium. Part II.1. De
transformatione serium} (Academia Imperialis Scientiarum Petropolitana,
1755). This book was reprinted as Vol. X of {\it Leonardi Euleri
Opera Omnia, Seria Prima} (Teubner, Leipzig and Berlin, 1913).

\item {[7]} G. H. Hardy, {\it Divergent series} (Oxford University
Press, Oxford, 1949).

\item {[8]} G. M. Petersen, {\it Regular matrix transformations}
(McGraw-Hill, London, 1966).

\item {[9]} A. Peyerimhoff, {\it Lectures on summability}
(Springer-Verlag, Berlin, 1969).

\item {[10]} K. Zeller and W. Beekmann, {\it Theorie der
Limitierungsverfahren} (Springer-Verlag, Berlin, 1970).

\item {[11]} R. E. Powell and S. M. Shah, {\it Summability theory and
its applications} (Van Nostrand Reinhold, London, 1972).

\item {[12]} A. C. Aitken, {\it On Bernoulli's numerical solution of
algebraic equations}, Proc. Roy. Soc. Edinburgh {\bf 46} (1926) , 289 -
305.

\item {[13]} J. Todd, {\it Motivation for working in numerical
analysis}, in J. Todd (ed.), {\it Survey of numerical analysis}
(McGraw-Hill, New York, 1962), pp. 1 -26.

\item {[14]} E. E. Kummer, {\it Eine neue Methode, die numerischen
Summen langsam convergirender Reihen zu berechnen}, J. Reine Angew.
Math. {\bf 16} (1837), 206 - 214.

\item {[15]} D. Shanks, {\it Non-linear transformations of divergent
and slowly convergent sequences}, J. Math. and Phys. (Cambridge, Mass.)
{\bf 34} (1955), 1 - 42.

\item {[16]} P. Wynn, {\it On a device for computing the $e_m (S_n)$
transformation}, Math. Tables Aids Comput. {\bf 10} (1956), 91 - 96.

\item {[17]} J. R. Schmidt, {\it On the numerical solution of linear
simultaneous equations by an iterative method}, Philos. Mag. {\bf 32}
(1941), 369 - 383.

\item {[18]} G. A. Baker, Jr., {\it Essentials of Pad\'e approximants}
(Academic Press, New York, 1975).

\item {[19]} C. Brezinski, {\it Acc\'el\'eration de la convergence en
analyse num\'erique} (Springer-Verlag, Berlin, 1977).

\item {[20]} C. Brezinski, {\it Algorithmes d'acc\'el\'eration de la
convergence -- \'Etude num\'erique} (Editions Technip, Paris, 1978).

\item {[21]} C. Brezinski, {\it Pad\'e-type approximation and general
orthogonal polynomials} (Birkh\"auser Verlag, Basel, 1980).

\item {[22]} G. A. Baker, Jr., and P. Graves-Morris, {\it Pad\'e
approximants. Part I: Basic theory. Part II: Extensions and
applications} (Addison-Wesley, Reading, Mass., 1981).

\item {[23]} J. Wimp, {\it Sequence transformations and their
applications} (Academic Press, New York, 1981).

\item {[24]} C. Brezinski, {\it Convergence acceleration methods: The
past decade}, J. Comput. Appl. Math. {\bf 12} \& {\bf 13} (1985), 19 -
36.

\item {[25]} P. Wynn, {\it On a Procrustean technique for the numerical
transformation of slowly convergent sequences and series}, Proc. Camb.
Phil. Soc. {\bf 52} (1956), 663 - 671.

\item {[26]} C. Brezinski, {\it Acc\'el\'eration de suites \`a convergence
logarithmique}, C. R. Acad. Sc. Paris {\bf 273} (1971), 727 - 730.

\item {[27]} C. Brezinski, {\it Some new convergence acceleration
methods}, Math. Comput. {\bf 39} (1982), 133 - 145.

\item {[28]} D. Levin, {\it Development of non-linear transformations
for improving convergence of sequences}, Int. J. Comput. Math. B {\bf
3} (1973), 371 - 388.

\item {[29]} D. A. Smith and W. F. Ford, {\it Acceleration of linear
and logarithmic convergence}, SIAM J. Numer. Anal. {\bf 16} (1979), 223
- 240.

\item {[30]} D. A. Smith and W. F. Ford, {\it Numerical comparisons of
nonlinear convergence accelerators}, Math. Comput. {\bf 38} (1982), 481
- 499.

\item {[31]} C. Brezinski, {\it A general extrapolation algorithm},
Numer. Math. {\bf 35} (1980), 175 - 180.

\item {[32]} T. H{\aa}vie, {\it Generalized Neville type extrapolation
schemes}, BIT {\bf 19} (1979), 204 - 213.

\item {[33]} B. Germain-Bonne, {\it Transformations de suites}, Rev.
Fran\c caise Automat. Informat. Rech. Operat. {\bf 7} (R-1) (1973), 84 -
90.

\item {[34]} W. Magnus, F. Oberhettinger, and R. P. Soni, {\it Formulas
and theorems for the special functions of mathematical physics}
(Springer-Verlag, New York, 1966).

\item {[35]} W. D. Clark, H. L. Gray, and J. E. Adams, {\it A note on
the T-transformation of Lubkin}, J. Res. Natl. Bur. Stand. B {\bf 73}
(1969), 25 - 29.

\item {[36]} C. Brezinski, {\it Algorithm 585: A subroutine for the
general interpolation and extrapolation problems}, ACM Trans. Math.
Software {\bf 8} (1982), 290 - 301.

\item {[37]} P. Wynn, {\it On the convergence and the stability of the
epsilon algorithm}, SIAM J. Numer. Anal. {\bf 3} (1966), 91 - 122.

\item {[38]} P. Wynn, {\it A note on programming repeated applications
of the $\epsilon$-algorithm}, R.F.T.I. - Chiffres {\bf 8} (1965), 23 -
62.

\item {[39]} P. Wynn, {\it Singular rules for certain non-linear
algorithms}, BIT {\bf 3} (1963), 175 - 195.

\item {[40]} S. Lubkin, {\it A method of summing infinite series}, J.
Res. Natl. Bur. Stand. {\bf 48} (1952), 228 - 254.

\item {[41]} R. R. Tucker, {\it The $\delta^2$ process and related
topics}, Pacif. J. Math. {\bf 22} (1967), 349 - 359.

\item {[42]} R. R. Tucker, {\it The $\delta^2$ process and related
topics II}, Pacif. J. Math. {\bf 28} (1969), 455 - 463.

\item {[43]} F. Cordellier, {\it Sur la r\'egularit\'e des proc\'ed\'es
$\delta^2$ d'Aitken et W de Lubkin}, in L. Wuytack (ed.), {\it Pad\'e
approximation and its applications} (Springer-Verlag, Berlin, 1979),
pp. 20 - 35.

\item {[44]} G. E. Bell and G. M. Phillips, {\it Aitken acceleration of
some alternating series}, BIT {\bf 24} (1984), 70 - 77.

\item {[45]} A. J. MacLeod, {\it Acceleration of vector sequences by
multi-dimensional $\Delta^2$-methods}, Commun. Appl. Numer. Meth. {\bf
2} (1986), 385 - 392.

\item {[46]} J. E. Drummond, {\it Summing a common type of slowly
convergent series of positive terms}, J. Austral. Math. Soc. {\bf B 19}
(1976), 416 - 421.

\item {[47]} P. Bj{\o}rstad, G. Dahlquist, and E. Grosse, {\it
Extrapolations of asymptotic expansions by a modified Aitken
$\delta^2$-formula}, BIT {\bf 21} (1981), 56 - 65.

\item {[48]} W. H. Press, B. P. Flannery, S. A. Teukolsky, and W. T.
Vetterling, {\it Numerical recipes} (Cambridge University Press,
Cambridge, 1986).

\item {[49]} P. Hillion, {\it M\'ethode d'Aitken it\'er\'ee pour les suites
oscillantes d'approximations}, C. R. Acad. Sc. Paris A {\bf 280}
(1975), 1701 - 1704.

\item {[50]} P. J. Davis, {\it Interpolation and approximation} (Dover,
New York, 1975).

\item {[51]} D. C. Joyce, {\it Survey of extrapolation processes in
numerical analysis}, SIAM Rev. {\bf 13} (1971), 435 - 490.

\item {[52]} E. H. Neville, {\it Iterative interpolation}, J. Indian
Math. Soc. {\bf 20} (1934), 87 - 120.

\item {[53]} L. F. Richardson, {\it The deferred approach to the limit.
I. Single lattice}, Phil. Trans. Roy. Soc. London A {\bf 226} (1927),
229 - 349.

\item {[54]} A. Cuyt and L. Wuytack, {\it Nonlinear methods in
numerical analysis} (North-Holland, Amsterdam, 1987).

\item {[55]} T. N. Thiele, {\it Interpolationsrechnung} (Teubner,
Leipzig, 1909).

\item {[56]} A. Sidi, {\it Convergence properties of some nonlinear
sequence transformations}, Math. Comput. {\bf 33} (1979), 315 - 326.

\item {[57]} T. Fessler, W. F. Ford, and D. A. Smith, {\it HURRY: An
acceleration algorithm for scalar sequences and series}, ACM Trans.
Math. Software {\bf 9 } (1983), 346 - 354.

\item {[58]} I. M. Longman, {\it Difficulties of convergence
acceleration}, in M. G. de Bruin and H. van Rossum (eds.), {\it Pad\'e
approximation and its applications Amsterdam 1980} (Springer-Verlag,
Berlin, 1981), pp. 273 - 289.

\item {[59]} Sister M. C. Fasenmyer, {\it Some generalized
hypergeometric polynomials}, Bull. Amer. Math. Soc. {\bf 53} (1947),
806 - 812. A good discussion of Sister Celine's technique can also be
found in chapter 14 of E. D. Rainville, {\it Special functions}
(Chelsea, New York, 1960).

\item {[60]} A. J. Thakkar, {\it A technique for increasing the utility
of the Wigner-Kirkwood expansion for the second virial coefficient},
Mol. Phys. {\bf 36} (1978), 887 - 892.

\item {[61]} A. C. Tanner and A. J. Thakkar, {\it Discrete and
continuum contributions to multipole polarizabilities and shielding
factors of hydrogen}, Int. J. Quantum Chem. {\bf 24} (1983), 345 - 352.

\item {[62]} E. J. Weniger, J. Grotendorst, and E. O. Steinborn, {\it
Some applications of nonlinear convergence accelerators}, Int. J.
Quantum Chem. Symp. {\bf 19} (1986), 181 - 191.

\item {[63]} J. Grotendorst and E. O. Steinborn, {\it Use of nonlinear
convergence accelerators for the efficient evaluation of GTO molecular
integrals}, J. Chem. Phys. {\bf 84} (1986), 5617 - 5623.

\item {[64]} J. Grotendorst, E. J. Weniger, and E. O. Steinborn, {\it
Efficient evaluation of infinite-series representations for overlap,
two-center nuclear attraction, and Coulomb integrals using nonlinear
convergence accelerators}, Phys. Rev. A {\bf 33} (1986), 3706 - 3726.

\item {[65]} I. S. Gradshteyn and I. M. Ryzhik, {\it Table of
integrals, series, and products} (Academic Press, New York, 1980).

\item {[66]} H. E. Salzer, {\it A Simple method for summing certain
slowly convergent series}, J. Math. and Phys. (Cambridge, Mass.) {\bf
33} (1954), 356 - 359.

\item {[67]} H. E. Salzer, {\it Formulas for the partial summation of
series}, Math. Tables Aids Comput. {\bf 10} (1956), 149 - 156.

\item {[68]} H. E. Salzer and G. M. Kimbro, {\it Improved formulas for
complete and partial summation of certain series}, Math. Comput. {\bf
15} (1961), 23 - 39.

\item {[69]} J. Wimp, {\it Some transformations of monotone sequences},
Math. Comput. {\bf 26} (1972), 251 - 254.

\item {[70]} A. Sidi, {\it An algorithm for a special case of a
generalization of the Richardson extrapolation process}, Numer. Math.
{\bf 38} (1982), 299 - 307.

\item {[71]} N. E. N\"orlund, {\it Vorlesungen \"uber Differenzenrechnung},
(Chelsea, New York, 1954).

\item {[72]} L. M. Milne-Thomson, {\it The calculus of finite
differences} (Chelsea, New York, 1981).

\item {[73]} A. Sidi, {\it Some properties of a generalization of the
Richardson extrapolation process}, J. Inst. Math. Appl. {\bf 24}
(1979), 327 - 346.

\item {[74]} W. F. Ford and A. Sidi, {\it An algorithm for a
generalization of the Richardson extrapolation process}, SIAM J. Numer.
Anal. {\bf 24} (1987), 1212 - 1232.

\item {[75]} D. Levin and A. Sidi, {\it Two new classes of nonlinear
transformations for accelerating the convergence of infinite integrals
and series}, Appl. Math. Comput. {\bf 9} (1981), 175 - 215.

\item {[76]} A. Sidi and D. Levin, {\it Rational approximations from
the d-transformation}, IMA J. Numer. Anal. {\bf 2} (1982), 153 - 167.

\item {[77]} N. Nielsen, {\it Die Gammafunktion} (Chelsea, New York,
1965).

\item {[78]} N. E. N\"orlund, {\it Le{\c c}ons sur les s\'eries
d'interpolation} (Gautier-Villars, Paris, 1926).

\item {[79]} E. Borel, {\it Le{\c c}ons sur les s\'eries divergentes}
(Gautier-Villars, Paris, 1928). Reprinted by \'Editions Jacques Gabay,
Paris, 1988.

\item {[80]} W. Wasow, {\it Asymptotic expansions for ordinary
differential equations} (Dover, New York, 1987).

\item {[81]} S. Iseki and Y. Iseki, {\it Asymptotic expansion for the
remainder of a factorial series}, Mem. Natl. Defense Acad. Japan {\bf
20} (1980), 1 - 6.

\item {[82]} E. Landau, {\it \"Uber die Grundlagen der Theorie der
Fakult\"atenreihen}, Sitzungsb. K\"onigl. Bay. Akad. Wissensch. M\"unchen,
math.-phys. Kl. {\bf 36} (1906), 151 - 218.

\item {[83]} G. N. Watson, {\it The transformation of an asymptotic
series into a convergent series of inverse factorials}, Rend. Circ.
Mat. Palermo {\bf 34} (1912), 41 - 88.

\item {[84]} A. Sidi, {\it A new method for deriving Pad\'e approximants
for some hypergeometric functions}, J. Comput. Appl. Math. {\bf 7}
(1981), 37 - 40.

\item {[85]} C. M. Bender and T. T. Wu, {\it Anharmonic oscillator},
Phys. Rev. {\bf 184}, (1969), 1231 - 1260.

\item {[86]} B. Simon, {\it The anharmonic oscillator: A singular
perturbation theory}, in D. Bessis (ed.), {\it Carg\`ese lectures in
physics} (Gordon and Breach, New York, 1972), Vol. 5, pp. 383 - 414.

\item {[87]} B. Simon, {\it Large orders and summability of eigenvalue
perturbation theory: A mathematical overview}, Int. J. Quantum Chem.
{\bf 21} (1982), 3 - 25.

\item {[88]} J. {\v C}{\' \i}{\v z}ek and E. R. Vrscay, {\it Large
order perturbation theory in the context of atomic and molecular
physics -- Interdisciplinary aspects}, Int. J. Quantum Chem. {\bf 21}
(1982), 27 - 68.

\item {[89]} H. J. Silverstone, J. G. Harris, J. {\v C}{\' \i}{\v z}ek,
and J. Paldus, {\it Asymptotics of high-order perturbation theory for
the one-dimensional anharmonic oscillator by quasisemiclassical
methods}, Phys. Rev. A {\bf 32} (1985), 1965 - 1980. See p. 1966, Eq.
(1), p. 1977, Eq. (69), and p. 1979, Eq. (71).

\item {[90]} J. {\v C}{\' \i}{\v z}ek, R. J. Damburg, S. Graffi, V.
Grecchi, E. M. Harrell II, J. G. Harris, S. Nakai, J. Paldus, R. Kh.
Propin, and H. J. Silverstone, {\it $1 / R$ expansion for $H_2^+$:
Calculation of exponentially small terms and asymptotics}, Phys. Rev. A
{\bf 33} (1986), 12 - 54. See p. 13, Eq. (2), p. 15, Eqs. (28) and
(29), p. 36, Eq. (229), p. 37, Table IV, pp. 38-39, Eq. (232), p. 43,
Eq. (236), and p. 45, Eq. (238).

\item {[91]} G. Alvarez, {\it Coupling-constant behavior of the cubic
anharmonic oscillator}, Phys. Rev. A {\bf 37} (1988), 4079 - 4083. See
p. 4079, Eq. (4), p. 4081, Eq. (24).

\item {[92]} J. E. Drummond, {\it A formula for accelerating the
convergence of a general series}, Bull. Austral. Math. Soc. {\bf 6}
(1972), 69 - 74.

\item {[93]} E. H. Moore and H. L. Smith, {\it A general theory of
limits}, Amer. J. Math. {\bf 44} (1922), 102 - 121.

\item {[94]} A. M. Gleason, {\it Fundamentals of abstract analysis}
(Addison-Wesley, Reading, Mass., 1966).

\item {[95]} O. Perron, {\it Die Lehre von den Kettenbr\"uchen, Band II:
Analytisch-funktionentheoretische Kettenbr\"uche}, (Teubner, Stuttgart,
1957).

\item {[96]} H. S. Wall, {\it Analytic theory of continued fractions},
(Chelsea, New York, 1973).

\item {[97]} B. Simon, {\it Coupling constant analyticity for the
anharmonic oscillator}, Ann. Phys. {\bf 58} (1970), 76 - 136.

\item {[98]} M. Reed and B. Simon, {\it Methods of modern mathematical
physics IV: Analysis of operators} (Academic Press, New York, 1978).

\item {[99]} P. Wynn, {\it Upon the Pad\'e table derived from a Stieltjes
series}, SIAM J. Numer. Anal. {\bf 5} (1968), 805 - 834.

\item {[100]} A. K. Common, {\it Pad\'e approximants and bounds to series
of Stieltjes}, J. Math. Phys. {\bf 9} (1968), 32 - 38.

\item {[101]} G. D. Allen, C. K. Chui, W. R. Madych, F. J. Narcowich,
and P. W. Smith, {\it Pad\'e approximation of Stieltjes series}, J.
Approx. Theor. {\bf 14} (1975), 302 - 316.

\item {[102]} J. Karlsson and B. von Sydow, {\it The convergence of
Pad\'e approximants to series of Stieltjes}, Ark. Matem. {\bf 14} (1976),
43 - 53.

\item {[103]} A. Sidi, {\it Borel summability and converging factors
for some everywhere divergent series}, SIAM J. Math. Anal. {\bf 17}
(1986), 1222 - 1231.

\item {[104]} L. J. Slater, {\it Generalized hypergeometric functions}
(Cambridge University Press, Cambridge, 1966).

\item {[105]} A. Sidi, {\it Analysis of convergence of the
T-transformation for power series}, Math. Comput. {\bf 35} (1980), 833
- 850.

\item {[106]} The NAG Library, Mark 5 (1975), Numerical Analysis Group,
NAG Central Office, Oxford, UK.

\item {[107]} E. J. Weniger and E. O. Steinborn, {\it Nonlinear
sequence transformations for the efficient evaluation of auxiliary
functions for GTO molecular integrals}, in M. Defranceschi and J.
Delhalle (eds.), {\it Numerical determination of the electronic
structure of atoms, diatomic and polyatomic molecules} (Kluwer,
Dordrecht, 1989), pp. 341 - 346.

\item {[108]} R. A. Levy, {\it Principles of solid state physics}
(Academic Press, New York, 1968).

\item {[109]} J. Killingbeck, {\it Quantum-mechanical perturbation
theory}, Rep. Prog. Phys. {\bf 40} (1977), 963 - 1031.

\item {[110]} J. P. Delahaye and B. Germain-Bonne, {\it The set of
logarithmically convergent sequences cannot be accelerated}, SIAM J.
Numer. Anal. {\bf 19} (1982), 840 - 844.

\item {[111]} J. Todd, {\it The lemniscate constants}, Commun. ACM {\bf
18} (1975), 14 - 19.

\item {[112]} E. Filter and E. O. Steinborn, {\it The three-dimensional
convolution of reduced Bessel functions and other functions of physical
interest}, J. Math. Phys. {\bf 19} (1978), 79 - 84.

\item {[113]} E. O. Steinborn and E. Filter, {\it Translations of
fields represented by spherical-harmonic expansions for molecular
calculations III. Translations of reduced Bessel functions, Slater-type
s-orbitals, and other functions}, Theor. Chim. Acta {\bf 38} (1975),
273 - 281.

\item {[114]} E. J. Weniger and E. O. Steinborn, {\it Numerical
properties of the convolution theorems of B functions}, Phys. Rev. A
{\bf 28} (1983), 2026 - 2041.

\item {[115]} E. Grosswald, {\it Bessel polynomials} (Springer-Verlag,
Berlin, 1978).

\item {[116]} H. E. Salzer, {\it Note on the Do{\v c}ev-Grosswald
asymptotic series for generalized Bessel polynomials}, J. Comput. Appl.
Math. {\bf 9} (1983), 131 - 135.

\item {[117]} F. Beleznay, {\it Estimations for asymptotic series using
a modified Romberg algorithm: I. Finite-size scaling calculations}, J.
Phys. A {\bf 19} (1986), 551 - 562.

\item {[118]} C. -M. Liegener, F. Beleznay, and J. Ladik, {\it
Application of a modified Romberg algorithm to Hartree-Fock
calculations on periodic chains}, Phys. Lett. A {\bf 123} (1987), 399 -
401.

\endAbschnittsebene

\bye